\documentclass[12pt]{amsbook}
\usepackage{amsthm}
\usepackage{amssymb}
\usepackage[dvips]{graphicx}

\makeatletter
\renewcommand\thesection{\thechapter.\arabic{section}}

\newcommand{\bmath}[1]{\mbox{\mathversion{bold}$#1$}}

\newcommand{\R}{\bmath{R}}

\newcommand{\SL}{\operatorname{SL}}
\newcommand{\SU}{\operatorname{SU}}
\newcommand{\PSL}{\operatorname{PSL}}

\newcommand{\slg}{\operatorname{sl}}
\newcommand{\su}{\operatorname{su}}

\newcommand{\ord}{\operatorname{ord}}
\newcommand{\id}{\operatorname{I}}

\renewcommand{\Re}{\operatorname{Re}}

\renewcommand{\Im}{\operatorname{Im}}

   \newtheorem{theorem}{Theorem}[section]
   \newtheorem{proposition}[theorem]{Proposition}
   \newtheorem{corollary}[theorem]{Corollary}
   \newtheorem{lemma}[theorem]{Lemma}
   \newtheorem{defn}[theorem]{Definition}

 \theoremstyle{remark}
   \newtheorem{example}[theorem]{Example}
   \newtheorem{remark}[theorem]{Remark}
   
\numberwithin{equation}{section}
\numberwithin{figure}{section}

\def\subsubsection{\@startsection{subsubsection}{3}%
  {\parindent}{.5\linespacing\@plus.7\linespacing}{-.5em}%
  {\bf}}
\def\paragraph{\@startsection{paragraph}{4}%
  {\z@}{\z@}{-\fontdimen2\font}%
  \normalfont\itshape}
\def\dfrac#1#2{{\displaystyle\frac{#1}{#2}}}
\def\eqnarray{%
 \stepcounter{equation}%
 \let\@currentlabel=\theequation
 \global\@eqnswtrue
 \global\@eqcnt\z@
 \tabskip\@centering
 \let\\=\@eqncr
 $$\halign to \displaywidth\bgroup\@eqnsel\hskip\@centering
 $\displaystyle\tabskip\z@{##}$&\global\@eqcnt\@ne
 \hfil$\displaystyle{{}##{}}$\hfil
 &\global\@eqcnt\tw@$\displaystyle\tabskip\z@{##}$\hfil
 \tabskip\@centering&\llap{##}\tabskip\z@\cr}
\def\appendix{\par
  \c@section\z@
  \let\sectionname\appendixname
  \def\thesection{\@Alph\c@section}}
\makeatother
\begin{document}
\begin{titlepage}
\vspace{24pt}
\begin{center}
{\Huge\bf Loop Group Methods for} \\  \vspace{10pt}
{\Huge\bf Constant Mean Curvature Surfaces}\\
\vspace{70pt}
{\LARGE Shoichi Fujimori}  \\ \vspace{12pt}
{\LARGE Shimpei Kobayashi} \\ \vspace{12pt}
{\LARGE Wayne Rossman} \\ 
\vspace{150pt}
\includegraphics[width=4.5in, height=-4.5in]{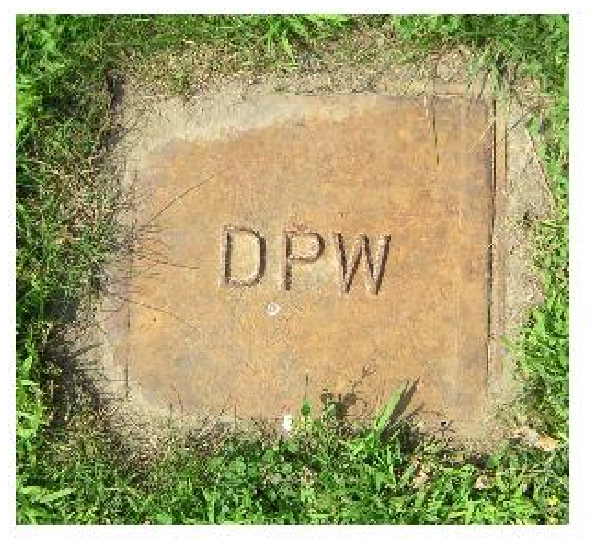}
\end{center}
\end{titlepage}

\pagenumbering{roman}
\setcounter{page}{2}

\chapter*{Introduction}

This is an elementary introduction to a method for studying harmonic 
maps into symmetric spaces, and for studying 
constant mean curvature (CMC) surfaces, that was developed by J. 
Dorfmeister, F. Pedit and H. Wu, and is often called 
the DPW method after them.  
There already exist a number of other introductions to this method, 
but all of them require a higher degree of mathematical sophistication from 
the reader than is needed here.  
The authors' goal was to create an exposition that would be readily accessible 
to a beginning graduate student, and even to a highly motivated undergraduate 
student.  The material here is elementary in the following ways: 
\begin{enumerate}
\item we include an introductory chapter to explain notations; 
\item we include many computations that 
are too trivial for inclusion in research papers, 
and hence are generally not found there, but nonetheless 
are often useful for newcomers to the field; 
\item essentially the only symmetric space we use is $\mathbb{S}^2$, the 
round unit sphere, allowing us to avoid any discussion 
of symmetric spaces; 
\item we consider the DPW method via its application to CMC surfaces, which 
are concrete objects in the sense that they (at least locally) model 
soap films; 
\item we first apply the DPW recipe to produce basic CMC surfaces 
(Chapter \ref{chapter1}), without explaining why it works, and only later 
(Chapters \ref{chapter2} and \ref{chapter3}) do we gradually explain the 
theory; 
\item we consider surfaces only in the three space forms $\mathbb{R}^3$ 
(Euclidean $3$-space), $\mathbb{S}^3$ (spherical $3$-space) and 
$\mathbb{H}^3$ (hyperbolic $3$-space), 
and we give the most emphasis to the least abstract case $\mathbb{R}^3$; 
\item we show how the theory leads to the Weierstrass representation for 
minimal surfaces in $\mathbb{R}^3$ and the related representation 
of Bryant for CMC $1$ surfaces in $\mathbb{H}^3$ and another related 
representation for flat surfaces in $\mathbb{H}^3$, in Sections 
\ref{Wrepsection}, \ref{bryantrep}, \ref{sectionGMM} 
(which are simpler than the DPW representation, 
as they do not require a spectral parameter); 
\item we assume only a minimal mathematical background from the reader: 
a basic knowledge of differential geometry (including submanifolds, 
fundamental forms and curvature) and just a bit of 
experience with Riemann surfaces, topology, matrix groups and Lie groups, and 
ordinary and partial differential equations.  
\end{enumerate}
This approach is suitable only for newcomers, who could read these 
notes as a precursor to reading research articles using or relating to 
the DPW method.  With this 
purpose in mind, we intentionally avoid the 
more theoretic and concise arguments already found in the 
literature.  After reading 
the utilitarian arguments here, we hope the reader would then find 
the arguments in research papers to be more transparent.  

These notes consists of five distinct parts: Chapter \ref{chapter0} 
explains the notations that appear 
throughout the rest of these notes.  Chapter \ref{chapter1} gives simple 
examples of how to make CMC surfaces in $\mathbb{R}^3$ using the DPW method.  
Chapters \ref{chapter2} and \ref{chapter3} describe the theory 
(Lax pairs, loop groups, and Iwasawa decomposition) behind the DPW recipe.  
Chapter \ref{chapter4} describes Lax pairs and the DPW method in the other 
two simply-connected space forms $\mathbb{H}^3$ and $\mathbb{S}^3$.  

Although we do include some computer graphics of surfaces in these notes, 
there are numerous places where one can find a wide variety of graphics.  
In particular, the web pages 
\begin{center}
{\tt http://www.gang.umass.edu/gallery/cmc/} $\;$, \\
{\tt http://www.math.tu-berlin.de/} $\;$, \\
{\tt http://www.iam.uni-bonn.de/sfb256/grape/examples.html} $\;$, \\
{\tt http://www.msri.org/publications/sgp/SGP/indexc.html} $\;$, \\
{\tt http://www.indiana.edu/\verb+~+minimal/index.html} $\;$, \\
{\tt http://www.ugr.es/\verb+~+surfaces/}
\end{center} 
are good resources.

We remark that an expanded and informal version of the notes here can be 
found in \cite{waynebook-maybe}, containing further details and extra topics.  

The authors are grateful for conversations with Alexander Bobenko, Fran Burstall, 
Joseph Dorfmeister, Martin Guest, Udo Hertrich-Jeromin, Jun-ichi Inoguchi, Martin 
Kilian, Ian McIntosh, Yoshihiro Ohnita, Franz Pedit, Pascal Romon, 
Takashi Sakai, Takeshi Sasaki, 
Nicolas Schmitt, Masaaki Umehara, Hongyu Wu and Kotaro Yamada, and much of the notes 
here is comprised of explanations we have received from them.  
The authors are especially grateful to Franz Pedit, who generously gave 
considerable amounts of his time for numerous insightful explanations.  
The authors also thank Yuji Morikawa, Nahid Sultana, Risa Takezaki, Koichi Shimose 
and Hiroya Shimizu for contributing computations and graphics, 
and Katsuhiro Moriya for finding and taking the picture on the front cover.  

\tableofcontents

\pagenumbering{arabic}

\chapter{CMC surfaces}
\label{chapter0}

\section{The ambient spaces (Riemannian, Lorentzian manifolds)}  
\label{ambientspaces}

CMC surfaces always sit in some larger ambient space.  
Although we will not encounter CMC surfaces in other 
non-Euclidean ambient spaces until 
we arrive at Chapter \ref{chapter4}, in this section we will first 
exhibit all of the ambient spaces that appear in these notes.  

{\bf Manifolds.} 
Let us begin by recalling the definition of a manifold.  

\begin{defn}\label{defn:differentiablemanifold}
An $n$-dimensional {\em differentiable manifold} of class $C^\infty$ 
(resp. of class $C^k$ for some $k \in \bmath{N}$, of real-analytic class)
is a Hausdorff topological space 
$M$ and a family of homeomorphisms $\phi_\alpha: U_\alpha \subseteq \mathbb{R}^n \to 
\phi_\alpha(U_\alpha) \subseteq M$ of open sets $U_\alpha$ of 
$\mathbb{R}^n$ to $\phi_\alpha(U_\alpha) \subseteq M$ such that 
\begin{enumerate}
\item $\cup_{\alpha} \phi_\alpha(U_\alpha) = M$, 
\item for any pair $\alpha,\beta$ with $W = \phi_\alpha (U_\alpha) \cap \phi_\beta 
(U_\beta) \neq \emptyset$, the sets $\phi_\alpha^{-1}(W)$ and $\phi_\beta^{-1}(W)$ are 
open sets in $\mathbb{R}^n$ and the mapping 
$f_{\beta\alpha}:=\phi_\beta^{-1} \circ \phi_\alpha$ 
from $\phi_\alpha^{-1}(W)$ to $\phi_\beta^{-1}(W)$ is $C^\infty$ differentiable 
(resp. $C^k$ differentiable, real-analytic), 
\item the family $\{(U_\alpha,\phi_\alpha)\}$ is maximal relative to the conditions 
(1) and (2) above.  
\end{enumerate}
The pair $(U_\alpha,\phi_\alpha)$ with $p \in \phi_\alpha(U_\alpha)$ is called a 
{\em coordinate chart} of $M$ at $p$, and $\phi_\alpha(U_\alpha)$ is called a 
{\em coordinate neighborhood} at $p$.  A family $\{(U_\alpha,\phi_\alpha)\}$ of 
coordinate 
charts satisfying (1) and (2) is called a 
{\em differentiable structure} on $M$.  The functions $f_{\beta\alpha}$ are called 
{\em transition functions}.  
\end{defn}

Any differentiable structure on $M$ can be extended uniquely to one that is 
maximal, i.e. 
to one that satisfies property (3) in the definition above.  Hence to define a manifold 
it is sufficient to give a differentiable structure on $M$.  

Because the maps $\phi_\alpha$ are all homeomorphisms, a set $A \subseteq M$ is open 
with respect to the topology of $M$ if and only if $\phi_\alpha^{-1}(A \cap 
\phi_\alpha(U_\alpha))$ is open in the usual topology of $\mathbb{R}^n$ for all 
$\alpha$.  

\begin{defn}
We say that a 
subset $M$ of an $m$-dimensional manifold $\hat M$ is an {\em $n$-dimensional 
submanifold} of $\hat M$ (with $n<m$) if there exist 
coordinate charts $(\phi_\alpha,U_\alpha)$ of $\hat M$ so that the restriction 
mappings $(\phi_\alpha,\{(x_1,...,x_n,0,...,0) \in U_\alpha \})$ form a differential 
structure for $M$.  
\end{defn}

{\bf Tangent spaces.} 
Every $n$-dimensional manifold $M$ has a tangent space $T_p M$ defined at 
each of its points $p \in 
M$.  Each tangent space is an $n$-dimensional vector space consisting 
of all linear differentials applied (at $p$) to functions defined on the manifold.  
We now describe the linear differentials in $T_p M$ using coordinate charts.  
If $\phi_\alpha : U_\alpha \subseteq \mathbb{R}^n \to M$ is a coordinate neighborhood at 
\[ p = \phi_\alpha(\hat x_1,...,\hat x_n) \in M \; , \] 
we can define values $x_j$ at points $q \in \phi_\alpha(U_\alpha)$ by 
taking the $j$'th coordinate 
$x_j$ of the point $\phi_\alpha^{-1}(q)=(x_1,...,x_n)$.  In this way, we 
have $n$ functions 
$x_j : \phi_\alpha (U_\alpha) \to \mathbb{R}$ called {\em coordinate functions}.  
So we have two different interpretations of $x_j$, one as a coordinate 
of $\mathbb{R}^n$ and 
the other as a function on $\phi_\alpha (U_\alpha)$; we will use both interpretations, 
and in each case, the interpretation we use can be determined from the context.  These 
functions $x_j$ are examples of smooth functions on $\phi_\alpha (U_\alpha)$, and we 
now define what we mean by "smooth": 

\begin{defn}
A function $f: M \to \mathbb{R}$ is {\em smooth} on a $C^\infty$ differentiable 
manifold $M$ if, for any coordinate chart 
$(U_\alpha,\phi_\alpha)$, $f \circ \phi_\alpha : U_\alpha \to \mathbb{R}$ is a 
$C^\infty$ function with respect to the coordinates $x_1,...,x_n$ of $U_\alpha$ coming 
from the differentiable structure of $M$.  
\end{defn}

Choosing 
one value for $j \in \{1,...,n\}$ and fixing all $x_i$ for $i \neq j$ to be 
constant, one 
has a curve 
\[ c_j(x_j) = \phi_\alpha(\hat x_1,...,\hat x_{j-1},x_j,\hat x_{j+1},...,\hat x_n) \] 
in $M$ parametrized by $x_j \in (\hat x_j-\epsilon_j,\hat x_j+\epsilon_j)$ for some 
sufficiently small $\epsilon_j > 0$.  
We then define the linear differential (the tangent vector) \[ 
\vec{X}_j = \frac{\partial}{\partial x_j} \] 
in $T_p M$ as the derivative of functions on $M$ along 
this curve in $M$ with respect to the parameter $x_j$.  That is, for a smooth function 
$f:M \to \mathbb{R}$, \begin{equation}\label{TpMA} \vec{X}_j (f) = 
\left. \frac{\partial (f \circ \phi_\alpha)}{\partial x_j} 
\right|_{\phi_\alpha^{-1}(p)} \; . \end{equation}  The tangent vectors 
\begin{equation}\label{TpMbasis} \vec{X}_1,...,\vec{X}_n \end{equation} then 
form a basis for 
the tangent space $T_p M$ at each point $p \in \phi_\alpha(U_\alpha)$, i.e. the linear 
combinations (with real scalars $a_j$) 
\[ \vec{X} = a_1 \vec{X}_1 + ... + a_n \vec{X}_n \]
of these $n$ vectors comprise the full tangent space $T_p M$.  
Furthermore, $\vec{X}$ can be described using a curve, analogously to the way 
$\vec{X}_j$ was described using $c_j(x_j)$, as follows: there exists a curve 
\[ c(t): [-\epsilon,\epsilon] \to \phi_\alpha(U_\alpha) \; , \;\;\; c(0) = p \] 
so that $x_j \circ c(t): [-\epsilon,\epsilon] \to \mathbb{R}$ is $C^\infty$ in $t$ for 
all $j$ and 
\begin{equation}\label{TpMD} a_j = \left. \tfrac{d}{dt} (x_j \circ c(t)) 
\right|_{t=0} \; . 
\end{equation}  In other words, 
by the chain rule we have 
\[ \vec{X}(f) = \left. \tfrac{d}{dt} (f \circ c(t)) \right|_{t=0} \]  
for any smooth function $f: M \to \mathbb{R}$.  

{\bf Tangent spaces for submanifolds of ambient vector spaces.} 
In all the cases we will consider in these notes, 
the manifold $M$ is either a vector space itself 
($\mathbb{R}^3$) or a subset of some larger vector space ($\mathbb{R}^4$).  This 
allows us to give a one-to-one correspondence between the linear 
differentials in $T_p M$ 
and actual vectors in the vector space ($\mathbb{R}^3$ or $\mathbb{R}^4$) that are 
placed at $p$ and tangent to the set $M$.  This correspondence can be made as follows: 
The curve $c_j : (\hat x_j - \epsilon_j,\hat x_j + \epsilon_j) \to M$ can be extended 
to the ambient vector space ($\mathbb{R}^3$ or $\mathbb{R}^4$) by composition 
with the inclusion map $\mathcal{I}: M \to \mathbb{R}^3 \text{ or } \mathbb{R}^4$ 
to the curve $\mathcal{I} \circ c_j : (\hat x_j - \epsilon_j,\hat x_j + \epsilon_j) \to 
\mathbb{R}^3 \text{ or } \mathbb{R}^4$.  Then, as this curve 
$\mathcal{I} \circ c_j$ lies in a vector space, we can compute its tangent 
vector at $p$ as 
\[ \left. \frac{\partial (\mathcal{I} \circ c_j)}{\partial x_j} \right|_{x_j = \hat x_j} 
\in \mathbb{R}^3 \text{ or } \mathbb{R}^4 \; . \]  
This vector is tangent to $M$ when placed at $p$, and is the vector in the ambient 
vector space that corresponds to the linear differential $\vec{X}_j$.  Then we can 
extend this correspondence linearly to 
\[ \sum_{j} a_j \vec{X}_j \leftrightarrow \sum_{j} a_j 
\left. \frac{\partial (\mathcal{I} \circ c_j)}{\partial x_j} \right|_{\hat x_j} \; . \]  
In fact, following \eqref{TpMD}, 
the one-to-one correspondence can be described explicitly as follows: 
For any curve \[ c(t):[-\epsilon,\epsilon] 
\to \phi_\alpha(U_\alpha) \subseteq M \; , \;\;\; c(0)=p \; , \]  so that 
$\mathcal{I} \circ c(t)$ is $C^\infty$, the vector 
\[ c^\prime(0) := 
\left. \tfrac{d}{dt} (\mathcal{I} \circ c(t)) \right|_{t=0} \] is tangent to $M$ at 
$p$, and the corresponding linear differential in $T_p M$ will be the operator 
\begin{equation}\label{TpMC} \vec{X}_c = 
\sum_{j=1}^n a_j \vec{X}_j \; , \;\;\; a_j = \left. 
\tfrac{d}{dt} (x_j \circ c(t)) \right|_{t=0} \; . \end{equation}  
This can be seen to be the correct correspondence, because 
it extends the definition \eqref{TpMA} to all of $T_p M$, i.e. 
for any smooth function $f: M \to \mathbb{R}$, we have 
\begin{equation}\label{TpMB} \vec{X}_c (f) = 
\sum_{j=1}^n \left( \left. \frac{d}{dt} (x_j \circ c(t)) \right|_{t=0} \right) 
\vec{X}_j (f) = \left. \frac{d(f \circ c(t))}{dt} \right|_{t=0} \; , \end{equation} 
by \eqref{TpMA} and 
the chain rule; when $c(t)$ is $c_j(t)$ with $t=x_j$, this is precisely \eqref{TpMA}.  

It is because of this correspondence that linear differentials in $T_p M$ are referred 
to as "tangent vectors".  Because this correspondence $\vec{X}_c \leftrightarrow 
c^\prime(0)$ is a linear bijection, we can at times allow ourselves to not 
distinguish between the two different types of objects.  However, we recommend the 
reader to keep the distinction between them in the back of his 
mind, in order to understand the meaning of the tangent space in the theory 
of abstract manifolds.  

Thus we have a one-to-one correspondence between linear differentials 
$\vec{X}_c$ in $T_p M$ and vectors $c^\prime(0)$ in 
the full space $\mathbb{R}^3$ (when $M = \mathbb{R}^3$) or the hyperplane in 
$\mathbb{R}^4$ tangent to a $3$-dimensional $M$ at $p$ (when $M \neq \mathbb{R}^3$).  
In the case of $M=\mathbb{R}^3$, 
this means that the tangent space $T_p M$ at each point $p \in \mathbb{R}^3$ is 
simply another 
copy of $\mathbb{R}^3$.  In the case of $M=\mathbb{S}^3$, the tangent 
space $T_p M$ at each 
point $p \in \mathbb{S}^3$ is a $3$-dimensional hyperplane of $\mathbb{R}^4$ containing 
$p$ and tangent to the sphere $\mathbb{S}^3$.  

{\bf Metrics.}
A metric on a manifold $M$ is a correspondence which associates to each 
point $p \in M$ a symmetric bilinear form $\langle , \rangle_p$ defined from 
$T_p M \times T_p M$ to $\mathbb{R}$ so 
that it varies smoothly in the following sense: If $\phi_\alpha : U_\alpha 
\subseteq \mathbb{R}^n \to M$ is a coordinate neighborhood at $p$, then each 
$g_{ij}$, defined at each $p \in \phi_\alpha(U_\alpha)$ by 
\[ g_{ij} = \left\langle \frac{\partial}{\partial x_i} , \frac{\partial}{\partial x_j} 
\right\rangle_p \; , \] is 
a $C^\infty$ function on $U_\alpha$ for any choices of $i,j \in \{1,...,n\}$.  
Because the inner product is symmetric, we have that $g_{ij}=g_{ji}$ for all $i,j$.  

Now we assume that $M$ is 
$3$-dimensional.  We take a point $p \in M$ and a coordinate chart $(U_\alpha,
\phi_\alpha)$ at $p$.  Then the coordinates $(x_1,x_2,x_3) \in U_\alpha$ produce 
a basis \[ \frac{\partial}{\partial x_1},\frac{\partial}{\partial x_2},
\frac{\partial}{\partial x_3} \] for $T_p M$.  Let us take two arbitrary vectors 
\[ \vec{v} = a_1 \frac{\partial}{\partial x_1}+a_2 \frac{\partial}{\partial x_2}+a_3 
\frac{\partial}{\partial x_3} \; , \;\;\; \vec{w} = b_1 \frac{\partial}{\partial x_1}+
b_2 \frac{\partial}{\partial x_2}+b_3 \frac{\partial}{\partial x_3} \] (with 
$a_j,b_j \in \mathbb{R}$) in $T_p M$ and consider two ways to write the inner product 
$\langle \vec{v} , \vec{w} \rangle_p$.  

Because the inner product is bilinear, we have that 
\[ \langle \vec{v} , \vec{w} \rangle_p = \sum_{i,j=1}^n a_i b_j g_{ij} \; . \]
This can be written as the product of one matrix and two vectors, giving us our first 
way to write the inner product, as follows:
\[ \langle \vec{v} , \vec{w} \rangle_p = 
\begin{pmatrix}
a_1 & a_2 & a_3
\end{pmatrix} 
g \begin{pmatrix}
b_1 \\ b_2 \\ b_3
\end{pmatrix} \; , \] 
where 
\begin{equation}\label{gmetricmatrixform} g := \begin{pmatrix}
g_{11} & g_{12} & g_{13} \\ 
g_{21} & g_{22} & g_{23} \\ 
g_{31} & g_{32} & g_{33} 
\end{pmatrix} \; . \end{equation} 
We can refer to this matrix $g$ as the metric of $M$, defined with respect to the 
basis \[ \frac{\partial}{\partial x_1},\frac{\partial}{\partial x_2},
\frac{\partial}{\partial x_3} \; . \]

We could also have chosen to write the metric as a 
symmetric $2$-form, as follows: First we define the $1$-form $dx_j$ by 
$dx_j(\frac{\partial}{\partial x_i}) = 0$ if $i \neq j$ and 
$dx_j(\frac{\partial}{\partial x_j}) = 1$ and then extend $dx_j$ 
linearly to all vectors in $T_p M$, i.e. 
\[ dx_j ( \vec{v} ) = a_j \] for $j=1,2,3$.  
We then define the symmetric 
product $dx_i dx_j$ on $T_p M \times T_p M$ by 
\[ dx_i dx_j (\vec{v},\vec{w})=\tfrac{1}{2}(dx_i(\vec{v}) dx_j(\vec{w}) + 
dx_i(\vec{w}) dx_j(\vec{v})) \] for the two arbitrary vectors 
$\vec{v},\vec{w} \in T_p M$.  So, for example, we have 
\[ dx_2^2(\vec{v},\vec{w})= \tfrac{1}{2} (dx_2(\vec{v}) dx_2(\vec{w}) + 
dx_2(\vec{w}) dx_2(\vec{v})) = a_2b_2 \] and 
\[ dx_1dx_3(\vec{v},\vec{w})= \tfrac{1}{2} (dx_1(\vec{v}) dx_3(\vec{w}) + 
dx_1(\vec{w}) dx_3(\vec{v})) 
= \frac{a_1b_3+a_3b_1}{2} \; . \]  
Then, defining a symmetric $2$-form by 
\begin{equation}\label{gmetrictwoformform} g = \sum_{i,j=1,2,3} g_{ij} dx_idx_j \; , 
\end{equation} we have 
\[ \langle \vec{v} , \vec{w} \rangle_p = g(\vec{v},\vec{w}) \; . \]

Note that in the definitions \eqref{gmetricmatrixform} and \eqref{gmetrictwoformform}, 
we have described the same metric in two different 
ways.  But because they both represent the same object, we have intentionally given 
both of them the same name "$g$".  

We now define when the metric on $M$ is Riemannian or Lorentzian.

\begin{defn}
If $M$ is a differentiable manifold of dimension $n$ with metric $\langle , \rangle$ so 
that $\langle \vec{v} , \vec{v} \rangle_p > 0$ for every point $p \in M$ and every 
nonzero vector $\vec{v} \in T_p M$, then $M$ is a {\em Riemannian manifold}.  

If $M$ is a differentiable manifold of dimension $n$ with metric $\langle , \rangle$ so 
that $\langle \vec{v} , \vec{v} \rangle_p > 0$ for every point $p \in M$ and every 
nonzero vector $\vec{v}$ in some $(n-1)$-dimensional subspace $\mathcal{V}_p$ of 
$T_p M$, and so that there exists a nonzero 
vector $\vec{v}$ in $T_p M \setminus \mathcal{V}_p$ 
so that $\langle \vec{v} , \vec{v} \rangle_p < 0$ for every $p \in M$, 
then $M$ is a {\em Lorentzian manifold}.  
\end{defn}

Another way of saying this is that when $g$ is written in matrix form with respect 
to some choice of coordinates (as in \eqref{gmetricmatrixform}), $M$ is a 
Riemannian manifold if all of the 
eigenvalues of $g$ are positive, and $M$ is a Lorentzian manifold one 
eigenvalue of $g$ is negative and all the others are positive, for all 
points $p \in M$.  
More briefly, it is often said that the metric $g$ is positive definite in the 
Riemannian case and of signature (-,+,...,+) in the Lorentzian case.  

{\bf Lengths of curves in a Riemannian manifold $M$.}  Given a smooth curve 
\[ c(t) : [a,b] \to M \; , \] the length of the curve can be approximated by 
summing up the lengths of the line segments from $c(\tfrac{(n-j)a+jb}{n})$ to 
$c(\tfrac{(n-j-1)a+(j+1)b}{n})$ for $j=0,1,2,...,n-1$.  Taking the limit as 
$n \to +\infty$, we arrive at 
\begin{equation}\label{curvelength} 
\text{Length}(c(t)) = \int_a^b \sqrt{g(c^\prime(t),c^\prime(t))} dt \; , \end{equation} 
which we define to be the length of the curve $c(t)$ on the interval $[a,b]$.  

{\bf Geodesics and sectional curvature (classical approach via coordinates).}  
From the metric, which we can now write as either $\langle , \rangle$ or 
in matrix form (as in \eqref{gmetricmatrixform}) or as a symmetric $2$-form 
(as in \eqref{gmetrictwoformform}), we can compute the 
sectional curvatures and the geodesics of $M$.  Let us begin with the very classical 
approach using coordinates, which begins with the definition of the Christoffel symbols 
and the components of the Riemannian curvature tensor.  

We define the Christoffel symbols 
\[ \Gamma_{ij}^m = \frac{1}{2} \sum_{k=1}^n 
\left( \frac{\partial g_{jk}}{\partial x_i} + 
\frac{\partial g_{ki}}{\partial x_j} - 
\frac{\partial g_{ij}}{\partial x_k} \right) g^{km}
\; , \] where the $g^{km}$ represent the entries of the inverse matrix $g^{-1}$ of 
the matrix $g = (g_{ij})_{i,j=1,...,n}$.  Then we next define 
\[ R_{ijk}^\ell = \sum_{s=1}^n \left( \Gamma_{ik}^s \Gamma_{js}^\ell - 
\Gamma_{jk}^s \Gamma_{is}^\ell \right) + 
\frac{\partial \Gamma_{ik}^\ell}{\partial x_j} - 
\frac{\partial \Gamma_{jk}^\ell}{\partial x_i} \; . 
\] (The terms $R_{ijk}^\ell$ are actually the components of the Riemannian curvature 
tensor that we will define a bit later.)  

Geodesics of $M$ are curves \[ c(t) : [a,b] \to \phi_\alpha(U_\alpha) \] from an 
interval $[a,b] \in \mathbb{R}$ to $M$, with tangent vector 
\[ c^\prime(t) = \frac{d(\mathcal{I} \circ c(t))}{dt} \] 
corresponding to \[ \vec{X}_c = \sum_{j=1}^n a_j \vec{X}_j \in 
T_{c(t)} M \; , \;\;\; 
a_j = a_j(t) = \tfrac{d}{dt} (x_j \circ c(t)) \] by \eqref{TpMC}, satisfying 
the system of equations 
\begin{equation}\label{geodEqn1} 
\frac{d a_k}{dt} + \sum_{i,j=1}^n a_ia_j \cdot \Gamma_{ij}^k |_{c(t)}
= 0 \end{equation} for $k=1,...,n$ and for all $t \in [a,b]$.  

Geodesics are the curves that are "straight", i.e. not bending, 
with respect to the manifold $M$ (with metric) that they lie in.  However, this is 
not immediately clear from the definition above.  So now let us explain in what sense 
geodesics are the "straight" curves.  We start by noting that geodesics have 
the following properties: 

\begin{enumerate}
\item Geodesics are curves of "constant speed".  More precisely, if $c(t)$ is a geodesic 
of $M$, then \[ \langle \vec{X}_c(t) , \vec{X}_c(t) \rangle_{c(t)} \] is a constant 
independent of $t \in [a,b]$.  This is not difficult to prove.  One can simply check 
that the derivative of this inner product with respect to $t$ is identically zero, 
by using the definitions of $\vec{X}_c(t)$ and $g_{ij}$ and $\Gamma_{ij}^k$ and 
Equation \eqref{geodEqn1}, and also that the chain rule gives $d(g_{ij})/dt = 
\sum_{k=1}^n (\partial (g_{ij})/ \partial x_k) \cdot (d(x_k)/dt)$.  When 
$\langle \vec{X}_c(t) , \vec{X}_c(t) \rangle_{c(t)}$ is identically $1$ and $0 \in 
(a,b)$, we say that $c(t)$ is parametrized by {\em arc-length} starting from $c(0)$.  
\item Let us consider the case that $M$ is a submanifold of some ambient vector space 
$\mathcal{V}$ with inner product $\langle \cdot , \cdot \rangle_{\mathcal{V}}$, 
and that $M$ is given the metric 
$g_M$ which is the restriction of $\langle \cdot , \cdot \rangle_{\mathcal{V}}$ to the 
tangent spaces of $M$.  Let 
$c(t):[a,b] \to M$ be a curve in $M$.  Using the inclusion map $\mathcal{I}$, we 
can consider 
$c^{\prime\prime}(t)=d^2(\mathcal{I} \circ c(t))/(dt^2)$, since $\mathcal{V}$ is a 
vector space.  If $c(t)$ is a geodesic of $M$, it can be shown that 
\[ \langle c^{\prime\prime}(t), \vec{v} \rangle_{\mathcal{V}} = 0 \] for all vectors 
$\vec{v}$ tangent to $M$ at $c(t)$, for all $t \in [a,b]$.  (Proving this, and also 
proving the next item, requires more machinery than we have given here, so we state 
these two items without proof.)  
\item For any two points $p,q$ in a manifold $M$, the path of shortest length from $p$ 
to $q$ (when it exists) is the image of a geodesic.  Furthermore, for any point $p \in 
M$, there exists an open neighborhood $U \subseteq M$ of $p$ so that for any two points 
$q_1,q_2 \in U$ there exists a unique geodesic parametrized by arc-length from $q_1$ to 
$q_2$ fully contained in $U$.  This unique geodesic is the shortest path from 
$q_1$ to $q_2$.  
\end{enumerate}

The first property tells us that geodesics are straight in the sense that they are 
not accelerating nor decelerating.  The second property tells us that, when the manifold 
lies in some ambient vector space, the "curvature direction" (the second derivative 
of $\mathcal{I} \circ c(t)$) of a geodesic is entirely perpendicular to the manifold.  
Thus a geodesic is straight in the sense that it does not "curve" in any direction 
tangent to the manifold.  The third property tells us that geodesics are straight in 
the sense that they are the shortest paths between any two points that are not too 
far apart.  In other words, to travel within the manifold from one point to another 
not-so-distant point, the most expedient route to take is to follow a geodesic.  

In the upcoming examples, we will describe the geodesics for specific examples of 
manifolds with metrics.  We will soon see that the geodesics of Euclidean $3$-space 
$\mathbb{R}^3$ are straight lines, and that the geodesics of spherical $3$-space 
$\mathbb{S}^3$ are great circles.  We will also give simple geometric descriptions 
for the geodesics of hyperbolic $3$-space $\mathbb{H}^3$.  

We now turn to a description of sectional curvature.  We start with 
a Riemannian (or Lorentzian) manifold $M$ with metric $g$. 
The sectional curvature is a value that is assigned to each 
$2$-dimensional subspace ${\mathcal V}_p$ of each tangent space $T_p M$ for any point 
$p \in M$, and geometrically 
it represents the intrinsic (Gaussian) curvature at the point $p$ 
of the $2$-dimensional manifold made by taking the union of short pieces of 
all geodesics through $p$ and 
tangent to ${\mathcal V}_p$.  The meaning of "intrinsic" is that the sectional curvature 
can be computed using only the metric $g$.  In the case of a $2$-dimensional manifold 
(a surface), this union of geodesics through $p$ is locally just the surface itself, 
so the sectional curvature is just the intrinsic Gaussian curvature in the case of 
surfaces.  We will be saying more about the Gaussian curvature of a surface later, but 
for now we simply note that a surface in $\mathbb{R}^3$ is locally convex when 
the Gaussian curvature is positive, but not when the Gaussian curvature is negative.  
(A surface is {\em locally convex} at a point $p$ is there exists an open neighborhood 
of $p$ in the surface that lies to one side of the surface's tangent plane at $p$.)  
For example, the 
$2$-dimensional upper sheet of the $2$-sheeted hyperboloid in $\mathbb{R}^3$ is locally 
convex, and has positive Gaussian curvature with respect to the metric induced by the 
usual inner product of $\mathbb{R}^3$; but the $2$-dimensional $1$-sheeted hyperboloid 
in $\mathbb{R}^3$ is not locally convex, and has negative Gaussian 
curvature with respect 
to the usual inner product of $\mathbb{R}^3$.  The Gaussian curvature of these 
examples (which, as noted above, is the same 
as the sectional curvature in the $2$-dimensional case) can be computed from the 
metric by the procedure we are about to give for computing sectional curvature.  

Although the sectional curvature can be 
computed using only the metric $g$, it is rather complicated to describe explicitly.  
Even though it is a bit long, let us describe once 
how the computation can be made: Take any basis $\vec{v}_1,\vec{v}_2$ of 
${\mathcal V}_p$.
(The resulting sectional curvature will be independent of the choice of this basis.)  
In terms of the coordinates $(x_1,...,x_n)$ determined by a coordinate chart 
$(U_\alpha,\phi_\alpha)$ at $p$, we can write $\vec{v}_1,\vec{v}_2$ as linear 
combinations 
\[ \vec{v}_1 = \sum_{i=1}^n a_i \vec{X}_i \; , \;\;\; 
\vec{v}_2 = \sum_{i=1}^n b_i \vec{X}_i \] of the basis \eqref{TpMbasis} 
of $T_p M$, for coefficients $a_i,b_j \in \mathbb{R}$.  

Then the sectional curvature of the $2$-dimensional subspace ${\mathcal V}_p$ is 
\begin{equation}\label{p14sectcurv} 
K({\mathcal V}_p) = \frac{\langle \sum_{i,j,k,\ell=1}^n R_{ijk}^\ell a_i b_j a_k 
\vec{X}_\ell , \vec{v}_2 \rangle}{\langle 
\vec{v}_1 , \vec{v}_1 \rangle \langle \vec{v}_2 , \vec{v}_2 
\rangle-\langle \vec{v}_1 , \vec{v}_2 \rangle^2} \; . \end{equation}  

Since computing the sectional curvature as above is a rather involved procedure, 
let us consider a trick that will work well for the ambient manifolds we will be 
using.  First we recall the meaning of the differential of a map between two manifolds.  
Given two manifolds $M_1$ and $M_2$ (lying in equal or larger vector 
spaces) and a mapping $\psi:M_1 \to M_2$, we say that $\psi$ is {\em differentiable} 
if, for any coordinate chart $(U_{\alpha,1},\phi_{\alpha,1})$ of $M_1$ and any 
coordinate chart $(U_{\beta,2},\phi_{\beta,2})$ of $M_2$, 
\[ \phi_{\beta,2}^{-1} \circ 
\psi \circ \phi_{\alpha,1} : \phi_{\alpha,1}^{-1}(W) \to 
\phi_{\beta,2}^{-1}(\phi(W)) \; , \;\;\; 
W = \{ p \in \phi_{\alpha,1}(U_{\alpha,1}) \, | \, \psi(p) \in 
\phi_{\beta,2}(U_{\beta,2}) \}
\] is a $C^1$-differentiable map with respect to the coordinates of 
$\phi_{\alpha,1}^{-1}(W)$ and $\phi_{\beta,2}^{-1}(\psi(W))$.  The differential 
$d\psi_p : T_p M_1 \to T_{\psi(p)} M_2$ of $\psi$ at each $p \in M_1$ can be defined as 
follows:  Given a point $p \in M_1$ and a vector $\vec{v}_1 \in T_p M_1$, choose a curve 
$c(t):[-\epsilon,\epsilon] \to M_1$ so that $\vec{v}_1 = \vec{X}_c$.  
Such a curve always exists.  We then have a linear differential $\vec{v}_2 = 
\vec{X}_{\psi \circ c} \in T_{\psi(p)} M_2$ corresponding to the curve $\psi \circ c(t)$ 
through $\psi(p)$ in $M_2$.  So we take 
$d\psi_p(\vec{v}_1)$ to be $\vec{v}_2$, i.e. 
\begin{equation}\label{thedifferential} 
d\psi_p(\vec{v}_1 = \vec{X}_c) = \vec{v}_2 = \vec{X}_{\psi \circ c} 
\in T_{\psi(p)} M_2 \; . 
\end{equation}  
One can check that this definition does not depend on the choice of $c$, provided that 
$c$ produces the desired $\vec{v}_1$ at $p$.  
It can also be shown that $d\psi_p$ is a linear 
map (see Proposition 2.7 of \cite{doCarmo2}, for example), i.e. 
\[ d\psi_p(a \vec{v}+b \vec{w}) = 
a \cdot d\psi_p(\vec{v})+b 
\cdot d\psi_p(\vec{w}) \] for any $\vec{v},\vec{w} \in T_p M_1$ and 
any $a,b \in \mathbb{R}$.  

\begin{defn}
Let $M_1$ and $M_2$ be Riemannian or Lorentzian manifolds with metrics 
$g_1=\langle , \rangle_1$ and $g_2=\langle , \rangle_2$, respectively.  
Consider a differentiable 
map $\psi:M_1 \to M_2$.  If \[ \langle \vec{v} , \vec{w} \rangle_{1,p} = 
\langle d\psi_p(\vec{v}) , d\psi_p(\vec{w}) \rangle_{2,\psi(p)} \] for any $p \in M$ and 
any $\vec{v},\vec{w} \in T_p M$, then $\psi$ is an {\em isometry} from $M_1$ to $M_2$.  
In the case that $M=M_1=M_2$ and $g=g_1=g_2$, we say that $\psi$ is an 
{\em isometry} of $M$ (or rather of $(M,g)$).  
\end{defn}

As we saw above, the sectional curvature of $M$ is determined by $g$, thus it is 
invariant under isometries, because isometries preserve the metric.  So 
if we fix one particular $2$-dimensional subspace ${\mathcal V}_p$ of $T_p M$ at 
one particular point $p \in M$, and 
if we fix another point $q \in M$ and another $2$-dimensional subspace 
${\mathcal V}_q$ of $T_q M$, and if there exists an isometry $\psi:M \to M$ such 
that $\psi(q)=p$ and $d\psi({\mathcal V}_q)={\mathcal V}_p$, then 
${\mathcal V}_p$ and ${\mathcal V}_q$ have the same sectional curvature.  
This allows us to consider the following method for determining if a manifold $M$ 
has constant sectional curvature, i.e. if the sectional curvature is the same 
for every $2$-dimensional subspace of every tangent space of $M$:  

\begin{lemma}\label{onconstantsectcurvature} 
Fix one particular $2$-dimensional subspace ${\mathcal V}_p$ of $T_p M$ at one 
particular point 
$p \in M$.  Suppose that for any point $q \in M$ and any $2$-dimensional subspace 
${\mathcal V}_q$ of $T_q M$ there exists an isometry $\psi:M \to M$ such that 
$\psi(q)=p$ and 
$d\psi({\mathcal V}_q)={\mathcal V}_p$, then every $2$-dimensional subspace 
${\mathcal V}_q$ has the same sectional curvature as the fixed $2$-dimensional subspace 
${\mathcal V}_p$.  Thus $M$ is a manifold of constant sectional curvature.  
\end{lemma}

By this method, without actually computing any sectional curvature, we conclude that 
$M$ has constant sectional curvature.  
Then to find the value of that constant, we need only compute 
the sectional curvature of the single fixed $2$-dimensional subspace ${\mathcal V}_p$ at 
a single fixed point $p \in M$.  We can apply this method below, when we describe the 
specific ambient spaces we will be using in these notes.  

{\bf Geodesics and sectional curvature (modern coordinate free approach).}  
The above classical description of geodesics and sectional curvature depended very 
essentially on coordinates.  But in fact, geodesics and sectional curvature are 
determined by the metric of $M$, and the metric is an object that is defined on 
the tangent spaces $T_p M$ without need of choosing particular coordinate charts.  
Thus we can expect that geodesics and sectional curvature should also be describable 
without making particular choices of coordinate charts.  We now give that 
description, which gives much more concise formulas for the definitions of geodesics 
and sectional curvature.  

Smooth vector fields on $M$ are maps from $\vec{X}=M \to T M = 
\cup_{p \in M} T_p M$ that 
are $C^\infty$ with respect to the coordinates $x_j$ of any coordinate charts.  
Let $\mathfrak{X}(M)$ denote the set of 
smooth vector fields on $M$.  Given a vector field 
$\vec{X} \in \mathfrak{X}(M)$, 
its evaluation $\vec{X}_p$ at any given point $p \in M$ is 
a linear differential; so for any smooth function $f : M \to \mathbb{R}$, 
$\vec{X}_p(f)$ gives us 
another $C^\infty$ function on $M$.  Also, we can define the Lie bracket 
$[\vec{X},\vec{Y}]$ 
by $[\vec{X},\vec{Y}](f) = \vec{X}(\vec{Y}(f))-\vec{Y}(\vec{X}(f))$, which also turns 
out to be a linear differential, and so $[\vec{X},\vec{Y}] \in T_p M$.  

\begin{defn}
\label{Riem-connection}
A $C^\infty$ {\em Riemannian connection} $\nabla$ on $M$ is a mapping $\nabla : 
\mathfrak{X}(M) \times \mathfrak{X}(M) \to \mathfrak{X}(M)$ denoted by $\nabla: 
(\vec{X},\vec{Y}) \to \nabla_{\vec{X}} \vec{Y}$ 
which has the following four properties: For all $f,g \in C^\infty (M)$ and all 
$\vec{X}, \vec{Y}, \vec{Z}, \vec{W} \in \mathfrak{X}(M)$, we have 
\begin{enumerate}
\item $\nabla_{f \vec{X} + g \vec{Y}} (\vec{Z}+\vec{W}) = f \nabla_{\vec{X}} \vec{Z} 
+ g \nabla_{\vec{Y}} \vec{Z} + f \nabla_{\vec{X}} \vec{W} + g \nabla_{\vec{Y}} \vec{W}$, 
\item $\nabla_{\vec{X}} (f \vec{Y}) = f \nabla_{\vec{X}} \vec{Y} 
+ (\vec{X}(f)) \vec{Y}$, 
\item $[\vec{X},\vec{Y}] = \nabla_{\vec{X}} \vec{Y} - \nabla_{\vec{Y}} \vec{X}$, 
\item $\vec{X}(\langle \vec{Y},\vec{Z} \rangle)=
\langle \nabla_{\vec{X}} \vec{Y},\vec{Z} \rangle 
+ \langle \vec{Y},\nabla_{\vec{X}} \vec{Z} \rangle$.  
\end{enumerate}
\end{defn}

The fundamental result in the theory of Riemannian geometry is this: 

\begin{theorem}
For a Riemannian manifold $M$, there exists a uniquely determined Riemannian connection 
on $M$.  
\end{theorem}

Since this result is so fundamental, its proof can be found in virtually any book on 
Riemannian geometry, so we will not repeat the proof here.  

A computation shows that, in local coordinates 
$x_1,...,x_n$, if we have two arbitrary vector fields 
\[ \vec{X} = \sum_{i=1}^n a_i \vec{X}_i \; , \;\;\; 
\vec{Y} = \sum_{j=1}^n b_j \vec{X}_j \] with $a_j$ and $b_j$ being 
smooth functions of the coordinates $x_1,...,x_n$, then the object defined by 
\begin{equation}\label{theconnection} 
\nabla_{\vec{X}} \vec{Y} = \sum_{k=1}^n \left(
\sum_{j=1}^n b_j \frac{\partial a_k}{\partial x_j} + \sum_{i,j} \Gamma_{ij}^k a_ib_j 
\right) \vec{X}_k \end{equation} 
satisfies all the conditions of the definition of a 
Riemannian connection.  So, although we just defined this object 
\eqref{theconnection} in terms of local 
coordinates, the fundamental theorem of Riemannian geometry tells us that it is the 
same object regardless of the choice of coordinates.  Thus this object 
$\nabla_{\vec{X}} \vec{Y}$ in 
\eqref{theconnection} is independent of the choice of particular coordinates.  
Furthermore, once specific coordinates are chosen, the Riemannian connection can be 
written as the object \eqref{theconnection} just above.  
Also, looking back at how we defined geodesics in 
\eqref{geodEqn1}, we can now say that geodesics are 
curves $c(t)$ that satisfy 
\begin{equation}\label{geodEqn2} \nabla_{\vec{X}_c} \vec{X}_c = 0 \; , \end{equation} 
where \[ \vec{X}_c(f) = \tfrac{d}{dt} (f \circ c(t)) \] for any 
smooth function $f : M \to \mathbb{R}$.  
This is a concise coordinate-free way to define geodesics.  

Now let us turn to a coordinate-free definition of the sectional curvature.  
First we define the Riemann curvature tensor: 

\begin{defn}
The {\em Riemann curvature tensor} is the multi-linear map $\mathfrak{X}(M) \times 
\mathfrak{X}(M) \times \mathfrak{X}(M) \to \mathfrak{X}(M)$ defined by 
\[ R(\vec{X},\vec{Y}) \cdot \vec{Z} = \nabla_{\vec{X}} (\nabla_{\vec{Y}} \vec{Z}) 
- \nabla_{\vec{Y}} (\nabla_{\vec{X}} \vec{Z}) + 
\nabla_{[\vec{X},\vec{Y}]} \vec{Z} \; . \]  
\end{defn}

Because the Riemann curvature tensor above was defined using only the Riemannian 
connection, it has been written in a coordinate-free form.  If we wish to write the 
Riemann curvature tensor in terms of local coordinates (i.e. in terms of the components 
$R_{ijk}^\ell$ that were defined using local coordinates), this can be done by using the 
relations \[ R(\vec{X}_k,\vec{X}_\ell) \cdot 
\vec{X}_i = \sum_{j=1}^n R_{\ell k i}^j \vec{X}_j \; , \] which follow from 
\eqref{theconnection} and the definition of the $R_{\ell k i}^j$.  

Looking back at our coordinate-based 
definition of the sectional curvature of a $2$-dimensional 
subspace $\mathcal{V}_p$ of $T_p M$ spanned by two vectors (i.e. linear 
differentials) $\vec{X}_p, \vec{Y}_p 
\in T_p M$, we can now write the sectional curvature in a coordinate-free form: 
\begin{equation}\label{p16sectcurv} 
K(\mathcal{V}_p)=\frac{-\langle R(\vec{X}_p,\vec{Y}_p) \cdot \vec{X}_p , \vec{Y}_p 
\rangle}{\langle \vec{X}_p , \vec{X}_p \rangle \langle \vec{Y}_p , \vec{Y}_p \rangle - 
\langle \vec{X}_p , \vec{Y}_p \rangle^2} \; . \end{equation}

We now describe the ambient spaces that we will be using in these notes.  The 
Euclidean $3$-space $\mathbb{R}^3$ is the primary ambient space that we will use, 
and the other ones do not appear again until Chapter \ref{chapter4}.  So if the 
reader wishes to, he could read only about $\mathbb{R}^3$ and save a reading of the 
other spaces for when he arrives at Chapter \ref{chapter4}.  

\subsection{Euclidean $3$-space $\mathbb{R}^3$}

As we saw before, Euclidean $3$-space is 
\[ \mathbb{R}^3 = \{(x_1,x_2,x_3) \, | \, x_j \in \mathbb{R} \} \; . \]  
The standard Euclidean metric $g$ with respect to these standard rectangular coordinates 
$(x_1,x_2,x_3)$ is \[ g = dx_1^2+dx_2^2+dx_3^2 \; , \]  
or equivalently, in matrix form $g$ is the $3 \times 3$ identity matrix.  
Because the functions $g_{ij}$ are all constant we can easily compute that 
the geodesics in $\mathbb{R}^3$ are straight lines parametrized linearly as 
\[ c(t) = (a_1 t + b_1,a_2 t + b_2,a_3 t + b_3) \; , \;\;\; t \in \mathbb{R} \; , \] 
for constants $a_j,b_j \in \mathbb{R}$.  Also, again because the functions $g_{ij}$ are 
all constant, the Christoffel symbols are all zero, and so the sectional curvature of 
$\mathbb{R}^3$ is identically zero.  So $\mathbb{R}^3$ is a simply-connected complete 
$3$-dimensional manifold with constant sectional curvature zero, and in fact it is the 
unique such manifold.  

The set of isometries of $\mathbb{R}^3$ is a group under the composition operation, and 
is generated by 
\begin{enumerate}
\item translations, 
\item reflections across planes containing the origin $(0,0,0)$, and 
\item rotations fixing the origin.  
\end{enumerate}  
Isometries of $\mathbb{R}^3$ are also called {\em rigid motions} of $\mathbb{R}^3$.  The 
reflections and rotations fixing the origin can be described by left multiplication 
(i.e. multiplication on the left) of 
$(x_1,x_2,x_3)^t$ (where the superscript "$t$" denotes transposition) 
by matrices in the orthogonal group 
\[ O(3) = \{A \in M_{3 \times 3} \; | \; A^t \cdot A = I \} = 
\{ A \in M_{3 \times 3} \; | \; \langle \vec{x},\vec{x} \rangle = 
\langle A\vec{x}, A\vec{x} \rangle \; \forall \vec{x} \in \mathbb{R}^3 \} \; . \]  
Here $M_{n \times n}$ denotes the set of all $n \times n$ matrices.  

There exists a rigid motion taking any $2$-dimensional subspace of any tangent space 
to any other.  By Lemma \ref{onconstantsectcurvature}, from this alone we can 
conclude that $\mathbb{R}^3$ is a manifold of constant sectional curvature.  

\subsection{Spherical $3$-space $\mathbb{S}^3$}
\label{spherical3space}

The spherical $3$-space $\mathbb{S}^3$ is the unique simply-connected 
$3$-dimensional complete Riemannian manifold with 
constant sectional curvature $+1$.  In fact, $\mathbb{S}^3$ is compact.  
As defined above, 
\[ \mathbb{S}^3 = \left\{ (x_1,x_2,x_3,x_4) \in \mathbb{R}^4 \, \left| \, 
\sum_{j=1}^4 x_j^2 = 1 \right. \right\} \; . \]  The metric $g$ on 
$\mathbb{S}^3$ can be defined by restricting 
the metric $g_{\mathbb{R}^4} = dx_1^2+dx_2^2+dx_3^2+dx_4^2$ to the $3$-dimensional 
tangent spaces of $\mathbb{S}^3$, and then the form of $g$ (but not $g$ itself) 
will depend on the choice of coordinates we make for $\mathbb{S}^3$.  

The isometries of $\mathbb{S}^3$ are the restrictions of the rotations and reflections 
of $\mathbb{R}^4$ fixing the origin $(0,0,0,0)$ to $\mathbb{S}^3$.  Thus the isometry 
group of $\mathbb{S}^3$ is represented by 
\[ O(4) = \{A \in M_{4 \times 4} \; | \; A^t \cdot A = I \} = 
\{ A \in M_{4 \times 4} \; | \; \langle \vec{x},\vec{x} \rangle = 
\langle (A\vec{x}^t)^t, (A\vec{x}^t)^t 
\rangle \; \forall \vec{x} \in \mathbb{R}^4 \} \; . \]  
These isometries can take any $2$-dimensional subspace of any tangent space of 
$\mathbb{S}^3$ to any other, so, by Lemma \ref{onconstantsectcurvature}, 
$\mathbb{S}^3$ is also a manifold of constant sectional curvature.  
Thus we conclude that $\mathbb{S}^3$ has constant sectional curvature.  Then to find 
the value of that constant sectional curvature, we can compute the sectional curvature 
of only one fixed $2$-dimensional subspace of any tangent space.  Making this 
computation 
using the Christoffel symbols as described above, we find that the sectional curvature 
$K$ is $+1$.  

The curves \[ c_\rho(t) = (\cos(\rho t),\sin(\rho t),0,0) \in \mathbb{S}^3 \; , \;\;\; 
t \in [0,2 \pi) \] for any constant $\rho \in \mathbb{R}$ are geodesics of 
$\mathbb{S}^3$, and all other geodesics of $\mathbb{S}^3$ can be produced by rotating 
or reflecting the $c_\rho(t)$ by elements of $O(4)$.  Thus the images of the 
geodesics are the great circles in $\mathbb{S}^3$.  

Since $\mathbb{S}^3$ lies in $\mathbb{R}^4$, we cannot simply take a picture of it 
and show it on a printed page.  When we show graphics of CMC surfaces in $\mathbb{S}^3$ 
later in these notes, we will need to choose a different model for $\mathbb{S}^3$.  
The model we will use is the stereographic projection in $\mathbb{R}^4$ from the north 
pole $(0,0,0,1)$, which is a one-to-one mapping of $\mathbb{S}^3 \setminus \{ (0,0,0,1) 
\}$ to the $3$-dimensional vector space 
\[ \{ (x_1,x_2,x_3,0) \in \mathbb{R}^4 \} \]
lying in $\mathbb{R}^4$.  Since the fourth coordinate is identically zero, we can remove 
it and then view this $3$-dimensional subspace of $\mathbb{R}^4$ 
simply as a $3$-dimensional space.  More explicitly, stereographic projection 
$\text{Pr} : \mathbb{S}^3 \setminus \{(0,0,0,1)\} \to \mathbb{R}^3$ is given by 
\[ \text{Pr}(x_1,x_2,x_3,x_4) = 
(\tfrac{x_1}{1-x_4},\tfrac{x_2}{1-x_4},\tfrac{x_3}{1-x_4})
\; , \] and its inverse map $\text{Pr}^{-1}$ is 
\[ \text{Pr}^{-1}(x,y,z) = (x^2+y^2+z^2+1)^{-1} \cdot (2 x, 2 y, 2 z, x^2+y^2+z^2-1)
\; . \]  
The metric $g$ is certainly not the Euclidean metric 
$dx_1^2+dx_2^2+dx_3^2$, as $g$ 
must be inherited from the metric of $\mathbb{R}^4$ restricted to the original 
$3$-sphere $\mathbb{S}^3$.  One can compute that 
\[ g = \left( \frac{2}{1 + x_1^2 + x_2^2 + x_3^2} \right)^2 
(dx_1^2+dx_2^2+dx_3^2) \] in the stereographic projection model of $\mathbb{S}^3$, 
and one can compute directly, with either \eqref{p14sectcurv} or \eqref{p16sectcurv}, 
from this $g$ that $K=+1$.  

The images of the geodesics (i.e. the great circles) of $S^3$ under the projection 
$\text{Pr}$ are a particular collection of circles and lines in $\mathbb{R}^3$ 
that can be 
described explicitly.  (This is done in \cite{Callahan}, for example.)  

\subsection{Minkowski $3$-space $\mathbb{R}^{2,1}$ and $4$-space $\mathbb{R}^{3,1}$}

To make the simplest example of a $3$-dimensional Lorentzian manifold, we could take 
the $3$-dimensional vector space 
\[ \{(x_1,x_2,x_0) \, | \, x_j \in \mathbb{R} \} \]  
and endow it with the non-Euclidean metric \[ g = dx_1^2+dx_2^2-dx_0^2 \; , \]  
or equivalently, in matrix form \[ g = 
\begin{pmatrix}
1 & 0 & 0 \\
0 & 1 & 0 \\
0 & 0 & -1
\end{pmatrix} \; . \]  
With this metric we have the Minkowski $3$-space $\mathbb{R}^{2,1}$, which clearly 
satisfies the conditions for being a Lorentzian manifold.  

In a similar manner, we can produce Minkowski spaces of higher dimensions.  For 
example, Minkowski $4$-space is the $4$-dimensional vector space 
\[ \{(x_1,x_2,x_3,x_0) \, | \, x_j \in \mathbb{R} \} \]  
endowed with the metric \[ g = dx_1^2+dx_2^2+dx_3^2-dx_0^2 \; , \]  
or equivalently, in matrix form \[ g = 
\begin{pmatrix}
1 & 0 & 0 & 0\\
0 & 1 & 0 & 0\\
0 & 0 & 1 & 0\\
0 & 0 & 0 & -1
\end{pmatrix} \; . \]  
This is then the Lorentzian manifold that is 
called Minkowski $4$-space $\mathbb{R}^{3,1}$.  

For these two manifolds, the entries $g_{ij}$ in the metrics $g$ are all constant, 
hence $\mathbb{R}^{2,1}$ and $\mathbb{R}^{3,1}$ have constant sectional curvature $0$, 
and the geodesics are the same as those in $\mathbb{R}^3$ and $\mathbb{R}^4$, that is, 
they are straight lines.  

The isometries of $\mathbb{R}^{2,1}$, respectively $\mathbb{R}^{3,1}$, that fix 
the origin $(0,0,0)$, respectively $(0,0,0,0)$, can be described by the matrix group 
\[ O(2,1) = \{ A \in M_{3 \times 3} \; | \; \langle 
\vec{x},\vec{x} \rangle =\langle (A \vec{x}^t)^t, (A \vec{x}^t)^t \rangle \; \forall 
\vec{x} \in \mathbb{R}^{2,1} \} \]\[ = 
\{ A \in M_{3 \times 3} \; | \; A^{\mathcal T} \cdot A = I \} \; , \] respectively, 
\[ O(3,1) = \{ A \in M_{4 \times 4} \; | \; \langle 
\vec{x},\vec{x} \rangle = \langle (A \vec{x}^t)^t, (A \vec{x}^t)^t 
\rangle \; \forall \vec{x} \in \mathbb{R}^{3,1} \} \]\[ = 
\{ A \in M_{4 \times 4} \; | \; A^{\mathcal T} \cdot A = I \} \; . \] 
The superscript $\mathcal T$ means the Lorentz transpose: for 
$A = (a_{ij})_{i,j=1,2,0} \in M_{3 \times 3}$, respectively 
$A = (a_{ij})_{i,j=1,2,3,0} \in M_{4 \times 4}$, 
the Lorentz transpose is the transformation $a_{ij} \rightarrow a_{ji}$ if
$i \neq 0$ and $j \neq 0$ or if $i = j = 0$, 
and $a_{ij} \rightarrow - a_{ji}$ if $i = 0$ or $j = 0$ but not $i = j = 0$.  
So, in particular, for a matrix $A \in M_{4 \times 4}$, if 
\[
  A=\begin{pmatrix}
    a_{11} & a_{12} & a_{13} & a_{10} \\
    a_{21} & a_{22} & a_{23} & a_{20} \\
    a_{31} & a_{32} & a_{33} & a_{30} \\
    a_{01} & a_{02} & a_{03} & a_{00} 
    \end{pmatrix} , \;\;\; \mbox{then} \;\;\; 
  A^{\mathcal T} = \begin{pmatrix}
    a_{11} & a_{21} & a_{31} & -a_{01} \\
    a_{12} & a_{22} & a_{32} & -a_{02} \\
    a_{13} & a_{23} & a_{33} & -a_{03} \\
    -a_{10} & -a_{20} & -a_{30} & a_{00} 
    \end{pmatrix} \; . 
\]

The group $O(2,1)$, respectively $O(3,1)$, is called 
the orthogonal group of $\mathbb{R}^{2,1}$, respectively $\mathbb{R}^{3,1}$.  
The isometry of $\mathbb{R}^{2,1}$ for each $A \in O(2,1)$ is the map 
\[ \vec{x} \in \mathbb{R}^{2,1} \to (A \vec{x}^t)^t \in \mathbb{R}^{2,1} \; , \]  
and the isometry of $\mathbb{R}^{3,1}$ for each $A \in O(3,1)$ is the map 
\begin{equation}\label{R31isometries} 
\vec{x} \in \mathbb{R}^{3,1} \to (A \vec{x}^t)^t \in \mathbb{R}^{3,1} \; . 
\end{equation}  
We note that these do not represent the full sets of isometries of 
$\mathbb{R}^{2,1}$ and $\mathbb{R}^{3,1}$, but only those isometries 
that preserve the origin.  We are interested in these particular isometries because 
they can be used to describe the isometries of 
hyperbolic $3$-space $\mathbb{H}^3$, as we are about to see.  

\subsection{Hyperbolic $3$-space $\mathbb{H}^3$}
\label{hyperbolic3space}

Hyperbolic $3$-space $\mathbb{H}^3$ is the unique simply-connected 
$3$-dimensional complete Riemannian manifold with 
constant sectional curvature $-1$.  However, it can be described by a 
variety of models, each with their own advantages.  Here we describe 
the two models we will need: the Poincare ball model, 
convenient for showing computer graphics; and the Hermitian 
matrix model, the one that is best suited for applying the DPW method and also the 
one that Bryant chose for his representation \cite{Br}.  We derive these models from 
the Minkowski model.  

We define $\mathbb{H}^3$ by way of the Minkowski $4$-space 
$\mathbb{R}^{3,1}$ with its Lorentzian metric $g_{\mathbb{R}^{3,1}}$.  
We take $\mathbb{H}^3$ to be the upper sheet of the two-sheeted hyperboloid 
\[ \left\{(x_1,x_2,x_3,x_0) \in \mathbb{R}^{3,1} \, \left| \, 
x_{0}^2-\sum_{j=1}^{3} x_j^2 = 1 \, , \; x_{0} > 0 \right. \right\} \; , \]  
where the metric $g$ is given by the restriction of $g_{\mathbb{R}^{3,1}}$ 
to the tangent spaces of this $3$-dimensional upper sheet.  
We call this the {\em Minkowski model for hyperbolic $3$-space}.  Although the 
metric $g_{\mathbb{R}^{3,1}}$ is Lorentzian and therefore not positive definite, 
the restriction $g$ to the upper sheet is actually positive 
definite, so $\mathbb{H}^3$ is a Riemannian manifold.  

As in the cases of $\mathbb{R}^3$ and $\mathbb{S}^3$, the isometry group of 
$\mathbb{H}^3$ can also be described using a matrix group, the orthogonal group 
\[ O_+(3,1) = \{ A = (a_{ij})_{i,j=1,2,3,0} \in O(3,1) \, | \, a_{00} > 0 \} \] 
of $\mathbb{R}^{3,1}$.  For $A \in O_+(3,1)$, the map in \eqref{R31isometries} 
is an isometry of $\mathbb{R}^{3,1}$ that preserves the Minkowski model for 
$\mathbb{H}^3$, hence it is an isometry of $\mathbb{H}^3$.  In fact, all isometries 
of $\mathbb{H}^3$ can be described this way.  

The above definition for the Minkowski model for hyperbolic $3$-space does not 
immediately imply that it has all the properties we wish it to have.  However, 
the following lemma tells us that the Minkowski model for hyperbolic $3$-space
is indeed the true hyperbolic $3$-space.  Lemma 
\ref{onconstantsectcurvature} is helpful for proving this result, and the proof 
can be found in the supplement \cite{waynebook-maybe} to these notes.  

\begin{lemma}\label{H3hascurvatureminus1}
The Minkowski model $\mathbb{H}^3$ for the hyperbolic $3$-space is a simply-connected 
$3$-dimensional complete Riemannian manifold with constant sectional curvature $-1$.  
\end{lemma}

Because the isometry group of $\mathbb{H}^3$ is the matrix group 
$O_+(3,1)$, the image of the geodesic $\alpha(t)=(0,0,\cosh t,\sinh t)$ under an 
isometry of $\mathbb{H}^3$ always lies in a $2$-dimensional plane of $\mathbb{R}^{3,1}$ 
containing the origin.  Thus we can conclude that the image of any 
geodesic in $\mathbb{H}^3$ is formed by the 
intersection of $\mathbb{H}^3$ with a $2$-dimensional plane in $\mathbb{R}^{3,1}$ 
which passes through the origin $(0,0,0,0)$ of $\mathbb{R}^{3,1}$.  

The Minkowski model is perhaps the best model of $\mathbb{H}^3$ for understanding 
the isometries and geodesics of $\mathbb{H}^3$.  However, since the Minkowski model 
lies in the $4$-dimensional space $\mathbb{R}^{3,1}$, we cannot use it to view graphics 
of surfaces in $\mathbb{H}^3$.  So we would like to have models that can be viewed 
on the printed page.  We would also like to have a model that uses $2 \times 2$ matrices 
to describe $\mathbb{H}^3$, as this is more compatible with the DPW method that lies 
at the heart of these notes.  With this in mind, 
we now give two other models for $\mathbb{H}^3$.  

{\bf The Poincare model:} 
Let $\mathcal P$ be the $3$-dimensional ball in $\mathbb{R}^{3,1}$ lying
in the hyperplane 
$\{ x_0 = 0 \}$ with Euclidean radius $1$ and center at the origin $(0,0,0,0)$.
By Euclidean stereographic projection from the point $(0,0,0,-1) \in 
\mathbb{R}^{3,1}$ of the Minkowski model for $\mathbb{H}^3$ to $\mathcal P$, 
one has the Poincare model $\mathcal P$ for $\mathbb{H}^3$.  This stereographic 
projection is 
\begin{equation}\label{H3toPoincar} (x_1,x_2,x_3,x_0) \in \mathbb{H}^3 \to 
\left( \frac{x_1}{1+x_0}, \frac{x_2}{1+x_0}, \frac{x_3}{1+x_0} \right) \in 
\mathcal{P} \; . \end{equation}  
$\mathcal P$ is given the metric $g$ that makes this stereographic 
projection an isometry.  We 
can view the Poincare model as the Euclidean unit ball 
\[ B^3 = \{(x_1,x_2,x_3) \in \mathbb{R}^3 \; | \; \; 
x_1^2+x_2^2+x_3^2 < 1\} \]  in 
$\mathbb{R}^3 = \{(x_1,x_2,x_3) \in \mathbb{R}^3 \}$.  One can compute that 
the hyperbolic metric $g$ is 
\begin{equation}\label{hypmet} g = \left( \frac{2}{1 - x_1^2-x_2^2-x_3^2} \right)^2 
(dx_1^2+dx_2^2+dx_3^2) . \end{equation}  
This metric \eqref{hypmet} is the one that will make the stereographic 
projection \eqref{H3toPoincar} an isometry.  By 
either Equation \eqref{p14sectcurv} or \eqref{p16sectcurv}, the sectional 
curvature is constantly $-1$.  This metric $g$ in 
\eqref{hypmet} is written as a function times the 
Euclidean metric $dx_1^2+dx_2^2+dx_3^2$, and we will see later that this means 
that the Poincare model's metric is conformal to the Euclidean metric.  From 
this it follows that angles between vectors in the tangent spaces are the 
same from the viewpoints of both the hyperbolic and Euclidean metrics, and this 
is why we prefer this model when showing graphics of surfaces in hyperbolic $3$-space.  
However, distances are clearly not Euclidean.  In fact, the boundary 
\[ \partial B^3 = \{(x_1,x_2,x_3) \in \mathbb{R}^3 \; | \; \; 
x_1^2+x_2^2+x_3^2 = 1 \} \] of the Poincare model is infinitely far from any 
point in $B^3$ with respect to the hyperbolic metric $g$ in 
\eqref{hypmet}.  For example, consider the curve 
\[ c(t) = (t,0,0,0) \; , \;\;\; t \in [0,1) \] in the Poincare model.  Its 
length is 
\[ \int_0^1 
\sqrt{g(c^\prime(t),c^\prime(t))} dt = \int_0^1 \tfrac{2 dt}{1-t^2} = + \infty 
\; . \]  Thus the point $(0,0,0,0)$ is infinitely far from the boundary point 
$(1,0,0,0)$ in the Poincare model.  For this reason, 
the boundary $\partial B^3$ is often called the {\em ideal boundary at infinity}.  

Geodesics in the Poincare model are not Euclidean straight 
lines.  Instead they are segments of Euclidean lines and circles that intersect the 
ideal boundary $\partial B^3$ at right angles.  

Let us define $\mathbb{H}^2$ analogously to the way we defined $\mathbb{H}^3$ but in 
one lower dimension (i.e. $\mathbb{H}^2$ is the unique complete simply-connected 
$2$-dimensional Riemannian manifold with constant sectional curvature $-1$). Then we 
have that the portions of Euclidean spheres and planes inside $B^3$ that 
intersect $\partial B^3$ 
orthogonally are isometric to $\mathbb{H}^2$ when they are given the restriction of 
the Poincare metric to their tangent spaces.  These surfaces in $\mathbb{H}^3$ are 
called {\em totally geodesic hypersurfaces} or {\em hyperbolic planes}.  
(Using terms to be defined later, hyperbolic planes are CMC 
surfaces with constant mean curvature $0$, so they are minimal surfaces in 
$\mathbb{H}^3$.)  

There are some other simply-described surfaces in the Poincare model that are also 
of interest to us.  The portions of Euclidean spheres and planes in 
$B^3$ that intersect $\partial B^3$ 
non-orthogonally but transversally turn out to be CMC surfaces in $\mathbb{H}^3$.  
These surfaces are often called {\em hyperspheres}.  (Without yet defining mean 
curvature, we note that the hyperspheres 
have constant mean curvature whose absolute value lies strictly between $0$ and $1$.)  
The Euclidean spheres that 
lie entirely in $B^3$ are simply called {\em spheres} and turn out to be CMC surfaces 
$\mathbb{H}^3$ (that have constant mean curvature whose absolute value is strictly 
greater than $1$).  The special case that the Euclidean spheres lie in 
$\overline{B^3} = B^3 \cup \partial B^3$ and are tangent to $\partial B^3$ at a 
single point give 
us the {\em horospheres}.  The horospheres turn out to be CMC surfaces in $\mathbb{H}^3$ 
(whose constant mean curvature has absolute value exactly $1$).  

\begin{figure}[phbt]
\begin{center}
\includegraphics[width=0.45\linewidth]{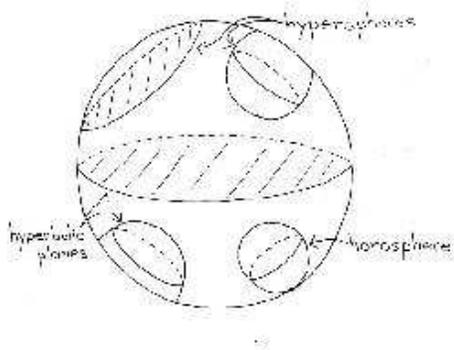}
\end{center}
\vspace{2in}
\caption{Various CMC surfaces in the Poincare model for hyperbolic $3$-space.}
\end{figure}

{\bf The Hermitian matrix model:}
Although we will use the Poincare model to show graphics of surfaces in 
$\mathbb{R}^3$, we will be using the Hermitian matrix model for all 
mathematical computations, because it is the most convenient model 
for making connections with the DPW method.  Unlike the other two 
models above, which can be used for hyperbolic spaces of any dimension, the 
Hermitian model can be used only when the hyperbolic space is $3$-dimensional.  

We first recall the following definitions: 
The group $\SL_2\!\mathbb{C}$ is all $2 \times 2$ matrices 
with complex 
entries and determinant $1$, with matrix multiplication as the group 
operation.  The vector space $\slg_2\!\mathbb{C}$ consists of all 
$2 \times 2$ complex matrices with trace $0$, with the vector space operations 
being matrix addition and scalar multiplication.  
(In Remark \ref{liegroups} we will note that $\SL_2\!\mathbb{C}$ is a Lie 
group.  $\SL_2\!\mathbb{C}$ is $6$-dimensional.  Also, 
$\slg_2\!\mathbb{C}$ is the associated Lie 
algebra, thus is the tangent space of $\SL_2\!\mathbb{C}$ at the identity
matrix.  So $\slg_2\!\mathbb{C}$ is also $6$-dimensional.)  

We also mention, as we use it just below, that 
$\SU_2$ is the subgroup of matrices 
$F \in \SL_2\!\mathbb{C}$ such that $F \cdot F^*$ is the identity matrix, where 
$F^* = \bar{F}^t$.  Equivalently, 
\[ F = \left(
\begin{array}{cc}
p & -\bar{q} \\ q & \bar{p}
\end{array} \right) \; , \]
for some $p$, $q \in \mathbb{C}$ with $|p|^2 + |q|^2 = 1$.  
(We will note that $\SU_2$ is a $3$-dimensional Lie subgroup in 
Remark \ref{liegroups}.)  

Finally, we define 
Hermitian symmetric matrices as matrices of the form 
\[ \left( \begin{array}{cc}
a_{11} & a_{12} \\ \overline{a_{12}} & a_{22} 
\end{array} \right) \; , \] where $a_{12} \in \mathbb{C}$ and $a_{11}, 
a_{22} \in \mathbb{R}$.  Hermitian symmetric matrices with determinant $1$ would then 
have the additional condition that $a_{11} a_{22} - a_{12} \overline{a_{12}} = 1$.  

The Lorentz $4$-space 
$\mathbb{R}^{3,1}$ can be mapped to the space of $2 \times 2$ Hermitian 
symmetric matrices by 
\[ \psi : \vec{x} = (x_1,x_2,x_3,x_0) \longrightarrow 
\psi (\vec{x}) = \left(
\begin{array}{cc}
x_0 + x_3 & x_1 + ix_2 \\
x_1 - ix_2 & x_0 - x_3
\end{array}
\right) \; . \]
For $\vec{x} \in \mathbb{R}^{3,1}$, the metric in the Hermitean matrix form 
is given by $\langle \vec{x} , 
\vec{x} \rangle_{\mathbb{R}^{3,1}} = -\mbox{det}(\psi(\vec{x}))$.  
Thus $\psi$ maps the Minkowski model for $\mathbb{H}^3$ to the set of 
Hermitian symmetric 
matrices with determinant $1$.  Any Hermitian symmetric 
matrix with determinant $1$ can be written as the product $F \bar{F}^t$ for some 
$F \in \SL_2\!\mathbb{C}$, and $F$ is determined uniquely up to right-multiplication 
by elements in $\SU_2$.  That is, for $F,\hat F \in \SL_2\!\mathbb{C}$, we have 
$F \bar{F}^t = \hat F \bar{ \hat F}^t$ if and only if $F = \hat F \cdot B$ for some 
$B \in \SU_2$.  
Therefore we have the Hermitian model 
\[ {\mathcal H} = \{ F F^* \; | \; F \in \SL_2\!\mathbb{C} \} \; , \;\;\; 
\text{where} \;\; F^* := \bar{F}^t \; , \] 
and $\mathcal H$ is given the metric so that $\psi$ is an isometry from 
the Minkowski model of $\mathbb{H}^3$ to $\mathcal H$.  

It follows that, when we compare the Hermitean and Poincare models $\mathcal{H}$ and 
$\mathcal{P}$ for $\mathbb{H}^3$, 
the mapping \[ \left( \begin{array}{cc}
a_{11} & a_{12} \\ \overline{a_{12}} & a_{22} 
\end{array} \right) \in {\mathcal H} \to 
\left( \frac{a_{12}+\overline{a_{12}}}{2+a_{11}+a_{22}}, 
\frac{i(\overline{a_{12}}-a_{12})}{2+a_{11}+a_{22}}, 
\frac{a_{11}-a_{22}}{2+a_{11}+a_{22}} \right) \in {\mathcal P} \]
is an isometry from $\mathcal H$ to $\mathcal P$.  

The Hermitian model is actually very convenient for describing the isometries of 
$\mathbb{H}^3$.  Up to scalar multiplication by $\pm 1$, the 
group $\SL_2\!\mathbb{C}$ represents the isometry group of $\mathbb{H}^3$ in the 
Hermitian model $\mathcal H$ in the following way: A matrix $h \in \SL_2\!\mathbb{C}$ 
acts isometrically on $\mathbb{H}^3$ in the model $\mathcal H$ by 
\[x \in \mathcal{H} \to h \cdot x := h \, x \, h^* \in \mathcal{H} \; , \]
where $h^* = \bar{h}^t$.  The kernel of this action is $\pm I$, 
hence $\PSL_2\!\mathbb{C} = \SL_2\!\mathbb{C}/\{\pm I\}$ is the isometry group of 
$\mathbb{H}^3$.  

\section{Immersions}

Now we wish to consider surfaces lying in these ambient spaces.  Before doing 
this, we first describe the notion of immersion in this 
section.  In these notes, we will be 
immersing surfaces into the ambient spaces described in the previous section.  

Given two manifolds $M_1$ and $M_2$ and a smooth mapping $\psi:M_1 \to M_2$, 
we now describe what it means for $\psi$ to be an immersion.  
(By {\em smooth} mapping, we mean that all of the component functions of $\psi$ 
are smooth with respect to the coordinates of the differentiable structures of 
$M_1$ and $M_2$.)  Recall that we defined the differential 
$d\psi$ of $\psi$ in Equation \eqref{thedifferential}.  

\begin{defn}
If the differential $d\psi_p$ is injective for all $p \in M_1$, then $\psi$ is 
an {\em immersion}.  
\end{defn}

We have already seen a number of immersions, because 
nonconstant geodesics are immersions.  For a geodesic $c(t)$, $M_1$ is an interval 
in $\mathbb{R}$ with parameter $t$, and $\psi$ is the map $c$, and $M_2$ is the 
manifold where the geodesic lies.  Surfaces lying in a larger space (of 
dimension at least three) will also be 
immersions, and we will soon be looking at those too.  
For surfaces, $M_1$ is a domain in the 
$(u,v)$-plane ($\mathbb{R}^2$), and $\psi$ is the map from that domain to a 
manifold $M_2$ where the surface lies.  

When $M_2$ is a Riemannian or Lorentzian manifold, its metric $\langle , \rangle_2$ can 
be restricted to the image $\psi(M_1)$.  We can pull this metric back to a metric 
$\langle , \rangle_1$ on $M_1$ by $d\psi$, i.e. 
\begin{equation}\label{inducedfirstfundform} 
\langle \vec{v},\vec{w} \rangle_1 = \langle d\psi(\vec{v}),d\psi(\vec{w}) 
\rangle_2 \end{equation} 
for any two vectors $\vec{v},\vec{w}$ for any tangent space $T_p M_1$ of $M_1$.  

It turns out that for all the surfaces we will be considering (which are always 
immersions of some $2$-dimensional manifold $M_1$ into a $3$-dimensional ambient space), 
we will be able to parametrize the manifold $M_1$ so that it becomes something called a 
{\em Riemann surface} with respect to the metric pulled back from the ambient space.
This is very useful to us, because Riemann surfaces have some very nice properties.  
We comment about Riemann surfaces in Remark \ref{Riemsurfs}.  

For using Riemann surfaces, we will find it useful to complexify 
the coordinates of surfaces; that is, we can rewrite the coordinates 
$u,v \in \mathbb{R}$ of a surface into 
\[ z=u+i v \; , \;\;\; \bar z=u-i v \; . \]  
Then surfaces can be considered as maps with respect to the variables $z$ and 
$\bar z$ instead of $u$ and $v$.  Also, we define the one-forms 
\[ dz=du+i dv \; , \;\;\; d\bar z=du-i dv \; . \]  
We will get our first hints about how we use Riemann surfaces, and how the coordinates 
$z$ and $\bar z$ appear, in Remark \ref{firststatementofconformality}.  

We end this section with a definition of {\em embeddings}: 

\begin{defn}
Let $\psi:M_1 \to M_2$ be an immersion.  If $\psi$ is also a homeomorphism from 
$M_1$ onto $\psi(M_1)$, where $\psi(M_1)$ has the subspace topology induced from 
$M_2$, then $\psi$ is an {\em embedding}.  
\end{defn}

If $\psi$ is an injective immersion, it does not follow in general that $\psi$ 
is an embedding.  However, all of the injectively immersed surfaces we will see 
in these notes will be embeddings.  

In the next section, we consider surfaces in more detail.  

\section{Surface theory}
\label{surface-theory-section}

{\bf The first fundamental form.} 
The first fundamental form is the metric induced by an immersion, as in 
Equation \eqref{inducedfirstfundform}.  So, although we have already discussed 
metrics, we now consider them again in terms of the theory of immersed surfaces.  
To start with, we assume the ambient space is $\mathbb{R}^3$.  
Let \[ f: \Sigma \to \mathbb{R}^3 \] be an immersion 
(with at least $C^2$ differentiability) of a $2$-dimensional domain 
$\Sigma$ in the $(u,v)$-plane into Euclidean 3-space $\mathbb{R}^3$.  
The standard metric $\langle \cdot , \cdot \rangle$ on 
$\mathbb{R}^3$, also written as the $3 \times 3$ identity matrix $g = \id_3$ in 
matrix form and as $g = dx_1^2+dx_2^2+dx_3^2$ in the form of a symmetric $2$-form, 
induces a metric (the first fundamental form) $ds^2: 
T_p \Sigma \times T_p \Sigma \rightarrow \mathbb{R}$ on $\Sigma$, which 
is a bilinear map, 
where $T_p(\Sigma)$ is the tangent space at $p \in \Sigma$.  Here we have introduced 
another notation $ds^2$ for the metric that we were calling $g$ before.  Both $g$ and 
$ds^2$ are commonly used notations and we will use both in these notes.  

Since $(u,v)$ is a coordinate for $\Sigma$ and $f$ is an immersion, 
a basis for $T_p \Sigma$ can be chosen as 
\[ f_u = \left( \frac{\partial f}{\partial u} \right)_p, f_v = \left( 
\frac{\partial f}{\partial v} \right)_p \; , \] then the 
metric $ds^2$ is represented by the matrix 
\[ g = \left( 
\begin{array}{cc}
g_{11} & g_{12} \\
g_{21} & g_{22} 
\end{array}
\right) = \left( 
\begin{array}{cc}
\langle f_u, f_u \rangle & \langle f_u, f_v \rangle \\
\langle f_v, f_u \rangle & \langle f_v, f_v \rangle 
\end{array}
\right) \; , \] that is, \[ ds^2(a_1 f_u+a_2 f_v,b_1 f_u+b_2 f_v)=
\begin{pmatrix} a_1 & a_2 \end{pmatrix} g \begin{pmatrix} b_1 \\ b_2 
\end{pmatrix} \] for any $a_1,a_2,b_1,b_2 \in \mathbb{R}$.  Note that 
here we are equating the vector $a_1 f_u + a_2 f_v \in T_p \Sigma$ with the 
$1 \times 2$ matrix $(a_1 \; a_2)$.  We could also equate $a_1 f_u + a_2 f_v \in 
T_p \Sigma$ with the $2 \times 1$ matrix $(a_1 \; a_2)^t$.  We will often make 
such "equating", without further comment.  

Clearly, $g$ is a symmetric matrix, since 
\[ 
g_{12} = \langle f_u , f_v \rangle = \langle f_v , f_u \rangle = g_{21} \; . 
\]

\begin{remark}\label{firststatementofconformality}
As we will see in Remark \ref{Riemsurfs}, we can choose the 
coordinates $(u,v)$ on $\Sigma$ so that $ds^2$ is a "conformal" metric.  This means 
that the vectors 
$f_u$ and $f_v$ are orthogonal and of equal positive length in $\mathbb{R}^3$ at every 
point $f(p)$.  Since now $f$ is a conformal immersion, $g_{12}=g_{21}=0$ and 
there exists some function $\hat u: U_\alpha \to \mathbb{R}$ so that $g_{11}=g_{22}=4 
e^{2 \hat u}$.  Then, as a symmetric $2$-form, 
$ds^2$ becomes \[ ds^2 = 4 e^{2 \hat u} dz d\bar z = 4 e^{2 \hat u} (du^2+dv^2) 
\; . \]  This is the way of writing 
the first fundamental form that we will most commonly use in these notes.  (We put a 
"hat" over the function $\hat u$ to distinguish it from the coordinate $u$.)  
\end{remark}

\begin{remark}
One can check that $f$ being an immersion is equivalent to $g$ having positive 
determinant.  
\end{remark}

We have already seen how the 
first fundamental form $g$ of a surface 
can be used to find the lengths of curves in the surface, in \eqref{curvelength}.  
We have also seen how to define geodesics in surfaces 
in terms of only $g$ and its derivatives, in \eqref{geodEqn1} and \eqref{geodEqn2}.  
The metric $g$ can also be used to 
find angles between intersecting curves in the surface: Let $c_1(t)$ and $c_2(t)$ 
be two curves in a surface intersecting at the point $p=c_1(t_1)=c_2(t_2)$.  The 
angle between the two curves at this point is the angle $\theta$ between the 
two tangent vectors $c_1^\prime(t_1)$ and $c_2^\prime(t_2)$ in the tangent plane at 
$p$, which can be computed by 
\[ \cos \theta = \frac{ds^2(c_1^\prime(t_1),c_2^\prime(t_2))}
{\sqrt{ds^2(c_1^\prime(t_1),c_1^\prime(t_1)) \cdot 
ds^2(c_2^\prime(t_2),c_2^\prime(t_2))}}
\; . \]  The area of a surface can be computed as well, from the metric $g$: 
The area is \[ \int_\Sigma \sqrt{\det g} dudv \; . \] There are many books 
explaining these uses of $g$, and \cite{doCarmo1}, \cite{doCarmo2}, \cite{kn} are 
just a few examples.  

Any notion that can fully described using only the first fundamental form $g$ 
is called "intrinsic".  Thus area, geodesics, lengths of curves and angles between 
curves are all intrinsic notions.  It is a fundamental result (found by Gauss) that 
the Gaussian curvature (to be defined in the next 
definition) of a surface immersed in $\mathbb{R}^3$ is intrinsic, and this can 
be seen from either Equation \eqref{p14sectcurv} or \eqref{p16sectcurv} (the Gaussian 
curvature will be the same as the sectional curvature in the case of a 
$2$-dimensional surface).  This is an interesting and surprising result, since the 
Gaussian curvature is not defined using only $g$, but also using 
the second fundamental form, as we will now see: 

{\bf The second fundamental form.} 
We can define a unit normal vector to the surface $f(\Sigma)$ on each coordinate chart 
by taking the cross product of $f_u$ and $f_v$, denoted by $f_u \times f_v$, and 
scaling it to have length $1$: 
\[ \vec{N} = \frac{f_u \times f_v}{|f_u \times f_v|} \; , \;\;\; 
|f_u \times f_v|^2 = \langle f_u \times f_v,f_u \times f_v \rangle \; . \]  
Since $f_u$ and $f_v$ are both perpendicular to their cross product $f_u \times f_v$, 
the vector $\vec{N}$ is perpendicular to the tangent plane $T_p \Sigma$ at 
$f(p)$ for each point $f(p)$.  Since the length of $\vec{N}$ is $1$, we conclude 
that $\vec{N}$ is a unit normal vector to the surface $f(\Sigma)$.  This vector is 
uniquely determined up to sign $\pm \vec{N}$, and the sign of $\vec{N}$ is determined 
by the orientation of the coordinate chart, i.e. a different coordinate chart at the 
same point in $f(\Sigma)$ would produce the same unit normal vector if and only if it 
is oriented in the same way as the original coordinate chart.  

Using the normal $\vec{N}$, we can define the second fundamental form of $f$, which 
(like the first fundamental form) is a symmetric bilinear map from $T_p \Sigma \times 
T_p \Sigma$ to $\mathbb{R}$.  Given two vectors $\vec{v},\vec{w} \in T_p \Sigma$, 
the second fundamental form return a value $b(\vec{v},\vec{w}) \in \mathbb{R}$.  
The second fundamental form can be represented by the matrix 
\begin{equation}\label{bmatrixdefined} b = \left( 
\begin{array}{cc}
b_{11} & b_{12} \\
b_{21} & b_{22} 
\end{array}
\right) = - 
\left( 
\begin{array}{cc}
\langle \vec{N}_u, f_u \rangle & \langle \vec{N}_v, f_u \rangle \\
\langle \vec{N}_u, f_v \rangle & \langle \vec{N}_v, f_v \rangle 
\end{array} \right) = 
\left( 
\begin{array}{cc}
\langle \vec{N}, f_{uu} \rangle & \langle \vec{N}, f_{uv} \rangle \\
\langle \vec{N}, f_{vu} \rangle & \langle \vec{N}, f_{vv} \rangle 
\end{array} \right) \; . \end{equation}  
Here $f_{uu}=\tfrac{\partial^2 f}{\partial u^2}$, 
$f_{uv}=\tfrac{\partial^2 f}{\partial u \partial v}$, 
$f_{vu}=\tfrac{\partial^2 f}{\partial v \partial u}$, 
$f_{vv}=\tfrac{\partial^2 f}{\partial v^2}$, and the right-hand equality in 
the matrix equation \eqref{bmatrixdefined} for $b$ follows from taking the partial 
derivatives of $\langle \vec{N},f_u \rangle = \langle \vec{N},f_v \rangle = 0$ 
with respect to $u$ and $v$, i.e. 
\[ 0 = \partial_u \langle \vec N , f_u \rangle = 
\langle \vec N_u , f_u \rangle + \langle \vec N , f_{uu} \rangle \; , \]  
\[ 0 = \partial_u \langle \vec N , f_v \rangle = 
\langle \vec N_u , f_v \rangle + \langle \vec N , f_{vu} \rangle \; , \]  
\[ 0 = \partial_v \langle \vec N , f_u \rangle = 
\langle \vec N_v , f_u \rangle + \langle \vec N , f_{uv} \rangle \; , \]  
\[ 0 = \partial_v \langle \vec N , f_v \rangle = 
\langle \vec N_v , f_v \rangle + \langle \vec N , f_{vv} \rangle \; . \]  
This matrix form for $b$ is written with respect to 
the basis $f_u$, $f_v$ of $T_p \Sigma$, so if 
\[ \vec{v} = a_1 \tfrac{\partial}{\partial u} + a_2 \tfrac{\partial}{\partial v} \; , 
\;\;\; \vec{w} = b_1 \tfrac{\partial}{\partial u} + b_2 \tfrac{\partial}{\partial v} \] 
for some $a_1,a_2,b_1,b_2 \in \mathbb{R}$, then 
\[ b(\vec{v},\vec{w}) = 
\begin{pmatrix}
a_1 & a_2 
\end{pmatrix}
\cdot b \cdot \begin{pmatrix}
b_1 \\ b_2 
\end{pmatrix} \; . \]  
Like the matrix $g$, the matrix $b$ is symmetric, because $f_{uv} = f_{vu}$.  

The second fundamental form can also be written using symmetric $2$-differentials as 
\[ b = b_{11}du^2 + b_{12}dudv + b_{21}dvdu + b_{22}dv^2 \; . \]
Converting, like we did for the first fundamental form in Remark 
\ref{firststatementofconformality}, to the 
complex coordinate $z=u+iv$, we can write the second fundamental form as 
\[ b = Qdz^2+
       \tilde H dzd\bar z+
       \bar Q d\bar z^2 \; , \]  
where $Q$ is the complex-valued function \[ Q:=\tfrac{1}{4} (b_{11}-b_{22}
-ib_{12}-ib_{21}) \] and $\tilde{H}$ is the real-valued function 
\[ \tilde{H} := \tfrac{1}{2} (b_{11}+b_{22}) \; . \]  

\begin{defn}\label{defn:Hopfdifferential}
We call the symmetric $2$-differential $Qdz^2$ the {\em Hopf differential} 
of the immersion $f$.  
\end{defn}

Sometimes, when the complex parameter $z$ is 
already clearly specified, we also simply refer to the function $Q$ as the 
Hopf differential (but strictly speaking, this is an improper use of the 
term).  

The second fundamental form $b$ is useful for computing directional derivatives of 
the normal vector $\vec N$, and for computing the principal and 
Gaussian and mean curvatures of a surface.  To explain this, we now discuss the 
shape operator.  

{\bf The shape operator.} 
For $g$ and $b$ in matrix form, the 
linear map \[ S : \begin{pmatrix} a_1 \\ a_2 \end{pmatrix} \to g^{-1} b \cdot 
\begin{pmatrix} a_1 \\ a_2 \end{pmatrix} \] is the 
{\em shape operator} and represents the map taking the vector 
$\vec{v} = a_1 f_u + a_2 f_v \in T_p \Sigma$ to 
$-D_{\vec{v}} \vec{N} \in T_p \Sigma$, where $D$ is the directional derivative of 
$\mathbb{R}^3$.  More precisely, 
\begin{equation}\label{directionalderiv} 
D_{\vec{v}} \vec N = a_1 \tfrac{\partial}{\partial u} \vec N + 
a_2 \tfrac{\partial}{\partial v} \vec N = a_1 \vec N_u + a_2 \vec N_v \; , 
\end{equation}  or equivalently, if 
\[ c(t) : [-\epsilon,\epsilon] \to f(\Sigma) \] 
is a curve in $f(\Sigma)$ such that $c(0) = p$ and $c^\prime(0) = \vec{v}$, 
then 
\[ D_{\vec{v}} \vec{N} = \left. \frac{d}{dt} \vec{N}_{c(t)} \right|_{t=0} \; . \]  
Let us prove this property of the shape operator: 

\begin{lemma}\label{shapeoperatorlemma}
The shape operator $S = g^{-1}b$ maps a vector $\vec{v} \in T_p \Sigma$ to the 
vector $-D_{\vec{v}} \vec N \in T_p \Sigma$.  
\end{lemma}

\begin{proof}
Since $\langle \vec N, \vec N \rangle=1$, we have 
\[ \partial_u \langle \vec N, \vec N \rangle = 2 \langle \vec N, \vec N_u \rangle 
= \partial_v \langle \vec N, \vec N \rangle = 
2 \langle \vec N, \vec N_v \rangle = 0 \; , \]  
so $\vec N_u, \vec N_v \in T_p \Sigma$.  Thus we can write 
\[ \vec N_u = n_{11} f_u + n_{21} f_v \; , \;\;\; \vec N_v = n_{12} f_u + n_{22} f_v \] 
for some functions $n_{11},n_{12},n_{21},n_{22}$ from the surface to $\mathbb{R}$.  
Writing $\vec v = a_1 f_u + a_2 f_v$, we have 
\[ -D_{\vec{v}} \vec N = -a_1 \vec N_u - a_2 \vec N_v = 
(- a_1 n_{11} - a_2 n_{12}) f_u + (- a_1 n_{21} - a_2 n_{22}) f_v \; . \]  
So when we write the vectors $\vec v$ and $-D_{\vec{v}} \vec N$ in $2 \times 1$ matrix 
form with respect to the basis $\{ f_u,f_v \}$ on $T_p \Sigma$, the map 
$\vec v \to -D_{\vec{v}} \vec N$ becomes 
\[ \begin{pmatrix} a_1 \\ a_2 \end{pmatrix} \to - 
n \cdot 
\begin{pmatrix} a_1 \\ a_2 \end{pmatrix} \; , \;\;\;\;\; \text{ where } n = 
\begin{pmatrix} n_{11} & n_{12} \\ n_{21} & n_{22} \end{pmatrix} \; . \]  
Thus our objective becomes to show that $-n = g^{-1} b$, or rather $-g n = b$.  This 
follows from the relations
\[ b_{11} = - \langle \vec N_u , f_u \rangle = 
- n_{11} \langle f_u , f_u \rangle - n_{21} \langle f_v , f_u \rangle = 
- n_{11} g_{11} - n_{21} g_{12} \; , \]  
\[ b_{21} = - \langle \vec N_u , f_v \rangle = 
- n_{11} \langle f_u , f_v \rangle - n_{21} \langle f_v , f_v \rangle = 
- n_{11} g_{21} - n_{21} g_{22} \; , \]  
\[ b_{12} = - \langle \vec N_v , f_u \rangle = 
- n_{12} \langle f_u , f_u \rangle - n_{22} \langle f_v , f_u \rangle = 
- n_{12} g_{11} - n_{22} g_{12} \; , \]  
\[ b_{22} = - \langle \vec N_v , f_v \rangle = 
- n_{12} \langle f_u , f_v \rangle - n_{22} \langle f_v , f_v \rangle = 
- n_{12} g_{21} - n_{22} g_{22} \; . \]  
\end{proof}

The eigenvalues $k_1, k_2$ and corresponding eigenvectors of 
the shape operator $g^{-1}b$ are the principal curvatures and principal 
curvature directions of the surface $f(\Sigma)$ at 
$f(p)$.  Although a bit imprecise, intuitively the principal curvatures 
tell us the maximum and minimum amounts of bending 
of the surface toward the normal $\vec{N}$ at $f(p)$, 
and the principal curvature directions tell us the directions in 
$T_p(\Sigma)$ of those maximal and minimal bendings.  Because $g$ 
and $b$ are symmetric, 
the principal curvatures are always real, and the principal curvature directions 
are always perpendicular wherever $k_1 \neq k_2$.  

The following definition tells us that the Gaussian and mean curvatures are the 
product and average of the principal curvatures, respectively.  

\begin{defn}\label{DefnOfKandH}
The determinant and the half-trace of the shape operator $g^{-1}b$ of $f : \Sigma \to 
\mathbb{R}^3$ are the {\em Gaussian curvature} $K$ and the {\em 
mean curvature} $H$, respectively.  
The immersion $f$ is {\em CMC} if $H$ is constant, 
and is {\em minimal} if $H$ is identically zero.  
\end{defn}

Thus \[ K = k_1k_2 = 
\frac{b_{11}b_{22}-b_{12}b_{21}}{g_{11}g_{22}-g_{12}g_{21}} \; , \] 
and \[ H = \frac{k_1 + k_2}{2} = \frac{g_{22}b_{11}-g_{12}b_{21}-g_{21}b_{12}+
g_{11}b_{22}}{2(g_{11}g_{22}-g_{12}g_{21})} \; . \]  

\begin{remark}\label{CMCareorientable}
If we would choose the unit normal vector $\vec N$ to be $-\vec N$ instead, i.e. 
if we would switch 
the direction of the normal vector, then the principal curvatures $k_1,k_2$ would change 
sign to $-k_1,-k_2$.  This would not affect the Gaussian curvature $G$, but it would 
switch 
the sign of the mean curvature $H$.  It follows that nonminimal CMC surfaces must have 
globally defined normal vectors.  In particular, nonminimal CMC surfaces are 
orientable.  
There actually do exist nonorientable minimal surfaces in $\mathbb{R}^3$, see 
\cite{MeeksinDuke}.  
\end{remark}

\begin{figure}[phbt]
\begin{center}
\includegraphics[width=0.9\linewidth]{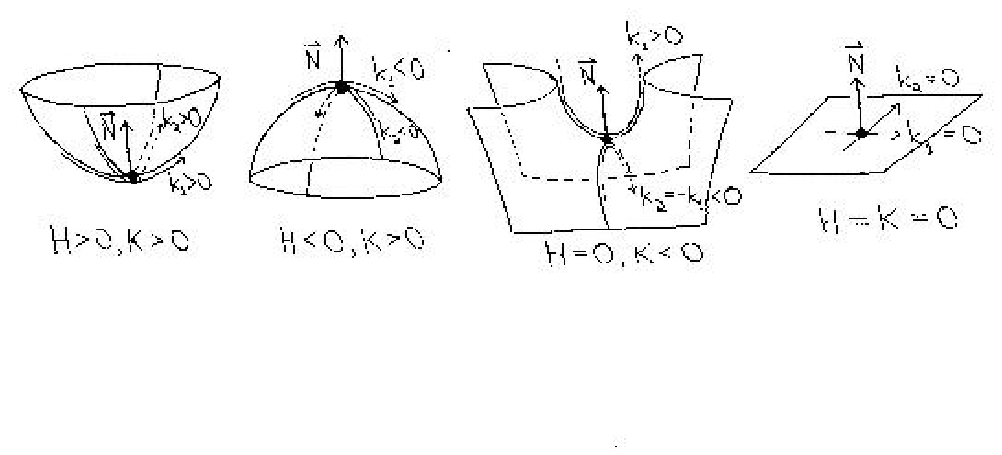}
\end{center}
\vspace{1.6in}
\caption{Various cases of $K>0$, $K<0$, $K=0$, $H>0$, $H<0$, $H=0$, for surfaces in 
$\mathbb{R}^3$.}
\end{figure}

{\bf Non-Euclidean ambient spaces.} 
In Chapter \ref{chapter4} we will consider immersions $f$ of surfaces into other 
$3$-dimensional ambient spaces $M$ (that lie in $4$-dimensional 
Riemannian or Lorentzian vector spaces $\hat M$ with metrics 
$\langle \cdot , \cdot \rangle_{\hat M}$).  To define the first fundamental form 
$g$ for immersions $f : \Sigma \to M \subseteq \hat M$, we can proceed in 
exactly the same way as we did for $\mathbb{R}^3$: 
\[ g = \left( 
\begin{array}{cc}
g_{11} & g_{12} \\
g_{21} & g_{22} 
\end{array}
\right) = \left( 
\begin{array}{cc}
\langle f_u, f_u \rangle_{\hat M} & \langle f_u, f_v \rangle_{\hat M} \\
\langle f_v, f_u \rangle_{\hat M} & \langle f_v, f_v \rangle_{\hat M} 
\end{array}
\right) \; . \]  The only difference is 
that we now use the metric $\langle , \rangle_{\hat{M}}$ rather than the metric of 
$\mathbb{R}^3$.  

To define the Gaussian and mean curvatures of immersions $f : \Sigma \to M$ into 
these other ambient spaces $M$, we need to describe what the shape operator is when 
the ambient space is not $\mathbb{R}^3$.  For this purpose, 
let us first give an alternate way to define the shape operator in $\mathbb{R}^3$.

Suppose that $f : \Sigma \to \mathbb{R}^3$ is a 
smooth immersion with normal vector $\vec N$, 
and that \[ \vec v = a_1 f_u + a_2 f_v \; , \;\;\; \vec w = b_1 f_u + b_2 f_v \] are 
smooth vector fields defined on a coordinate chart 
of $\Sigma$ that lie in the tangent 
space of $f$.  In other words $a_1=a_1(u,v),a_2=a_2(u,v),b_1=b_1(u,v),
b_2=b_2(u,v)$ are smooth real-valued functions of $u$ and $v$.  We can then define the 
directional derivative of $\vec w$ in the direction $\vec v$ as 
\begin{equation}\label{directderivforUandV}
D_{\vec{v}} \vec w = (a_1 \tfrac{\partial}{\partial u} + 
a_2 \tfrac{\partial}{\partial v})(b_1 f_u + b_2 f_v) = 
a_1 \partial_u b_1 f_u + a_1 b_1 f_{uu} + \end{equation} \[ 
a_1 \partial_u b_2 f_v + a_1 b_2 f_{vu} + 
a_2 \partial_v b_1 f_u + a_2 b_1 f_{uv} + 
a_2 \partial_v b_2 f_v + a_2 b_2 f_{vv} \; . \]  
Let $\vec N \in TM$ be the normal to the immersion $f$.  
We can now see that the inner product 
\[ \langle D_{\vec{v}} \vec{w} , \vec{N} \rangle \] at a point $p \in \Sigma$ 
depends only on the values of $a_1,a_2,b_1,b_2$ at $p$ and not on how they extend to a 
neighborhood of $p$ in $\Sigma$.  This is because the only places where derivatives of 
$a_j,b_j$ appear in $D_{\vec{v}} \vec w$ are in the coefficients of $f_u$ and $f_v$, and 
these terms will disappear when we take the inner product with $\vec N$.  We now have 
the relation, by Equation \eqref{directionalderiv} and Lemma \ref{shapeoperatorlemma}, 
\begin{equation}\label{alternateformforS}
\langle \vec w , S(\vec v) \rangle = \langle D_{\vec v} \vec w , \vec N \rangle
\end{equation}
at each $p \in \Sigma$ for the shape operator $S=g^{-1}b$.  Because the 
$\vec v_p,\vec w_p \in T_p \Sigma$ are arbitrary and $S$ is a linear map, 
Equation \eqref{alternateformforS} uniquely determines $S$.  So 
Equation \eqref{alternateformforS} actually gives us an alternate definition of $S$.  

Now, for other 
$3$-dimensional ambient spaces $M$ with Riemannian 
metric $\langle , \rangle_M$ and immersions $f : \Sigma \to M$, 
we can define the shape operator as follows: We take arbitrary 
$\vec v, \vec w \in T_p \Sigma$ (now they must be viewed as linear differentials) 
for $p \in f(\Sigma)$ and insert them into Equation 
\eqref{alternateformforS}.  We then define $\vec N$ so that it has length $1$ and is 
perpendicular to $T_p \Sigma$, and insert this $\vec N$ into 
Equation \eqref{alternateformforS} as well.  Finally, we define the 
Riemannian connection of $M$ as in Definition 
\ref{Riem-connection}, denoted by $\nabla : T_p M \times T_p M \to 
T_p M$, and extend $\vec v,\vec w$ to smooth vector fields of $M$ so that an 
object $\nabla_{\vec{v}} \vec{w} \in T_p M$ is defined, and then replace 
$D_{\vec{v}} \vec{w}$ with $\nabla_{\vec{v}} \vec{w}$ in 
Equation \eqref{alternateformforS}.  Then Equation \eqref{alternateformforS} will 
uniquely determine the linear shape operator $S$ for $f : \Sigma \to M$.  

However, we wish to avoid here a discussion of the linear differentials that comprise 
$T_p M$ and of the general Riemannian connection $\nabla$, and we can get away 
with this because all of the non-Euclidean ambient spaces $M$ we consider 
($\mathbb{S}^3$, $\mathbb{H}^3$) lie in some larger 
vector space $V$ ($\mathbb{R}^4$, $\mathbb{R}^{3,1}$ respectively).  Because 
of this, we can think of $\vec v,\vec w$ as 
actual vectors in the larger vector space $V$, and we do not need to replace 
$D_{\vec{v}} \vec{w}$ by $\nabla_{\vec{v}} \vec{w}$ and so can simply use 
$D_{\vec{v}} \vec{w}$ (as defined in \eqref{directderivforUandV}) in 
Equation \eqref{alternateformforS} to define the shape operator $S$.  Now 
Equation \eqref{alternateformforS} implies that 
\begin{equation}\label{Sinotherspaces}
S = g^{-1} b \; , \;\;\;\; b = \left( \begin{array}{cc}
-\langle \vec{N}_u, f_u \rangle_V & -\langle \vec{N}_v, f_u \rangle_V \\
-\langle \vec{N}_u, f_v \rangle_V & -\langle \vec{N}_v, f_v \rangle_V 
\end{array} \right) \; , 
\end{equation} which is exactly the same as the definition used in $\mathbb{R}^3$, 
except 
that the metric of $\mathbb{R}^3$ has been replaced by the metric $\langle , \rangle_V$ 
of $V$.  Then, like as in $\mathbb{R}^3$: 

\begin{defn}
The Gaussian curvature $K$ and mean curvature $H$ of the 
immersion $f : \Sigma \to M$ are the determinant and half-trace of the shape operator 
$S$ in \eqref{Sinotherspaces}.  
\end{defn}

Also, like as in $\mathbb{R}^3$, we can call the matrix 
$b$ the second fundamental form.  

\begin{remark}
We note that, for general $M$, this $K$ is not the same as the intrinsic curvature.  
\end{remark}

{\bf The fundamental theorem of surface theory.} 
Suppose we define two symmetric $2$-form valued functions 
$g$ and $b$ on a simply-connected open domain $U \subseteq \Sigma$.  Suppose that 
$g$ is positive definite at every point of $U$.  Again, $U$ has a single 
coordinate chart with coordinates 
$(u,v)$, and suppose that $g$ and $b$ are smooth with respect to the variables 
$u$ and $v$.  The fundamental theorem of surface theory 
(Bonnet theorem) tell us there exists an 
immersion $f : U \to \mathbb{R}^3$ with first fundamental form $g$ and second 
fundamental form $b$ (with respect to the coordinates $u$ and $v$) if and only if 
$g$ and $b$ satisfy a pair of equations called the Gauss and 
Codazzi equations.  Furthermore, the fundamental theorem tells us that when such an 
$f$ exists, it is uniquely determined by $g$ and $b$ up to rigid motions (that is, up to 
combinations of rotations and translations) of $\mathbb{R}^3$.  

When such an $f$ exists, we mentioned in Remark \ref{firststatementofconformality} 
that we can choose coordinates $(u,v)$ so that $g$ is conformal.  In this case, 
in the $1$-form formulation for $g$ and $b$, we have 
\begin{equation}\label{gb-fundthm}
g = 4 e^{2\hat u} dz d\bar z \; , \;\;\; b= Qdz^2+\tilde H dzd\bar z+\bar Q d\bar z^2 
\; , \end{equation}
where $z=u+iv$ and $\hat u : U \to \mathbb{R}$, $\tilde H : U \to \mathbb{R}$, 
$Q : U \to \mathbb{C}$ are functions.  The symmetric $2$-form $Q dz^2$ is the 
Hopf differential as in Definition 
\ref{defn:Hopfdifferential}, and as noted before, 
after a specific conformal coordinate $z$ is chosen, 
sometimes just the function $Q$ is referred to as the Hopf differential.  

In the conformal situation, the Gauss and 
Codazzi equations can then be written in 
terms of these functions $\hat u$, $\tilde H$, and $Q$.  
From \eqref{gb-fundthm}, we have that the half-trace of the shape operator $g^{-1}b$ is 
the mean curvature 
\begin{equation}\label{HandHtilde} 
H = \frac{\tilde H}{4 e^{2\hat u}} \; . \end{equation}  Hence the 
three functions $\hat u$, $Q$ and $H$ in \eqref{gb-fundthm} determine the two 
fundamental forms, and 
in turn determine whether the surface $f$ exists.  Then if $f$ does exist, it is 
uniquely determined up to rigid motions of $\mathbb{R}^3$.  

When we consider the theory behind the 
DPW method in Chapter \ref{chapter2}, we will see the exact forms that the Gauss and 
Codazzi equations take.  These two 
equations appear naturally in the matrix formulations 
that we will use to study surface theory and the DPW method in these notes, so we will 
wait until we have discussed the matrix formulations before computing them 
(in Chapter \ref{chapter2}).  For the moment, we only state the 
Gauss and Codazzi equations (for the ambient space $\mathbb{R}^3$), which are 
\[ 4 \partial_{\bar z} \partial_z \hat u-Q \bar Q e^{-2 \hat u}+4 H^2 e^{2 \hat u} 
= 0 \; , \;\;\; \partial_{\bar z} Q = 2 e^{2 \hat u} \partial_z H \; , \]  
respectively.  Here the partial derivatives are defined by 
\[ \partial_{z} = \frac{1}{2} 
\left( \frac{\partial}{\partial u}-i\frac{\partial}{\partial v} \right) \; , \;\;\; 
\partial_{\bar z} = \frac{1}{2} 
\left( \frac{\partial}{\partial u}+i\frac{\partial}{\partial v} \right) \; , \] and they 
are defined this way to ensure that the relations 
\[ dz(\partial_z)=d\bar z(\partial_{\bar z})=1 \; , \;\;\; 
dz(\partial_{\bar z})=d\bar z(\partial_z)=0 \] will hold.  
(Later in these notes, we will often abbreviate notations like 
$\partial_{\bar z} \partial_z 
\hat u$, $\partial_{\bar z} Q$ and $\partial_z H$ to $\hat u_{z\bar z}$, $Q_{\bar z}$ 
and $H_z$.)  

We can draw the following two conclusions from the Gauss and Codazzi equations: 
\begin{enumerate}
\item If $f$ is CMC with $H=1/2$ and if the Hopf differential function 
$Q$ is identically $1$, then the Gauss equation becomes the sinh-Gordon equation 
\[ \partial_{\bar z} \partial_z (2 \hat u) + \sinh(2 \hat u) = 0 \; . \]  
The case $H=1/2$ and $Q=1$ is certainly a special case, but it is an important special 
case, for the following reason: By applying a homothety of $\mathbb{R}^3$ and 
possibly switching the orientation of the surface (see Remark \ref{CMCareorientable}), 
any nonminimal CMC surface can be converted, changing only its size and not its 
shape, into a CMC surface with $H=1/2$.  So the restriction $H=1/2$ is no restriction 
at all.  Then, away from umbilic points (i.e. points where the two principal curvatures 
are equal, and we say more about umbilic points in Section 
\ref{hopf-diff-hopf-thm}), which are usually 
only isolated points of the surface, one can actually 
choose the local coordinates $z$ so that $Q$ is identically $1$.  
Thus, at generic points 
of nonminimal CMC surfaces, the Gauss equation can be locally viewed as the sinh-Gordon 
equation.  
\item If the surface $f$ is CMC, then $\partial_z H=0$ and the Codazzi equation implies 
that the Hopf differential function $Q$ 
is holomorphic, i.e. \[ Q_{\bar z}=0 \; . \]  Furthermore, as we will prove in 
Section \ref{harmonic}, the normal vector $\vec N$, viewed as a map from the surface to 
the $2$-sphere $\mathbb{S}^2$, is a harmonic map if and only if $Q$ is holomorphic.  
Thus $f$ is CMC if and only if the normal map \[ \vec N : \Sigma \to \mathbb{S}^2 \] is 
harmonic.  This fact is what gives 
us a strong connection between the theory of CMC surfaces and the theory of harmonic 
maps.  And this is why the original paper \cite{DPW} on the 
DPW method was formulated in terms of harmonic maps.  
\end{enumerate}

The fundamental theorem applies in the same  way for other $3$-dimensional Riemannian 
manifolds as well, and it is applicable to all of the $3$-dimensional ambient spaces 
$M$ we consider in these notes.  In other $3$-dimensional ambient spaces $M$, 
here again the fundamental theorem for surface theory 
tells us that $f : \Sigma \to M$ is uniquely determined up to rigid motions by its two 
fundamental forms, which again must satisfy the Gauss and Codazzi equations for $f$ 
to even exist at all.  The only changes are: 
\begin{itemize}
\item the Gauss and Codazzi equations will change with respect to the choice of $M$, and 
\item $f$ is uniquely determined up to rigid motions, where the rigid motions are now 
the isometries of the space $M$.  
Hence the rigid motions will change with respect to the choice of $M$.  
\end{itemize}
The fundamental theorem is actually a much more general statement about immersions of 
arbitrary dimension and codimension, but as we do not need that here, we will 
not say anything further about the more general theorem.  

\section{Other preliminaries}
\label{hopf-diff-hopf-thm}

\begin{remark}\label{section1.6}
A discussion on the variational properties of CMC surfaces can be found in 
the supplement \cite{waynebook-maybe} to these notes.
\end{remark}

{\bf Existence of conformal coordinates.}  

\begin{remark}\label{Riemsurfs}
When the dimension of a differentiable manifold $M$ is $2$, then we have some 
special properties.  This is because the coordinate charts are maps from $\mathbb{R}^2$, 
and $\mathbb{R}^2$ can be thought of as the complex plane $\mathbb{C} \approx 
\mathbb{R}^2$.  Thus we can consider the notion of holomorphic functions on 
$M$.  This leads to the theory of Riemann surfaces, which contains a result 
(Theorem \ref{conformalityispossible} just below) that will be 
very important for us here.  
In these notes, we will always be considering CMC surfaces as immersions of 
$2$-dimensional differentiable manifolds $M$.  
Each immersion will determine an induced metric $g$ on $M$ that makes it into a 
Riemannian manifold.  Theorem \ref{conformalityispossible} tells us that 
we can choose coordinates 
on $\Sigma$ so that $g$ is conformal.  Thus without loss of generality we can restrict 
ourselves to those immersions that have conformal induced metric, and we will do this 
on every occasion that we can.  
\end{remark}

To distinguish $2$-dimensional manifolds from other manifolds, we will often denote 
them by $\Sigma$ instead of $M$.  

\begin{theorem}\label{conformalityispossible}
Let $\Sigma$ be a $2$-dimensional orientable manifold with an 
oriented $C^2$ 
family $\{(U_\alpha,\phi_\alpha)\}$ of coordinate charts that determines the 
differentiable structure and with a positive definite 
$C^1$ metric $g$.  Then there exists another family of coordinate charts 
with respect to which the metric $g$ is conformal.  
\end{theorem}

In fact, the $C^2$ and $C^1$ differentiability 
conditions can be weakened to Holder continuity 
conditions, see \cite{chern2}, \cite{bers}, \cite{UY8}.  But as 
will later consider constant 
mean curvature immersions, for which the mean curvature 
cannot even be defined unless the immersion is at least $C^2$, it does us no harm 
to make the stronger assumptions here.  

{\bf The Hopf differential and Hopf theorem.} 
The Hopf differential $Qdz^2$ will be of central importance to us.  We have already 
seen that the Hopf differential can be used to measure whether or not a conformal 
immersion parametrized by a complex coordinate $z$ has 
constant mean curvature, because the surface will have 
constant mean curvature if and only if $Q$ is holomorphic.  The Hopf differential 
can also be used to determine the umbilic points of a surface, as we now see: 

\begin{defn}
Let $\Sigma$ be a $2$-dimensional manifold.  
The {\em umbilic points} of an immersion $f : \Sigma \to \mathbb{R}^3$ are the points 
where the two principal curvatures are equal.  
\end{defn}

The following lemmas and theorem are elementary (proofs can be found in 
the supplement \cite{waynebook-maybe}).  

\begin{lemma}\label{umbIsQis0}
If $\Sigma$ is a Riemann surface and $f : \Sigma \to \mathbb{R}^3$ is a conformal 
immersion, then $p \in \Sigma$ is an umbilic point if and only if $Q=0$.  
\end{lemma}

Thus the Hopf differential tells us where the umbilic points are.  
When $Q$ is holomorphic, it follows that $Q$ is either identically zero or is zero 
only at isolated points.  So, in the case of a 
CMC surface, if there are any points that are not umbilics, then 
all the umbilic points must be isolated.  

If every point is an umbilic, we say that 
the surface is {\em totally umbilic}, and then the 
surface must be a round sphere.  This is proven in \cite{doCarmo1}, for example.  
A proof is also included in \cite{waynebook-maybe}:  

\begin{lemma}\label{lemmaonumbilics}
Let $\Sigma$ be a Riemann surface and $f : \Sigma \to \mathbb{R}^3$ a 
totally umbilic conformal immersion.  Then $f(\Sigma)$ is part of a plane or 
sphere.  
\end{lemma}

\begin{theorem}\label{thmofHopf} {\bf (The Hopf theorem.)}  
If $\Sigma$ is a closed $2$-dimensional manifold of genus zero and if 
$f : \Sigma \to \mathbb{R}^3$ 
is a nonminimal CMC immersion, then $f(\Sigma)$ is a round sphere.  
\end{theorem}

\begin{remark}\label{no-minimal-compact-surfaces}
In fact, there do not exist any compact minimal surfaces without 
boundary in $\mathbb{R}^3$, 
and this can be proved using the maximum principle (see the supplement 
\cite{waynebook-maybe} to these notes).  
Therefore, without assuming that $f$ in the above theorem is nonminimal, the result 
would still be true.  
\end{remark}

Now let us consider the case that $\Sigma$ is a closed Riemann surface of genus 
$\frak{g} \geq 1$ and $f : \Sigma \to \mathbb{R}^3$ 
is a conformal CMC immersion (by Remark 
\ref{no-minimal-compact-surfaces}, because there do not exist any closed compact 
minimal surfaces in $\mathbb{R}^3$, $f$ is guaranteed to be nonminimal).  
In this case, $f(\Sigma)$ certainly cannot be a sphere (because $\frak{g} \geq 1$), 
so $Q$ is not identically zero (by Lemma 
\ref{lemmaonumbilics}).  We then have that (see \cite{fark-kra}, for example)
\[ \sum_{p \in \Sigma} \ord_p(Qdz^2) = 4 \frak{g} - 4 \;\; \text{ and } \;\; 
\ord_p(Qdz^2) \geq 0 \;\;\; \forall \, p \in \Sigma \; . \]  It follows that, 
counted with 
multiplicity, there are exactly $4 \frak{g} -4$ umbilic points on the surface.
We conclude the following: 

\begin{corollary}
A closed CMC surface in $\mathbb{R}^3$ of genus $1$ has no umbilic points, and a closed 
CMC surface in $\mathbb{R}^3$ of genus strictly greater than $1$ must 
have umbilic points.  
\end{corollary}

This corollary is closely related to the question of existence of CMC surfaces 
of "finite type", which is a topic that is explored in the supplement 
\cite{waynebook-maybe} to these notes.  Also, further motivations for studying 
CMC surfaces are given in this supplement.  

\bigskip

{\bf Lie groups and algebras.}  

\begin{remark}\label{liegroups}
Another tool we will need in these notes is Lie groups and algebras.  We say a bit 
more about them in the supplement \cite{waynebook-maybe}, and here we only comment 
that the particular Lie groups and algebras we will use are: 
\[ \text{SO}_3 = \{A \in M_{3 \times 3}(\mathbb{R}) \, | A \cdot A^t = I , 
\det A = 1 \} \; , \;\;\; \text{and} \]  
\[ \SL_2\!\mathbb{C} = \{A \in M_{2 \times 2}(\mathbb{C}) 
\, | \det A = 1 \} \; . \]  In fact, 
$\SL_2\!\mathbb{C}$ is a double cover of $\text{SO}_3$, and has 
corresponding Lie algebra $\slg_2\!\mathbb{C}$, and has the subgroup 
\[ \SU_2 = \{A \in \SL_2\!\mathbb{C} \, | A \cdot \bar{A}^t = I \} 
= \left\{ \left. \begin{pmatrix} p & q \\ -\bar q & \bar p \end{pmatrix} \, 
\right| p,q \in \mathbb{C} , p \bar p+q \bar q = 1 \right\} \; , \] which is also a 
Lie group, with corresponding Lie algebra 
\[ \su_2 = \left\{ \left.\frac{-i}{2}\begin{pmatrix}
        -x_3 & x_1+i x_2 \\
        x_1-i x_2 & x_3
     \end{pmatrix} \; \right| \; x_1,x_2,x_3 \in \mathbb{R} \right\} \; . \] 
\end{remark}

We will have occasion to use the following lemma, which can be proven by 
direct computation:  

\begin{lemma}\label{matrix-identity-lemma}
For any square matrix $\mathcal{A} \in M_{n \times n}$ with $\det A \neq 0$ 
that depends smoothly on some parameter 
$t$, we have \begin{equation}\label{matrix-identity-lemma-eqn} 
\tfrac{d}{dt} (\det \mathcal{A})=\text{trace} 
((\tfrac{d}{dt} \mathcal{A}) \cdot 
\mathcal{A}^{-1}) \det \mathcal{A} \; . \end{equation}  
\end{lemma}

\chapter{Basic examples via the DPW method}
\label{chapter1}

In this chapter we give a description of the DPW method (omitting some 
important technical points that we will return to in Chapters \ref{chapter2} 
and \ref{chapter3}), and then we show how the method can be used to construct 
some basic examples of CMC surfaces.  
The examples we include are the sphere, cylinder, Delaunay surfaces and Smyth 
surfaces.  

\section{The DPW recipe} 
\label{section0}

{\bf $2 \times 2$ matrix notation.}  
The $2 \times 2$ matrix representation begins with writing $\mathbb{R}^3$ as a 
$3$-dimensional matrix vector space.  We will rewrite points 
\[ (x_1,x_2,x_3) \in \mathbb{R}^3 \] into the matrix form 
\begin{equation}\label{first-eqn-in-chap2} \frac{-i}{2}\begin{pmatrix}
        -x_3 & x_1+i x_2 \\
        x_1-i x_2 & x_3
     \end{pmatrix} \; . \end{equation}  There is certainly no harm in doing this, 
as the $3$-dimensional vector space 
\[ \su_2 = \left\{ \left.\frac{-i}{2}\begin{pmatrix}
        -x_3 & x_1+i x_2 \\
        x_1-i x_2 & x_3
     \end{pmatrix} \; \right| \; x_1,x_2,x_3 \in \mathbb{R} \right\} \] is 
isomorphic to $\mathbb{R}^3$.  This way of notating $\mathbb{R}^3$ 
is useful for considering integrable systems techniques.  

{\bf The metric of $\mathbb{R}^3$ in the $\su_2$ matrix formulation.}  
Let \[ \vec x = \frac{-i}{2}\begin{pmatrix}
        -x_3 & x_1+i x_2 \\
        x_1-i x_2 & x_3
     \end{pmatrix} \; , \;\;\; \vec y = \frac{-i}{2}\begin{pmatrix}
        -y_3 & y_1+i y_2 \\
        y_1-i y_2 & y_3
     \end{pmatrix} \] be two points 
in the $\su_2$ matrix formulation for $\mathbb{R}^3$.  (Later we will simply 
say that $\vec x, \vec y$ are contained in $\mathbb{R}^3$, even though they are 
actually matrices contained in a model that is merely equivalent to $\mathbb{R}^3$.)  
Then we find that $-2$ times the trace of the product of $\vec x$ and $\vec y$ 
is equal to $x_1y_1+x_2y_2+x_3y_3$.  Since $x_1y_1+x_2y_2+x_3y_3$ is the standard 
inner product of $\mathbb{R}^3$, we conclude that 
\begin{equation}\label{metricinmatrices1}
\langle \vec x , \vec y \rangle_{\mathbb{R}^3} = -2 \cdot \text{trace} 
(\vec x \cdot \vec y) \; . \end{equation}
So this is how we can compute the $\mathbb{R}^3$ inner product in the 
$\su_2$ matrix model for $\mathbb{R}^3$.  

{\bf The DPW recipe.}
Now we will 
begin to consider the DPW recipe, which uses $2 \times 2$ matrices right from 
the outset.  

The DPW algorithm for producing CMC surfaces starts with a potential 
$\xi$, called the {\em holomorphic potential}, of the form 
\[ \xi = \sum_{j=-1}^\infty \xi_j(z) \lambda^j dz \; , \] where 
$\lambda$ (which we call the {\em spectral parameter}) is the standard 
complex coordinate for the complex plane, and the 
$\xi_j(z)$ are $2\times 2$ matrices that are independent of $\lambda$ and are 
holomorphic functions of a coordinate $z$ of some Riemann surface $\Sigma$.  
Furthermore, the $\xi_j(z)$ are traceless and 
diagonal (resp. off-diagonal) when $j$ is even 
(resp. odd).  We then solve the equation 
\[ d\phi = \phi \xi \] for 
\[ \phi \in \SL_2\!\mathbb{C} \] so that the diagonal (resp. off-diagonal) 
of $\phi$ is a pair of even (resp. odd) functions 
in $\lambda$ (we will call such a $\phi$ ``twisted''), and 
so that $\phi=\phi(z,\lambda)$ is holomorphic in $\lambda$ for $\lambda \neq 
0,\infty$.  ($\phi(z,\lambda)$ is automatically holomorphic in $z$ for 
$z \in \Sigma$.)  Then we split $\phi$ (via Iwasawa splitting, to be 
explained later) into a product 
\[ \phi = F \cdot B \; , \] where \[ F=F(z,\bar{z},\lambda) \in \SU_2 \] for all 
$\lambda \in \mathbb{S}^1:=\{\lambda \in \mathbb{C} \, | \, |\lambda |=1\}$, 
and $B=B(z,\bar{z},\lambda)$ has no negative powers of 
$\lambda$ in its series expansion with respect 
to $\lambda$ and is in fact holomorphic in $\lambda$ for $\lambda \in 
D_1:=\{\lambda \in \mathbb{C} \, | \, |\lambda|<1\}$.  
Furthermore, we choose the splitting so that $B(z,\bar z,0)$ is real 
and diagonal with positive diagonal for all $z \in \Sigma$.  
Such an $F$ and $B$ are uniquely determined.  ($F(z,\bar{z},\lambda)$ and 
$B(z,\bar{z},\lambda)$ are holomorphic in $\lambda$, 
but not in $z$, which is why we write that they are dependent on both $z$ and 
$\bar{z}$.)  Finally, we compute the 
Sym-Bobenko formula 
\[ \left. \left[ \frac{-1}{2} F \begin{pmatrix}
        i & 0 \\
        0 & -i
     \end{pmatrix} F^{-1} - i \lambda (\partial_\lambda F) F^{-1} 
     \right] \right|_{\lambda=1} \; , \]
which is of the form 
\[ \frac{-i}{2}\begin{pmatrix}
        -x_3 & x_1+i x_2 \\
        x_1-i x_2 & x_3
     \end{pmatrix} \; , \] for real-valued functions 
$x_1=x_1(z,\bar{z}),x_2=x_2(z,\bar{z}),x_3=x_3(z,\bar{z})$ and 
$i=\sqrt{-1}$, and then the theory behind the DPW method implies that 
\[ \Sigma \ni z \mapsto (x_1,x_2,x_3) \in \mathbb{R}^3 \] 
becomes a conformal parametrization of a CMC $H=1/2$ surface.  

We have obviously omitted much in this initial brief description of the DPW recipe, 
but we have already said enough to be able to show how this method produces some simple 
well-known CMC surfaces, and we do this in the next few sections.  

\section{Cylinders}
\label{section1}

Here we show how the DPW method makes a cylinder.  Define 
\[ \xi = \lambda^{-1} \begin{pmatrix}
        0 & 1 \\
        1 & 0
     \end{pmatrix} dz \; , \] 
for the complex variable $z \in \Sigma = \mathbb{C}$ 
and $\lambda \in \mathbb{S}^1$.  
Checking that a solution to \[ d\phi = \phi \xi \; , \] is 
\[ \phi = \exp \left( z \begin{pmatrix}
        0 & 1 \\
        1 & 0
     \end{pmatrix} \lambda^{-1} \right) \] and then using that 
$\exp A = I+A+(1/2!)A^2+\cdots$, we have 
\[\phi = \begin{pmatrix}
        \cosh(\lambda^{-1} z) & \sinh(\lambda^{-1} z) \\
        \sinh(\lambda^{-1} z) & \cosh(\lambda^{-1} z)
     \end{pmatrix} \; . \]

Because \[ z \begin{pmatrix}
        0 & 1 \\
        1 & 0
     \end{pmatrix} \lambda^{-1} = \hat{F}+\hat{B} \; , \] where 
\[ \hat{F}=(z \lambda^{-1}-\bar{z} \lambda) \begin{pmatrix}
        0 & 1 \\
        1 & 0
     \end{pmatrix} \; , \;\;\; \hat{B}= \bar{z}\lambda \begin{pmatrix}
        0 & 1 \\
        1 & 0
     \end{pmatrix} \; , \] 
and because $\hat{F}\cdot\hat{B}=\hat{B}\cdot\hat{F}$ 
implies that $\exp(\hat{F}+\hat{B})=\exp \hat{F} \exp \hat{B}$, 
we have that $\phi$ can be split in the following way: 
\[ \phi = F \cdot B \; , \;\;\; F=\exp \left( (z \lambda^{-1}-\bar{z} 
\lambda) \begin{pmatrix}
        0 & 1 \\
        1 & 0
     \end{pmatrix} \right) \; , 
\; \;\; B= \exp \left( \bar{z} \lambda \begin{pmatrix}
        0 & 1 \\
        1 & 0
     \end{pmatrix} \right) \; . \] 
Note here that $\hat{F}$ is in the Lie algebra $\su_2$, and so $F$ is in the 
Lie group $\SU_2$, by the following lemma: 

\begin{lemma}
\label{simpleSU2lemma}
$A \in \su_2$ implies 
$\exp A \in \SU_2$.  \end{lemma} 

\begin{proof} 
Note that in general $(\det \mathcal{A})^\prime=\text{trace} 
(\mathcal{A}^\prime \cdot \mathcal{A}^{-1}) \det \mathcal{A}$, or 
equivalently, 
\begin{equation}
\left(\log (\det \mathcal{A})\right)^\prime=\text{trace} 
(\mathcal{A}^\prime \cdot \mathcal{A}^{-1}) 
\label{eq:matrixlemma}
\end{equation}
($^\prime$ means derivative with respect to some parameter that $\mathcal{A}$ 
depends on, see 
Lemma \ref{matrix-identity-lemma}).  Replacing $\mathcal{A}$ by $e^{t A}$ and 
taking the derivative with respect to $t$, we have that $\det e^{t A}$ is constant, 
because the trace of $A$ is zero.  Since $\det e^{t A}$ is $1$ when $t=0$, 
it follows that $\det e^{A}$ is also $1$.  
Then, using that $A \in \su_2$, one can check that $e^A$ is of the form 
\begin{equation}
\begin{pmatrix}
         p & q \\ -\bar{q} & \bar{p} \end{pmatrix} \; ,
\label{eq:SU2form}
\end{equation}
so $e^A$ is in $\SU_2$.  
\end{proof}

\begin{remark}
In fact, any matrix in $\SU_2$ is of the form \eqref{eq:SU2form} 
for some $p,q\in\mathbb{C}$ such that $p \bar p + q \bar q = 1$.  
\end{remark}

Note also that $\hat{B}$, and hence $B$ as well, 
has no negative powers of $\lambda$.  In fact, 
this $F$ and $B$ satisfy all the properties described in Section \ref{section0}.  
We can compute that 
\[ F = \begin{pmatrix}
        \cosh(\lambda^{-1} z-\bar{z} \lambda) & \sinh(\lambda^{-1} z
        -\bar{z} \lambda) \\
        \sinh(\lambda^{-1} z-\bar{z} \lambda) & \cosh(\lambda^{-1} z
        -\bar{z} \lambda)
     \end{pmatrix} \; . 
\] We now define the $2 \times 2$ matrix via the Sym-Bobenko formula 
\begin{equation}\label{firstSym} 
  f (z,\bar{z},\lambda) = \frac{-1}{2} F \begin{pmatrix}
        i & 0 \\
        0 & -i
  \end{pmatrix} F^{-1} - i \lambda (\partial_\lambda F) F^{-1} \; , \end{equation}
and let $f(z,\bar{z})=f(z,\bar{z},1)$.  We can then show that 
$f(z,\bar{z})$ equals 
\[ \frac{-i}{2}\begin{pmatrix}
        \cosh^2(z-\bar{z})+\sinh^2(z-\bar{z}) & 
        -2 \cosh(z-\bar{z}) \sinh(z-\bar{z})-2 (z+\bar{z}) \\
        2 \cosh(z-\bar{z}) \sinh(z-\bar{z})-2 (z+\bar{z}) & 
        -\cosh^2(z-\bar{z})-\sinh^2(z-\bar{z})
     \end{pmatrix} \; . \]  Noting that $f(z,\bar{z})$ is of the form 
\begin{equation}\label{prob3-a} \frac{-i}{2}\begin{pmatrix}
        -x_3 & x_1+i x_2 \\
        x_1-i x_2 & x_3
     \end{pmatrix} \; , \end{equation} 
     for real-valued functions $x_1,x_2,x_3$ of $z,\bar{z}$, 
     we see that (with $z=x+i y$) 
\[ (x_1,x_2,x_3)= (-4x,-\sin(4y),-\cos(4y)) \; . \] 
This is a parametrization of a cylinder of radius $1$ 
with axis the $x_1$-axis.  One can easily check that 
this parametrization is conformal and has CMC $H=1/2$.

\section{Spheres}
\label{section2}

Here we show how the DPW method makes a sphere.  Define 
\[ \xi = \lambda^{-1} \begin{pmatrix}
        0 & 1 \\
        0 & 0
     \end{pmatrix} dz \; , \] for $z \in \Sigma = \mathbb{C}$ and $\lambda \in \mathbb{S}^1$.  
Check that a solution to \[ d\phi = \phi \xi \; , \] 
is 
\[ \phi = \exp \left( z \begin{pmatrix}
        0 & 1 \\
        0 & 0
     \end{pmatrix} \lambda^{-1} \right) = \begin{pmatrix}
        1 & z\lambda^{-1} \\
        0 & 1
     \end{pmatrix} \; . \]
Note that $\phi$ can be split in the following way: $\phi = F \cdot B$, where
\[ F=\frac{1}{\sqrt{1+|z|^2}} 
\begin{pmatrix}
        1 & z \lambda^{-1} \\
        -\bar{z} \lambda & 1 \end{pmatrix} \; , 
\; \;\; B=     \begin{pmatrix}
        1 & 0 \\
        \bar{z} \lambda & 1
     \end{pmatrix}\begin{pmatrix}
        \tfrac{1}{\sqrt{1+|z|^2}} & 0 \\
        0 & \sqrt{1+|z|^2}
     \end{pmatrix} \; . \] 
Note here that $F$ is in the 
Lie group $\SU_2$ (because $\lambda \in \mathbb{S}^1$), and $B$ has no negative 
powers of $\lambda$.  In fact, 
this $F$ and $B$ satisfy all the properties described in Section \ref{section0}.  
Now define the $2 \times 2$ matrix $f(z,\bar{z},\lambda)$ as in \eqref{firstSym}, 
and let $f(z,\bar{z})=f(z,\bar{z},1)$.  We can then show that 
\[ f(z,\bar{z})= \frac{-i}{2(1+z\bar{z})}\begin{pmatrix}
        1-3z\bar{z} & -4 z \\
        -4 \bar{z} & -1+3z\bar{z}
     \end{pmatrix} \; . \]  Noting that $f(z,\bar{z})$ is of the form 
\eqref{prob3-a}, we have that 
\[ (x_1,x_2,x_3)=(1+x^2+y^2)^{-1} (-4x,-4y,-1+3x^2+3y^2) \; . \]
This is a parametrization of a sphere of radius $2$ centered at $(0,0,1)$.  
It is easily shown to be a conformal parametrization with CMC $H=1/2$.  

\section{Rigid motions} 
\label{section4}

Given some holomorphic potential $\xi$ and some solution $\phi$ of 
$d\phi=\phi \xi$, as in Section \ref{section0}, we could also have chosen 
any other solution 
\[ A(\lambda) \cdot \phi \; , \] 
where $A(\lambda) \in \SL_2\!\mathbb{C}$ is a twisted matrix independent 
of $z$ that is holomorphic in $\lambda$ for $\lambda \neq 0,\infty$.  
Let us assume that $A(\lambda) \in \SU_2$ for all $\lambda \in \mathbb{S}^1$, 
and then the new Iwasawa splitting must be 
\[ A(\lambda) \cdot \phi = \hat{F}\cdot B \; , 
   \;\;\;\; \hat F = A(\lambda) \cdot F 
\; . \]  (This follows from the uniqueness of Iwasawa splitting.)  
With this new $\hat F$, the Sym-Bobenko formula \eqref{firstSym} changes as 
\begin{equation}\label{eq:ftofhat}
 f(z,\bar z) \to \hat f(z,\bar z) = 
 A(1) \cdot f(z,\bar z) \cdot A^{-1}(1) 
-i (\partial_\lambda A(\lambda))|_{\lambda=1} 
A^{-1}(1) \; . \end{equation}  
We will now show that: 

\begin{lemma}\label{lm:rigid_motion}
The mapping \eqref{eq:ftofhat} is a rigid motion of $f(z,\bar z)$ in $\mathbb{R}^3$, 
and the rigid motion has an elementary and explicit description 
in terms of the entries of $A(\lambda)$.  
\end{lemma}

\begin{proof}
We first show that $f(z,\bar z) \to A(1) \cdot f(z,\bar z) \cdot A^{-1}(1)$ 
amounts to a rotation of $f(z,\bar z)$: If 
\[ A(1) = \begin{pmatrix} a & b \\ -\bar b & \bar a \end{pmatrix} \] 
with $a\bar a+b\bar b=1$, 
then 
\[
A(1) \cdot \frac{-i}{2}\begin{pmatrix}-x_3 & x_1+i x_2 \\x_1-i x_2 & x_3
     \end{pmatrix} \cdot A^{-1}(1)
\] 
changes points $\vec{x}=(x_1,x_2,x_3)^t \in \mathbb{R}^3$ to 
points $R \vec{x}$, where 
$R$ is the $3 \times 3$ matrix (with $a=a_1+ia_2$ and $b=b_1+ib_2$, 
$a_1,\,a_2,\,b_1,\,b_2\in\mathbb{R}$) 
\[
R = \begin{pmatrix}
       a_1^2-a_2^2-b_1^2+b_2^2 & -2a_1a_2-2b_1b_2 & 2a_1b_1-2a_2b_2 \\ 
      2a_1a_2-2b_1b_2 & a_1^2-a_2^2+b_1^2-b_2^2 & 2a_2b_1+2a_1b_2 \\
     -2a_1b_1-2a_2b_2 & 2a_2b_1-2a_1b_2 & a_1^2+a_2^2-b_1^2-b_2^2
    \end{pmatrix} \; . \] 
$R$ is a matrix in $\text{SO}_3$ and hence represents a 
rotation of $\mathbb{R}^3$ fixing 
the origin $\vec O \in \mathbb{R}^3$.  
And computing the eigenvalues and eigenvectors of $R$, 
we can confirm that the axis of this rotation contains the vector 
\[ (-b_2, b_1, a_2) 
\; . \] To determine the angle and orientation of rotation, 
let us assume $a_1,a_2,b$ are all nonzero (other cases are easier to 
examine).  
Note that $(b_1,b_2,0)$ is perpendicular to $(-b_2, b_1, a_2)$, so 
\[ \theta = \arccos \left( \frac{\langle (b_1,b_2,0),
   (R \cdot (b_1, b_2, 0)^t)^t\rangle}{|b|^2} \right) \] 
is the angle of rotation.  
One can compute that 
\[ \theta = \arccos (a_1^2 - a_2^2 - |b|^2) \; . \]
To determine the orientation of the rotation, 
one can see that the matrix with first row $(-b_2, b_1, a_2)$ and second row 
$(b_1,b_2,0)$ and third row $(R \cdot (b_1,b_2,0)^t)^t$ has 
determinant $2 a_1 |b|^2 (a_2^2 + |b|^2)$, 
which is positive if and only if $a_1>0$.  
So $R$ is a clockwise (resp. counterclockwise) rotation of angle $\theta$ 
about the axis (oriented in the direction $(-b_2,b_1,a_2)$) when 
$a_1<0$ (resp. $a_1>0$).

Secondly, we show that the term
$-(i \lambda \partial_\lambda A(\lambda))|_{\lambda=1} A^{-1}(1)$ in 
\eqref{eq:ftofhat} amounts to a translation of the original surface.  
To see that, we show that it is of the form 
\eqref{prob3-a}, for some real-valued functions $x_1,x_2,x_3$.  
$A:=A(\lambda)$ is holomorphic in $\lambda$, 
and $i\lambda\partial_\lambda = \partial_t$, where $t \in \mathbb{R}$ is defined by 
$\lambda =e^{it}$, so $[i\lambda\partial_\lambda A \cdot A^{-1}]_{\lambda =1}
 = [\partial_tA \cdot A^{-1}]_{t=0}$, and $\lambda = e^{it} \in \mathbb{S}^1$ implies 
$A|_{\mathbb{S}^1}$ is of the form
\[ A|_{\mathbb{S}^1}=\begin{pmatrix}
              p(t) &           q(t) \\
   -\overline{q(t)}& \overline{p(t)}
               \end{pmatrix} \]
with $|p(t)|^2+|q(t)|^2\equiv 1$, by the assumption that $A \in \SU_2$ 
for all $\lambda \in \mathbb{S}^1$.  
So, setting $p:=p(t)$ and $q:=q(t)$, 
\begin{eqnarray*}
\left[\partial_tA \cdot A^{-1}\right]_{t=0}
&=& \left[\begin{pmatrix}
              \partial_t p &           \partial_t q \\
   -\overline{\partial_t q}& \overline{\partial_t p}
          \end{pmatrix}
          \begin{pmatrix}
                    \bar p &-q \\
                    \bar q & p
          \end{pmatrix}\right]_{t=0} \\
&=& \left[\begin{pmatrix}
              \bar p \cdot \partial_t p + \bar q \cdot \partial_t q &
             -q \cdot \partial_t p + p \cdot \partial_t q \\
 -\bar p \cdot \overline{\partial_t q} +\bar q \cdot \overline{\partial_t p} &
    p \cdot \overline{\partial_t p} + q \cdot \overline{\partial_t q}
          \end{pmatrix} \right]_{t=0} \; , 
\end{eqnarray*}
which implies that the off-diagonal terms are negatives of the conjugates 
of each other, and the diagonal terms are conjugate to each other.  
Furthermore, $p\bar p+q\bar q\equiv 1$ implies that 
$[\partial_t A \cdot A^{-1}]_{t=0}$ is traceless, hence 
$[i\lambda\partial_\lambda A \cdot A^{-1}]_{\lambda =1}$ 
is of the desired form.  
\end{proof}

\begin{remark}
If we had not assumed that $A(\lambda) \in \SU_2$ for all $\lambda \in \mathbb{S}^1$, 
but only that $A(\lambda) \in \SL_2\!\mathbb{C}$, then the resulting transformation 
of the surface $f$ would not be merely a rigid motion, in general.  This 
more general transformation is what we will later be calling a {\em dressing} 
in Section \ref{section14}.  
\end{remark}

\section{Delaunay surfaces}
\label{section5}

We define \[ \xi=D \frac{dz}{z} \; , \;\;\; \text{where} \;\;\; 
D = \begin{pmatrix}
        r & s \lambda^{-1}+\bar t \lambda \\
        \bar s \lambda+t \lambda^{-1} & -r
     \end{pmatrix} \; , \]
with $s,t \in \mathbb{C}$ and $r\in \mathbb{R}$ and $st \in \mathbb{R}$ and 
$r^2+|s+\bar t|^2=1/4$, for $z$ in $\Sigma =\mathbb{C}\setminus\{0\}$.  
One can check that one solution to 
\[ d\phi = \phi \xi \] is 
\[ \phi = \exp \left( \log z  \cdot D  \right) \; . \]
We define \[ \hat F=\exp \left( i \theta D \right) \; , 
\; \;\; \hat B= \exp \left( \log \rho \cdot 
D \right) \; , \] where $z=\rho e^{i \theta}$.  
We have that $i \theta D$ and $\log \rho \cdot D$ commute, 
so, as seen in Section \ref{section1}, we know that 
$\exp(i \theta D + \log \rho \cdot D) = \hat F \cdot \hat B$, 
and so $\phi$ can be split (this is not yet 
Iwasawa splitting) in the following way: 
\[ \phi = \hat F \cdot \hat B \; . \] 
We also have that $\hat F \in \SU_2$ 
for all $\lambda \in \mathbb{S}^1$.  (This follows from showing 
that $i \theta D \in \su_2$ for all $\lambda \in \mathbb{S}^1$, and then 
applying Lemma \ref{simpleSU2lemma}.)  

Since 
\begin{equation}
D^2 = \kappa^2 \cdot \id \; , \;\;\; \text{where} \;\;\; 
\kappa = \sqrt{1/4 +s t (\lambda-\lambda^{-1})^2} \; , 
\label{eq:D^2=kappa^2id}
\end{equation}
we see that 
\[ \hat F = \begin{pmatrix} \cos (\kappa \theta) + i r \kappa^{-1} 
\sin (\kappa \theta) & i \kappa^{-1} \sin 
(\kappa \theta) (s \lambda^{-1} + \bar t \lambda) \\
i \kappa^{-1} \sin 
(\kappa \theta) (\bar s \lambda + t \lambda^{-1}) & 
\cos (\kappa \theta) - i r \kappa^{-1} \sin (\kappa \theta)
\end{pmatrix} \; , \]\[ 
\hat B = \begin{pmatrix} \cosh (\log \rho \cdot \kappa) + r \kappa^{-1} 
\sinh (\log \rho \cdot \kappa) & \kappa^{-1} \sinh 
(\log \rho \cdot \kappa) (s \lambda^{-1} + \bar t \lambda) \\
\kappa^{-1} \sinh 
(\log \rho \cdot \kappa) (\bar s \lambda + t \lambda^{-1}) & 
\cosh (\log \rho \cdot \kappa) - r \kappa^{-1} \sinh (\log \rho \cdot 
\kappa)
\end{pmatrix} \; . \] 
Note that $\hat F$ and $\hat B$ are both even in $\lambda$ along the diagonal and 
odd in $\lambda$ on the off-diagonal.  We can now do Iwasawa splitting on 
$\hat B$, that is, $\hat B=\tilde F \cdot B$, where $\tilde F \in \SU_2$ for all 
$\lambda \in \mathbb{S}^1$ and 
the power series expansion with respect to $\lambda$ of $B$ has no terms with 
negative powers of $\lambda$.  We define $F = \hat F \cdot \tilde F$.  Then 
$\phi = F \cdot B$ is the Iwasawa splitting of $\phi$.  

For each fixed $\lambda$, $\hat B$ depends only on $|z|=\rho$, and so also 
$\tilde F$ and $B$ depend only on $\rho$.  
Also, $\hat F$ depends only on $\theta$.  
So, under the rotation of the domain 
\[ z \to e^{i \theta_0} \cdot z \; , \]  
the following transformations occur: 
\[ F \to \exp (i \theta_0 D) \cdot F \; , \] 
\[ B \to B \; , \]
and by the map \eqref{eq:ftofhat}, we have
\begin{equation}\label{delaunay}
 f(z,\bar z) \to \left[\exp(i \theta_0 D) f(z,\bar z) \exp(-i \theta_0 D) - i 
(\partial_\lambda (\exp(i \theta_0 D) )) \exp(-i \theta_0 D) 
\right]_{\lambda=1} \; . \end{equation}
Using $D^2 = \kappa^2 \id$ and $\partial_\lambda \kappa|_{\lambda =1}=0$, 
we see that 
\[ \partial_\lambda (\exp(2\pi i D) )|_{\lambda=1} 
=\partial_\lambda \left( 
\cos(2\pi \kappa) \cdot \id +i \kappa^{-1} \sin(2\pi \kappa) \cdot D 
\right) |_{\lambda=1} = 0 \; , \]
and $\exp(2\pi i D)|_{\lambda=1}=-\id$, 
so $f \to f$ is unchanged in \eqref{delaunay} when $\theta_0=2\pi$, 
and hence the surface $f$ is well defined on $\mathbb{C}^*=\mathbb{C}\setminus\{0\}$.  

The computations in Section \ref{section4}, 
with $A(\lambda):=\exp (i\theta_0D)$, imply that 
\[ f(z,\bar z) \to \left[ \exp(i \theta_0 D) f(z,\bar z) \exp(-i \theta_0 
D) \right]_{\lambda=1} \] represents a rotation of the surface about the 
line through the origin in the direction 
\[ \hat L=(-\Re (s+t),-\Im (s-t),r) \] of angle $\theta_0$, 
and the orientation of the rotation can be determined as described in 
Section \ref{section4}.  

The term $- i(\partial_\lambda  (\exp(i \theta_0 D) )) \exp(-i \theta_0 D)$ 
added into the map \eqref{delaunay} represents a translation by the vector 
\begin{eqnarray*}
\hat T
 &=& 2 \cdot \sin \theta_0 \cdot
     \big( \Im (s+t),-\Re (s-t),0 \big) \\
 & & \hspace{70pt}
    +4 \cdot (1-\cos \theta_0) \cdot
     \big( r \cdot \Re (s-t),r \cdot \Im (s+t),|s|^2-|t|^2 \big)
 \; . \end{eqnarray*}
Noting that $\langle \hat L,\hat T\rangle =0$ for all $\theta_0\in\mathbb{R}$, 
the map \eqref{delaunay} must represent a rotation of angle $\theta_0$ 
about a line parallel to $\hat L$.  
Because the map \eqref{delaunay} has a fixed point 
\[
\hat P= (\hat x_1, \hat x_2, \hat x_3)
\]
satisfying
\[
\Im \left[\left(\hat x_1-i\hat x_2\right)\left(s+\bar t\right)\right]=0 \; , 
\;\;\; \hat x_3(s+\bar t)+r(\hat x_1+i\hat x_2) = s-\bar t \; ,
\]
the point $\hat P$ can be chosen independent of $\theta_0$, and so
the axis of rotation is also independent of $\theta_0$.  
(In fact, when $i(s+t)\notin\mathbb{R}$, we may simply choose 
$\hat P=(0,0,\frac{\Re (s-t)}{\Re (s+t)})$.)
It follows that the map \eqref{delaunay} represents a rotation 
of angle $\theta_0$ about the line 
\[ L = \{x \cdot \hat L + \hat P \, | \, x \in \mathbb{R} \} \; . \] 

The conclusion is that the CMC surface we have just made is a surface 
of revolution, so it must be a Delaunay surface 
(see \cite{Eel, HY, KorKS}).  Which Delaunay surface 
one gets depends on the choice of $r,s,t$.  An unduloid is produced when 
$s t > 0$, 
a nodoid is produced when $s t < 0$, and for the limiting singular case of 
a chain of spheres, $s t = 0$.  
A cylinder is produced when $|s| = |t| =1/4$ and $r = 0$.  

We will give more details about Delaunay surfaces in $\mathbb{R}^3$ 
(including a description of the weight, or flux, of a Delaunay surface), 
and in $\mathbb{H}^3$ and $\mathbb{S}^3$ as well, in Section 
\ref{spaceformdelaunay}.  
Delaunay surfaces in $\mathbb{R}^3$ via the DPW method are also described in 
detail in \cite{k}, \cite{DorKS} and \cite{k2}.  
Also (after describing the associated family of 
a CMC surface) it is shown in the supplement 
\cite{waynebook-maybe} to these notes the 
existence of certain surfaces that we call ``twizzlers''.  Twizzlers 
are special surfaces in the associated families of Delaunay surfaces.  

\begin{figure}[phbt]
\begin{center}
\includegraphics[width=0.9\linewidth]{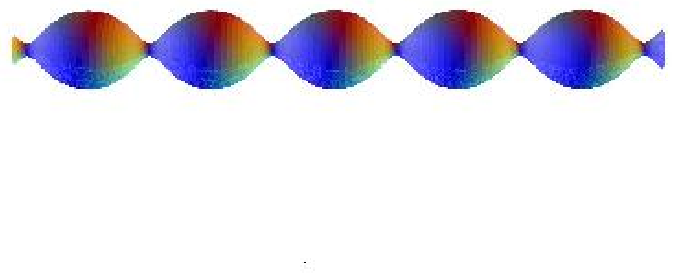}
\end{center}
\vspace{0.5in}
\begin{center}
\includegraphics[width=0.9\linewidth]{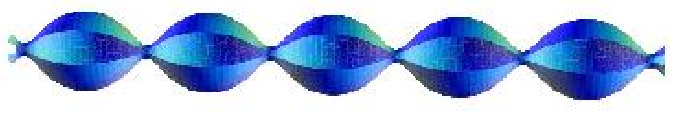}
\end{center}
\vspace{0.8in}
\begin{center}
\includegraphics[width=0.25\linewidth]{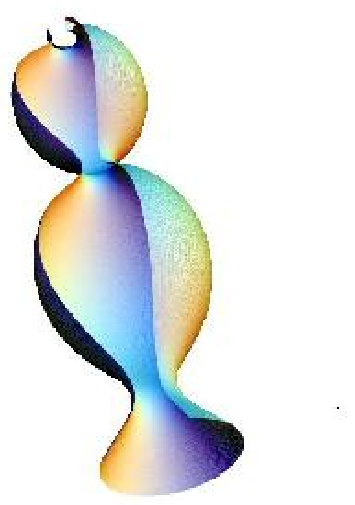}
\end{center}
\vspace{3in}
\caption{Various portions of an unduloid.  (These graphics were made by 
Koichi Shimose using Mathematica.)}
\end{figure}

\begin{figure}[phbt]
\begin{center}
\includegraphics[width=0.6\linewidth]{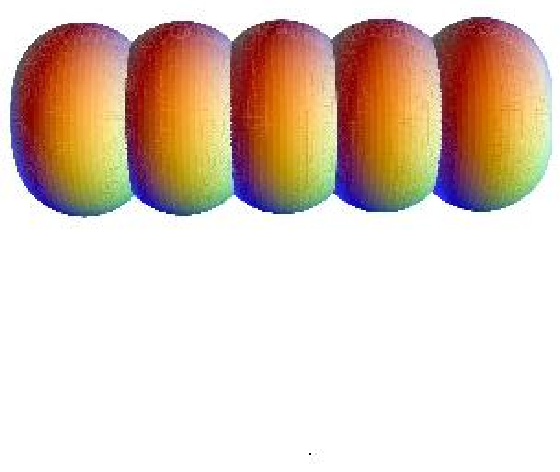}
\end{center}
\vspace{1.2in}
\begin{center}
\includegraphics[width=0.6\linewidth]{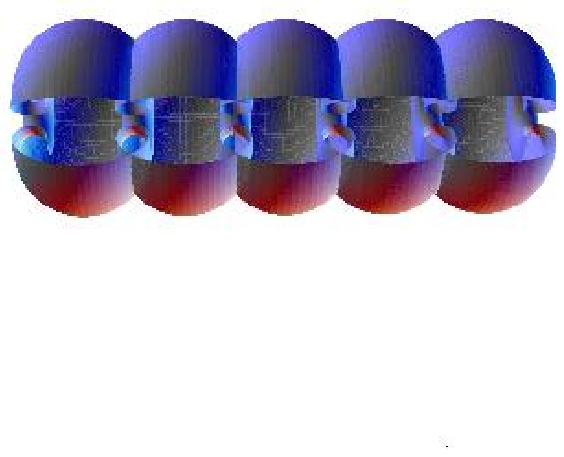}
\end{center}
\vspace{1.6in}
\begin{center}
\includegraphics[width=0.45\linewidth]{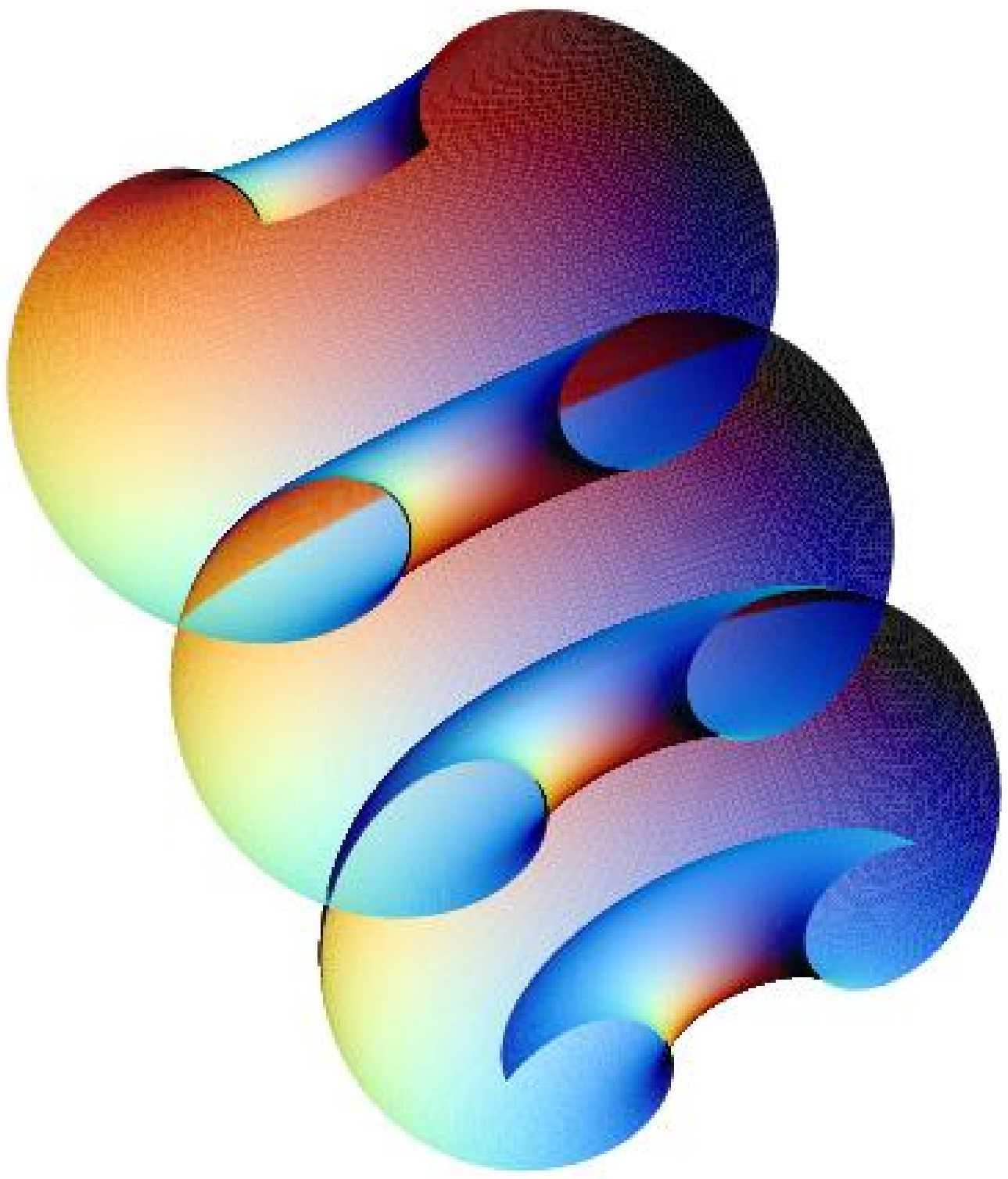}
\end{center}
\vspace{0.2in}
\caption{Various portions of a nodoid.  (These graphics were made by 
Koichi Shimose using Mathematica.)}
\end{figure}

\begin{figure}[phbt]
\begin{center}
\includegraphics[width=0.4\linewidth]{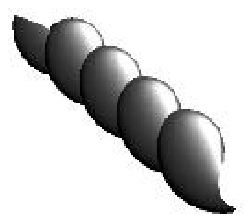}
\includegraphics[width=0.4\linewidth]{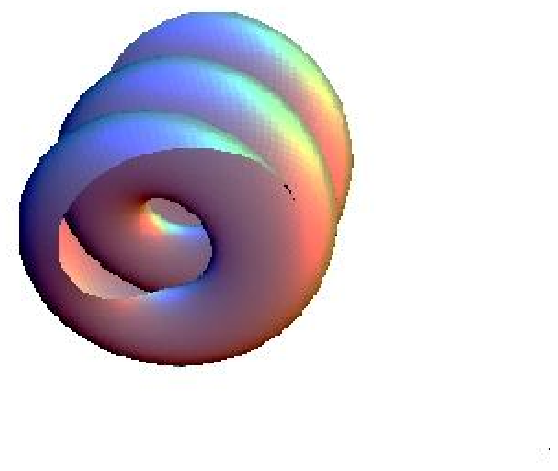}
\end{center}
\vspace{2.6in}
\caption{Two different twizzlers in $\mathbb{R}^3$.  (These graphics were made by 
Yuji Morikawa using Mathematica.)}
\end{figure}

\section{Smyth surfaces}
\label{section3}

Here we show how the DPW method makes Smyth surfaces.  Define 
\[ \xi = \lambda^{-1} \begin{pmatrix}
        0 & 1 \\
        c & 0
     \end{pmatrix} dz \; , \] 
for the complex variable $z \in \Sigma = \mathbb{C}$ 
and $\lambda \in \mathbb{S}^1$ and some $c \in \mathbb{C}$.  
One can check that a solution to 
\begin{equation}\label{df=fxi}
d\phi = \phi \xi \end{equation} is 
\[ \phi = \exp \left( \lambda^{-1} z \begin{pmatrix}
        0 & 1 \\
        c & 0
     \end{pmatrix} \right) = \begin{pmatrix}
        \cosh (\lambda^{-1} z \sqrt{c}) & \sqrt{c}^{-1} 
                         \sinh (\lambda^{-1} z \sqrt{c}) \\
  \sqrt{c} \sinh (\lambda^{-1} z \sqrt{c}) & \cosh (\lambda^{-1} z \sqrt{c})
     \end{pmatrix} \; . \] 

When $c \in \mathbb{S}^1$, we can do Iwasawa splitting of $\phi$ into $F \cdot B$ explicitly, 
in exactly the same way as we did in Section \ref{section1}, that is 
\begin{eqnarray*}
F &=& \exp \left( (\lambda^{-1} z-\lambda \bar{z} c^{-1}) \begin{pmatrix}
        0 & 1 \\
        c & 0
     \end{pmatrix} \right) \\
  &=& \begin{pmatrix}
        \cosh (\lambda^{-1} z \sqrt{c}-\lambda \bar{z} \sqrt{c}^{-1}) & \sqrt{c}^{-1} 
        \sinh (\lambda^{-1} z \sqrt{c}-\lambda \bar{z} \sqrt{c}^{-1}) \\
  \sqrt{c} \sinh (\lambda^{-1} z \sqrt{c}-\lambda \bar{z} \sqrt{c}^{-1}) & 
  \cosh (\lambda^{-1} z \sqrt{c}-\lambda \bar{z} \sqrt{c}^{-1})
     \end{pmatrix} \; , \\  
B &=& \exp \left( \lambda \bar{z} c^{-1} \begin{pmatrix}
        0 & 1 \\
        c & 0
     \end{pmatrix} \right) \; .  
\end{eqnarray*}
One can check that the resulting surface is again a cylinder of radius $1$, 
following exactly the same procedure as in Section \ref{section1}.  The axis 
of the cylinder depends on $c$.  

When $c \not\in \mathbb{S}^1 \cup \{ 0 \}$, 
there is no simple way to do Iwasawa splitting of $\phi$.  
With $\phi |_{z=0}=\id$, the resulting surface is a 2-legged Smyth 
surface, also called a 2-legged Mister Bubble.  More generally, 
\[ \xi = \lambda^{-1} \begin{pmatrix}
        0 & 1 \\
        c z^k & 0
     \end{pmatrix} dz \] 
gives a solution $\phi$ with initial condition $\phi |_{z=0}=\id$, 
and the resulting surface is a ($k+2$)-legged Mister Bubble.  
The reason for using ``$(k+2)$-legged'' in the naming of these surfaces is that
pictures of compact portions of the surfaces appear to have $k+2$ ``legs'', 
although in fact the surfaces always have only a single end.  

When $k \geq 1$, there is again no simple way to 
do Iwasawa splitting $\phi = F \cdot B$ even when $c\in \mathbb{S}^1$.  

One can check that changing 
$c$ to $c e^{i \theta_0}$ for any $\theta_0 \in \mathbb{R}$ only changes the surface 
by a rigid motion and a reparametrization $z\to ze^{-\frac{i\theta_0}{k+2}}$.  
To prove this, note that with these changes, 
\[
\xi \to A\xi A^{-1} \; , \;\;\; \text{where} \;\;\;
A=\begin{pmatrix} e^{-\frac{ i \theta_0}{2(k+2)}} & 0 \\ 0 & 
e^{\frac{ i \theta_0}{2(k+2)}} \end{pmatrix} \in \SU_2 \] and so $F \to AFA^{-1}$, 
by an argument similar to that in the proof of Lemma \ref{lm:smyth}.  
So $f \to AfA^{-1}$.  
Then Lemma \ref{lm:rigid_motion} implies this is a rigid motion.  
So without loss of generality we may assume and do assume that $c \in \mathbb{R}$.  

It is shown in \cite{bi} that $u=u(z,\bar z):\mathbb{C} \to \mathbb{R}$ 
in the metric (first fundamental form) $4e^{2u}dzd\bar z$ 
of the surface resulting from the frame $F$ 
is constant on each circle of radius $r$ centered at 
the origin in the $z$-plane, that is, 
$u=u(r)$ is independent of $\theta$ in $z=r e^{i \theta}$.  
Having this internal rotational symmetry of the metric (without actually 
having a surface of revolution) is the definition of a Smyth surface, 
hence this potential $\xi$ (with initial condition $\phi|_{z=0}=\id$) produces 
CMC conformal immersions $f$:$\mathbb{C} \to \mathbb{R}^3$ that are Smyth surfaces. 
Furthermore, Smyth surfaces are related to the Painleve III equation, 
see \cite{bi}.  We will also prove those facts in the supplement \cite{waynebook-maybe} 
to these notes.  Furthermore, although we will not prove it in these notes, 
Timmreck, Pinkall and Ferus \cite{FerPT} showed that the Smyth surfaces are 
proper immersions.  

\begin{lemma}\label{lm:smyth}
The Smyth surfaces have reflective symmetry with respect to $k+2$ geodesic 
planes that meet equiangularly along a geodesic line.  
\end{lemma}

\begin{proof}
Consider the reflections 
\[ R_\ell(z)=e^{\frac{2 \pi i \ell}{k+2}} \bar z \] of the domain $\mathbb{C}$, for 
$\ell \in \{0,1,...,k+1\}$.  Note that 
\[ \xi(R_\ell(z),\lambda)= A_\ell \xi(\bar z,\lambda) A_\ell^{-1} \; , \] 
where 
\[ A_\ell = \begin{pmatrix} e^{\frac{\pi i \ell}{k+2}} & 0 \\ 0 & 
e^{\frac{-\pi i \ell}{k+2}} \end{pmatrix} \; . \] 
Since 
\[ d(\phi(R_\ell(z),\lambda) A_\ell) = \phi(R_\ell(z),\lambda) A_\ell 
\cdot\xi(\bar z,\lambda) \; , \] 
and since any two solutions of this equation differ 
by a factor that is constant in $z$, we have 
\[ \phi(R_\ell(z),\lambda) \cdot A_\ell 
 = \mathcal{A} \cdot \phi(\bar z,\lambda) \]  
for some $\mathcal A$ depending only on $\lambda$.  However the initial 
condition $\phi(z,\lambda)|_{z=0}=\id$ for the solution $\phi$ of 
\eqref{df=fxi} implies that $\mathcal A = A_\ell$, so 
\[ \phi(R_\ell(z),\lambda) = A_\ell \phi(\bar z,\lambda) A_\ell^{-1} \] 
and so the Iwasawa splitting $\phi=F \cdot B$, 
with $F=F(z,\bar z,\lambda )$ and $B=B(z,\bar z,\lambda )$ 
satisfying all the conditions in Section \ref{section0}, satisfies 
\[ F(R_\ell(z),\overline{R_\ell(z)},\lambda) \cdot 
B(R_\ell(z),\overline{R_\ell(z)},\lambda) = 
(A_\ell F(\bar z,z,\lambda) A_\ell^{-1}) 
\cdot (A_\ell B(\bar z,z,\lambda) A_\ell^{-1}) \; . \] 
Since $A_\ell F(\bar z,z,\lambda) A_\ell^{-1}$ and 
$A_\ell B(\bar z,z,\lambda) A_\ell^{-1}$ also satisfy all the properties
of the $F$ and $B$ in Section \ref{section0}, respectively, 
uniqueness of the Iwasawa decomposition implies 
\[ F(R_\ell(z),\overline{R_\ell(z)},\lambda) 
 = A_\ell F(\bar z,z,\lambda) A_\ell^{-1} \; . \] 
Note that $c \in \mathbb{R}$ implies that 
$\xi (\bar z,\lambda) = \overline{\xi (z,\bar \lambda)}$ and that 
$\phi |_{z=0}=\id$ then implies that 
$\phi (\bar z,\lambda) = \overline{\phi (z,\bar \lambda)}$, therefore 
\[ F(\bar z,z,\lambda) = \overline{F(z,\bar z,\bar \lambda)} \; . \] 
Then, we have	
\[ F(R_\ell(z),\overline{R_\ell(z)},\lambda) 
 = A_\ell \overline{F(z,\bar z,\bar \lambda)} A_\ell^{-1} \; . \] 
Inserting this into \eqref{firstSym}, we have
\[ f(R_\ell(z),\overline{R_\ell(z)},\lambda) 
 = -A_\ell \overline{f(z,\bar z,\bar \lambda)} A_\ell^{-1} \; . \] 
Recalling that we defined $f(z,\bar z)=f(z,\bar z,1)$, 
the transformation $f(z,\bar z) \to -\overline{f(z,\bar z)}$ 
represents reflection across the plane $\{x_2=0\}$ of $\mathbb{R}^3$, 
and conjugation by $A_\ell$ represents a rotation by angle $2\pi i\ell /(k+2)$ 
about the $x_3$-axis of $\mathbb{R}^3$ (see Lemma \ref{lm:rigid_motion}), 
hence the proof is completed.
\end{proof}

\begin{figure}[phbt]
\begin{center}
\includegraphics[width=0.4\linewidth]{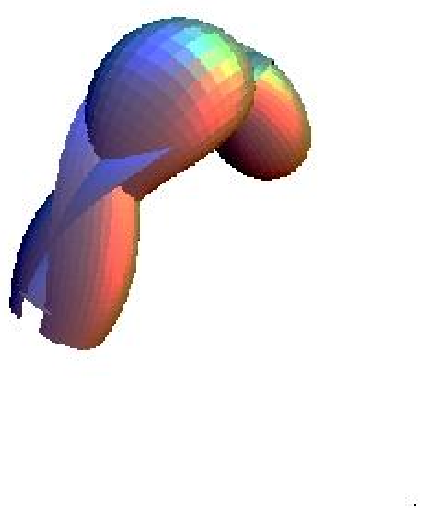}
\includegraphics[width=0.4\linewidth]{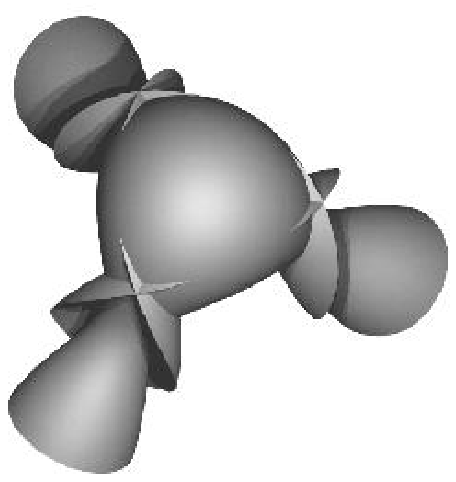}
\end{center}
\vspace{2.6in}
\caption{A 2-legged Smyth surface and a 3-legged Smyth surface. (The graphics were made 
by Yuji Morikawa using Mathematica and CMCLab \cite{Sch:cmclab}.)}
\end{figure}

\chapter{The Lax pair for CMC surfaces}
\label{chapter2}

In this chapter we describe CMC surfaces in terms of a pair of linear 
first-order matrix partial differential equations called a Lax pair.  
The Lax pair formulation is fundamental to the DPW method and will be used to 
prove results about the DPW method in Chapter \ref{chapter3}.  

At the end of this chapter we give some other applications of the Lax pair 
formulation: we prove the Weierstrass representation for minimal surfaces; and 
we prove the harmonicity of the Gauss map of CMC surfaces.  

\section{The $3 \times 3$ Lax pair for CMC surfaces}
\label{section6}

Before giving the proof of the DPW recipe in the next chapter, we show 
in this chapter that any CMC $H$ surface gives a solution to 
a certain kind of Lax pair, and that any solution of this certain 
kind of Lax pair determines a unique CMC $H$ surface.  (Then proving that the 
DPW recipe gives all CMC $H$ surfaces just amounts to showing that it gives 
all possible solutions for this certain kind of Lax pair.)  
Here we give arguments similar to those in \cite{b1, DW, ps}.  

We first note that any $C^2$ immersion of a simply-connected domain 
$\Sigma\subseteq\mathbb{C}$ into $\mathbb{R}^3$ can be made conformal by a 
reparametrization. (In fact, more generally any $2$-dimensional Riemannian 
manifold can be given a covering of coordinate charts with respect to which 
the Riemannian metric is conformal, see Theorem \ref{conformalityispossible}.) 
So, without loss of generality, we will always assume that CMC immersions are 
conformal.  

For the moment, we will not assume $H$ is constant.  
Let $\Sigma$ be a simply-connected domain in $\mathbb{C}$ and 
\begin{equation}
 f=f(z,\bar z): \Sigma \to \mathbb{R}^3
\label{eq:f:SigmaToR3}
\end{equation}
be a conformal immersion (with at least $C^2$ differentiability), and so, 
with $z=x+iy$ and $\partial_z=(\partial_x-i\partial_y)/2$ and 
$\partial_{\bar z}=(\partial_x+i\partial_y)/2$, 
\begin{equation}
\langle f_z , f_z \rangle=\langle f_{\bar{z}}, f_{\bar{z}}\rangle
=0 \; , \;\; \langle f_z, f_{\bar{z}}\rangle =2e^{2 u} \; , 
\label{eq:metric}
\end{equation}
where $u:\Sigma \to \mathbb{R}$ is defined this way.  
$\langle \cdot , \cdot \rangle$ denotes the complex 
bilinear extension of the usual $\mathbb{R}^3$ inner product (which 
is no longer an actual inner product).  
The unit normal vector $N$ to $f$ is defined by 
$N=(f_x \times f_y)/|f_x \times f_y|$.  
Then the {\em shape operator} $S$ of the surface $f$ 
(see Section \ref{surface-theory-section}) is defined to be the 
inverse of the first fundamental form $g=(g_{ij})_{i,j=1,2}$ times 
the second fundamental form, that is, 
\[
S=\begin{pmatrix} 4e^{2u} & 0 \\ 0 & 4e^{2u} \end{pmatrix}^{-1}
  \begin{pmatrix} \langle f_{xx},N\rangle &
                  \langle f_{xy},N\rangle \\
                  \langle f_{yx},N\rangle &
                  \langle f_{yy},N\rangle \end{pmatrix} \; . \]
the mean curvature $H$ is defined as 
\[
H=\frac{1}{2}\text{trace}(S) \; ,
\]
we have
\[ H=\frac{1}{8e^{2u}} \langle f_{xx}+f_{yy}, N \rangle
=\frac{1}{2} \langle \triangle f, N \rangle \; . 
\]  
(Recall that $\triangle = \sqrt{\det g}^{-1} \sum_{i,j} 
\partial_{u_i} (\sqrt{\det g} g^{ij} \partial_{u_j})$ for $(x,y)=(u_1,u_2)$, 
where $g^{-1}=(g_{ij})^{-1}=(g^{ij})$.)  
So $\langle f_{z\bar{z}},N\rangle =2He^{2u}$.  Thus we have 
\begin{eqnarray*}
& \langle f_z, f_{zz}\rangle
  = \langle f_z, f_{z\bar{z}}\rangle = 0 \; , \;\;\; 
  \langle f_{\bar{z}}, f_{zz}\rangle = 4 u_z e^{2u} \; , & \\ 
& \langle N_z, N\rangle = 0 \; , \;\;\; 
   \langle N_z, f_{\bar{z}}\rangle = -\langle N, f_{z\bar{z}}\rangle = 
- 2He^{2u} \; , \;\;\;
  \langle N_z, f_{z}\rangle = -\langle N, f_{zz}\rangle = - Q \; ,
\end{eqnarray*}
where this is the definition of the Hopf differential $Q$ as in 
Section \ref{surface-theory-section}.  
This implies the following Gauss-Weingarten equations:
\begin{eqnarray*}
& f_{zz} = 2 u_z f_z + Q N \; , \quad f_{z\bar{z}} = 2He^{2u} N \; , 
\quad f_{\bar{z}\bar{z}} = 2 u_{\bar{z}} f_{\bar{z}} + \bar{Q} N , & \\ 
& N_z = \frac{1}{2} (-2Hf_z-Q e^{-2u} f_{\bar{z}}) \; , \quad 
N_{\bar{z}} = \frac{1}{2} (-2Hf_{\bar{z}}-\bar{Q} e^{-2u} f_{z}) \; . &
\end{eqnarray*}
The Gaussian curvature $K$ was defined as $K=\det (S)$ 
in Section \ref{surface-theory-section}, and we can compute that 
\begin{equation}\label{eq:K=detS}
K=H^2-\frac{1}{4}Q\bar Qe^{-4u} \; .
\end{equation}
We define $e_1=f_x/|f_x|=e^{-u}(f_z+f_{\bar z})/2$ and 
$e_2=f_y/|f_y|=ie^{-u}(f_z-f_{\bar z})/2$ and then 
\[\mathcal{F}=(e_1,e_2,N)\] is an orthonormal frame of the surface, that is, 
$\mathcal{F}:\Sigma\to\text{SO}_3$, and 
\[ \mathcal{F}_z = \mathcal{F} (\Theta +\Upsilon_z) \; , \;\;\;
\mathcal{F}_{\bar{z}} = \mathcal{F} (\bar{\Theta}-\Upsilon_{\bar{z}}) \; , \]  
where 
\begin{eqnarray*}
 \Upsilon &=& \begin{pmatrix} 0 & i u & 0 \\
-i u & 0 & 0 \\
0 & 0 & 0 \end{pmatrix} \; , \\ 
\Theta &=& \frac{1}{2} \begin{pmatrix} 0 & 0 & -Q e^{-u}-2He^u \\
0 & 0 & -i Q e^{-u} + 2H i e^u \\
Q e^{-u}+2He^u & i Q e^{-u}-2Hi e^u & 0 \end{pmatrix} \; . 
\end{eqnarray*}
Note that the compatibility condition for the existence of a solution 
$\mathcal{F}$ is that $\mathcal{F}_{z\bar{z}}=\mathcal{F}_{\bar{z}z}$, 
in other words, 
\begin{equation}
   (\Theta +\Upsilon_z)_{\bar{z}}
  -(\bar{\Theta}-\Upsilon_{\bar{z}})_z
  - [ \Theta +\Upsilon_z ,\bar{\Theta}-\Upsilon_{\bar{z}} ] = 0 \; , 
\label{eq:3by3MCeqn}
\end{equation}
which is also called the {\em Maurer-Cartan equation}.  This implies 
\begin{equation}\label{firstMC} 
4u_{z\bar{z}}-Q\bar{Q}e^{-2u}+4H^2e^{2u}=0 \; , \;\;\; 
Q_{\bar{z}}=2 H_z e^{2 u} \; , 
\end{equation} 
which are the Gauss equation and Codazzi equation for $f$, respectively.  

\begin{remark}\label{re:GCeqns}
(1) The Gauss equation implies that 
\begin{equation}\label{eq:egregium}
K=-e^{-2u}u_{z\bar z} \; ,
\end{equation}
since $K$ satisfies \eqref{eq:K=detS}. Hence $K$ depends only on $u$, that is, 
$K$ is determined by only the first fundamental form, 
without requiring any knowledge of the second fundamental form, 
even though the definition of $K$ uses the second fundamental form.  
This fact was first found by Gauss \cite{G1828, G}.  

Thus, for a Riemannian manifold $\Sigma$ with the Riemannian metric 
$g=4e^{2u}dzd\bar z$, we could have chosen to define the Gaussian curvature of 
$(\Sigma ,g)$ by \eqref{eq:egregium}.  
Hereafter we will use \eqref{eq:egregium} as our definition of $K$.  
(We can use either \eqref{eq:K=detS} or \eqref{eq:egregium} when the 
surface is immersed in $\mathbb{R}^3$, but when the ambient space is not 
$\mathbb{R}^3$, then \eqref{eq:egregium} and \eqref{eq:K=detS} are not 
equivalent, and only \eqref{eq:egregium} can be the correct definition of the 
intrinsic curvature of the surface, as we will see in Chapter \ref{chapter4}.  
``Intrinsic'' means that the object in question is determined by only the 
first fundamental form.)  

(2) The Codazzi equation implies that $H$ is constant if and only if $Q$ is 
holomorphic, since $H$ is a real valued function.  
\end{remark}

Now we assume $H$ is constant.  Then by part (2) of Remark \ref{re:GCeqns}, 
$Q$ is holomorphic.  
Note that the Gauss and Codazzi equations still hold if $Q$ is replaced with 
$\lambda^{-2} Q$ for any $\lambda \in \mathbb{S}^1$, so the fundamental theorem of 
surface theory (Bonnet theorem) tells us we have a $1$-parameter 
family of surfaces $f(z,\bar z,\lambda)$ for $\lambda\in \mathbb{S}^1$, 
all with CMC $H$ 
and the same metric, and with Hopf differential $\lambda^{-2}Q$ (see, 
for example, \cite{doCarmo1, doCarmo2}).  
We call $f(z,\bar z,\lambda)$ for $\lambda\in \mathbb{S}^1$ the {\em associated family} 
of $f(z,\bar z)=f(z,\bar z,1)$.  

We now replace $Q$ with $\lambda^{-2} Q$, that is, we will add a 
``spectral'' parameter $\lambda$.  
The Gauss-Weingarten equations then become (now $f=f(z,\bar z,\lambda)$ depends on 
$\lambda$) 
\begin{eqnarray*}
& f_{zz} = 2 u_z f_z + \lambda^{-2} Q N \; , \quad 
  f_{z\bar{z}} = 2He^{2u} N \; , \quad 
  f_{\bar{z}\bar{z}} 
 = 2 u_{\bar{z}} f_{\bar{z}} + \lambda^{2} \bar{Q} N , & \\ 
& N_z
 = \frac{1}{2} (-2Hf_z-\lambda^{-2}Q e^{-2u} f_{\bar{z}}) \; , \quad 
  N_{\bar{z}}
 = \frac{1}{2} (-2Hf_{\bar{z}}-\lambda^{2}\bar{Q} e^{-2u} f_{z}) \; , &
\end{eqnarray*}
with
\[ e_1=\frac{f_x}{|f_x|} \; , \quad
   e_2=\frac{f_y}{|f_y|} \; , \quad \text{and} \quad
    N =e_1 \times e_2 \; . \] Then
\[ \mathcal{F}=(e_1,e_2,N) \; , \quad
   \mathcal{F}_z
  = \mathcal{F} (\Theta + \Upsilon_z) \; , \quad
   \mathcal{F}_{\bar{z}}
  = \mathcal{F} (\bar{\Theta} - \Upsilon_{\bar{z}}) \; , \]
with 
\[ \Upsilon = \begin{pmatrix} 0 & i u & 0 \\
  -i u & 0 & 0 \\ 0 & 0 & 0 \end{pmatrix} \]
as before, and with $\Theta$ now being
\[ \Theta = \frac{1}{2} 
                 \begin{pmatrix} 0 & 0 & -\lambda^{-2}Q e^{-u}-2He^u \\
0 & 0 & -i \lambda^{-2}Q e^{-u} + 2H i e^u \\
\lambda^{-2}Q e^{-u}+2He^u & i \lambda^{-2}Q e^{-u}-2Hi e^u & 0 \end{pmatrix} \; . 
\]
The Maurer-Cartan equation is again of the form \eqref{firstMC}, where 
$H$ is now constant.  

Here we prove a basic result about the Maurer-Cartan equation for a Lax pair, 
that this equation holds if and only if there 
exists a solution to the Lax pair.  
This result is fundamental for the discussion here, as it tells us that the 
Lax pair has solutions if and only if the Maurer-Cartan equation holds.  
This proposition can be applied to the above Lax pair for each fixed 
$\lambda$, where $z$ and $w=\bar z$ are considered as independent variables 
and $U=\Theta +\Upsilon_z$, $V=\bar \Theta -\Upsilon_{\bar z}$, 
$F=\mathcal{F}$ and $n=3$.  

\begin{proposition}\label{prop:MCeqn}
Let $\mathcal{U}\subset \mathbb{C}\times\mathbb{C}$ be an open set containing $(0,0)$.  
For $U,V: \mathcal{U} \to \slg_n\!\mathbb{C}$, 
there exists a solution $F=F(z,w):\mathcal{U}\to\SL_n\!\mathbb{C}$ of the Lax pair 
\begin{equation}\label{LaxPair-in-a-proposition} 
F_z = F U \; , \;\;\; F_w = F V \end{equation} for any initial condition 
$F(0,0)\in\SL_n\!\mathbb{C}$ if and only if 
\[ U_w-V_z+[V,U]=0 \; . \]  Furthermore, for any two solutions $F$ and 
$\tilde F$, $\tilde F = G F$ for some constant matrix $G$.  
\end{proposition}

\begin{proof}
Assume there exists an invertible solution $F(z,w)$.  
Then the fact that $(F_z)_w=(F_w)_z$ quickly implies 
\[ U_w-V_z+[V,U]=0 \; . \]  

To prove the converse, suppose that $U_w-V_z+[V,U]=0$.  Then solve 
\[ (F(z,0))_z = F(z,0) U(z,0) \] for variable $z$ 
with initial condition $F(0,0)$.  So now 
$F(z,0)$ is defined.  Now, for each fixed $z_0$, solve 
\[ (F(z_0,w))_w = F(z_0,w) V(z_0,w) \] for the variable $w$ 
with initial condition $F(z_0,0)$.  So now 
$F(z,w)$ is defined and $F_w = F V$ for all $z,w$.  It remains to show $F_z=FU$ for 
all $z,w$.  Note that 
\[ (F_z - F U)(z,w) = 0 \] if $w=0$ (for all $z$).  Using 
$F_{zw}=F_{wz}$, we have 
\begin{equation}\label{difference-eqn} 
(F_z - F U)_w = F_z V - F (U_w-V_z+[V,U]) - F U V = 
(F_z-F U) V \; , \end{equation} since $U_w-V_z+[V,U]=0$ and $F_w = F V$ for 
all $z,w$.  The desired relation $F_z=FU$ was never used in computing 
\eqref{difference-eqn}.  (Of course, we cannot use it, as this is what we wish to 
prove.)  Hence $\mathcal{G} =F_z - F U$ is a solution of $\mathcal{G}_w = 
\mathcal{G} V$ with initial condition $0$.  Hence $\mathcal{G}$ is identically zero, 
by the uniqueness of solutions of $\mathcal{G}_w = 
\mathcal{G} V$ with given initial condition.  So $F_z - F U$ is identically 
zero and hence $F$ is a solution to the Lax pair \eqref{LaxPair-in-a-proposition}.  
Since $\det F\cdot\text{trace}(F^{-1}F_z)=(\det F)_z$ and 
$\det F\cdot\text{trace}(F^{-1}F_w)=(\det F)_w$ by 
Lemma \ref{matrix-identity-lemma}, and $U,\,V\in\slg_n\!\mathbb{C}$, 
we have $(\det F)_z=(\det F)_w=0$.  
So $\det F=1$ because $\det (F(0,0))=1$.

Now, suppose $F$ and $\tilde F$ both solve the Lax pair 
\eqref{LaxPair-in-a-proposition}.  Then 
\[ (\tilde F F^{-1})_z = 
\tilde F_z F^{-1}-\tilde F F^{-1} F_z F^{-1} = 
\tilde F U F^{-1}-\tilde F U F^{-1} = 0 \; . \]  
Similarly, $(\tilde F F^{-1})_w = 0$.  So $\tilde F = G F$ for 
some constant matrix $G$.  
\end{proof}

\section{The $2 \times 2$ Lax pair for CMC surfaces}
\label{2by2laxpair}

Now we will rework the $3 \times 3$ frame $\mathcal{F}$ into a $2 \times 2$ frame $F$.  
Fix a starting point $z_* \in \Sigma$, 
and fix $\hat F(z_*,\overline{z_*},\lambda)= \id$ 
for all $\lambda \in \mathbb{S}^1$.  Consider the Pauli matrices 
\[ \sigma_1=\begin{pmatrix} 0 & 1 \\ 1 & 0 \end{pmatrix} \; , \;\;\; 
\sigma_2=\begin{pmatrix} 0 & -i \\ i & 0 \end{pmatrix} \; , \;\;\; 
\sigma_3=\begin{pmatrix} 1 & 0 \\ 0 & -1 \end{pmatrix} \; . \] 
These matrices, along with the $2 \times 2$ identity matrix $I$, 
can be used to form a basis \[ \{ I,-i \sigma_1, -i \sigma_2, -i \sigma_3 \} \] for 
the ring of quaternions, since we have the quaternionic relations 
\[ (-i\sigma_1)(-i\sigma_2)=(-i\sigma_3)=-(-i\sigma_2)(-i\sigma_1) \; , \] 
\[ (-i\sigma_2)(-i\sigma_3)=(-i\sigma_1)=-(-i\sigma_3)(-i\sigma_2) \; , \]
\[ (-i\sigma_3)(-i\sigma_1)=(-i\sigma_2)=-(-i\sigma_1)(-i\sigma_3) \] and 
$(-i\sigma_j)^2=-\id$ for $j=1,2,3$.  

We identify $\mathbb{R}^3$ with $\su_2$, 
identifying $(x_1,x_2,x_3)\in\mathbb{R}^3$ with the matrix
\[
-x_1\frac{i}{2}\sigma_1+x_2\frac{i}{2}\sigma_2+x_3\frac{i}{2}\sigma_3
=\frac{-i}{2}
\begin{pmatrix} -x_3 & x_1+i x_2 \\ x_1-i x_2 & x_3\end{pmatrix} \; . 
\]
Let $\hat F=\hat F(z,\bar z,\lambda) 
\in \SU_2$ be the matrix that rotates $(-i\sigma_1)/2$, 
$(-i\sigma_2)/2$ and $(-i\sigma_3)/2$ to the $2 \times 2$ matrix forms of $e_1$, 
$e_2$, and $N$, respectively, that is, 
\begin{equation}\label{e1e2}
e_1 = \hat F \left(\frac{-i}{2}\sigma_1\right) \hat F^{-1} \; , \;\;\; 
e_2 = \hat F \left(\frac{-i}{2}\sigma_2\right) \hat F^{-1} \; , \;\;\; 
N = \hat F \left(\frac{-i}{2}\sigma_3\right) \hat F^{-1} \; . \end{equation}  
These relations determine $\hat F$ uniquely, up to sign; and 
because $\hat F \in \SU_2$, conjugation by $\hat F$ as in \eqref{e1e2} really is 
a rotation of $\mathbb{R}^3$ (see Sections \ref{section0} and \ref{section4}).  

Note that $\hat F(z_*,\overline{z_*},\lambda)=\id$ for all $\lambda \in 
\mathbb{S}^1$ means 
we assumed that, at the point $z_*$, the vectors $e_1$, $e_2$, and $N$ are 
the unit vectors in the positive $x_1$-axis, $x_2$-axis, and $x_3$-axis directions, 
respectively, for all $\lambda \in \mathbb{S}^1$.  This can be assumed, without loss 
of generality, just by applying a rigid motion to the surface $f$, for each 
$\lambda \in \mathbb{S}^1$.  

Define 
\[ \hat U = \begin{pmatrix} U_{11} & U_{12} \\ U_{21} & U_{22} \end{pmatrix}
         := \hat F^{-1} \hat F_z \; , \;\;\;
   \hat V = \begin{pmatrix} V_{11} & V_{12} \\ V_{21} & V_{22} \end{pmatrix} 
         := \hat F^{-1} \hat F_{\bar{z}} \; . \] 
Now we compute $\hat U$ and $\hat V$ explicitly in terms of the metric factor term 
$u$, the mean curvature $H$ and the Hopf differential $Q$.  
Using 
\[ e_1 = \frac{f_x}{|f_x|} = 
\frac{f_x}{2e^u} = \hat F \frac{1}{2} 
\begin{pmatrix} 0 & -i \\ -i & 0 \end{pmatrix} \hat F^{-1} \; , \;\;\; 
e_2 = \frac{f_y}{|f_y|} = 
\frac{f_y}{2e^u} = \hat F \frac{1}{2} 
\begin{pmatrix} 0 & -1 \\ 1 & 0 \end{pmatrix} \hat F^{-1} \; , \]
we have 
\begin{equation}\label{phiform} 
f_z = -i e^{u} \hat F \begin{pmatrix} 0 & 0 \\ 1 & 0 
\end{pmatrix} \hat F^{-1} \; , \;\;\; 
f_{\bar{z}} = -i e^{u} \hat F \begin{pmatrix} 0 & 1 \\ 0 & 0 
\end{pmatrix} \hat F^{-1} \; . \end{equation}
Then $f_{z\bar{z}}=f_{\bar{z}z}$ implies that 
$V_{11}-V_{22}=u_{\bar{z}}$ and $U_{22}-U_{11}=u_{z}$ and $V_{12}=-U_{21}$.  
Then $f_{z\bar{z}}=2He^{2u} N$ implies $V_{12}=He^u$ and 
$f_{zz}=2u_z f_z+\lambda^{-2} Q N$ implies 
$U_{12}=e^{-u} \lambda^{-2} Q/2$ and 
$f_{\bar z \bar z}=2u_{\bar z} f_{\bar z}+ 
\lambda^{2} \bar Q N$ implies $V_{21}=-e^{-u} \lambda^{2} \bar Q/2$.  
Since the determinant of $\hat F$ is $1$, \eqref{eq:matrixlemma} implies that 
$\hat U$ and $\hat V$ have trace zero (see also Lemma 
\ref{matrix-identity-lemma}), so it is determined that 
\begin{equation}
 \hat U = \frac{1}{2} \begin{pmatrix} -u_z & e^{-u} \lambda^{-2} Q \\ 
-2He^u & u_z \end{pmatrix}  \; , \;\;\; \hat V = 
\frac{1}{2} \begin{pmatrix} u_{\bar z} & 2He^{u} \\ 
-e^{-u} \lambda^{2} \bar Q & -u_{\bar z} \end{pmatrix} \; .  
\label{eq:UhatVhat}
\end{equation}
Now the Lax pair for $f$ is 
\begin{equation}
 \hat F_z = \hat F \hat U \; , \;\;\; \hat F_{\bar{z}} = \hat F \hat V 
\label{eq:laxpairforf}
\end{equation}
with $\hat U$ and $\hat V$ as in Equation \eqref{eq:UhatVhat}.  
Take any solution $\hat F$ of this Lax pair such that $\hat F\in\SU_2$ for all 
$\lambda\in\mathbb{S}^1$, and define the Sym-Bobenko formula
\begin{equation}\label{firstSymBobenko}
\hat f(z,\bar z,\lambda) =\frac{1}{2H}\cdot\left[ \hat F 
\begin{pmatrix} i & 0 \\ 0 & -i \end{pmatrix}\hat F^{-1}
 -i \lambda (\partial_\lambda \hat F) \cdot \hat F^{-1}\right]  \; . 
\end{equation}

\begin{remark}
The lemma just below tells us that the Sym-Bobenko formula \eqref{firstSymBobenko} 
retrieves the CMC immersion $f$ for us.  There are advantages to using this 
formula, and we note two of them here: 1) Given the frame $\hat F$, we then know 
$f_z$ and $f_{\bar z}$, so we would expect to need to integrate in order to find 
$f$ itself.  The Sym-Bobenko formula gives us a way to avoid integration, using 
the derivative of $\hat F$ with respect to the spectral parameter $\lambda$ instead.  
2) When the domain (locally parametrized by $z$) is not simply connected, it is 
not always clear that the surface $f$ resulting from the frame $\hat F$ will be 
well defined on that domain, because the frame $\hat F$ produced as a solution 
to the Lax pair might not be well defined on the domain.  
The Sym-Bobenko formula gives us a means for 
understanding when $f$ will be well defined.  We had some hints of this 
in Sections \ref{section4} and \ref{section5}, and we will see more explicitly 
why this is so in Section \ref{section13}.  
\end{remark}

\begin{lemma}\label{lm:2by2Laxpair}
The CMC $H$ surfaces $f(z,\bar z,\lambda)$ with $H \ne 0$ as in 
Section $\ref{section6}$ and the surfaces $\hat f(z,\bar z,\lambda)$ 
as in \eqref{firstSymBobenko} differ only by rigid motions of $\mathbb{R}^3$.  
Thus Equation \eqref{firstSymBobenko} produces the associated family of 
any CMC $H$ surface $f$ from a frame $\hat F$ solving 
\eqref{eq:laxpairforf}-\eqref{eq:UhatVhat}.  

Conversely, for any $u$ and $Q$ satisfying the Gauss and Codazzi equations 
\eqref{firstMC}, and any solution $\hat F$ of the Lax pair 
\eqref{eq:laxpairforf} satisfying \eqref{eq:UhatVhat} so that $\hat F\in\SU_2$ 
for all $\lambda\in\mathbb{S}^1$, 
$\hat f(z,\bar z,\lambda)$ as defined in \eqref{firstSymBobenko} 
is a conformal CMC $H$ immersion into $\mathbb{R}^3$ with metric $4e^{2u}(dx^2+dy^2)$ 
and Hopf differential $\lambda^{-2}Q$.  
\end{lemma}

\begin{proof}
Using the previous forms for $f_z$ and $f_{\bar z}$ in 
Equation \eqref{phiform} and computing $\hat f_z$ and 
$\hat f_{\bar z}$, we see that $\hat f_{z}=f_{z}$ and 
$\hat f_{\bar z}=f_{\bar z}$.  Thus 
$f$ and $\hat f$ are the same surfaces, up to translations.  

To prove the converse, we see that 
\[ \hat f_x = e^u \hat F \begin{pmatrix} 0 & -i \\ -i & 0 \end{pmatrix} 
\hat F^{-1} \; , \;\;\; \hat f_y = e^u \hat F \begin{pmatrix} 0 & -1 \\ 
1 & 0 \end{pmatrix} \hat F^{-1} \; , \] 
because $\hat f_x=\hat f_z+\hat f_{\bar z}$ and 
$\hat f_y=i(\hat f_z-\hat f_{\bar z})$.  And so 
\[ \hat N:=\frac{\hat f_x \times \hat f_y}{|\hat f_x \times \hat f_y|}
= \hat F \frac{1}{2} \begin{pmatrix} -i & 0 \\ 0 & i \end{pmatrix}
\hat F^{-1} \; . \]  
Then, using that $\langle A,B\rangle$ for vectors $A,B \in \mathbb{R}^3$ 
becomes 
$-2 \cdot \text{trace}(AB)$ when $A,B$ are in $2 \times 2$ matrix form 
(see Equation \eqref{metricinmatrices1}), we have 
\[ \langle\hat f_x, \hat f_x\rangle = \langle\hat f_y, \hat f_y\rangle = 4 
e^{2u} \; , \;\;\; 
\langle \hat f_x, \hat f_y\rangle = 0 \; . \] 
So the metric of $\hat f$ is conformal and is $4 e^{2 u} (dx^2+dy^2)$.  
Let $\hat Q$ and $\hat H$ be the Hopf differential 
and mean curvature of $\hat f$.  Then 
\[ \hat Q = \langle \hat f_{zz}, \hat N\rangle = i e^u \cdot \text{trace} 
\left[ \left( \begin{pmatrix} 0 & 0 \\ u_z & 0 \end{pmatrix} + 
\left[ \hat U, \begin{pmatrix} 0 & 0 \\ 1 & 0 \end{pmatrix} \right] 
\right) \begin{pmatrix} -i & 0 \\ 0 & i \end{pmatrix} 
\right] =\lambda^{-2}Q \; , \]  
using the form of $\hat U$, 
so $\lambda^{-2}Q$ is the Hopf differential of $\hat f$.  Also, 
\[
\hat H = \frac{1}{2} \langle\hat f_{xx}+\hat f_{yy},\hat N\rangle \cdot 
                     \langle\hat f_x , \hat f_x \rangle ^{-1} 
       =  \frac{1}{2} e^{-2u} \langle\hat f_{z \bar z},\hat N\rangle 
       = \]\[ = i e^{-u} \cdot \text{trace} 
\left[ \left( \begin{pmatrix} 0 & 0 \\ u_{\bar z} & 0 \end{pmatrix} + 
\left[ \hat V, \begin{pmatrix} 0 & 0 \\ 1 & 0 \end{pmatrix} \right] 
\right) \begin{pmatrix} 1 & 0 \\ 0 & -1 \end{pmatrix} 
\right] = H \; . \]
So $H$ is the mean curvature of $\hat f$.  This completes the proof.  
\end{proof}

So we have shown that solutions of the integrable Lax pair 
\eqref{eq:laxpairforf} satisfying \eqref{eq:UhatVhat} and the 
Sym-Bobenko formula \eqref{firstSymBobenko} locally produce 
all CMC $H$ surfaces, and then that any solution $\hat F\in\SU_2$ of 
\eqref{eq:laxpairforf} produces a CMC $H$ surface 
via \eqref{firstSymBobenko}.  Therefore, finding arbitrary CMC $H$ surfaces is 
equivalent to finding arbitrary solutions $\hat F\in\SU_2$ of the 
Lax pair \eqref{eq:laxpairforf}-\eqref{eq:UhatVhat}.  

\section{The $2 \times 2$ Lax pair in twisted form}
\label{laxpairintwisted}

Now we make some changes in the formulation of Section \ref{2by2laxpair}.  
These changes are not essential, in that they change the resulting CMC 
immersion only by a rigid motion, but we make them so that the formulation 
here will be the same as in \cite{DW} and in the Sym-Bobenko formula 
\eqref{firstSym}.  

Let $F$ be defined by  
\[ \hat F = -\sigma_3\left( F^{-1}\right) ^t 
\begin{pmatrix} \sqrt{\lambda} & 0 \\ 0 & 1/\sqrt{\lambda} 
\end{pmatrix} \sigma_3 \; . \]
We find that the Lax pair \eqref{eq:laxpairforf}-\eqref{eq:UhatVhat} is 
written in the form 
\begin{equation}\label{LaxPair}
F_z = F U \; , \;\;\; F_{\bar z} = F V \; , 
\end{equation}  where
\begin{equation}
U = \frac{1}{2} \begin{pmatrix} u_z & -2He^{u} \lambda^{-1} \\ 
Q e^{-u} \lambda^{-1} & -u_z \end{pmatrix} \; , \;\;\; 
V = \frac{1}{2} \begin{pmatrix} -u_{\bar z} & -\bar Q e^{-u} \lambda \\ 
2He^{u} \lambda & u_{\bar z} \end{pmatrix} \; , 
\label{eq:UandV}
\end{equation}
and that 
\begin{equation}\label{fhatrelationwithsigma3} \hat f(z,\bar z,\lambda)
 = -\sigma_3 \left( f(z,\bar z,\lambda)\right) ^t \sigma_3 \; , \end{equation} where
\begin{equation}\label{Sym-Bobenko} 
f(z,\bar z,\lambda) = 
    \frac{1}{2H}\cdot\left[ \frac{-1}{2}F\begin{pmatrix}i&0\\0&-i\end{pmatrix}
    F^{-1} - i\lambda (\partial_{\lambda} F) \cdot F^{-1}\right] \; , 
\end{equation}
and using the $\mathbb{R}^3$ metric \eqref{metricinmatrices1} we can check 
that the normal vector to $f(z,\bar z,\lambda)$ is 
\begin{equation}\label{eq:normalR3}
N=F\frac{1}{2}\begin{pmatrix}i&0\\0&-i\end{pmatrix}F^{-1} \; .
\end{equation}
Note that the transformation $\hat f(z,\bar z,\lambda) \to f(z,\bar z,\lambda)$ 
in \eqref{fhatrelationwithsigma3} just 
represents the rotation by angle $\pi$ about the $x_1$-axis 
in $\mathbb{R}^3$ (we are using \eqref{first-eqn-in-chap2} here), 
so $f(z,\bar z,\lambda)$ in \eqref{Sym-Bobenko} and $\hat f(z,\bar z,\lambda)$ 
in \eqref{firstSymBobenko} are the same surface (up to a rigid motion).  
Note also that Lemma \ref{lm:2by2Laxpair} implies that then 
$f(z,\bar z,\lambda)$ in \eqref{Sym-Bobenko} and the previous $f(z,\bar z,\lambda)$ 
defined in Section \ref{section6} differ by only a rigid motion.  
Thus Lemma \ref{lm:2by2Laxpair} holds for the $f(z,\bar z,\lambda)$ 
in \eqref{Sym-Bobenko} as well: 

\begin{corollary}\label{cor:twistedlaxpair}
The CMC $H$ surfaces $f(z,\bar z,\lambda)$ with $H \ne 0$ as in 
Section $\ref{section6}$ and the surfaces $f(z,\bar z,\lambda)$ 
as in \eqref{Sym-Bobenko} differ only by a rigid motion of $\mathbb{R}^3$.  
Thus Equation \eqref{Sym-Bobenko} produces the associated family of 
any CMC $H$ surface $f$ from a frame $F$ solving 
\eqref{eq:UandV}-\eqref{LaxPair}.  

Conversely, for any $u$ and $Q$ satisfying the Gauss and Codazzi equations 
\eqref{firstMC}, and any solution $F$ of the Lax pair 
\eqref{LaxPair} satisfying \eqref{eq:UandV}, 
$f(z,\bar z,\lambda)$ as defined in \eqref{Sym-Bobenko} 
is a conformal CMC $H$ immersion into $\mathbb{R}^3$ with metric $4e^{2u}(dx^2+dy^2)$ 
and Hopf differential $\lambda^{-2}Q$.  
\end{corollary}

As the formulation in \eqref{LaxPair} and \eqref{eq:UandV} and 
\eqref{Sym-Bobenko} is used by Dorfmeister and Wu \cite{DW}, 
we will often use this formulation from here on.  

\begin{remark}
One often sees the matrices $U$ and $V$ in \eqref{eq:UandV} without $H$ 
included, that is, as 
\begin{equation}
    \frac{1}{2} \begin{pmatrix} u_z & -e^{u} \lambda^{-1} \\ 
Q e^{-u} \lambda^{-1} & -u_z \end{pmatrix} \; , \;\;\text{and}\;\; 
    \frac{1}{2} \begin{pmatrix} -u_{\bar z} & -\bar Q e^{-u} \lambda \\ 
e^{u} \lambda & u_{\bar z} \end{pmatrix} \; .  
\end{equation}
However, this is not an essential difference.  One can easily ``remove'' the 
$H$ by either fixing $H=1/2$, or changing the parameter $z$ to 
$2Hz$ and renaming $Q/(2H)$ as $Q$.  
\end{remark}

\section{The Weierstrass representation}
\label{Wrepsection}

We can now give two applications of the Lax pairs in Sections \ref{section6}, 
\ref{2by2laxpair} and \ref{laxpairintwisted}.  The first, in this section, 
describes a representation for minimal surfaces in $\mathbb{R}^3$.  
The second, in Section \ref{harmonic}, shows that the Gauss map of a conformal 
CMC immersion into $\mathbb{R}^3$ is harmonic.  
A third application is to show existence 
of CMC surfaces that are topologically cylindrical in the associated 
families of Delaunay surfaces, and is included in the supplement 
\cite{waynebook-maybe} to these notes.  

Let $f: \Sigma \to \mathbb{R}^3 $ be a minimal ($H=0$) conformal immersion, 
where $\Sigma$ is a simply-connected domain in $\mathbb{C}$.  
So $\langle f_z,f_z\rangle =\langle f_{\bar{z}},f_{\bar{z}}\rangle =0$, 
$\langle f_z,f_{\bar{z}}\rangle =2e^{2u}$
(this is the definition of $u:\Sigma\to\mathbb{R}$).  
By the Gauss-Weingarten equations (with $H=0$ and $\lambda =1$) 
in Section \ref{section6}, we have
\[ f_{zz}=2u_zf_z+QN \; , \;\;\; f_{z{\bar{z}}}=0 \; , \;\;\; 
   f_{\bar{z}\bar{z}}=2u_{\bar{z}}f_{\bar{z}}+\bar{Q}N \; . \]
Recall that the Hopf differential $Q$ is defined as 
$Q=\langle f_{zz},N\rangle$, and $N:\Sigma\to \mathbb{S}^2$ is the unit normal vector 
of $f$, oriented as in Section \ref{section6}.  

\begin{remark}\label{re:KformininR3}
Equation \eqref{eq:K=detS} implies that the Gaussian curvature for a minimal surface is 
always non-positive.
\end{remark}

As in Section \ref{2by2laxpair}, we can take $\hat F \in \SU_2$ satisfying
\[ e_1=\frac{f_x}{|f_x|}=\hat F \frac{-i}{2}\sigma_1\hat F^{-1} \; , \;\;\; 
   e_2=\frac{f_y}{|f_y|}=\hat F \frac{-i}{2}\sigma_2\hat F^{-1} \; , \;\;\; 
                       N=\hat F \frac{-i}{2}\sigma_3\hat F^{-1} \; . \]
Then, as in Section \ref{2by2laxpair}, now with $H=0$ and $\lambda =1$, 
the Lax pair for $f$ is 
\[ \hat U=\hat F^{-1}\hat F_z
=\frac{1}{2}\begin{pmatrix}-u_z&e^{-u}Q\\0&u_z \end{pmatrix} \; , \;\;\; 
\hat V=\hat F^{-1}\hat F_{\bar z}=\frac{1}{2}
\begin{pmatrix}u_{\bar z}&0\\-e^{-u}\bar{Q}&-u_{\bar z}\end{pmatrix} \; . \]
Defining functions $a,b:\Sigma\to\mathbb{C}$ so that
\[ \hat F= \begin{pmatrix} e^{-u/2} \bar a & e^{-u/2} b \\ 
-e^{-u/2} \bar b & e^{-u/2} a \end{pmatrix} \] 
holds, then $a\bar{a}+b\bar{b}=e^u$ because $\det \hat F = 1$.  We have 
\[ \hat V=
\frac{1}{2}\begin{pmatrix}u_{\bar z}&0\\-e^{-u}\bar Q&-u_{\bar z} \end{pmatrix}
=\begin{pmatrix} -u_{\bar z}/2+e^{-u}(a\bar{a}_{\bar{z}}+b\bar{b}_{\bar{z}})&
e^{-u}(ab_{\bar{z}}-ba_{\bar{z}})\\e^{-u}(\bar{b}\bar{a}_{\bar{z}}-
\bar{a}\bar{b}_{\bar{z}})&-u_{\bar{z}}/2+e^{-u}(\bar{a}a_{\bar{z}}+
\bar{b}b_{\bar{z}})
\end{pmatrix} \; .  \]
It follows that 
\[ \begin{pmatrix} -b & a \\ \bar a & \bar b \end{pmatrix}
   \begin{pmatrix} a_{\bar z} \\ b_{\bar z} \end{pmatrix}
  =\begin{pmatrix} 0 \\ 0 \end{pmatrix}
\]
and so $a_{\bar{z}}=b_{\bar{z}}=0$; that is, $a$ and $b$ are 
holomorphic.  

\begin{remark}\label{re:hol-fct}
For any real valued function $\psi :\Sigma\to\mathbb{R}$ and holomorphic 
function $\Psi :\Sigma\to\mathbb{C}$, we have the general property, from the 
Cauchy-Riemann equations for $\Psi$, that 
$\psi_z = \Psi_z/2$ if and only if $\psi = \Re \Psi + c$ 
for some real constant $c$.  We will use this property below.  
\end{remark}

Note that \[ f_z=-ie^u\hat F \begin{pmatrix}0&0\\1&0\end{pmatrix}\hat F^{-1}
= i \begin{pmatrix} -ab&b^2\\-a^2&ab \end{pmatrix} \; , \] 
which is holomorphic and is written as \[ f_z=(a^2-b^2,i(a^2+b^2),-2ab) \] 
in the standard $\mathbb{R}^3$ coordinates via the identification 
\eqref{first-eqn-in-chap2} (also found at the beginning of Section \ref{2by2laxpair}).  
Since $f$ is real-valued, by Remark \ref{re:hol-fct}, we have 
\[ \Re \int f_z dz = \frac{1}{2} f + (c_1, c_2, c_3) \]
for some constant $(c_1,c_2,c_3)\in\mathbb{R}^3$.  So up to a translation, 
\begin{equation}\label{eq:W-rep}
\begin{array}{r@{=}l}
f=2 \Re \displaystyle\int f_z dz
 & 2\Re \displaystyle\int\limits_{\!\!\!} (a^2-b^2,i(a^2+b^2),-2ab) dz \\
 &  \Re \displaystyle\int (1-g^2,i(1+g^2),2g) \eta \; , 
\end{array}
\end{equation}
where $g=-b/a$ and $\eta=2a^2 dz$.
This is the Weierstrass representation for minimal surfaces: 
\begin{quote}
{\bf The Weierstrass representation}: 
Any simply-connected minimal immersion from 
$\Sigma$ into $\mathbb{R}^3$ can be parametrized as in \eqref{eq:W-rep}, 
using a meromorphic function $g:\Sigma\to\mathbb{C}$ and 
holomorphic 1-form $\eta$ on $\Sigma$.  
\end{quote}

The term "meromorphic" refers to any function that is holomorphic where it 
is finite, and that has only isolated singularities, which are always poles (not 
essential singularities).  

Also, 
\[
(1+g\bar g)^2\eta\bar\eta =4e^{2u}dzd\bar z 
\]
is the metric for the surface $f$ in \eqref{eq:W-rep}.  Furthermore, 
\[
N=\hat F\frac{-i}{2}\sigma_3 \hat F^{-1}=\frac{-ie^{-u}}{2}
\begin{pmatrix}a\bar a-b\bar b & -2\bar ab \\ -2a\bar b & b\bar b -a\bar a
\end{pmatrix} \; , 
\]
which is written as 
\begin{eqnarray*}
N &=& 
e^{-u}\left( -(\bar a b+a\bar b),i(\bar a b-a\bar b),b\bar b-a\bar a\right) \\
  &=& 
\left( \frac{g+\bar g}{g\bar g +1},i\frac{\bar g-g}{g\bar g +1},
       \frac{g\bar g -1}{g\bar g +1}\right)
\end{eqnarray*}
in the standard $\mathbb{R}^3$ coordinates via the identification of the 
beginning of Section \ref{2by2laxpair}.  
Thus the function $g$ is the composition of the Gauss map with stereographic 
projection from $\mathbb{S}^2$ to the Riemann sphere $\mathbb{C}\cup\{\infty\}$.  

Now we consider some examples.

\begin{example}
The plane can be given by the Weierstrass representation with 
$(g,\eta)=(c_1,c_2dz)$ on $\Sigma=\mathbb{C}$, where $c_1\in\mathbb{C}$, 
$c_2\in\mathbb{C}\setminus\{0\}$ are constants.  
The plane is the most trivial example of a minimal surface.  
\end{example}

\begin{figure}[phbt]
\begin{center}
\includegraphics[width=0.25\linewidth]{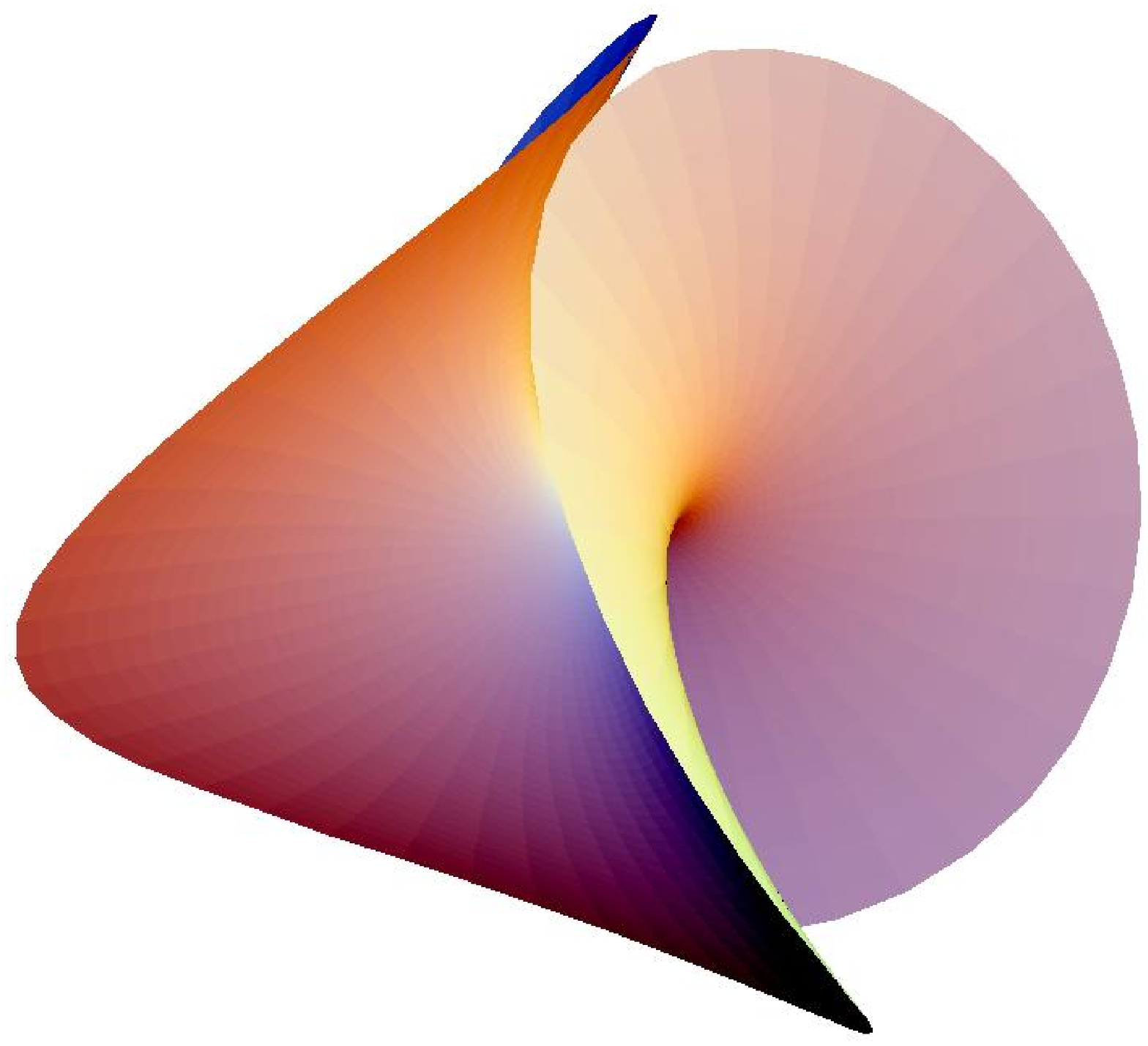}
\includegraphics[width=0.3\linewidth]{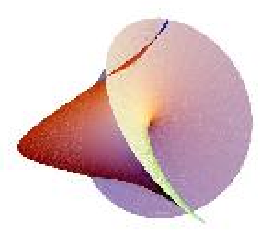}
\end{center}
\vspace{2in}
\caption{A compact portion of Enneper's surface.  
A smaller portion is shown on the left, 
and a larger portion is shown on the right, 
to indicate how the surface intersects itself.}
\label{fg:Enn}
\end{figure}

\begin{figure}[phbt]
\begin{center}
\includegraphics[width=0.4\linewidth]{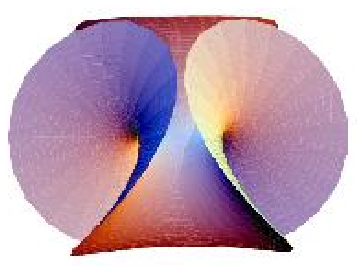}
\includegraphics[width=0.4\linewidth]{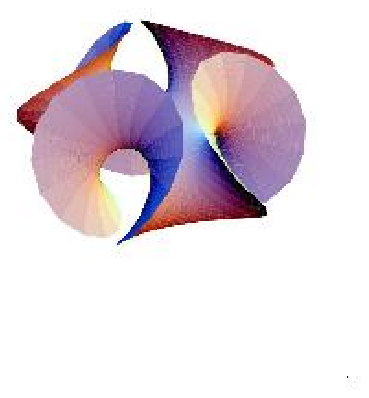}
\end{center}
\vspace{2in}
\caption{Higher order versions of Enneper's surface, for $n=2$ (on the left), 
and $n=3$ (on the right).}
\label{fg:H-Enn}
\end{figure}

\begin{example}
The minimal Enneper's surface can be given by the Weierstrass representation 
with $(g,\eta)=(z,cdz)$ on $\Sigma=\mathbb{C}$, 
where $c\in\mathbb{R}\setminus\{0\}$ is a constant.  
An explicit parametrization of Enneper's surface is 
\[
f(x,y)=c\left( x+xy^2-\frac{x^3}{3},-y-x^2y+\frac{y^3}{3},x^2-y^2\right) \; , 
\]
where $z=x+iy$.  
See Figure \ref{fg:Enn}.  
This surface was first found by A. Enneper in 1864 (see \cite{Enn}).  

If we replace $g=z$ by $z^n$ ($n\in\mathbb{N}\setminus\{1\}$), then we have 
higher order versions of Enneper's surface.  See Figure \ref{fg:H-Enn}.
\end{example}

\begin{figure}[phbt]
\begin{center}
\includegraphics[width=0.3\linewidth]{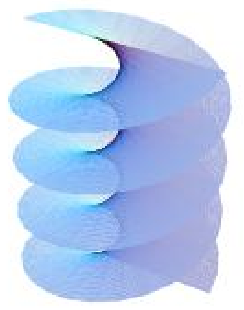}
\raisebox{5ex}{\includegraphics[width=0.3\linewidth]{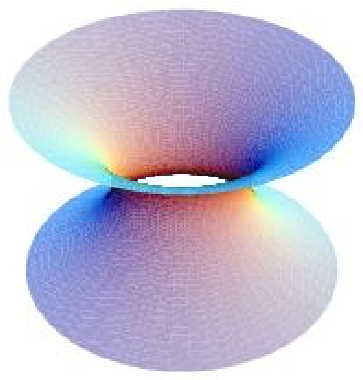}}
\end{center}
\vspace{2in}
\caption{A helicoid and a catenoid.}
\label{fg:cat-heli}
\end{figure}

\begin{example}\label{ex:helicoid}
The minimal helicoid is given by the Weierstrass representation with 
$(g,\eta)=(e^z,cie^{-z}dz)$ on $\Sigma=\mathbb{C}$, 
where $c\in\mathbb{R}\setminus\{0\}$ is a constant.  
The resulting explicit parametrization of the helicoid is 
\[
f(x,y)=2c(\sinh x\sin y, -\sinh x\cos y,-y) \; , 
\]
where $z=x+iy$.  
Setting $u=-2\sinh x$ and $v=-y$, we have
\[
f(u,v)=c(u\sin v,u\cos v,v) \; .
\]
So one can easily see that this surface has screw-motion symmetry in the 
direction of the $x_3$-axis.  
See the left-hand side of Figure \ref{fg:cat-heli}.
\end{example}

When the Riemann surface $\Sigma$ is not simply connected, the integration in 
Equation \eqref{eq:W-rep} might not be independent of the choice of path, 
so $f$ might not be well-defined on $\Sigma$ in general.
The immersion $f$ being well-defined on $\Sigma$ is equivalent to the vector 
\begin{equation}\label{eq:period_for_min_surf}
\oint_\delta(1-g^2,i(1+g^2),2g) \eta 
\end{equation}
being purely imaginary for any closed loop in $\Sigma$.  

\begin{example}
The minimal catenoid is given by the Weierstrass representation with 
$(g,\eta)=(z,cz^{-2}dz)$ on $\Sigma=\mathbb{C}\setminus\{0\}$, 
where $c\in\mathbb{R}\setminus\{0\}$ is a constant.  
To verify that $f$ is well-defined on $\Sigma$, it suffices to check that the 
vector \eqref{eq:period_for_min_surf} is purely imaginary for the loop
\[
\delta :[0,2\pi]\ni\theta\mapsto e^{i\theta}\in\Sigma \; ,
\]
and this is indeed the case.  
The resulting explicit parametrization of the catenoid is 
\[
f(x,y)=c\left( -x\frac{x^2+y^2+1}{x^2+y^2},
               -y\frac{x^2+y^2+1}{x^2+y^2},\log(x^2+y^2)\right) \; , 
\]
where $z=x+iy$.  
Setting $x=-e^v\cos u$ and $y=-e^v\sin u$, we have
\[
f(u,v)=2c(\cos u\cosh v,\sin u\cosh v,v) \; .
\]
So one can easily see that the catenoid is a surface of revolution about the 
$x_3$-axis generated by the catenary $x_1=\cosh x_3$.  
See the right-hand side of Figure \ref{fg:cat-heli}.  
\end{example}

\begin{example}
Consider the surface given by the Weierstrass representation with 
$(g,\eta)=(z,ciz^{-2}dz)$ on $\Sigma=\mathbb{C}\setminus\{0\}$, 
where $c\in\mathbb{R}\setminus\{0\}$ is a constant.  
Computing the vector in \eqref{eq:period_for_min_surf} similarly to the 
previous example, we find that this vector is $(0,0,4c\pi)$, which is not 
purely imaginary, so $f$ is not well-defined on $\Sigma$.  The resulting 
surface is now singly-periodic in the direction of the $x_3$-axis.  
If we define $\widetilde{\Sigma}=\mathbb{C}$ and 
$\rho :\widetilde{\Sigma}\to\Sigma$ by $\rho (w)=e^w$ for 
$w\in\widetilde{\Sigma}$, then $\rho$ is surjective and $\widetilde{\Sigma}$ 
is the universal cover of $\Sigma$.  
Furthermore, setting $\tilde{\eta}=\eta\circ\rho$ and 
$\tilde{g}=g\circ\rho$, then by Example \ref{ex:helicoid}, 
$\tilde{f}=f\circ\rho$ is a helicoid that is singly-periodic in the $x_3$-axis 
direction.  
\end{example}

\begin{figure}[phbt]
\begin{center}
\includegraphics[width=0.4\linewidth]{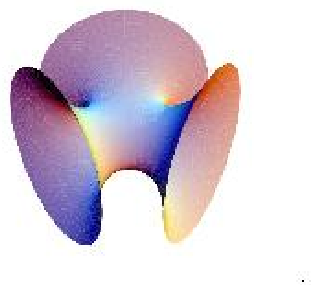} 
\end{center}
\vspace{2in}
\begin{center}
\includegraphics[width=0.4\linewidth]{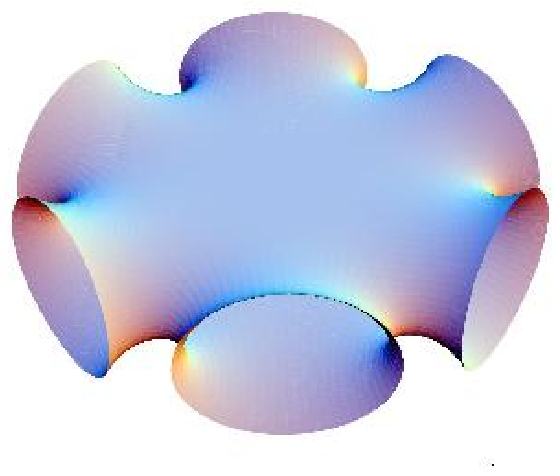} 
\end{center}
\vspace{2in}
\begin{center}
\includegraphics[width=0.4\linewidth]{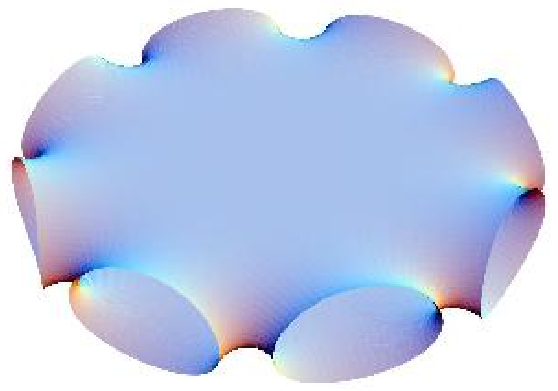} 
\end{center}
\vspace{2in}
\begin{center}
\includegraphics[width=0.4\linewidth]{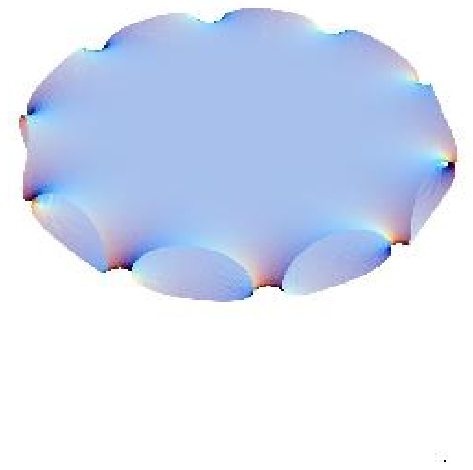}
\end{center}
\vspace{2in}
\caption{Jorge-Meeks $n$-noids for $n=3$, $6$, $8$ and $10$.}
\label{fg:JM3-10}
\end{figure}

\begin{example}
Minimal Jorge-Meeks $n$-noids ($n\in\mathbb{N}\setminus\{ 1,2 \}$) are given by 
the Weierstrass representation with $(g,\eta)=(z^{n-1},c(z^n-1)^{-2}dz)$ on 
$\Sigma=(\mathbb{C}\cup\{\infty\})\setminus\{z\in\mathbb{C}|z^n=1\}$, 
where $c\in\mathbb{R}\setminus\{0\}$ is a constant.  
Computing residues of each component, one can verify that this surface is 
well-defined on $\Sigma$ for any $n\in\mathbb{N}\setminus\{1\}$.  
This surface was first found by L. P. Jorge and W. H. Meeks III in 1983 
(see \cite{JM, BaCo}).  When $n=2$, this surface is a catenoid.  
See Figure \ref{fg:JM3-10}.
\end{example}

\begin{figure}[phbt]
\begin{center}
\includegraphics[width=0.4\linewidth]{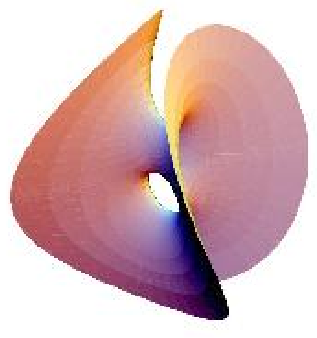}
\includegraphics[width=0.4\linewidth]{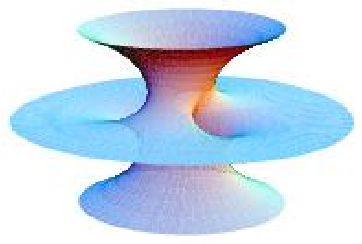}
\end{center}
\vspace{2in}
\caption{Genus 1 minimal surfaces: 
the Chen-Gackstatter surface and the Costa surface.}
\label{fg:ChenGack}
\end{figure}

\begin{example}
The minimal Chen-Gackstatter surface is given by the Weierstrass representation 
with $(g,\eta)=(Bw/z,zdz/w)$ on the hyperelliptic Riemann surface 
$\Sigma=\{(z,w)\in(\mathbb{C}\cup\{\infty\})^2\;|\;w^2=z(z^2-1)\}\setminus
        \{(\infty ,\infty )\}$, 
where $B$ is the constant  
\[
B=\sqrt{\int_0^1\frac{tdt}{\sqrt{t(1-t^2)}}
\left/\int_0^1\frac{(1-t^2)dt}{\sqrt{t(1-t^2)}}\right.} \; .
\]
This surface was first found by C.~C.~Chen, F.~Gackstatter in 1982 
(see \cite{cg, BaCo}).  See the surface on the left-hand side of 
Figure \ref{fg:ChenGack}.
\end{example}

\begin{example}
The Costa-Hoffman-Meeks minimal surface of genus $k$ is given by the 
Weierstrass representation with $(g,\eta)=(B/w,wdz/(z^2-1))$ on the 
Riemann surface 
$\Sigma=\{(z,w)\in(\mathbb{C}\cup\{\infty\})^2\;|\;w^{k+1}=z^k(z^2-1)\}
        \setminus\{(-1,0),(1,0),(\infty ,\infty )\}$, 
where $B$ is the constant  
\[
B=\sqrt{2\int_0^1\left(\frac{t}{1-t^2}\right)^{k/(k+1)} dt \left/
         \int_0^1\frac{dt}{t^{k/(k+1)}(1-t^2)^{1/(k+1)}}\right.} \; ,
\]
and $k$ is an arbitrary positive integer.  
The surface for $k=1$ was first found by C.~J.~Costa in 1982 in his Ph.D. 
thesis (see \cite{Co, BaCo}).  See the right-hand side of Figure \ref{fg:ChenGack}.
In 1985, D.~Hoffman and W.~Meeks \cite{HM1} showed that these surfaces are 
embedded.  So these surfaces are counterexamples to a longstanding conjecture 
that the plane, catenoid and helicoid are the only complete embedded minimal 
surfaces with finite topology.  See also \cite{HM2}.  
\end{example}
\begin{figure}[phbt]
\begin{center}
\begin{tabular}{cc}
\includegraphics[width=0.4\linewidth]{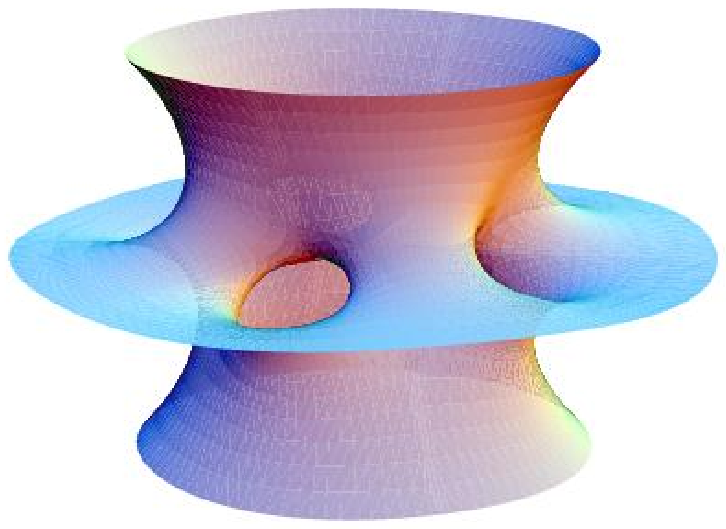} &
\includegraphics[width=0.4\linewidth]{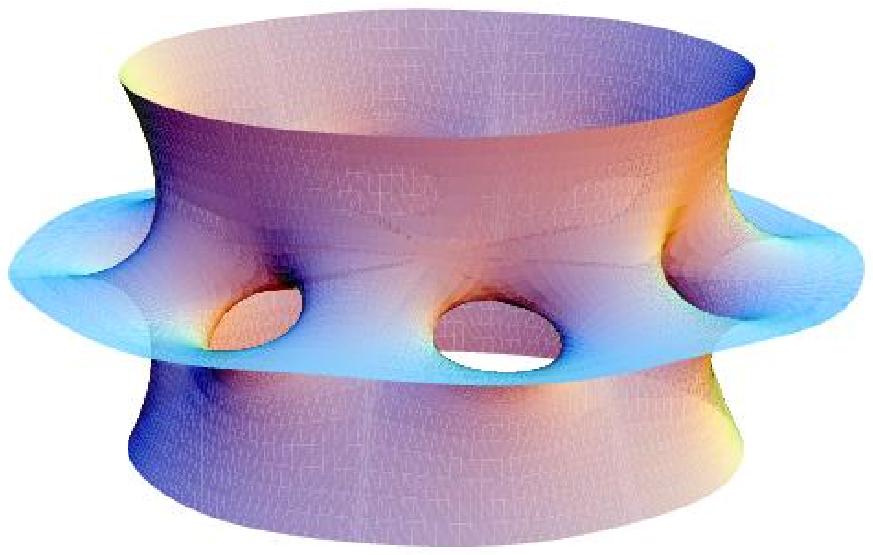} \\
\includegraphics[width=0.4\linewidth]{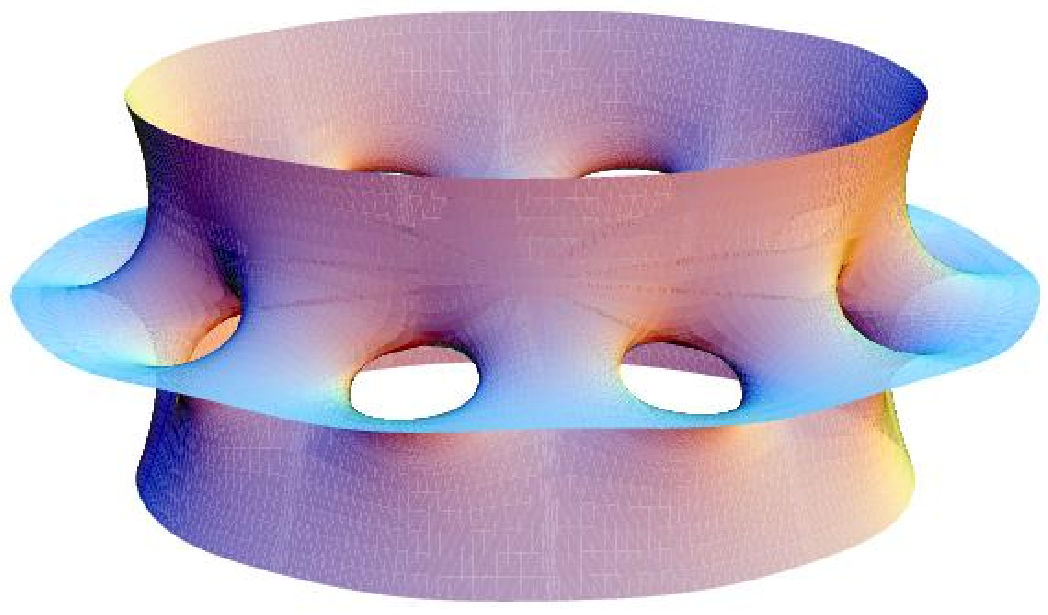} &
\includegraphics[width=0.4\linewidth]{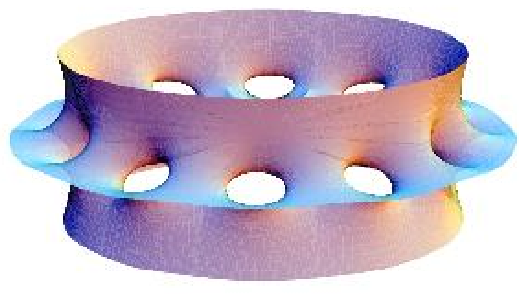}
\end{tabular}
\end{center}
\vspace{2in}
\caption{The Costa-Hoffman-Meeks minimal surface of genus $k$ for 
         $k=3$, $6$, $8$ and $10$.  Like the Costa surface, these surfaces have
         3 ends and are embedded.}
\label{fg:CHM3-10}
\end{figure}
\begin{figure}[phbt]
\begin{center}
\includegraphics[width=0.4\linewidth]{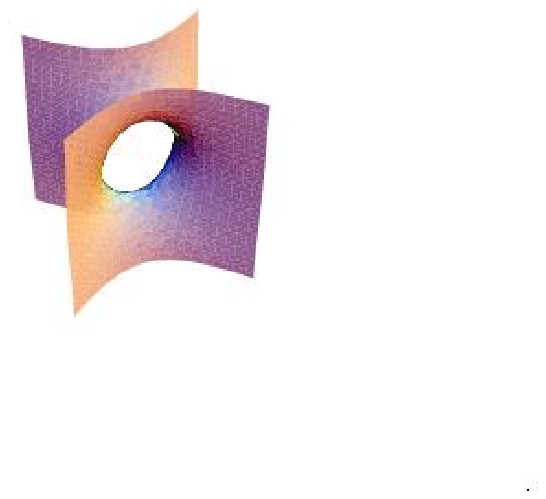}
\includegraphics[width=0.2\linewidth]{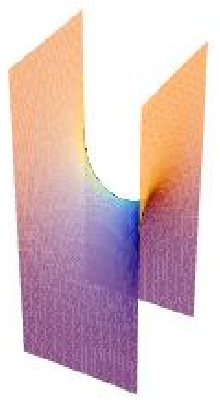}
\end{center}
\vspace{2.8in}
\caption{The singly-periodic Scherk surface (on the left) and the
doubly-periodic Scherk surface (on the right).
The singly-periodic Scherk surface is given by
the Weierstrass representation with $(g,\eta)=(z,ic(z^4-1)^{-1}dz)$ on
$\Sigma=(\mathbb{C}\cup\{\infty\})\setminus\{\pm 1,\pm i\}$,
where $c\in\mathbb{R}\setminus\{0\}$ is a constant.  If we replace $\eta$
by $c(z^4-1)^{-1}dz$, then we have the doubly-periodic Scherk surface.}
\label{fg:Scherk}
\end{figure}
\begin{figure}[phbt]
\begin{center}
\includegraphics[width=0.3\linewidth]{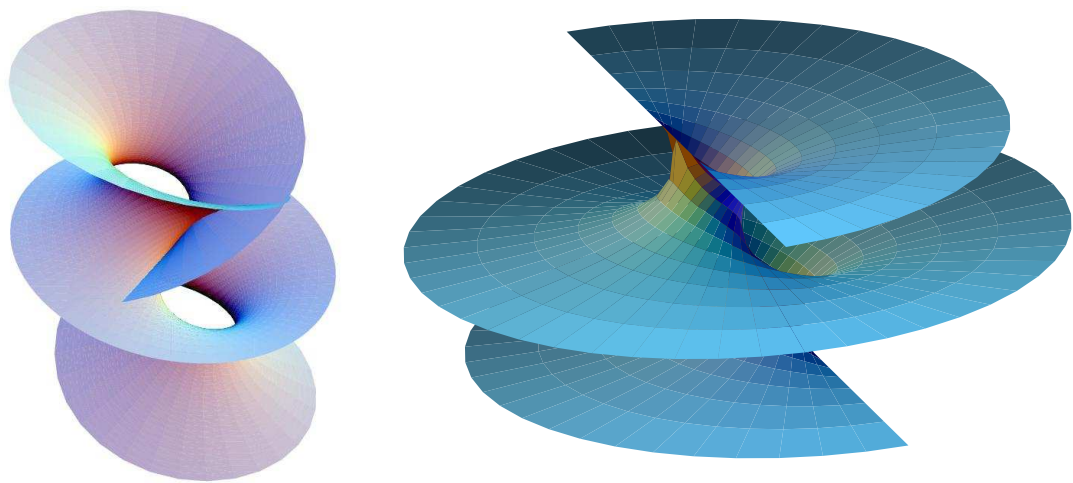}
\includegraphics[width=0.65\linewidth]{riemann_1200.eps}
\end{center}
\caption{Richmond's minimal surface (on the left) and the
singly-periodic Riemann's staircase minimal surface (on the right).
Richmond's surface is given by
the Weierstrass representation with $(g,\eta)=(z^2,cz^{-2}dz)$ on
$\Sigma=\mathbb{C}\setminus\{0\}$,
where $c\in\mathbb{R}\setminus\{0\}$ is a constant.
The singly-periodic Riemann's staircase is given by the
Weierstrass representation with $(g,\eta)=(z,dz/(zw))$ on the
hyperelliptic Riemann surface
$\Sigma=\{(z,w)\in(\mathbb{C}\cup\{\infty\})^2\;|\;w^2=z(z^2-1)\}
        \setminus\{(0,0),(\infty ,\infty )\}$.}
\label{fg:Richmond-Riemann}
\end{figure}
\begin{figure}[phbt]
\begin{center}
\includegraphics[width=0.45\linewidth]{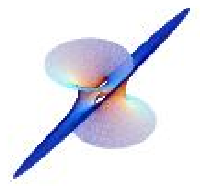}
\raisebox{5ex}{\includegraphics[width=0.45\linewidth]{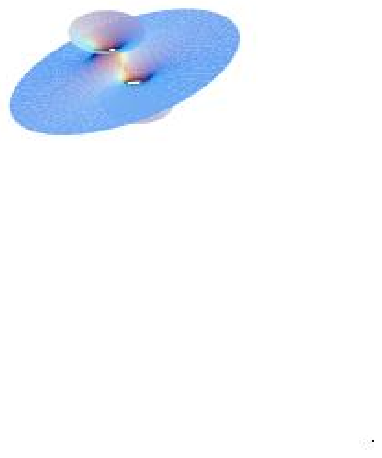}}
\end{center}
\vspace{2.6in}
\caption{Minimal Lopez-Ros surfaces.  These surfaces are given by 
the Weierstrass representation with 
$(g,\eta)=(\rho(z^2+3)/(z^2-1),\rho^{-1}dz)$ on 
$\Sigma=\mathbb{C}\setminus\{\pm 1\}$, 
where $\rho\in (0,\infty)$ is a constant.  
They each have two ends asymptotic to catenoid ends and one end asymptotic to 
a plane, just like the Costa and Costa-Hoffman-Meeks surfaces, 
but the surfaces here are not embedded.}
\label{fg:Lopez-Ros}
\end{figure}

\begin{figure}[phbt]
\begin{center}
\includegraphics[width=0.35\linewidth]{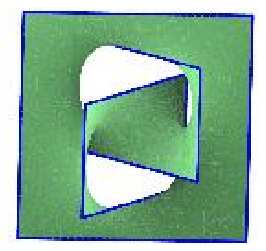}
\includegraphics[width=0.4\linewidth]{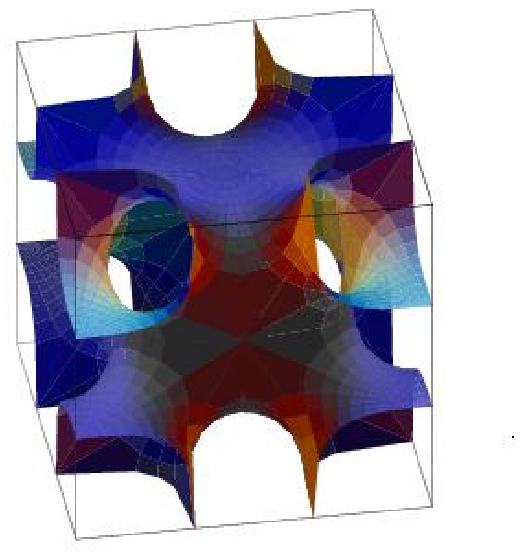}
\end{center}
\vspace{2.6in}
\caption{A triply-periodic embedded minimal surface with noncubic symmetry found by 
Fischer and Koch \cite{FischerKoch}, and A.~Schoen's triply-periodic I-Wp surface.  
The I-Wp surface is given by the Weierstrass representation with
$(g,\eta)=(z^3,dz/w^2)$ on the Riemann surface
$\Sigma=\{(z,w)\in(\mathbb{C}\cup\{\infty\})^2\;|\;w^3=z^{12}-1\}$.}
\label{fg:Fischer-Koch}
\end{figure}

\begin{figure}[phbt]
\begin{center}
\includegraphics[width=0.4\linewidth]{SchwarzD.eps}
\includegraphics[width=0.4\linewidth]{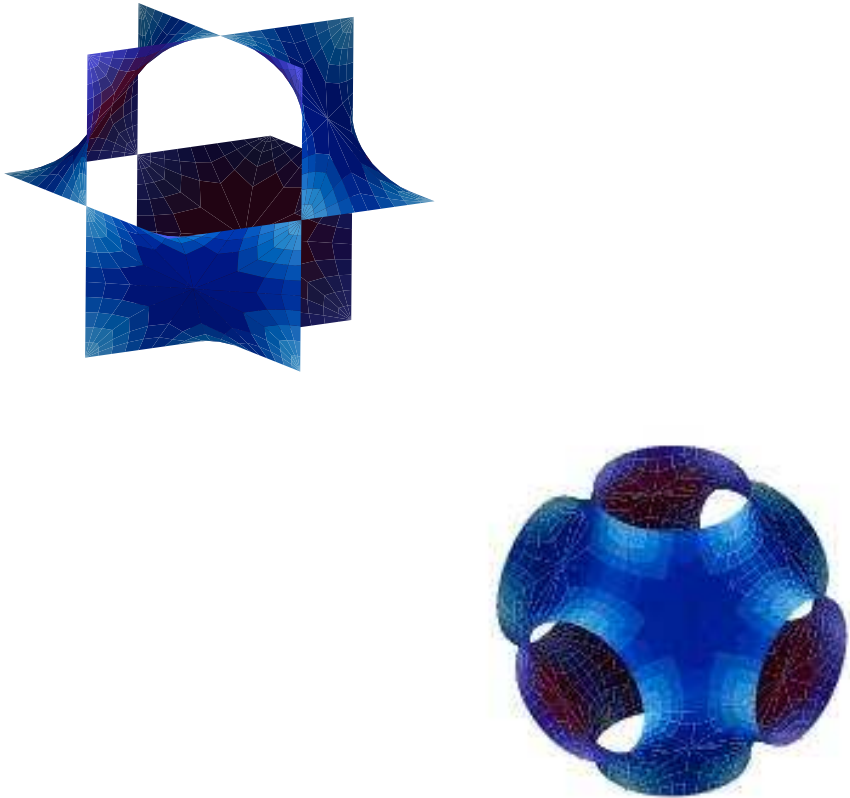}
\end{center}
\vspace{2in}
\caption{The triply-periodic Schwarz D surface (on the left) and Schwarz P
surface (on the right).
The Schwarz D surface is given by
the Weierstrass representation with $(g,\eta)=(z,dz/w)$ on the
hyperelliptic Riemann surface
$\Sigma=\{(z,w)\in(\mathbb{C}\cup\{\infty\})^2\;|\;w^2=z^8-14z^4+1\}$.
The Schwarz P surface is the conjugate surface of the Schwarz D surface, that is,
the Schwarz D surface is given by the Weierstrass data 
$(g,\eta)=(z,idz/w)$ on the same Riemann surface.}
\label{fg:SchwarzD-P}
\end{figure}

\section{Harmonicity of the Gauss map}
\label{harmonic}

Before discussing harmonicity of the Gauss map, we consider the meaning 
of harmonicity.  
Let $(M,g)$ and $(\mathcal{N},h)$ be Riemannian manifolds with $\dim M=m$ 
and $\dim \mathcal{N}=n$, and Riemannian metrics $g$ and $h$, respectively.  
Suppose that $M$ is compact, possibly with boundary.  Let 
\[ f : (M,g) \to (\mathcal{N},h) \] be a smooth map.  
The energy functional of $f$ is 
\[ E(f) := \frac{1}{2} \int_M |df|^2 dM \; , \]
where $df = (\partial_{x_1}f) dx_1 + \dots + (\partial_{x_m}f) dx_m$ 
for a local coordinate $(x_1,\dots ,x_m)$ on $M$.  
The volume element of $M$ with respect to $g$ is 
\[ dM=\sqrt{\det g}\,dx_1 \wedge \dots \wedge dx_m \; , \]
and the energy density of $f$ is
\[
|df|^2=\sum_{i,j=1}^m\sum_{\alpha ,\beta =1}^n g^{ij}h_{\alpha \beta}
\frac{\partial f_\alpha}{\partial x_i}\frac{\partial f_\beta}{\partial x_j}
\in C^\infty (M) \; , \]
where
\[
f(x_1,\dots x_m)=\left( f_1(x_1,\dots ,x_m),\dots ,f_n(x_1,\dots ,x_m) \right)
\]
in the local coordinates of $M$ and $\mathcal{N}$.  
For a bounded domain $U$ in $M$, we define a {\em smooth boundary-fixing 
variation} of the map $f$ to be a smooth map 
$\mathfrak{F}:(-\varepsilon ,\varepsilon )\times \bar U\to \mathcal{N}$ 
with these three properties: 

 (1) $f_t:=\mathfrak{F}(t,\cdot):\bar U\to \mathcal{N}$ is a smooth map 
for all $t\in (-\varepsilon ,\varepsilon )$, 

 (2) $f_0=f|_{\bar U}$ on $\bar U$, 

 (3) $f_t|_{\partial U}=f|_{\partial U}$ for all 
$t\in (-\varepsilon ,\varepsilon )$, where $\partial U=\bar U\setminus U$ is 
the boundary of $U$.  \\
We say that $f$ is {\em harmonic} if $f$ is an extremal value for the energy 
functional, that is, if for any smooth boundary-fixing 
variation $f_t$ of $f$ we have 
\[ \left. \frac{d}{dt} E(f_t) \right|_{t=0} = 0 \; . \]  
Let us now describe two simple examples:

\begin{example}
Let $M$ be an open bounded domain of $\mathbb{R}^2$ with the standard Euclidean 
coordinates $(x,y)$ and the standard metric.  
Let $\mathcal{N}=\mathbb{R}$ with the standard coordinate $x$ and standard metric.  
Consider a smooth map $f$ from $M$ to $\mathcal{N}$.  Then 
$f$ is harmonic if and only if $\triangle f = f_{xx}+f_{yy} = 0$, 
and this can be shown as follows:

The energy functional of $f$ is 
\[ E(f) := \frac{1}{2} \int_M |df|^2 dM = 
\frac{1}{2} \int_M (f_x^2+f_y^2) \,dxdy \; . \]  Let $f_t = f + 
t \cdot v + O(t^2)$, $v=v(x,y)$, $v|_{\partial M}=0$, be an arbitrary smooth 
boundary-fixing variation of $f$.  Then $(d/dt) E(f_t) |_{t=0} = 
0$ if and only if 
\begin{eqnarray*}
0 &=& \frac{1}{2} \int_M \partial_t 
((f_t)_x^2+(f_t)_y^2) |_{t=0} \, dxdy \\
  &=& \frac{1}{2} \int_M \partial_t (f_x^2+2 t f_x v_x + 
      f_y^2+2 t f_y v_y + O(t^2)) |_{t=0} \, dxdy \\
  &=& \frac{1}{2} \int_M (2 f_x v_x + 2 f_y v_y) \, dxdy = 
     -\int_M \triangle f \cdot v \, dxdy 
     +\int_M ((f_x v)_x+(f_y v)_y) \, dxdy \\
  &=&-\int_M \triangle f \cdot v \, dxdy + \int_{\partial M} f_x v \, dy 
     -\int_{\partial M} f_y v \, dx = - \int_M \triangle f \cdot v \, dxdy
\end{eqnarray*}
for all $v$, since $v$ is zero on $\partial M$ -- we used Green's theorem 
here.  And this holds if and only if $\triangle f = 0$.  
\end{example}

\begin{example}
Consider a smooth map $N$ from an open region $(\Sigma ,g)$ of a compact 
Riemann surface with Riemannian metric $g$ and 
with holomorphic coordinate $z=x+i y$ to an open domain of $\mathbb{S}^2$, 
where $\mathbb{S}^2$ has the metric induced from $\mathbb{R}^3$.  Then 
$N$ is harmonic if and only if $N_{xx}+N_{yy} = c \cdot N$ 
for some scalar function $c=c(x,y)$, and this can be shown as follows: 

Let $V:\Sigma\to\mathbb{R}^3$ be an arbitrary smooth vector field such that $V$ is 
always perpendicular to $N$ and $V|_{\partial \Sigma}\equiv 0$.  
An arbitrary smooth boundary-fixing variation of $N$ is of the form 
$N_t=N+t\cdot V +O(t^2)$.  
The metric $g$ on $\Sigma$ is conformal with respect to the coordinates 
$(x,y)$, that is, $g_{ij} = 4e^{2u} \delta_{ij}$.  
\[ |dN_t|^2 = \frac{1}{4}e^{-2u} (
\langle \partial_x N_t, \partial_x N_t \rangle + 
\langle \partial_y N_t, \partial_y N_t \rangle) \; , \] 
where $\langle\cdot ,\cdot\rangle$ is the standard inner product of 
$\mathbb{R}^3$, so 
\[ \partial_t |dN_t|^2 = \frac{1}{2}e^{-2u} (
\langle \partial_t \partial_x N_t, \partial_x N_t \rangle + 
\langle \partial_t \partial_y N_t, \partial_y N_t \rangle) \; . \] 
We have 
\begin{eqnarray*}
\frac{d}{dt} E(N_t)|_{t=0}
 &=& \frac{1}{2} \int_\Sigma \partial_t |dN_t|^2 |_{t=0} \, d\Sigma \\
 &=& \int_\Sigma \frac{e^{-2u}}{4} (\langle N_x, V_x \rangle
                           + \langle N_y, V_y \rangle) 4e^{2u} \, dxdy \\
 &=& \int_\Sigma (\langle N_x,V_x \rangle+ \langle N_y, V_y \rangle) \, dxdy \\
 &=&-\int_\Sigma (\langle N_{xx},V \rangle+ \langle N_{yy}, V \rangle) \, dxdy 
    +\int_\Sigma (\langle N_x,V\rangle_x + \langle N_y, V \rangle_y) \, dxdy \\
 &=&-\int_\Sigma \langle N_{xx}+N_{yy}, V \rangle \, dxdy \; , 
\end{eqnarray*}
because 
$\int_\Sigma (\langle N_x,V \rangle_x + \langle N_y, V \rangle_y) \, dxdy = 0$ 
by the same arguments as in the previous example.  So 
\[ \frac{d}{dt} E(N_t)|_{t=0} = 0 \]
for all smooth boundary-fixing variations $N_t$ if and only if $N_{xx}+N_{yy}$ 
is perpendicular to $V$ for all $V$ tangent to $\mathbb{S}^2$, 
and this holds if and only if $N_{xx}+N_{yy}$ is parallel to $N$.  
\end{example}

So we conclude that the following definition is correct:

\begin{defn}\label{df:harmonic}{\rm
$N:\Sigma \to \mathbb{S}^2$ is {\em harmonic} if $N_{z\bar z}$ is parallel to $N$.  
}\end{defn}

With Definition \ref{df:harmonic} in mind, we now 
discuss the harmonicity of the Gauss map of a surface.  
Consider a Riemann surface $\Sigma$ with holomorphic coordinate $z$ and a 
conformal immersion 
\[ f=f(z,\bar z) : \Sigma \to \mathbb{R}^3 \] with unit normal 
$N: \Sigma \to \mathbb{S}^2$.  
We have a function $u:\Sigma\to\mathbb{R}$ so that \eqref{eq:metric} holds, 
and $Q$ and $H$ are as defined in Section \ref{section6}.  
Then, since $(f_{zz})_{\bar z}=(f_{z\bar z})_z$ and 
$(N_z)_{\bar z}=(N_{\bar z})_z$, computations as in 
Section \ref{section6} (without assuming $H$ is constant here) imply
\begin{eqnarray}
 & Q_{\bar z} = 2 H_z e^{2 u} \; , &
\label{harm1} \\
 & \langle N_{z\bar z},N\rangle =-  \frac{1}{2} Q \bar Q e^{-2 u}
                                 -2 H^2 e^{2 u} \; , &
\label{harm2} \\
& N_{z\bar z}=- H_{\bar z} f_z- H_z f_{\bar z}
              + \langle N_{z\bar z},N\rangle  N \; . &
\label{harm3}
\end{eqnarray}
Equation \eqref{harm1} is the Codazzi equation, and it 
tells us that $H$ is constant if and only if $Q$ is holomorphic.  
Furthermore, Equation \eqref{harm3} tells us that $N$ is harmonic if and only 
if $H$ is constant.  

We note that we can make an extended $3\times 3$ frame $\mathcal{F}$ for the 
associated family $f(z,\bar z,\lambda)$ of $f(z,\bar z)$, 
as in Section \ref{section6}, so that the Maurer-Cartan 1-form is 
\[
\mathcal{A}:= \mathcal{F}^{-1}d\mathcal{F}
            = \mathcal{F}^{-1}\mathcal{F}_z dz
             +\mathcal{F}^{-1} \mathcal{F}_{\bar z}d\bar z
            = (\Theta +\Upsilon _z)dz 
             +(\bar{\Theta}-\Upsilon_{\bar{z}})d\bar z \; , 
\]
with the compatibility condition \eqref{eq:3by3MCeqn} holding for all 
$\lambda\in \mathbb{S}^1$.  
This compatibility condition can be equivalently written as 
\[ d\mathcal{A}+\mathcal{A}\wedge\mathcal{A} = 0 \;\;\; 
  \text{for all}\;\;\; \lambda \in \mathbb{S}^1 \; . \] 

Similarly, $f(z,\bar z,\lambda)$ has an extended $2\times 2$ frame $F$, 
as in Section \ref{laxpairintwisted}, 
so that its Maurer-Cartan 1-form 
$ \alpha := F^{-1} d F = U dz + V d\bar z $ satisfies 
\[ d\alpha + \alpha \wedge \alpha = 0 \] for all $\lambda \in \mathbb{S}^1$.  
Here $U$ and $V$ are as in Section \ref{laxpairintwisted}.  

If $H$ is not constant, then $Q_{\bar z} \neq 0$, so the Codazzi equation will 
no longer hold if $Q$ is replaced by $\lambda^{-2} Q$, and the associated 
family $f(z, \bar z, \lambda)$ does not exist.  
Hence the extended frame $\mathcal{F}$, or $F$, cannot be made.  
Thus, we have proven:  

\begin{proposition}\label{prop:harmonic}
Let $f:\Sigma \to \mathbb{R}^3$ be a conformal immersion with metric 
$4 e^{2 u} dzd\bar z$ and Hopf differential $Q$ and mean curvature $H$ 
and unit normal $N:\Sigma \to \mathbb{S}^2$.  
Then the following are equivalent: 

 $(1)$ $H$ is constant.  

 $(2)$ $Q$ is holomorphic.  

 $(3)$ $N$ is harmonic.  

 $(4)$ There exists an extended $3\times 3$ frame $\mathcal{F}$, 
and the Maurer-Cartan $1$-form $\mathcal{A}=\mathcal{F}^{-1}d\mathcal{F}$ 
satisfies $d\mathcal{A}+\mathcal{A}\wedge\mathcal{A}=0$ 
for all $\lambda \in \mathbb{S}^1$.  

 $(5)$ There exists an extended $2\times 2$ frame $F$, and the Maurer-Cartan 
$1$-form $\alpha=F^{-1} d F$ satisfies $d\alpha + \alpha \wedge 
\alpha = 0$ for all $\lambda \in \mathbb{S}^1$.  
\end{proposition}

That (1) implies (3) in Proposition \ref{prop:harmonic} is known as a 
theorem by Ruh-Vilms \cite{rv}.  

We note that 
$\mathcal{F}$, and also $F$, can be viewed either as the extended frame of $N$ 
or the extended frame of $f$, depending on whether your primary interest is in 
harmonic maps or in CMC surfaces.  

\chapter{Theory for the DPW method}
\label{chapter3}

In this chapter we provide the technical materials needed to prove the DPW 
results (loop groups, form of the potentials, Iwasawa and Birkhoff 
decompositions).  We also prove those results of DPW, and we describe 
what dressings and gaugings are.  

Our naming of the Iwasawa decomposition here is historically somewhat inaccurate.  
Iwasawa did actually consider the splitting described in Section \ref{section9}, 
which we do not call {\em Iwasawa splitting} here, and he did not consider the 
infinite dimensional splitting in Theorem \ref{th:iwasawa-sp}, which we do 
actually call {\em Iwasawa splitting} here.  However, this misnomer is 
now commonly used amongst practitioners of the DPW method, so we leave it this 
way here.  

\section{Gram-Schmidt orthogonalization}
\label{section9}

Iwasawa splitting (also called Iwasawa decomposition) is done on infinite 
dimensional spaces.  
So, to explain Iwasawa splitting, we must first define the infinite 
dimensional spaces we need.  
However, before doing this in Section \ref{section15}, 
we briefly describe a finite-dimensional analog of Iwasawa splitting, 
to provide some intuition for those 
who have never seen Iwasawa splitting before.  

Let us consider a decomposition of $\SL_2\!\mathbb{C}$, 
that is, the decomposition $\SL_2\!\mathbb{C}=\SU_2 \times \Delta_2$, 
where $\Delta_2$ is the set of upper triangular 
$2 \times 2$ matrices with positive real numbers on the diagonal.  
This decomposition of an arbitrary matrix in $\SL_2\!\mathbb{C}$ is 
unique and is 
\[ \begin{pmatrix} a & b \\ c & d 
\end{pmatrix} = \begin{pmatrix} \dfrac{a}{\sqrt{a\bar a+c\bar c}} & 
\dfrac{-\bar c}{\sqrt{a\bar a+c\bar c}} \\ \dfrac{c}{\sqrt{a\bar a+c\bar c}} & 
\dfrac{\bar a}{\sqrt{a\bar a+c\bar c}} 
\end{pmatrix} \begin{pmatrix} \sqrt{a\bar a+c\bar c} & 
\dfrac{\bar a b+\bar c d}{\sqrt{a\bar a+c\bar c}} \\ 0 & 
\dfrac{1}{\sqrt{a\bar a+c\bar c}} 
\end{pmatrix} \; , \] where we have used that $a d - b c = 1$.  
Also, $a\bar a+c\bar c>0$, so this splitting is always defined.  

Note that Gram-Schmidt orthogonalization replaces an arbitrary basis of 
$\mathbb{C}^2$ with an orthonormal basis via the following map: 
\[ \begin{pmatrix}a\\c\end{pmatrix} \; , \; 
   \begin{pmatrix}b\\d\end{pmatrix} \; \to \; 
   \begin{pmatrix}\dfrac{a}{\sqrt{a\bar a+c\bar c}}\\
                  \dfrac{c}{\sqrt{a\bar a+c\bar c}}\end{pmatrix} \; , \; 
   \begin{pmatrix}\dfrac{-\bar c}{\sqrt{a\bar a+c\bar c}}\\
                  \dfrac{ \bar a}{\sqrt{a\bar a+c\bar c}}\end{pmatrix} \; . 
\]
These resulting two vectors are precisely the columns of the $\SU_2$ 
matrix in the splitting of the $\SL_2\!\mathbb{C}$ matrix above.  
So the splitting above is equivalent to Gram-Schmidt orthogonalization.  

Iwasawa decomposition is similar in spirit to the above splitting 
$\SL_2\!\mathbb{C}=\SU_2 \times \Delta_2$, but is done on infinite-dimensional 
spaces called loop groups, and we consider this in the next section.

We will also consider Birkhoff splitting on infinite-dimensional loop groups 
in the next section.  Again for the sake of intuition, let us briefly look at 
the following finite-dimensional analog of Birkhoff splitting:

If the matrix (holomorphic in $z$) 
\[ \begin{pmatrix} a & b \\ c & d 
\end{pmatrix} = \begin{pmatrix} a(z) & b(z) \\ c(z) & d(z) 
\end{pmatrix}
\in \SL_2\!\mathbb{C}\;\;\text{for all}\;\; z\in\mathbb{C}\; , \;\; |z|<1 \; , 
\] and $a$ has only isolated zeros, then one Birkhoff splitting is 
\[ \begin{pmatrix} a & b \\ c & d 
\end{pmatrix} = \begin{pmatrix} 1 & 0 \\ c/a & 1 
\end{pmatrix} \cdot \begin{pmatrix} a & b \\ 0 & 1/a 
\end{pmatrix} \; . 
\]
Unlike the finite-dimensional analog of Iwasawa splitting above, 
this splitting has the advantage that the two component matrices are 
holomorphic in $z$ (at points where $a\ne 0$) when the original matrix is 
holomorphic in $z$.  However, it has the weakness of not being globally 
defined (that is, it is not defined at points where $a=0$), unlike the Iwasawa 
case.  So if the original matrix is holomorphic in $z$ (and $a$ is not 
identically zero), we can conclude only 
that the two component matrices of the splitting are meromorphic.  

\section{Loop groups and algebras, Iwasawa and Birkhoff decomposition}
\label{section15}

We first define the loop groups that we will need, and their loop algebras.  
The matrix $\sigma_3$ is as in Section \ref{2by2laxpair}, 
and $\mathbb{S}^1$ and $D_1$ are as in Section \ref{section0}.  
\begin{defn}\label{df:loop}{\rm 
\begin{eqnarray*}
\Lambda \SL_2\!\mathbb{C} &=& 
\left\{ \left. \phi :\mathbb{S}^1 \stackrel{C^\infty}{\longrightarrow} 
\SL_2\!\mathbb{C} \; \right| \; 
\phi (-\lambda) = \sigma_3 \phi (\lambda) \sigma_3 \right\} \; , \\
\Lambda \slg_2\!\mathbb{C} &=& 
\left\{ \left. 
A:\mathbb{S}^1 \stackrel{C^\infty}{\longrightarrow} \slg_2\!\mathbb{C} 
\; \right| \; A(-\lambda) = \sigma_3 A(\lambda) \sigma_3 \right\} \; , \\ 
\Lambda \SU_2 &=& 
\{ F\in\Lambda \SL_2\!\mathbb{C} \; | \; F :\mathbb{S}^1 \to \SU_2 \} \; , \\
\Lambda \su_2 &=& 
\{ A\in\Lambda \slg_2\!\mathbb{C} \; | \; A:\mathbb{S}^1 \to \su_2 \} \; , \\
\Lambda_+ \SL_2\!\mathbb{C} &=& 
\left\{ B \in \Lambda \SL_2\!\mathbb{C} \; | \; 
\text{$B$ extends holomorphically to $D_1$}\right\} \; , \\
\Lambda_+^\mathbb{R} \SL_2\!\mathbb{C} &=& 
\left\{ B \in \Lambda_+ \SL_2\!\mathbb{C} \; \left| \;
B|_{\lambda =0}=\begin{pmatrix}\rho & 0 \\ 0 & \rho^{-1}\end{pmatrix}
\text{ with $\rho>0$}\right.\!\right\} \; , \\
\Lambda_-^* \SL_2\!\mathbb{C} &=& 
\left\{ B \in \Lambda \SL_2\!\mathbb{C} \; \left| \!
\begin{array}{c}
\text{$B$ extends holomorphically to} \\
\text{$\mathbb{C}\cup\{\infty\}\setminus\overline{D_1}$ and $B|_{\lambda =\infty}=\id$}
\end{array}\right.\!\!\!\!\right\} \; .
\end{eqnarray*}
}\end{defn}

The conditions 
\[
\phi (-\lambda )=\sigma_3\phi (\lambda )\sigma_3 \; , \;\;\;
A(-\lambda )=\sigma_3A(\lambda )\sigma_3
\]
are the reason that these loop groups and algebras are referred to as 
``twisted''.  These conditions imply that 
$\phi (\lambda )\in\Lambda\SL_2\!\mathbb{C}$ and 
$A(\lambda )\in\Lambda\slg_2\!\mathbb{C}$ are even functions of $\lambda$ on 
their diagonals and odd functions on their off-diagonals.  

Note that $B \in \Lambda_+ \SL_2\!\mathbb{C}$ satisfies the following extension 
property: not only is $B$ holomorphic on the open disk $D_1$, it extends continuously 
to a matrix-valued function defined on the closed disk $\overline{D_1}$.  

$\Lambda\SL_2\!\mathbb{C}$ is the infinite-dimensional analog of 
$\SL_2\!\mathbb{C}$, and is the object that will be split by Iwasawa 
decomposition, analogous to the way $\SL_2\!\mathbb{C}$ was split in 
Section \ref{section9}.  
$\Lambda\SU_2$ is the infinite-dimensional analog to the $\SU_2$ part of the 
splitting in Section \ref{section9}.  
$\Lambda_+^\mathbb{R} \SL_2\!\mathbb{C}$ is the infinite-dimensional analog of 
$\Delta_2$, and can be thought of as an infinite-dimensional version of being 
upper triangular (that is, the power series expansion of 
$B\in\Lambda_+^\mathbb{R} \SL_2\!\mathbb{C}$ at $\lambda =0$ has no negative 
powers of $\lambda$).  

We gave the definitions of $\Lambda_+ \SL_2\!\mathbb{C}$ and 
$\Lambda_-^* \SL_2\!\mathbb{C}$ because they are required for stating the 
Birkhoff decomposition.  
We will be using the Birkhoff splitting as well.  

The definition for $\Lambda\SU_2$ is commonly stated in this equivalent form: 
$\Lambda\SU_2$ is the set of $F\in\Lambda\SL_2\!\mathbb{C}$ such that 
\begin{equation}
\overline{F(1/\bar \lambda)}^t=\left(F(\lambda)\right)^{-1}
\;\;\;\text{for all}\;\;\lambda\in\mathbb{S}^1 \; . 
\label{eq:realitycondition}
\end{equation}
This relation \eqref{eq:realitycondition} is called 
the {\it reality condition}, and is equivalent to saying that $F\in\SU_2$ for 
all $\lambda\in\mathbb{S}^1$.  
Of course, $1/\bar \lambda =\lambda$ when $\lambda\in\mathbb{S}^1$, 
but in \eqref{eq:realitycondition} we cannot replace $1/\bar \lambda$ by 
$\lambda$, since there will be occasions when we are not 
restricting $\lambda$ to $\mathbb{S}^1$, as we will see later.  

Before stating the Iwasawa and Birkhoff decompositions, we state a theorem 
about the relations between the above loop groups and loop algebras.  
This result is basic and central to Lie group theory, but we state it here, 
to emphasize that it still holds even though the groups and algebras 
under consideration are not finite dimensional.  

\begin{theorem}\label{th:loopgroup}
Suppose that $\phi: \Sigma \times \mathbb{S}^1 \to M_{2 \times 2}$ and 
$F: \Sigma \times \mathbb{S}^1 \to M_{2 \times 2}$ are matrix-valued functions 
depending on a complex coordinate $z$ of a Riemann surface $\Sigma$ and the 
parameter $\lambda \in \mathbb{S}^1$; 
that is, $\phi=\phi (z,\bar z,\lambda)=\phi (x,y,\lambda)$ and 
$F=F(z,\bar z,\lambda)=F(x,y,\lambda)$, where $x=\Re (z)$ and $y=\Im (z)$.  
Then $\phi \in \Lambda \SL_2\!\mathbb{C}$ (resp. $F\in\Lambda \SU_2$) for all 
$z=x+iy$ if and only if 

$(1)$ $\phi^{-1} \phi_x$, $\phi^{-1} \phi_y \in \Lambda \slg_2\!\mathbb{C}$ 
(resp. $F^{-1}F_x, \;\; F^{-1}F_y \in \Lambda \su_2$) at all points $z$, and 

$(2)$ there exists a fixed $z_*$ such that $\phi(z_*,\overline{z_*},\lambda) 
\in \Lambda \SL_2\!\mathbb{C}$ (resp. $F(z_*,\overline{z_*},\lambda) 
\in \Lambda \SU_2$).  
\end{theorem}

\begin{proof}
First we will show that 
{\it $F \in \SU_2$ for all $\lambda \in \mathbb{S}^1$ if and only if 
\[ F^{-1} F_x, F^{-1} F_y \in \su_2 \] and there exists a $z_*$ such 
that $F(z_*,\overline{z_*},\lambda) 
\in \SU_2$ for all $\lambda \in \mathbb{S}^1$}.  

$F \in \SU_2$ for all $\lambda \in \mathbb{S}^1$ means that $F\bar F^t=\id$ 
and $\det F=1$ for all $\lambda \in \mathbb{S}^1$.  
Equivalently, we can say that $F\bar F^t$ and $\det F$ are constant in 
$\lambda$ and there exists a fixed $z_*$ such that 
$F(z_*,\overline{z_*},\lambda)\in\SU_2$ for all $\lambda \in \mathbb{S}^1$, 
that is, condition (2) of the theorem holds for $F$.  
(Any choice of $z_*$ will suffice.) 
Again, equivalently, 
\[ (F \bar F^t)_x=(F \bar F^t)_y = 0 \; , \;\;\; \partial_x 
\log (\det F) = \partial_y \log (\det F) = 0 \; , \] 
and condition (2) holds for $F$.  
Note that $(F \bar F^t)_x=(F \bar F^t)_y = 0$ if and only if 
$F^{-1} F_x+\overline{F^{-1} F_x}^t=
F^{-1} F_y+\overline{F^{-1} F_y}^t= 0$.  Hence 
$F \in \SU_2$ for all $\lambda \in \mathbb{S}^1$ if and only if 
\[ F^{-1} F_x = \begin{pmatrix} i r_1 & w_1 \\ 
-\overline{w_1} & i s_1 \end{pmatrix} \; , \;\;\; 
F^{-1} F_y = \begin{pmatrix} i r_2 & w_2 \\ 
-\overline{w_2} & i s_2 \end{pmatrix} \; , \;\; \text{for some} \;\;
r_j,\;s_j \in \mathbb{R} \; , \; w_j \in \mathbb{C} \] and 
$\text{trace}(F^{-1}F_x)=\text{trace}(F_xF^{-1})
=\text{trace}(F^{-1}F_y)=\text{trace}(F_yF^{-1})=0$ 
(by Equation \eqref{eq:matrixlemma}) and condition (2) holds for $F$.  
Finally, this is equivalent to saying that 
$F^{-1} F_x , F^{-1} F_y \in \su_2$ and condition (2) holds for $F$.  

A shorter version of this argument shows that {\it $\phi \in \SL_2\!\mathbb{C}$ 
for all $\lambda \in \mathbb{S}^1$ if and only if $\phi^{-1} \phi_x 
\in \slg_2\!\mathbb{C}$ and 
$\phi^{-1} \phi_y \in \slg_2\!\mathbb{C}$ and there exists a $z_*$ such that 
$\phi(z_*,\overline{z_*},\lambda) \in \SL_2\!\mathbb{C}$ for all 
$\lambda \in \mathbb{S}^1$} (that is, condition (2) holds for $\phi$).  

It remains only to confirm the twistedness property to complete the proof 
of the theorem.  We wish to show the following: 
{\it Assume $\phi \in \SL_2\!\mathbb{C}$ for all $\lambda \in \mathbb{S}^1$.  
Then $$\sigma_3 \phi(z,\bar z,-\lambda) \sigma_3 = \phi(z,\bar z,\lambda)$$ 
for all $\lambda \in \mathbb{S}^1$ if and only if, 
for $\bmath{U} := \phi^{-1} \phi_x$ and $\bmath{V} := \phi^{-1} \phi_y$, 
$\sigma_3 \bmath{U}(z,\bar z,-\lambda) \sigma_3 = \bmath{U}(z,\bar z,\lambda)$ and 
$\sigma_3 \bmath{V}(z,\bar z,-\lambda) \sigma_3 = \bmath{V}(z,\bar z,\lambda)$ and 
condition (2) holds for $\phi$}.  

One direction is obvious.  To prove the other direction, note that 
\[ (\phi(z,\bar z,\lambda))_x=\phi(z,\bar z,\lambda) \cdot 
\bmath{U}(z,\bar z,\lambda)=\phi(z,\bar z,\lambda) \sigma_3 \bmath{U}(z,\bar z,-\lambda) 
\sigma_3 \;\; \text{for all}\;\; \lambda \in \mathbb{S}^1 \] 
implies \[ (\sigma_3 \phi(z,\bar z,\lambda) \sigma_3)_x=
(\sigma_3 \phi(z,\bar z,\lambda) \sigma_3) \cdot 
\bmath{U}(z,\bar z,-\lambda) \;\; \text{for all}\;\; \lambda \in \mathbb{S}^1 \; . 
\]  Similarly, 
\[ (\sigma_3 \phi(z,\bar z,\lambda) \sigma_3)_y=
(\sigma_3 \phi(z,\bar z,\lambda) \sigma_3) \cdot 
\bmath{V}(z,\bar z,-\lambda) \;\; \text{for all}\;\; \lambda \in \mathbb{S}^1 \; . 
\]  Also, for all $\lambda \in \mathbb{S}^1$ we have 
\[ (\phi(z,\bar z,-\lambda))_x=\phi(z,\bar z,-\lambda) 
\bmath{U}(z,\bar z,-\lambda) \; , \;\; (\phi(z,\bar z,-\lambda))_y=
\phi(z,\bar z,-\lambda) \bmath{V}(z,\bar z,-\lambda) \; , 
\] and at the initial point $z_*$, 
\[ \sigma_3 \phi(z_*,\overline{z_*},\lambda) \sigma_3 
 = \phi(z_*,\overline{z_*},-\lambda) 
\;\; \text{for all}\;\; \lambda \in \mathbb{S}^1 \; , \]  by condition (2).  
So uniqueness of the solution with respect to the initial condition 
implies that
\[ \sigma_3 \phi(z,\bar z,-\lambda) \sigma_3 = \phi(z,\bar z,\lambda) 
\;\; \text{for all} \;\; \lambda \in \mathbb{S}^1 \; . \]
Since $\Lambda\SU_2$ is a subgroup of $\Lambda\SL_2\!\mathbb{C}$, 
the corresponding statement about twistedness for $F\in\Lambda\SU_2$ must also 
hold.  
\end{proof}

Now we come to the splitting theorems:

\begin{theorem}\label{th:iwasawa-sp}
{\rm (Iwasawa decomposition)} 
The multiplication $\Lambda\SU_2\times\Lambda_+^\mathbb{R}\SL_2\!\mathbb{C}$ to 
$\Lambda\SL_2\!\mathbb{C}$ is a real-analytic bijective diffeomorphism 
(with respect to the natural smooth manifold structure, 
as in Chapter 3 of \cite{prse}).  
The unique splitting of an element $\phi\in\Lambda\SL_2\!\mathbb{C}$ 
\begin{equation}\label{eq:iwasawa-sp} 
\phi = FB \; , 
\end{equation}
with $F\in\Lambda\SU_2$ and $B\in\Lambda_+^\mathbb{R}\SL_2\!\mathbb{C}$, will 
be called {\em Iwasawa decomposition}.
Because the diffeomorphism is real-analytic, if $\phi$ depends 
real-analytically (resp. smoothly) on some parameter $z$, then $F$ and $B$ do 
as well.  
\end{theorem}

\begin{theorem}\label{th:birkhoff-sp}
{\rm (Birkhoff decomposition)} 
Multiplication $\Lambda_-^*\SL_2\!\mathbb{C}\times\Lambda_+\SL_2\!\mathbb{C}$ 
is a complex-analytic bijective diffeomorphism to an open dense subset 
$\mathcal{U}$ of $\Lambda\SL_2\!\mathbb{C}$.  
The unique splitting of an element $\phi\in\mathcal{U}$ 
\begin{equation}\label{eq:birkhoff-sp} 
\phi = B_-B_+ \; , 
\end{equation}
with $B_-\in\Lambda_-^*\SL_2\!\mathbb{C}$ and 
$B_+\in\Lambda_+\SL_2\!\mathbb{C}$, will be called {\em Birkhoff decomposition}.
Because the diffeomorphism is complex-analytic, if $\phi$ depends 
complex-analytically (resp. real-analytically, smoothly) on some parameter 
$z$, then $B_-$ and $B_+$ do as well.  
\end{theorem}

Here we point out two differences between these two theorems.

(1) When $\phi =\phi (z,\lambda)$ depends complex-analytically on $z$, 
it does not follow that $F$ and $B$ in \eqref{eq:iwasawa-sp} will as well, 
and any of the explicit Iwasawa splittings in Sections \ref{section1}, 
\ref{section2}, \ref{section5} and \ref{section3} provide counterexamples.  
However, the $B_-$ and $B_+$ in \eqref{eq:birkhoff-sp} will depend 
complex-analytically on $z$, when $\phi$ does.  

(2) Iwasawa decomposition can be applied to {\em any} 
$\phi\in\Lambda\SL_2\!\mathbb{C}$.  However, Birkhoff decomposition can only 
be applied to those $\phi$ lying in $\mathcal{U}$.  
$\mathcal{U}$ is referred to as the {\em big cell}.  

The standard place to find a proof of Iwasawa splitting is 
\cite{prse}.  There is also a proof in \cite{g}, using flag manifolds.  

\begin{remark}\label{re:rIwasawa}
The condition on $\Lambda_+^{\mathbb{R}} \SL_2\!\mathbb{C}$ that 
$\rho \in \mathbb{R}^+$ makes the Iwasawa splitting $\phi =F\cdot B$ unique.  
But actually one could more generally allow that $\rho \in \mathbb{C}$ is not 
real, that is, that $B$ lies only in $\Lambda_+ \SL_2\!\mathbb{C}$.  
Then the splitting is no longer unique -- let us choose one such splitting 
$\phi = \tilde F \cdot \tilde B$, where 
\[
\tilde B|_{\lambda =0}
=\begin{pmatrix} \tilde\rho & 0 \\ 0 & \tilde\rho^{-1} \end{pmatrix} \; , 
\] 
for some $\tilde\rho = r e^{i \theta} \in \mathbb{C}$ with 
$\theta\in\mathbb{R}$ and $r\in\mathbb{R}\setminus\{0\}$.  
However, it follows that 
\[ B = \begin{pmatrix}
e^{-i \theta} & 0 \\
0 & e^{i \theta}
\end{pmatrix} \tilde B \; , \;\;\; 
F = \tilde F \begin{pmatrix}
e^{i \theta} & 0 \\
0 & e^{-i \theta}
\end{pmatrix} \; , \] since the matrix 
\[
\begin{pmatrix}e^{i \theta} & 0 \\ 0 & e^{-i \theta}\end{pmatrix}
\]
is in $\SU_2$.  Thus $F$ and $\tilde F$ are not the same, 
but this will not create any problems, as using either one in 
Equation \eqref{Sym-Bobenko} results in the same immersion $f$.  
\end{remark}

\begin{remark}\label{twisted}
Via the transformation 
\begin{equation}\label{eq:twisted}\begin{pmatrix}
a(\lambda) & b(\lambda) \\ c(\lambda) & d(\lambda) 
\end{pmatrix} \to \begin{pmatrix}
a(\sqrt{\lambda}) & \sqrt{\lambda}^{-1} b(\sqrt{\lambda}) \\ 
\sqrt{\lambda} c(\sqrt{\lambda}) & d(\sqrt{\lambda}) 
\end{pmatrix} \; , 
\end{equation}
one can shift to an equivalent equation where the properties 
$\phi (-\lambda )=\sigma_3\phi (\lambda )\sigma_3$ and
$A(-\lambda )=\sigma_3A(\lambda )\sigma_3$ are no longer required in 
Definition \ref{df:loop}.  
This changes some other parts of Definition \ref{df:loop} as well, 
and is a shift from the ``twisted'' setting to the ``untwisted'' setting.  
The inverse shift is 
\[ \begin{pmatrix}
a(\lambda) & b(\lambda) \\ c(\lambda) & d(\lambda) 
\end{pmatrix} \to \begin{pmatrix}
a(\lambda^2) & \lambda b(\lambda^2) \\ 
\lambda^{-1} c(\lambda^2) & d(\lambda^2) 
\end{pmatrix} \; .  
\] 
For example, the coefficient matrix in the potential of a Delaunay surface 
would change by \eqref{eq:twisted} as follows: 
\[ \begin{pmatrix}
c & a \lambda^{-1} + \bar b \lambda \\ b \lambda^{-1} + \bar a \lambda & -c 
\end{pmatrix} \to \begin{pmatrix}
c & a \lambda^{-1} + \bar b \\ b + \bar a \lambda & -c 
\end{pmatrix} \; .  
\] 
This shift is the same mapping for all of $\phi$, $\xi$, $F$ and $B$ in the 
DPW recipe.  There is also a corresponding change in the Sym-Bobenko formula.  

We will not use the untwisted form in these notes.  We mention it here only 
because it does appear in some papers related to the DPW method, such as in 
\cite{KilKRS}, \cite{KSS} and \cite{s}.  
\end{remark}

\section{Proof of one direction of DPW, with holomorphic potentials}
\label{section7}

In Chapter \ref{chapter2} we saw that finding CMC $H$ surfaces is 
equivalent to finding integrable Lax pairs of the form 
\eqref{LaxPair}-\eqref{eq:UandV}, and 
then the surface is found by using the Sym-Bobenko formula 
\eqref{Sym-Bobenko}.  So to prove that the DPW recipe finds all CMC 
$H$ surfaces, we want to prove that the DPW recipe produces all integrable 
Lax pairs of the form \eqref{LaxPair}-\eqref{eq:UandV} and all their solutions $F$.  
(By Theorem \ref{conformalityispossible}, without loss of generality we may assume 
that the CMC $H$ surfaces are conformally parametrized.)  We consider only 
simply-connected domains $\Sigma$ here, saving global considerations for 
Chapter \ref{chapter4}.  

In this section we will show that for a given holomorphic potential $\xi$ 
and solution $\phi$ to $d\phi = \phi \xi$, a solution to a Lax pair of the 
form \eqref{LaxPair}-\eqref{eq:UandV} is produced.  
In Section \ref{section8} will we show the converse: 
that any CMC $H$ surface comes from some potential $\xi$ and solution $\phi$.  
So the full proof of DPW for CMC surfaces is in this section and Section 
\ref{section8}.  
First we must define which kinds of potentials $\xi$ we will use: 

\begin{defn}\label{df:potential}
The {\em holomorphic potentials} $\xi$ we will use in the DPW recipe are of the form 
\[ \xi =Adz \; , \;\;\; 
A=A(z,\lambda)=\sum_{j=-1}^\infty A_{j}(z)\lambda^{j} \; , \] 
where $z$ is a complex coordinate of a simply-connected Riemann surface 
$\Sigma$, and where the 
$A_j=A_j(z)$ are $2 \times 2$ matrices that depend on $z$ and not on 
$\lambda$.  The matrix $A=A(z,\lambda)$ is holomorphic in both $z\in\Sigma$ 
and $\lambda\in\mathbb{C}\setminus\{0\}$ and lies 
in the loop algebra $\Lambda\slg_2\!\mathbb{C}$.  Since 
$A\in\Lambda \slg_2\!\mathbb{C}$, $A_j$ is off-diagonal (resp. diagonal) 
when $j$ is odd (resp. even), and $\text{trace}(A_j)=0$ for all $j$.  
Finally, we require that the upper-right entry of $A_{-1}(z)$ is never zero on 
$\Sigma$.  
\end{defn}

\begin{remark}\label{re:332}
The upper-right entry of $A_{-1}(z)$ being nonzero is needed to make the 
resulting CMC surface an immersion (that is, free of branch points), 
as we will see in Lemma \ref{lm:ProofofDPW1}.  
\end{remark}

Now we solve 
\[ d\phi = \phi \xi \; , \;\;\; \phi(z_*)=\id \] for some base point $z_*$.  
It follows that $\phi$ is holomorphic, 
and so $\phi^{-1}\phi_{\bar z}$ is identically zero.  Then, since also 
$\text{trace}(\phi^{-1}\phi_z)=\text{trace}(A)=0$, Theorem \ref{th:loopgroup} 
implies that $\phi \in \Lambda \SL_2\!\mathbb{C}$.  
By Theorem \ref{th:iwasawa-sp}, $\phi$ can be Iwasawa split into 
\[ \phi = F \cdot B \; , \;\;\; F \in \Lambda\SU_2 \; , \;\; 
B \in \Lambda_+^{\mathbb{R}} \SL_2\!\mathbb{C}. \]  
We have that 
\begin{equation}\label{eq:BAB^-1dz}
F^{-1} dF = B A B^{-1}dz - dB \cdot B^{-1} \; , 
\end{equation}
and that $dB \cdot B^{-1}$ has a power series expansion at $\lambda =0$ of 
the form 
\[ dB\cdot B^{-1}=\sum_{j=0}^\infty B_j \lambda^j \; , \] 
where the $B_j$ are matrix-valued 1-forms.  Furthermore, 
$B A B^{-1}$ has a power series expansion at $\lambda =0$ of the form 
\[ BAB^{-1}dz=\sum_{j=-1}^\infty\gamma_j \lambda^j \; , \] 
where the $\gamma_j$ are matrix-valued 1-forms.  
Therefore, by \eqref{eq:BAB^-1dz}, 
\[ F^{-1} dF = \sum_{j=-1}^\infty \alpha_j \lambda^j \] 
for some matrix-valued 1-form coefficients $\alpha_j=\alpha_j(z,\bar z)$.  
Now, the fact that $F \in \Lambda\SU_2$ implies that $F^{-1} dF$ is skew 
Hermitian, that is, 
\begin{equation}\label{eq:skewhermitean}
\overline{F^{-1} dF}^t = - F^{-1} dF \; , 
\end{equation}
and by Theorem \ref{th:loopgroup} the trace of $F^{-1} dF$ is zero.  
We conclude that 
\[ F^{-1} dF=\alpha_{-1} \lambda^{-1} + \alpha_0 + \alpha_1 \lambda \; , \] 
that is, $\alpha_j=0$ for $j \geq 2$.  Furthermore, Equation 
\eqref{eq:skewhermitean} implies that 
$\overline{\alpha_0}^t=-\alpha_0$ and $\overline{\alpha_1}^t = - \alpha_{-1}$ 
(using $\bar{\lambda}=\lambda^{-1}$ for $\lambda \in \mathbb{S}^1$).  So if, 
using dashes and double dashes to differentiate between the $(1,0)$ and 
$(0,1)$ parts of a $1$-form, we write 
\[ \alpha_0 =\alpha_{0}' dz+\alpha_{0}'' d\bar{z}, \quad 
      \alpha_1=\alpha_{1}'dz+\alpha_{1}''d\bar{z}, \quad 
      \alpha_{-1} = \alpha_{-1}'dz+\alpha_{-1}''d\bar{z} \; , \] then 
\[ \alpha_{ 0}''=-\overline{\alpha_{0}' }^t , \quad 
   \alpha_{-1}''=-\overline{\alpha_{1}' }^t , \quad 
   \alpha_{-1}' =-\overline{\alpha_{1}''}^t \; . \]  
Now, writing 
\[ A_{-1}(z) = \begin{pmatrix} 0 & a_{-1} \\ b_{-1} & 0 \end{pmatrix} 
\; , \;\;\;a_{-1}=a_{-1}(z),\;\;b_{-1}=b_{-1}(z) \] and 
\begin{equation}\label{forlater1} 
B|_{\lambda =0}=\begin{pmatrix} \rho & 0 \\ 0 & \rho^{-1} \end{pmatrix} 
\; , \;\;\; \rho :\Sigma \to \mathbb{R}^+ \; , \end{equation}
we have 
\[ 
\gamma_{-1}
=\begin{pmatrix}0&\rho^2a_{-1}\\\rho^{-2}b_{-1}&0\end{pmatrix} dz \; . 
\] 
This implies, by \eqref{eq:BAB^-1dz}, that
\[ F^{-1} dF = \lambda^{-1} \begin{pmatrix} 
0 & \rho^2 a_{-1} \\ \rho^{-2} b_{-1} & 0 \end{pmatrix} dz
+\alpha_0 +\alpha_1\lambda \; . \]
Therefore $\alpha_{-1}$ has no $d\bar{z}$ term, and hence $\alpha_{-1}''=0$ 
and $\alpha_{-1}'dz=\gamma_{-1}$.  
We conclude that 
\begin{equation}\label{forlater2}
F^{-1} dF = \lambda^{-1} \alpha_{-1}'dz+ \alpha_{0}'dz-
\overline{\alpha_{0}'}^td\bar{z}
+\lambda \left( -\overline{\alpha_{-1}'}^t\right) 
d\bar{z} \; , \end{equation} where, using \eqref{eq:BAB^-1dz}, we have that 
\[
-\overline{\alpha_{0}'}^t+\lambda \left( -\overline{\alpha_{-1}'}^t\right)
=-B_{\bar z}\cdot B^{-1} \; ,
\]
which evaluated at $\lambda =0$ implies 
\[
-\overline{\alpha_{0}'}^t = -\begin{pmatrix} \rho_{\bar z} & 0 \\ 
0 & (\rho^{-1})_{\bar z} \end{pmatrix} \cdot 
\begin{pmatrix} \rho^{-1} & 0 \\ 0 & \rho \end{pmatrix} \; , 
\] and therefore 
\begin{equation}\label{u1u0} 
\alpha_{-1}' = \begin{pmatrix} 0 & \rho^{2} a_{-1} \\ 
\rho^{-2} b_{-1} & 0 \end{pmatrix} \; , \;\;\; \alpha_{0}' = 
\begin{pmatrix} \rho_z/\rho & 0 \\ 
0 & -\rho_z/\rho \end{pmatrix} \; , 
\end{equation} with $\rho :\Sigma\to \mathbb{R}^+$ and 
$a_{-1}(z), b_{-1}(z)$ holomorphic.  

Now, because we required that $a_{-1}(z)$ is never zero on $\Sigma$ (see 
Definition \ref{df:potential}), we can make the conformal change of parameter 
\[ w=\frac{-1}{H}\int a_{-1}(z)dz \] 
and the substitutions $u=2\log( \rho )$ and $Q=-2Hb_{-1}/a_{-1}$, and then 
$F$ satisfies the Lax pair \eqref{LaxPair}-\eqref{eq:UandV} with respect to 
the parameter $w$.  

This is the Lax pair for a general CMC $H$ surface, so the above argument and 
Corollary \ref{cor:twistedlaxpair} now imply the following lemma: 

\begin{lemma}\label{lm:ProofofDPW1}
The DPW recipe always produces a CMC $H$ conformal immersion, with metric 
\[ 4e^{2u}dwd\bar w=4\rho^4dwd\bar w
  =4\rho^4H^{-2}|a_{-1}|^2dzd\bar z \; , \] 
where $\rho >0$ is the upper-left entry of $B|_{\lambda =0}$.  The surface 
is an immersion because $|a_{-1}|$ is strictly positive on $\Sigma$.  
\end{lemma}

\section{Proof of the other direction of DPW, with holomorphic potentials}
\label{section8}

We saw in Corollary \ref{cor:twistedlaxpair} that all simply-connected CMC 
$H\ne 0$ surfaces can be constructed from loops $F\in\Lambda\SU_2$ 
satisfying a Lax pair of the form \eqref{LaxPair}-\eqref{eq:UandV}.  
So, to show that the DPW recipe constructs all CMC $H$ surfaces, we need only 
prove the following lemma, which in this sense is the converse 
of Lemma \ref{lm:ProofofDPW1}.  

\begin{lemma}\label{lm:ProofofDPW2}
Given a frame $F\in\Lambda\SU_2$ defined on a simply-connected Riemann surface 
$\Sigma$ with coordinate $z$ so that $F$ is real-analytic in $z$ and $\bar z$ 
and solves a Lax pair of the form \eqref{LaxPair}-\eqref{eq:UandV}, there 
exists a holomorphic potential $\xi$ as in Definition \ref{df:potential}, and 
a holomorphic solution $\phi\in\Lambda\SL_2\!\mathbb{C}$ of $d\phi =\phi\xi$ 
with Iwasawa splitting $\phi =F \cdot B$, that is, with $F$ as given and with 
$B \in \Lambda_+ \SL_2\!\mathbb{C}$. 
\end{lemma}

\begin{remark}\label{re:352}
For any CMC $H$ surface as in \eqref{Sym-Bobenko} with frame $F$ satisfying 
the Lax pair \eqref{eq:laxpairforf}-\eqref{eq:UhatVhat}, we know that $F$ is a 
real-analytic function of $z$ and $\bar z$.  
This follows from the well-known fact that 
the metric function $u$ is real-analytic in $z$ and $\bar z$, as it is a 
solution of the Gauss equation in \eqref{firstMC}, which is elliptic with 
real-analytic coefficients (see \cite{masuda}, or see \cite{kazdan} and 
references therein).  Hence we do not lose 
anything by assuming $F$ is real-analytic in $z$ and $\bar z$, in Lemma 
\ref{lm:ProofofDPW2}.  

In fact, we can assume without loss of generality that any CMC 
immersion is both real-analytic and conformal, because 
of this: results in \cite{masuda} and \cite{kazdan} imply that 
any CMC surface is real-analytic as a graph over a plane with the usual pair of 
Euclidean coordinates, so the surface can be described as a real-analytic 
immersion.  Then a result in \cite{Sp} (volume IV, Addendum 1 of Chapter 9) 
implies the immersion can be transformed into another immersion for the same 
surface that is now conformal, without losing real-analyticity -- and this is a 
fact that was first shown by Gauss.  
\end{remark}

The proof follows this procedure:

(1) find a $B \in \Lambda_+ \SL_2\!\mathbb{C}$ so that 
$(F B)_{\bar{z}} = 0$, 

(2) define $\phi :=F B$, 

(3) define $A=A(z,\lambda )$ by $\phi_z = \phi A(z,\lambda )$, 

(4) show that $A=\sum_{j=-1}^\infty A_j(z)\lambda^j$ with the $A_j(z)$ 
holomorphic, 

(5) define the needed holomorphic potential as $\xi =Adz$.  

Noting that $F$ satisfies \eqref{LaxPair} with $U$ and $V$ as in 
\eqref{eq:UandV}, if the first step is completed, then 
\[
(F B)_{\bar{z}} = F_{\bar z} B + F B_{\bar z} = F V B + F B_{\bar z} = 0 \; , 
\]
and then 
\begin{equation}\label{eq:B_barz=-VB}
B_{\bar z}=-V B \; . 
\end{equation}
So to complete the first step, we need to find a solution $B$ of Equation 
\eqref{eq:B_barz=-VB} so that $B \in \Lambda_+ \SL_2\!\mathbb{C}$.  
For this purpose we need to consider the dell bar problem.  

\begin{remark}
If we simply noted that the compatibility condition for the Lax pair 
\[ B_{\bar z} = -V B \; , \;\;\; B_z=-U B \] 
is satisfied, we can conclude existence of a solution $B$ to Equation 
\eqref{eq:B_barz=-VB}.  However, this does not necessarily imply that 
$B \in \Lambda_+ \SL_2\!\mathbb{C}$.  In fact, we would then have $\phi$ 
constant and $A\equiv 0$, which is clearly not what we want.  
So we cannot do something as simple as this to find $B$.  
We need to do more to find an appropriate $B$, and this is the purpose of 
using the dell bar problem.  
\end{remark}

\begin{proof}
We claim that we can solve step (1), that is, that there exists a 
$B \in \Lambda_+ \SL_2\!\mathbb{C}$ globally defined on $\Sigma$ such that 
$B_{\bar z}=-VB$.  Since the frame $F$ is real-analytic, 
$V=F^{-1}F_{\bar z}$ is also real-analytic.  So we can write 
\[ V = \sum_{j,k\in\mathbb{Z}} V_{j,k} z^j\bar z^k \; . \] 
Now we replace $\bar z$ by the free variable $w$, that is, we define 
\[ \tilde V(z,w) = \sum_{j,k\in\mathbb{Z}} V_{j,k} z^jw^k \; . \]
Fix a point $z_0$.  Since $V$ is real-analytic, 
there exists $\epsilon >0$ so that $\tilde V(z,w)$ is defined for all 
$(z,w)\in B_\epsilon (z_0)\times B_\epsilon (\overline{z_0})$ (where 
$B_\epsilon (z_0):=\{z\in\Sigma\;|\;|z-z_0|<\epsilon\}$) and all 
$\lambda\in\overline{D_1}$, because $V$ contains no negative powers of 
$\lambda$.  Now $\partial _w\tilde B(z,w)=-\tilde V(z,w)\tilde B(z,w)$ is a 
single ordinary differential equation for each fixed $z$, so we solve it with 
initial condition $\tilde B(z_0,\overline{z_0})=\id$.  Note that 
$\det\tilde B(z,w)\equiv 1$, 
since $\text{trace}(\tilde V(z,w))\equiv 0$,  like in Equation 
\eqref{eq:matrixlemma}.  Note as well that $\tilde B(z,w)$ is defined on 
$B_\epsilon (z_0)\times B_\epsilon (\overline{z_0})$ because the differential 
equation is linear.  
Then $\hat B:=\tilde B(z,\bar z)$ is defined for $z\in B_\epsilon (z_0)$.  
Note that $\hat B$ is twisted, that is, 
$\hat B(-\lambda )=\sigma_3\hat B(\lambda )\sigma_3$, 
because $\hat B|_{z=z_0}=\id$ and $V$ is twisted, 
by an argument as in the proof of Theorem \ref{th:loopgroup}. 
We have $\det \hat B\equiv 1$, and that $\hat B$ is defined for all 
$\lambda\in\overline{D_1}$.  

So we have local solutions $\hat B_\alpha$ of Equation \eqref{eq:B_barz=-VB} 
on simply-connected open sets $U_\alpha\subseteq\Sigma$ that cover $\Sigma$.  
We can choose the $U_\alpha$ to be a locally finite covering of $\Sigma$, 
so that all intersections $U_\alpha \cap U_\beta$ and 
$U_\alpha \cap U_\beta \cap U_\gamma$ are also simply-connected 
whenever they are nonempty.  On $U_\alpha \cap U_\beta$, we define 
$h_{\alpha\beta}:=\hat B_\alpha^{-1} \hat B_\beta$, 
and it is easily checked that $(h_{\alpha\beta})_{\bar z}=0$, that is, 
$h_{\alpha\beta}$ is holomorphic.  Clearly 
$h_{\alpha\beta}h_{\beta\gamma}=h_{\alpha\gamma}$, so $h_{\alpha\beta}$ 
represents a holomorphic bundle.  Since $\Sigma$ is a domain in $\mathbb{C}$ 
and hence is a Stein manifold (see Section 5.1.5 of \cite{GR}), 
a generalization of Grauert's theorem to $\infty$-dimensional holomorphic 
bundles implies that any $\infty$-dimensional holomorphic bundle over $\Sigma$ 
is isomorphic to a trivial holomorphic bundle (see also Chapter 3 of 
\cite{forster}).  In other words, there exist nonzero holomorphic 
$h_\alpha\in \Lambda_+ \SL_2\!\mathbb{C}$ 
defined on $U_\alpha$ so that $h_{\alpha\beta}=h_\alpha h_\beta^{-1}$ for all 
$\alpha,\beta$.  
Then defining $B:=\hat B_\alpha h_\alpha\in\Lambda_+\SL_2\!\mathbb{C}$ for all 
$z\in\Sigma$, we have that $B$ is a globally well-defined solution of 
\eqref{eq:B_barz=-VB} on $\Sigma$, that is, 
\[
B_{\bar z}=(\hat B_\alpha h_\alpha)_{\bar z}=(\hat B_\alpha)_{\bar z} h_\alpha
          = -V \hat B_\alpha h_\alpha = -V B \; . 
\]
Now defining $\phi=F B$, we see that $\phi$ is holomorphic, because 
\[ (FB)_{\bar z} = (F_{\bar z}-F V)B = 0 \; . \]

Define $A$ by $A = \phi^{-1} \phi_z$.  Then $A$ is holomorphic in $z$.  
Writing $A$ as a power series $A = \sum_j \lambda^j A_j(z)$, 
the $A_j=A_j(z)$ are all holomorphic in $z$:  
Since $\phi_{\bar z}=0$ and so $\phi_{z\bar z}=0$, we have 
\[ 0=\phi_{\bar z} \sum_j \lambda^j A_j + \phi \sum_j\lambda^j (A_j)_{\bar z} 
= \phi \sum_j \lambda^j (A_j)_{\bar z} \; , \] and thus 
$\sum_j \lambda^j (A_j)_{\bar z} = 0$ and so $(A_j)_{\bar z}=0$ for all 
$j$.  Finally 
\begin{equation}\label{eq:phi^(-1)phi_z} 
A = B^{-1} U B + B^{-1} B_z \; . 
\end{equation}
The lowest power of $\lambda$ in the power series expansion of the right-hand 
side at $\lambda =0$ of Equation \eqref{eq:phi^(-1)phi_z} is $\lambda^{-1}$, 
coming from $U$, since $B\in\Lambda_+\SL_2\!\mathbb{C}$.  
Therefore $A_j=0$ for any $j \leq -2$. 
Also, defining $\xi=\phi^{-1}d\phi =Adz$, we have that $\xi$ is a globally 
well-defined holomorphic potential on $\Sigma$.  This completes the proof.  
\end{proof}

\begin{remark}
In Lemma \ref{lm:ProofofDPW2}, we have $B\in\Lambda_+\SL_2\!\mathbb{C}$, but 
we do not that $B\in\Lambda_+^{\mathbb{R}}\SL_2\!\mathbb{C}$, so we are using 
a non-unique Iwasawa splitting in this lemma.  However, this is irrelevant to 
the construction of CMC surfaces, 
as noted in Remark \ref{re:rIwasawa}.  
\end{remark}

\section{Normalized potentials}

The DPW method is sometimes used with respect to another type of potential, 
called the {\em normalized potential} (see, for example, \cite{D, Wu1, Wu2, 
Wu3, Wu4}).  
Therefore, in this section we prove that the DPW construction produces all CMC 
surfaces with respect to normalized potentials as well.  
As we use normalized potentials only in this section, this section is 
independent of the other sections of these notes.  

Let $\Sigma$ be a simply-connected Riemann surface with coordinate $z$.  

\begin{defn}\label{df:normalizedpot}
A {\em normalized potential} on $\Sigma$ is of the form
\[
\xi =\lambda^{-1}\begin{pmatrix}0&a(z)\\b(z)&0\end{pmatrix} dz \; , 
\]
where $a(z)$ and $b(z)$ are meromorphic on $\Sigma$.  
\end{defn}

\begin{remark}
In \cite{DH} and \cite{DPW}, the potentials in Definition 
\ref{df:normalizedpot} were called {\em meromorphic potentials}, but nowadays 
they are referred to as normalized potentials, to avoid confusion with other 
uses of the phrase ``meromorphic potential''.  
\end{remark}

Away from points where $a(z)$ or $b(z)$ have poles, or where $a(z)$ is zero, 
$\xi$ is actually a holomorphic potential.  
So away from these isolated points, the DPW recipe applied to this $\xi$ 
produces a CMC immersion, by Lemma \ref{lm:ProofofDPW1}.  
Therefore, the DPW recipe applied to a normalized potential $\xi$ on $\Sigma$ 
produces a CMC surface possibly with isolated singularities.  

Conversely, we wish to show that any simply-connected portion of any CMC 
immersion can be constructed by applying the DPW recipe to some normalized 
potential.  This is the goal of the next two lemmas.  

\begin{lemma}\label{lm:normalized1}
{\rm (DPW for normalized potentials, local version)} 
Let $\Sigma$ be a simply-connected Riemann surface with coordinate $z$.  
Let $z_*$ be a fixed point in $\Sigma$.  
Let $f:\Sigma\to\mathbb{R}^3$ be a conformal CMC $H(\ne 0)$ immersion with 
associated extended frame $F:\Sigma\to\Lambda\SU_2$, 
and suppose that $F|_{z=z_*}=\id$ for all $\lambda\in\mathbb{S}^1$.  
Then there exist an open set $\mathcal{V}\subseteq\Sigma$ with 
$z_*\in\mathcal{V}$ and holomorphic functions $a(z)$, $b(z)$ on $\mathcal{V}$, 
and there exists a holomorphic matrix $\phi\in\Lambda\SL_2\!\mathbb{C}$ 
defined on $\mathcal{V}$, so that 
\[
d\phi =\phi\xi \; , \;\;\; 
\xi =\lambda ^{-1}\begin{pmatrix}0&a(z)\\b(z)&0\end{pmatrix}dz \; , 
\]
and there exists a $B\in\Lambda_+ \SL_2\!\mathbb{C}$ such that 
$\phi =F \cdot B$ $($this is nonunique Iwasawa splitting$)$, where $F$ is the 
frame of $f$ as above.  
\end{lemma}

\begin{proof}
By Lemma \ref{lm:ProofofDPW2}, there exist 
$\tilde G\in\Lambda\SL_2\!\mathbb{C}$ and 
$\tilde B\in\Lambda_+ \SL_2\!\mathbb{C}$ such that $\tilde G =F\cdot\tilde B$ 
and $\tilde G^{-1}d\tilde G$ is a holomorphic potential and 
$\partial_{\bar z}\tilde G=0$.  
Note that $F$ satisfies $F|_{z=z_*}=\id$.  
Define $B\in\Lambda_+ \SL_2\!\mathbb{C}$ by 
$B=\tilde B\cdot (\tilde B|_{z=z_*})^{-1}$.  
Then also $B|_{z=z_*}=\id$.  Define $G\in\Lambda\SL_2\!\mathbb{C}$ by 
$G=\tilde G\cdot (\tilde B|_{z=z_*})^{-1}$.  
Then $G=F\cdot B$ and $G|_{z=z_*}=\id$ and $G^{-1}dG$ is a holomorphic 
potential and $\partial_{\bar z}G=0$.  
Since $\id$ is contained in the big cell 
for Birkhoff splitting, there exists an open set $\mathcal{V}$ such that 
$z_*\in\mathcal{V}$ and $G$ has a well-defined Birkhoff splitting on 
$\mathcal{V}$ as follows: $G =\phi\cdot G_+$, where 
$\phi\in\Lambda_-^*\SL_2\!\mathbb{C}$ and $G_+\in\Lambda_+\SL_2\!\mathbb{C}$ 
(we do not need to require that $G_+\in\Lambda_+^\mathbb{R}\SL_2\!\mathbb{C}$ 
as explained in Remark \ref{re:rIwasawa}) and $\phi$ and $G_+$ have no 
singularities on $\mathcal{V}$.  Furthermore, $\phi$ and $G_+$ are holomorphic 
in $z$ on $\mathcal{V}$, since $G$ is holomorphic on $\mathcal{V}$ and 
Birkhoff decomposition is complex analytic.  Then 
\begin{equation}\label{eq:phi^-1dphi1}
\phi^{-1}d\phi =G_+(B^{-1}F^{-1}dFB+B^{-1}dB-G_+^{-1}dG_+)G_+^{-1} \; . 
\end{equation}
Equation \eqref{LaxPair} implies that 
\begin{equation}\label{eq:F^-1dF}
F^{-1}dF = \lambda^{-1}\alpha_{-1}dz+\alpha_0dz
 -\bar \alpha_0^t d\bar z - \lambda \bar \alpha_{-1}^t d\bar z \; , 
\end{equation}
where 
\[ 
\alpha_0   =\frac{1}{2}\begin{pmatrix}u_z&0\\0&-u_z\end{pmatrix} \text{ and } 
\alpha_{-1}=\frac{1}{2}\begin{pmatrix}0&-2He^u\\Qe^{-u}&0\end{pmatrix} \; .  
\]
Writing 
\[
\phi=\id+\sum_{j=-1}^{-\infty}\lambda^{j} \phi_j(z,\bar z)
\]
for some functions $\phi_j(z,\bar z)$, we have 
$d\phi=\sum_{j=-1}^{-\infty} \lambda^{j}d(\phi_j(z,\bar z))$ and so 
\begin{equation}\label{eq:phi^-1dphi2}
\phi^{-1} d\phi = \sum_{j=-1}^{-\infty}\lambda^j\left(
\hat\phi_{j,1}dz+\hat\phi_{j,2}d\bar z \right) \; , 
\end{equation}
for some functions $\hat\phi_{j,1}$ and $\hat\phi_{j,2}$.  Also note that 
\begin{equation}\label{eq:G_+^-1dG_+}
G_+^{-1} dG_+ = \sum_{j=0}^\infty \lambda^j\left(
\hat G_{j,1}dz+\hat G_{j,2}d\bar z \right) \; , \;\;\; 
B^{-1} dB = \sum_{j=0}^\infty \lambda^j\left(
\hat B_{j,1}dz+\hat B_{j,2}d\bar z \right) \; , 
\end{equation}
for some functions $\hat G_{j,1}$, $\hat G_{j,2}$, $\hat B_{j,1}$ and 
$\hat B_{j,2}$.  Therefore, with $G_+=\sum_{j=0}^\infty \lambda^j G_{+,j}$ and 
$B=\sum_{j=0}^\infty \lambda^j B_j$ for some functions $G_{+,j}$ and $B_j$, 
and using \eqref{eq:F^-1dF} and comparing the coefficients in the expansions 
\eqref{eq:phi^-1dphi2} and \eqref{eq:G_+^-1dG_+} with respect to $\lambda$, 
we have 
\[
\phi^{-1}d\phi =\lambda^{-1}G_{+,0}B_0^{-1}\alpha_{-1}B_0 G_{+,0}^{-1}dz \; .  
\]
Since $G_{+,0} B_0^{-1}$ is diagonal (by twistedness), 
$G_{+,0} B_0^{-1} \alpha_{-1} B_0 G_{+,0}^{-1}$ is off-diagonal; and since no 
$d\bar z$ term appears on the right hand side, $\phi_{\bar z}=0$, so 
\[ \phi^{-1} d\phi = \lambda^{-1} \begin{pmatrix}
0 & a \\ b & 0
\end{pmatrix} dz \; , \] where $a_{\bar z}=0$ and $b_{\bar z}=0$, 
and $a$ and $b$ have no singularities on $\mathcal{V}$.  
\end{proof}

The above lemma can be extended to a global version, that is, 
we may extend $\mathcal{V}$ to be all of $\Sigma$.  
However, we must then allow $a(z)$ and $b(z)$ to have poles.  
We do not prove the following lemma (global version) here, 
and simply note that it follows directly from Theorem 4.10 in \cite{DPW}.  

In the proof of Lemma \ref{lm:normalized1}, $G$ is nonsingular on $\Sigma$.  
The proof of Theorem 4.10 in \cite{DPW} shows that the elements 
$\phi$ and $G_+$ have at most pole singularities at fixed isolated values of 
$z$ where $G$ leaves the big cell $\mathcal{U}$.  

\begin{lemma}
{\rm (DPW for normalized potentials, global version)} 
Let $\Sigma$ be a simply-connected Riemann surface with coordinate $z$.  
Let $z_*$ be a fixed point in $\Sigma$.  
Let $f:\Sigma\to\mathbb{R}^3$ be a conformal CMC $H(\ne 0)$ immersion with 
associated extended frame $F:\Sigma\to\Lambda\SU_2$, 
and suppose that $F|_{z=z_*}=\id$ for all $\lambda\in\mathbb{S}^1$.  
Then there exist meromorphic functions $a(z)$, $b(z)$ on $\Sigma$, 
and there exists a meromorphic matrix $\phi$ on $\Sigma$ which is in 
$\Lambda\SL_2\!\mathbb{C}$ away from its poles, so that 
\[
d\phi =\phi\xi \; , \;\;\; 
\xi =\lambda ^{-1}\begin{pmatrix}0&a(z)\\b(z)&0\end{pmatrix}dz \; , 
\]
and there exists a $B\in\Lambda_+ \SL_2\!\mathbb{C}$, 
defined wherever $\phi$ is defined, such that $\phi =F \cdot B$ 
$($this is nonunique Iwasawa splitting$)$.  
\end{lemma}

\section{Dressing and gauging}
\label{section14}

Let $\Sigma$ be a Riemann surface with coordinate $z$.  
Here we do {\em not} assume that $\Sigma$ is simply-connected.  
Let $\xi$ be a holomorphic potential on $\Sigma$.  
Given a solution $\phi\in\Lambda\SL_2\!\mathbb{C}$ to $d\phi =\phi\xi$, 
if we define 
\[ \hat\phi = h_+(\lambda)\cdot\phi\cdot p_+(z,\bar z,\lambda) \; , 
\qquad h_+ , \; p_+ \in \Lambda_+ \SL_2\!\mathbb{C} \; , \] 
then the multiplication on the left by $h_+$ is a dressing, and the 
multiplication on the right by $p_+$ is a gauging.  The matrix $h_+$ is not 
allowed to depend on $z$.  The matrix $p_+$ can depend on $z$, 
but must be a single-valued matrix function on $\Sigma$.  

Let $f:\widetilde{\Sigma}\to\mathbb{R}^3$ be the CMC $H$ surface resulting 
from applying the DPW recipe to $\phi$, where $\widetilde{\Sigma}$ is the 
universal cover of $\Sigma$.  

Note that $\hat\phi$ satisfies $d\hat\phi = \hat\phi\hat\xi$, where 
\[
\hat\xi = \hat\phi^{-1}d\hat\phi 
        = (h_+\phi p_+)^{-1}d(h_+\phi p_+) 
        = p_+^{-1}\xi p_+ + p_+^{-1} dp_+ \; . 
\]
Hence, 
\begin{quote}
the dressing 
$h_+$ does not change the potential $\xi$, 
and changes only the resulting surface $f$.  
\end{quote}
Furthermore, if we look at the Iwasawa splittings $\phi =FB$ and 
$\hat\phi = \hat F\hat B$, then the change $F\to\hat F$ is 
affected only by $h_+$ and is independent of $p_+$, 
since $p_+\in\Lambda_+\SL_2\!\mathbb{C}$ and we are multiplying on the right 
by $p_+$.  Hence, 
\begin{quote}
the gauging $p_+$ does not change the surface $f$, 
and changes only the potential $\xi$.  
\end{quote}
To see how the surface $f$ is changed by $h_+$, one must Iwasawa split 
$h_+\cdot F$ into $h_+\cdot F= \tilde F\tilde B$ (and then $\hat F$ equals 
$\tilde F$) and then study the change in the frame $F\to\hat F$.  
This change is not trivial to understand, hence the change in the surface $f$ 
is also not trivial to understand.  

\begin{figure}[phbt]
\begin{center}
\includegraphics[width=0.4\linewidth]{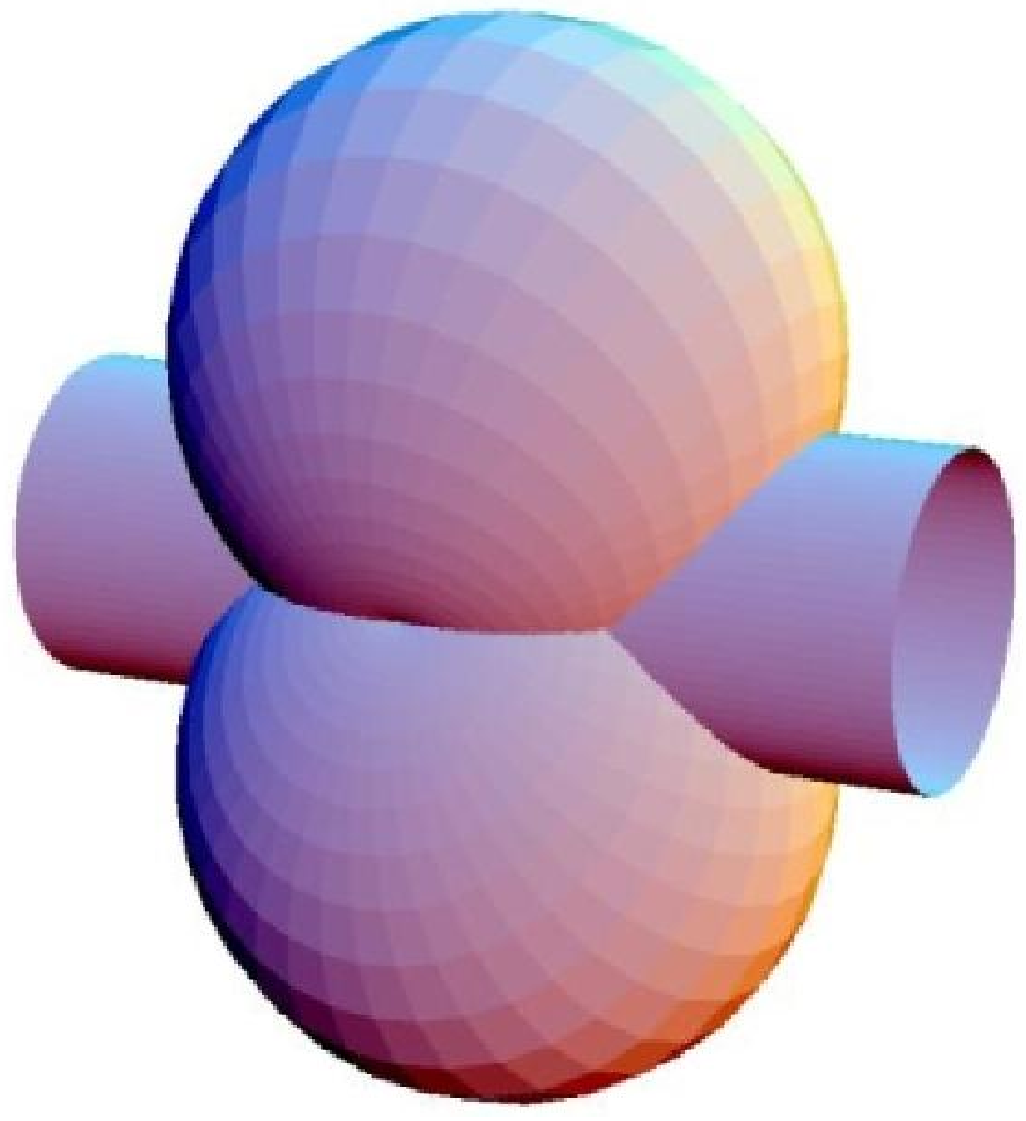}
\includegraphics[width=0.6\linewidth]{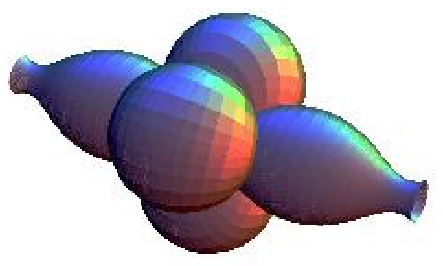}
\end{center}
\vspace{2in}
\caption{A cylinder bubbleton and a Delaunay bubbleton.  (The right-hand graphic 
was made by Yuji Morikawa using Mathematica.)}
\label{fig:bubbletons}
\end{figure}

It is precisely the above way of dressing (and the resulting nontrivial change 
$F\to\hat F$ of the frames) that produces bubbletons (with frames $\hat F$) 
from cylinders and Delaunay surfaces (with frames $F$), see Figure 
\ref{fig:bubbletons}.  Producing bubbletons requires the use of a 
generalization of Iwasawa splitting to something called $r$-Iwasawa splitting, 
so we do not discuss bubbletons here.  
For more on bubbletons, see \cite{Koba}.  

\begin{remark}\label{re:monodromy}
It is easier to understand how objects called the monodromy matrices of $\phi$ 
and $\hat\phi$ are related by $h_+$ (easier than understanding the change 
$F\to\hat F$), and this is often just the information we 
need to show that the dressed CMC surface resulting from $\hat F$ is 
well-defined on $\Sigma$ if $f$ is.  For example,  this is how one can see 
that dressing of Delaunay surfaces results in bubbletons that are 
homeomorphically annuli \cite{k, Koba}.  We will say more about monodromy 
matrices and well-definedness of CMC immersions of non-simply-connected 
Riemann surfaces in Section \ref{section13}.  
\end{remark}

\begin{figure}[phbt]
\begin{center}
\includegraphics[width=0.4\linewidth]{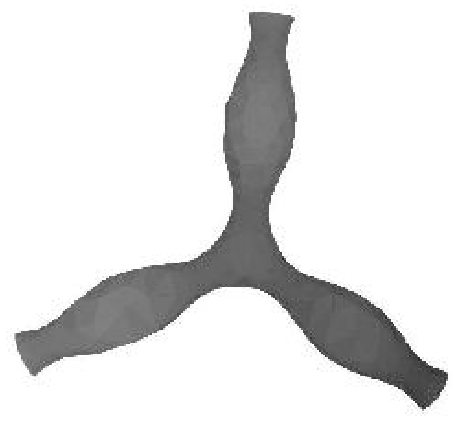}
\includegraphics[width=0.32\linewidth]{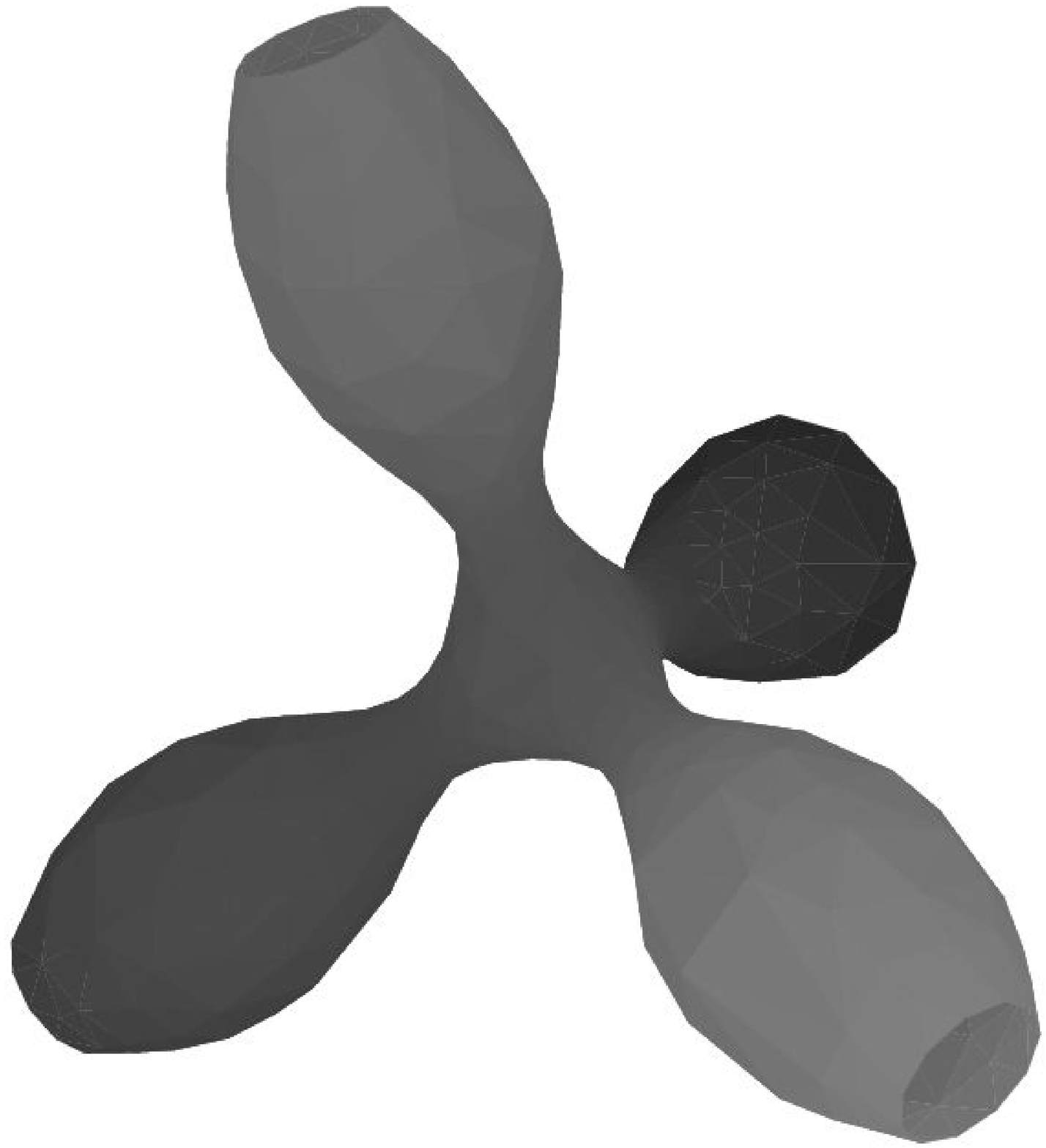}
\end{center}
\vspace{2in}
\caption{A CMC trinoid and CMC tetranoid in $\mathbb{R}^3$.  
         These are surfaces with $3$ and $4$ ends, respectively, 
         that all converge to ends of Delaunay surfaces.  (These graphics were made 
using N. Schmitt's software CMCLab \cite{Sch:cmclab}.)}
\label{fg:tri-tetra}
\end{figure}

\section{Properties of holomorphic potentials}
\label{section10}

Let $\Sigma$ be a Riemann surface with local coordinate $z$, 
which might not be simply-connected.  
Let $f:\Sigma\to\mathbb{R}^3$ be a CMC $H$ surface, 
and let $\xi$ be a holomorphic potential, defined on the universal cover 
$\widetilde\Sigma$ of $\Sigma$, that produces $f$ via the DPW recipe.  
Here we describe two assumptions that can be made about the potentials $\xi$ 
that will still allow us to produce all CMC surfaces $f$.  Hence either of 
these two assumptions can be made without loss of generality.  

\begin{figure}[phbt]
\begin{center}
\includegraphics[width=0.4\linewidth]{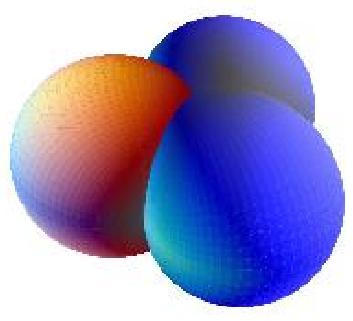}
\includegraphics[width=0.5\linewidth]{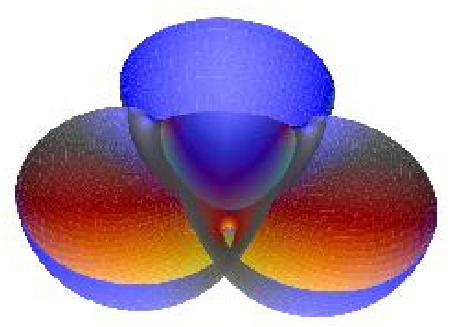}
\end{center}
\vspace{2in}
\caption{A Wente torus, and a half cut-away of it.  (These graphics were made by 
Koichi Shimose using Mathematica.)}
\label{fig:Wente}
\end{figure}

\begin{lemma}\label{lm:DorPW}
To produce $f$ via the DPW recipe, 
we may assume that the holomorphic potential $\xi =Adz$ is off-diagonal.  
\end{lemma}  

This lemma is proven by finding an appropriate gauge $p_+$ that changes $\xi$ 
into an off-diagonal potential (and does not change $f$).  
This lemma is proven in \cite{DPW}.  

\begin{proof}
Write $A=(A_{ij})_{i,j=1,2}=A_{\rm off}+A_{\rm diag}$ in off-diagonal and 
diagonal parts, that is, 
\[
A_{\rm off} =\begin{pmatrix}0&A_{12}\\A_{21}&0\end{pmatrix} 
\;\;\;\text{and}\;\;\; 
A_{\rm diag}=\begin{pmatrix}A_{11}&0\\0&A_{22}\end{pmatrix} \; .  
\]
Let $z_*$ be a fixed point in $\Sigma$.  Solve 
\begin{equation}\label{eq:dp+=p+Adiagdz}
dp_+=p_+ \cdot A_{\rm diag} dz \; , \;\;\; p_+(z_*)=\id \; , \;\;\; 
p_+=(p_{ij})_{i,j=1,2} \; .  
\end{equation}
Then 
$p_{11}^\prime = p_{11} A_{11}$ and $p_{12}^\prime = - p_{12} A_{11}$ and 
$p_{21}^\prime = p_{21} A_{11}$ and $p_{22}^\prime = - p_{22} A_{11}$, since 
$A_{11}=-A_{22}$.  
Since $p_{12}$ and $p_{21}$ are initially zero at $z_*$, $p_{12}$ and $p_{21}$ 
are identically zero, so $p_+$ is diagonal, which implies that 
$p_+ A_{\rm off} p_+^{-1}$ is off-diagonal.  
Also $A_{\rm diag} \in \Lambda_+ \slg_2\!\mathbb{C}$, 
so $p_+ \in \Lambda_+ \SL_2\!\mathbb{C}$.  
Note that generally $p_+$ is only defined on $\widetilde\Sigma$.  

For a solution $\phi$ of $d\phi=\phi \cdot A dz$ that produces $f$, 
set $\tilde \phi:=\phi p_+^{-1}$.  Then 
\begin{eqnarray*}
     d\tilde \phi
 &=& d\phi \cdot p_+^{-1} - \phi p_+^{-1} \cdot dp_+ \cdot p_+^{-1} \\
 &=& \tilde \phi p_+ \left( (A_{\rm off} + A_{\rm diag}) dz
                         - p_+^{-1}dp_+ \right) p_+^{-1}
 = \tilde \phi p_+ A_{\rm off} p_+^{-1} dz \; , 
\end{eqnarray*}
since \eqref{eq:dp+=p+Adiagdz} holds and since 
$p_+^{-1} \cdot dp_+=dp_+ \cdot p_+^{-1}$ as $p_+$ is diagonal.  
The new potential is $\tilde\xi = \tilde Adz = p_+ A_{\rm off} p_+^{-1}dz$, 
which is off-diagonal.  Note that $\phi$ and $\tilde \phi$ produce the same 
CMC $H$ surface $f$, as shown in Section \ref{section14}, that is, 
if $\phi =F \cdot B$ then $\tilde\phi =F(Bp_+^{-1})$, 
so $\phi$ and $\tilde\phi$ have the same unitary part $F$ and hence produce 
the same surface $f$.  
\end{proof}

		\begin{figure}[phbt]
			\begin{center}
			\includegraphics[width=0.38\linewidth]{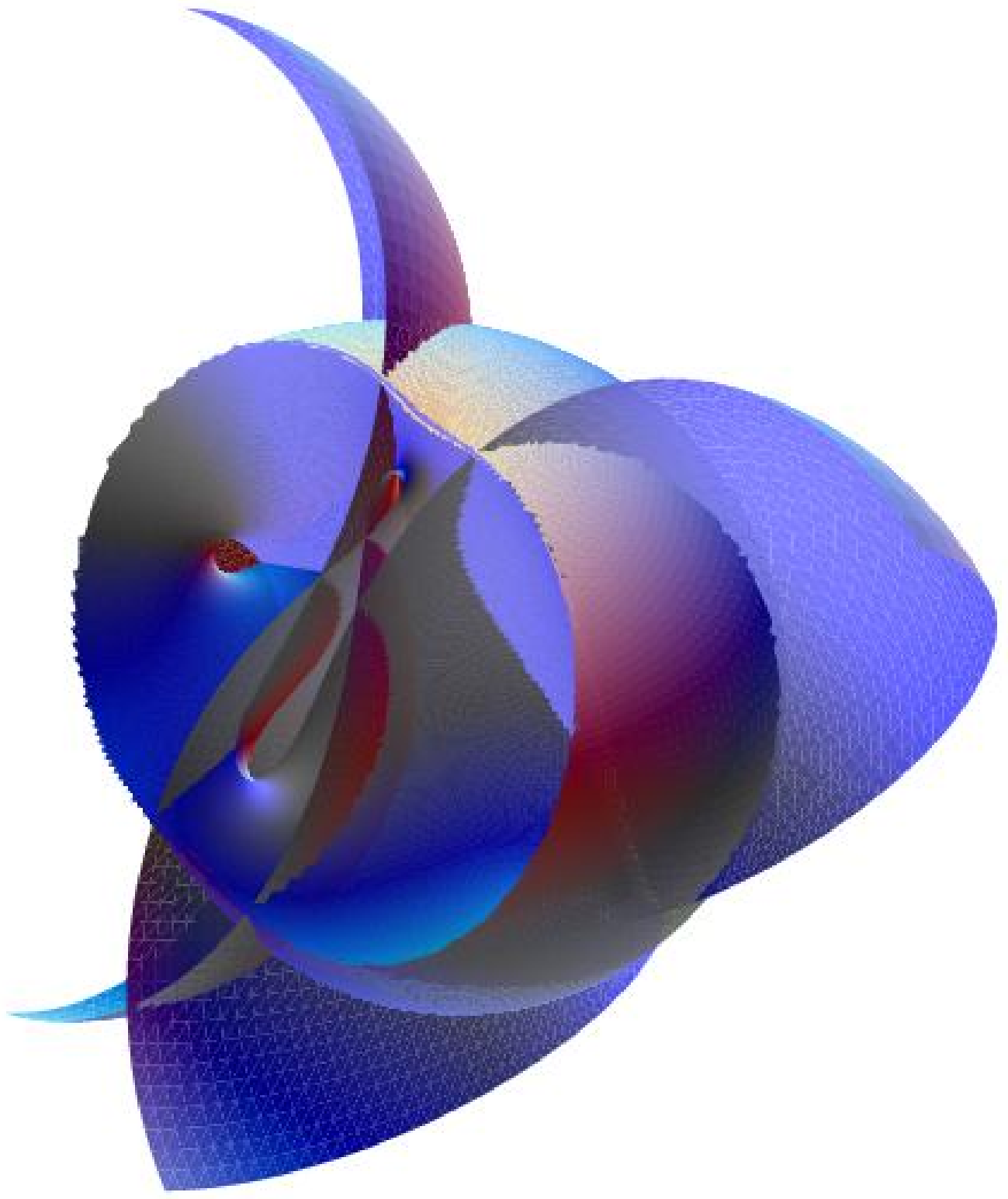}
			\hspace{1cm}
			\includegraphics[width=0.38\linewidth]{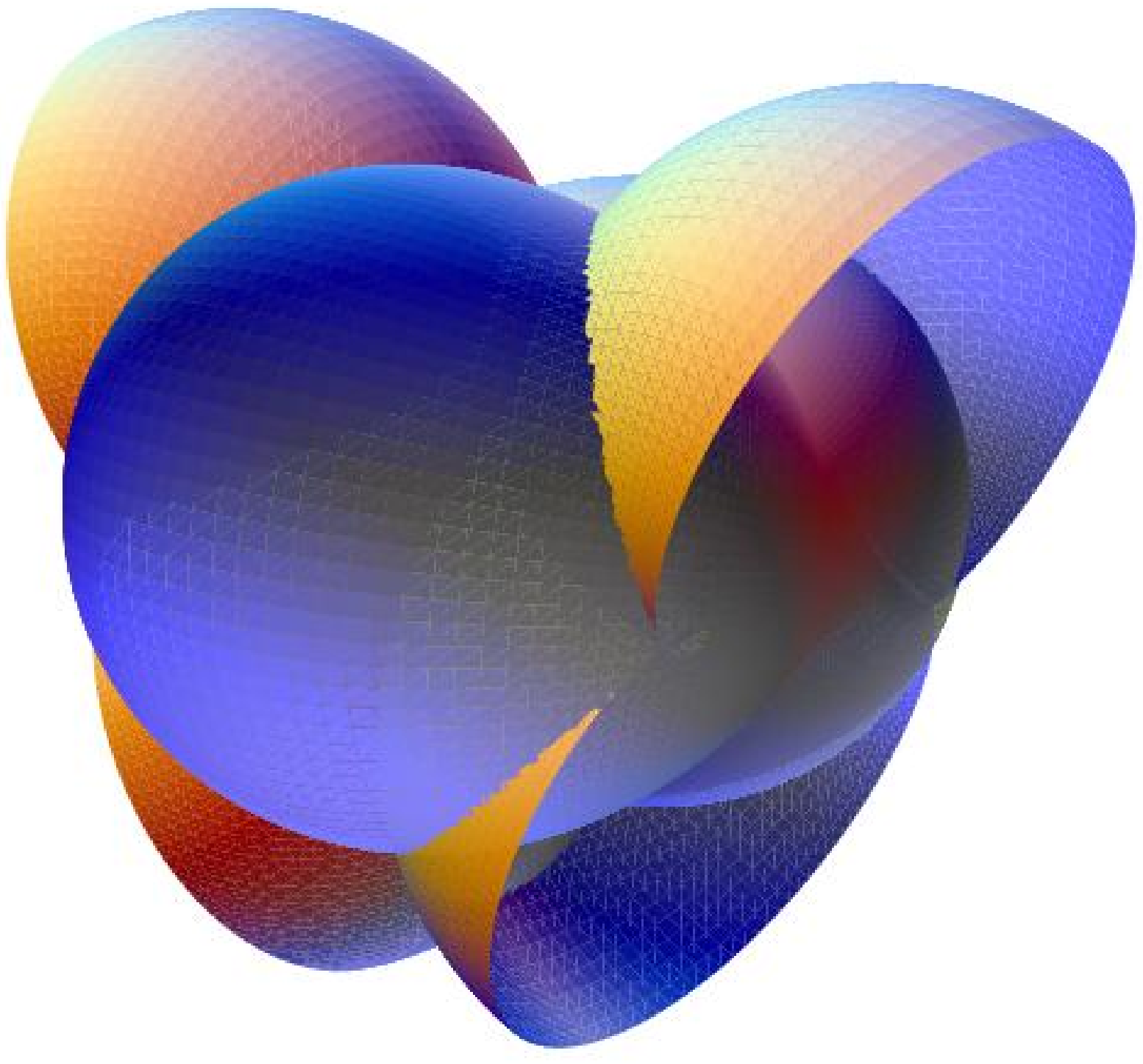}
			\end{center}
	\caption{A small and somewhat larger portion of a Dobriner torus.  
			(These graphics were made by Koichi Shimose using Mathematica.)}
			\label{dob1}
		\end{figure}

		\begin{figure}[phbt]
			\begin{center}
			\includegraphics[width=0.45\linewidth]{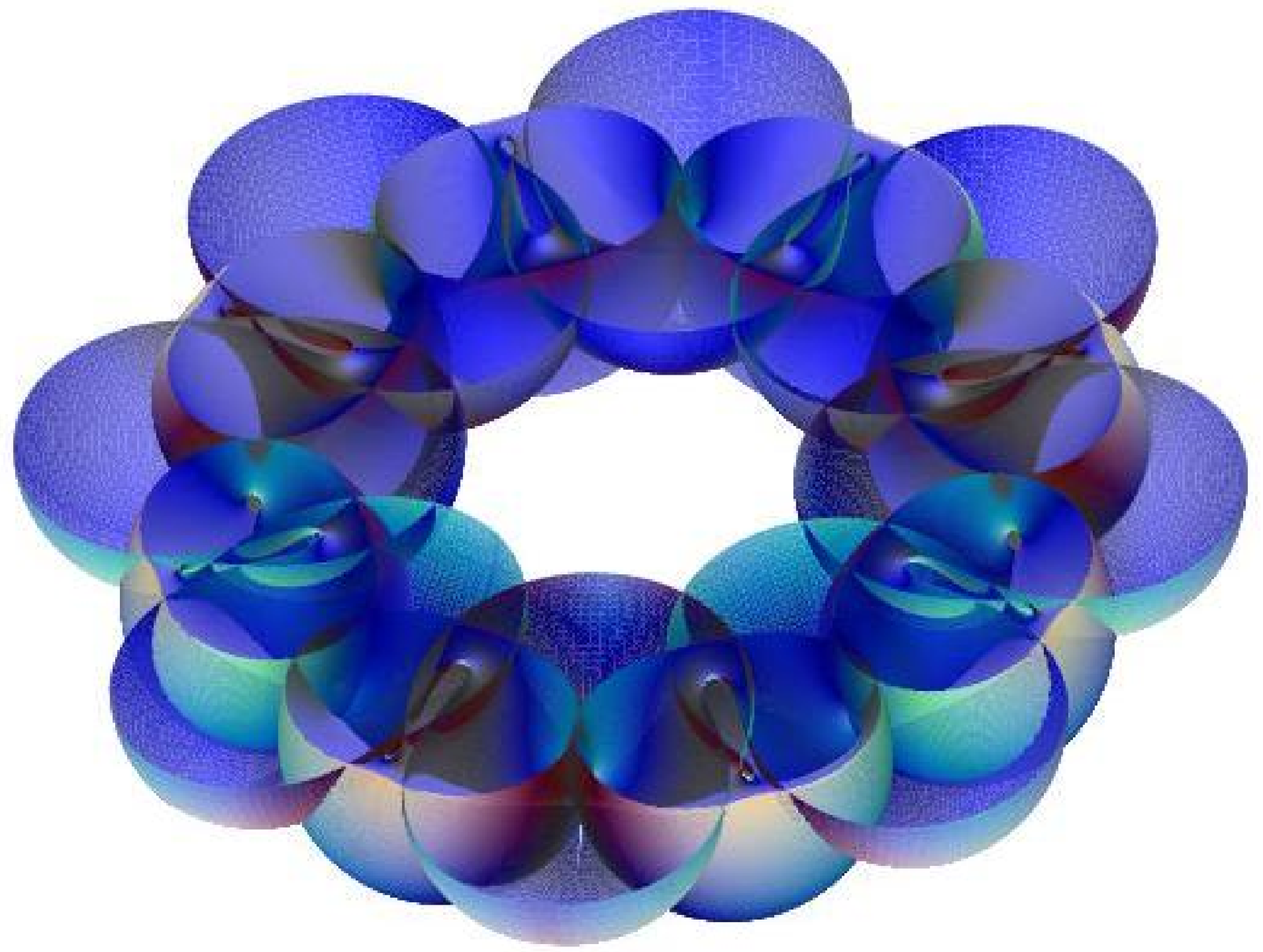}
			\includegraphics[width=0.45\linewidth]{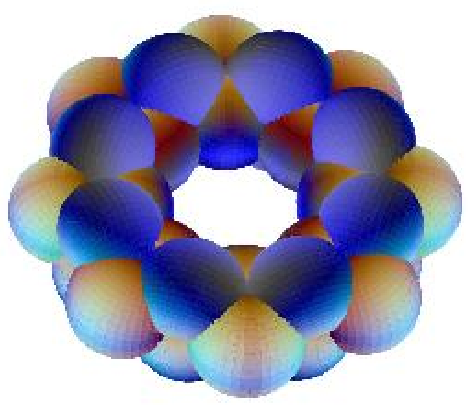}
			\end{center}
\vspace{2in}
			\caption{Half and all of a Dobriner torus.  
			(These graphics were made by Koichi Shimose using Mathematica.)}
			\label{dob2}
		\end{figure}

Note that the off-diagonal potential in Lemma \ref{lm:DorPW} is only defined 
on $\widetilde\Sigma$ in general, 
as $p_+$ is only defined on $\widetilde\Sigma$.  

The next lemma shows that we can always choose $\xi$ to be well-defined on 
$\Sigma$ (but then not necessarily off-diagonal).  
This lemma is proven in \cite{DPW} and \cite{k}.  

\begin{lemma}\label{lm:well-definedness}
Assume that $\Sigma$ is noncompact.  To produce $f$ via the DPW recipe, we may 
assume that the holomorphic potential $\xi$ is well-defined on $\Sigma$.  
\end{lemma}

\begin{proof}
Let $F$ be the extended frame of $f$.  Because $f$ is well-defined on 
$\Sigma$, $F$ is as well.  Although Lemma \ref{lm:ProofofDPW2} is stated only 
for simply-connected Riemann surfaces $\Sigma$, we can see from its proof that 
it is valid whenever $\Sigma$ is a Stein manifold.  
Furthermore, because $\Sigma$ is noncompact, it is a Stein manifold (see 
Section 5.1.5 of \cite{GR}, see also Chapter 3 of \cite{forster}).  We 
conclude that there exists a well-defined 
holomorphic potential that produces the extended frame $F$, 
and hence the CMC $H$ immersion $f$, via the DPW recipe.  
\end{proof}

Thus, when a CMC surface in $\mathbb{R}^3$ is noncompact, 
we have just seen that $\xi$ can be chosen to be well-defined on the Riemann 
surface $\Sigma$ itself.  
However, it is not always true that a complete CMC surface is noncompact.  
The sphere is a trivial example of a complete compact CMC surface, 
and there are nontrivial known examples as well, such as the Wente tori 
(see Figure \ref{fig:Wente}) and Dobriner tori (see Figures \ref{dob1} and 
\ref{dob2}).  More on CMC tori can be found in 
\cite{Wen:hopf, W, b1, Bob:tor, Bob:2x2, ps, Wal, Ab}.  
Furthermore, Kapouleas \cite{Kap2, Kap3} has found complete compact examples 
of genus two or more.  

\chapter{Lax pairs and the DPW method in $\mathbb{S}^3$ and $\mathbb{H}^3$} 
\label{chapter4} 

In this chapter we consider how Lax pairs and the DPW method can be used 
to study surfaces in some spaces other than $\mathbb{R}^3$, those other spaces 
being $\mathbb{S}^3$ (spherical $3$-space) and $\mathbb{H}^3$ 
(hyperbolic $3$-space).  
We use the DPW method to produce CMC surfaces of revolution (analogs of 
Delaunay surfaces) in $\mathbb{S}^3$ and $\mathbb{H}^3$.  
We use Lax pairs to prove a representation by Bryant for CMC $1$ surfaces in 
$\mathbb{H}^3$, and to study flat surfaces in $\mathbb{H}^3$.  

\section{The Lax pair for CMC surfaces in $\mathbb{S}^3$}
\label{section11}

Let $\mathbb{S}^3$ denote the unit $3$-sphere in $\mathbb{R}^4$, as in 
Subsection \ref{spherical3space}.  
($\mathbb{S}^3$ is a space form because it is the unique simply-connected 
complete $3$-dimensional Riemannian manifold with constant sectional curvature $+1$.) 
Consider a conformal immersion 
\[ f:\Sigma\to\mathbb{S}^3 \subset \mathbb{R}^4 \] 
of a simply-connected Riemann surface $\Sigma$ with complex coordinate $w$.  
(By Theorem \ref{conformalityispossible}, any immersion $f:\Sigma\to\mathbb{S}^3$ can be 
reparametrized so that it is 
conformal, so without loss of generality we assume here that $f$ is conformal.)
Then 
\[ \langle f, f \rangle_{\mathbb{R}^4} = 1 \, , \, \text{ where } 
\langle (x_1,x_2,x_3,x_4),(y_1,y_2,y_3,y_4) \rangle_{\mathbb{R}^4} = 
\sum_{k=1}^4 x_k y_k \; , \] and 
\[ \langle f_w,N\rangle =\langle f_{\bar w},N\rangle =\langle f,N\rangle =0 
\; , \;\; \langle N,N\rangle =1 \; , \] 
where $N$ is a unit normal of $f$ in $\mathbb{S}^3$ and $\langle \cdot , \cdot 
\rangle$ denotes the complex bilinear extension of the inner product 
$\langle \cdot , \cdot \rangle_{\mathbb{R}^4}$ above.  
The fact that $f$ is conformal implies 
\[ \langle f_w,f_w\rangle =\langle f_{\bar w},f_{\bar w}\rangle =0 \; , 
\;\; \langle f_w,f_{\bar w}\rangle =2 e^{2 u} \; , \] 
where $u:\Sigma\to\mathbb{R}$ is defined by this equation.  
We now refer to the Hopf differential as $A:\Sigma\to\mathbb{C}$ 
(we reserve the symbol $Q$ for a slightly different purpose later) 
and define it and the mean curvature $H:\Sigma\to\mathbb{R}$ by 
\[ A = \langle f_{ww},N\rangle \; , \;\;\; 
   2 H e^{2 u} =\langle f_{w\bar w},N\rangle \; . \] 
For the frame $\mathcal{F} = (f,f_w,f_{\bar w},N)$, 
similar to the computation in Section \ref{section6}, we have that 
$\mathcal{F}$ satisfies 
\[ \mathcal{F}_w = \mathcal{F}{\mathcal U} \; , \;\;\; 
   \mathcal{F}_{\bar w} = \mathcal{F}{\mathcal V} \; , \] where 
\[ {\mathcal U}=\begin{pmatrix}
0 & 0    & -2e^{2u} & 0 \\
1 & 2u_w & 0        & -H \\
0 & 0    & 0        & -Ae^{-2 u}/2 \\
0 & A    & 2He^{2u} & 0 
\end{pmatrix} \; , \;\;\; 
{\mathcal V}=\begin{pmatrix}
0 & -2e^{2u} & 0           & 0 \\
0 & 0        & 0           & -\bar Ae^{-2u}/2 \\
1 & 0        & 2u_{\bar w} & -H \\
0 & 2He^{2u} & \bar A      & 0 
\end{pmatrix} \; . \] 
The compatibility condition is 
\[ {\mathcal U}_{\bar w}-{\mathcal V}_w-[{\mathcal U},{\mathcal V}]=0 \; , \] 
which implies that 
\begin{equation}
\label{originalS3MC} 2 u_{w\bar w}+2 e^{2 u} (1+H^2) -\frac{1}{2} A 
\bar A e^{-2 u} = 0 \; , \;\;\; A_{\bar w}= 2H_we^{2u} \; . \end{equation} 
Now suppose that $H$ is constant, and thus $A_{\bar w=0}$.  
Making the change of variable 
$z=2\sqrt{1+H^2} w$, then choosing a real constant 
$\psi$ and defining a complex function $Q:\Sigma\to\mathbb{C}$ by 
$A=2\sqrt{1+H^2}e^{i\psi} Q$, we have 
\begin{equation}\label{maurer-cartan}
4 u_{z\bar z}+e^{2 u} - Q \bar Q e^{-2 u} = 0 \; , \;\;\; 
Q_{\bar z}= 0 \; , 
\end{equation}
which is the same as Equation \eqref{firstMC} (if $H$ is chosen to be $1/2$ 
there).  
We conclude that a CMC $H$ conformal immersion into $\mathbb{S}^3$ has a 
conformal parameter $z$ such that $Q= 2 \sqrt{1+H^2} e^{-i \psi} 
\langle f_{zz},N\rangle$ is holomorphic and 
$u$ satisfies the first equation in \eqref{maurer-cartan}.
The factor $e^{-i\psi}$ in $Q$ will play the role of the spectral parameter.  

We now consider how to get a surface from given $u$, $Q$, and $H$.  
With the Pauli matrices 
$\sigma_1$, $\sigma_2$, $\sigma_3$ as in Section \ref{2by2laxpair}, 
we first note that for any $2 \times 2$ matrix $X$, 
\begin{equation} \label{R4form}
X = \sigma_2 \bar X \sigma_2 \Rightarrow
X = \begin{pmatrix} a & b \\ -\bar b & \bar a \end{pmatrix} \; , 
\end{equation}
and that such a matrix $X$ represents a point in 
$\mathbb{R}^4$ via $X \equiv (a_1,b_2,b_1,a_2) \in \mathbb{R}^4$, 
where $a=a_1+i a_2$, $b= b_1+i b_2$.  
(This is a slightly different convention than we used to notate points in 
$\mathbb{R}^3$, in \eqref{prob3-a}.)  So we may consider $\mathbb{R}^4$ to 
be the set of matrices $X$ satisfying Equation \eqref{R4form}.  The $3$-sphere 
$\mathbb{S}^3$ is then the set of those $X$ in $\mathbb{R}^4$ such that 
$|a|^2+|b|^2=1$, that is, $\mathbb{S}^3$ is identified with $\SU_2$.  

For $X,Y\in\mathbb{R}^4$, that is, for $X$ and $Y$ of the form \eqref{R4form}, 
we have
\begin{equation}
\langle X,Y\rangle = \frac{1}{2}\text{trace}(X\sigma_2 Y^t\sigma_2) \, .  
\label{eq:S3metric}
\end{equation}

\begin{remark}
In fact, $\langle\cdot ,\cdot\rangle$ is the complex bilinear extension of 
$\langle\cdot ,\cdot\rangle_{\mathbb{R}^4}$, so \eqref{eq:S3metric} applies 
even when $X,Y\in\mathbb{C}^4$.  
Because $Y\in\mathbb{R}^4$ implies that $\bar Y^t=\sigma_2Y^t\sigma_2$,
we may use $\langle X,Y\rangle = (1/2)\text{trace}(X\bar Y^t)$ instead of 
\eqref{eq:S3metric} when $X,Y\in\mathbb{R}^4$, but 
$\langle X,Y\rangle \ne (1/2)\text{trace}(X\bar Y^t)$ for general 
$X,Y\in\mathbb{C}^4$.  
\end{remark}

Theorem \ref{S3formulation} now gives us a method for constructing 
CMC $H$ surfaces in $\mathbb{S}^3$ from given data 
$u$ and $Q$.  We must first choose a solution $F$ of 
the Lax pair system \eqref{LaxPairS3} with matrices as in 
\eqref{LaxPairMatricesS3}.  We then insert $F$ into the 
Sym-Bobenko-type formula in Equation \eqref{SymS3} and a CMC 
$H$ surface in $\mathbb{S}^3$, conformally parametrized by $z$ and 
with metric and Hopf differential that are constant nonzero multiples of 
$4 e^{2 u}dzd\bar z$ and $Q$ respectively, is produced.  

The advantage of the formulation in Theorem \ref{S3formulation} is 
that the Lax pair 
\eqref{LaxPairS3}-\eqref{LaxPairMatricesS3} is precisely the same as that 
in \eqref{LaxPair}-\eqref{eq:UandV} (with $H$ chosen to be $1/2$ there).  
Therefore, with no modification at all, we 
can use the DPW method to produce general solutions $F$, and then 
Theorem \ref{S3formulation} gives all CMC 
$H$ surfaces in $\mathbb{S}^3$.  We will later see how this 
approach can be applied to make Delaunay surfaces in $\mathbb{S}^3$ and 
also to formulate period problems for non-simply-connected CMC surfaces 
in $\mathbb{S}^3$.  

This theorem (and the theorem in the next section as well) 
are proven in \cite{b1}, but we include proofs here because we are using 
different notation than that of \cite{b1}.  

\begin{theorem}\label{S3formulation}
Let $\Sigma$ be a simply-connected domain in $\mathbb{C}$ with complex 
coordinate $z$.  
Choose $\gamma_1,\gamma_2 \in \mathbb{R}$ so that $(\pi^{-1}) (\gamma_2 - 
\gamma_1)$ is not an integer.  
Let $u$ and $Q$ solve \eqref{maurer-cartan} and let $F=F(z,\bar z,
\lambda)$ be a solution of the system 
\begin{equation}\label{LaxPairS3} F_z = FU \; , \;\;\; F_{\bar z}= FV \; , 
\end{equation} with 
\begin{equation}\label{LaxPairMatricesS3} 
U = \frac{1}{2} \begin{pmatrix} u_z & - \lambda^{-1} e^{u} \\
\lambda^{-1} Q e^{-u} & -u_z \end{pmatrix} \; , \;\;\; 
V = \frac{1}{2} \begin{pmatrix} -u_{\bar z} & - \lambda {\bar Q} 
e^{-u} \\ \lambda e^{u} & u_{\bar z} \end{pmatrix} \; , \end{equation}
Suppose that $det (F) = 1$ for all $\lambda$ and $z$.  
Defining $F_j = F|_{\lambda=e^{i\gamma_j}}$, $j=1,2$, 
suppose also that $F_j = \sigma_2 \bar F_j \sigma_2$ (that is, 
``$F_j \in \SU_2 = \mathbb{S}^3$'').  Define 
\begin{equation}\label{SymS3} f= 
F_1 \begin{pmatrix} e^{i(\gamma_1-\gamma_2)/2} 
& 0 \\ 0 & e^{i(\gamma_2-\gamma_1)/2} \end{pmatrix} F_2^{-1} 
\end{equation} 
and 
\begin{equation}\label{eq:normalS3}
N=iF_1 \begin{pmatrix} e^{i(\gamma_1-\gamma_2)/2} 
& 0 \\ 0 & -e^{i(\gamma_2-\gamma_1)/2} \end{pmatrix} F_2^{-1} \; . 
\end{equation}
Then $f$ satisfies $\mathcal{F}_w=\mathcal{F}{\mathcal U}$, 
$\mathcal{F}_{\bar w}=\mathcal{F}{\mathcal V}$ and is a surface 
in $\mathbb{S}^3$ with CMC $H=\cot (\gamma_2-\gamma_1)$ and normal $N$.  
\end{theorem}

\begin{proof}
First we note that $f \in \SU_2 = \mathbb{S}^3$ for all $z$.  
Then we note that 
\[ \langle f_z,N\rangle = \langle f_{\bar z},N\rangle = \langle f,N\rangle 
  = \langle f_z,f_z\rangle = \langle f_{\bar z},f_{\bar z}\rangle = 0 
 \; , \;\;\; \langle N,N\rangle =1 \; , \] 
and 
\[ 
\langle f_z,f_{\bar z}\rangle = \frac{1}{2} e^{2 u} \sin^2(\gamma_2-\gamma_1) 
\; , \;\;\; 
\langle f_{zz},N\rangle = \frac{1}{2} Q e^{-i(\gamma_1+\gamma_2)} 
\sin (\gamma_2-\gamma_1) \; , \]
and 
\[ f_{z\bar z}= \frac{-1}{2} e^{2 u} \sin^2 (\gamma_2-\gamma_1) f + 
\frac{1}{2} e^{2 u} \sin (\gamma_2-\gamma_1) \cos (\gamma_2-\gamma_1) N \; .\]
If $H$ is given by $H=\cot (\gamma_2-\gamma_1)$, $z=2\sqrt{1+H^2} w = 
2 w / \sin (\gamma_2-\gamma_1)$, $A=\langle f_{ww},N\rangle 
= (2 \sqrt{1+H^2})^2 \langle f_{zz},N\rangle 
= 2 Q e^{-i(\gamma_1+\gamma_2)}/\sin (\gamma_2-\gamma_1)$, then 
$\mathcal{F}_w=\mathcal{F}{\mathcal U}$, 
$\mathcal{F}_{\bar w}=\mathcal{F}{\mathcal V}$ hold, proving the theorem.  
\end{proof}

\subsection{Rotations in $\mathbb{S}^3$}
\label{subsectionS3}

In preparation for Sections \ref{section13} and \ref{spaceformdelaunay}, 
we briefly give a description of rotations of $\mathbb{S}^3$ in terms of 
$2 \times 2$ matrices here.
For points $X \in \mathbb{R}^4$ as in \eqref{R4form}, we claim that 
\[X \rightarrow A_1 \cdot X \cdot A_2^{-1}\] 
represents a general rotation of $\mathbb{R}^4$ fixing the origin, 
equivalently a general rotation of $\mathbb{S}^3$, where
 \[A_1=
     \begin{pmatrix}
     c & d\\
     -\bar{d}& \bar{c}
     \end{pmatrix}
\;,\;\;\;A_2^{-1}=
         \begin{pmatrix}
         e & f \\
         -\bar{f} & \bar{e}
         \end{pmatrix}
\in \SU_2\;.\]
for $c=c_1+i c_2$ , $d=d_1+i d_2$ , $e=e_1+i e_2$ , $f=f_1+i f_2$ with $c_j,d_j,e_j,f_j 
\in \mathbb{R}$.  
To see this, writing the point $ X=(a_1,b_2,b_1,a_2)^t\in \mathbb{R}^4$ in vector form, 
this map  $X \rightarrow A_1 \cdot X \cdot A_2^{-1}$ translates into $ X \rightarrow 
R\cdot X $ in the vector formulation for $\mathbb{R}^4$, where
\[
  R=\begin{pmatrix}
    \Re (ce-\bar{d}f) & -\Im (de+\bar{c}f) & -\Re(de+\bar{c}f) & \Im(\bar{d}f-ce) \\    
    \Im (cf-\bar{d}e)& \Re (\bar{c}e+df) & -\Im(\bar{c}e+df) & \Re(cf-\bar{d}e) \\
    \Re(\bar{d}e+cf) & \Im(\bar{c}e-df) & \Re(\bar{c}e-df) & -\Im(\bar{d}e+cf) \\
    \Im (ce+\bar{d}f) & \Re(de-\bar{c}f) & \Im(\bar{c}f-de) & \Re(ce+\bar{d}f)
    \end{pmatrix}\;.
\]
One can check that $R \in \text{SO}_4$, implying that $X \rightarrow A_1 \cdot X 
\cdot A_2^{-1}\ $ does in fact represent a rotation of $\mathbb{S}^3$. 

As an example, in the special case that $ A_1=A_2$, we have
\[
  R=\begin{pmatrix}
    1 & 0 & 0 & 0 \\
    0 & c_1^2-c_2^2-d_1^2+d_2^2 & 2(c_1c_2+d_1d_2) & -2c_1d_1+2c_2d_2 \\ 
    0 & -2c_1c_2+2d_1d_2 & c_1^2-c_2^2+d_1^2-d_2^2 & 2(c_2d_1+c_1d_2) \\ 
    0 & 2(c_1d_1+c_2d_2) & 2c_2d_1-2c_1d_2 & c_1^2+c_2^2-d_1^2-d_2^2
    \end{pmatrix}\;.
\]
Thus the lower-right $3 \times 3$ cofactor matrix is in $\text{SO}_3$ and is the same as 
the matrix in the proof of Lemma \ref{lm:rigid_motion} (up to a sign change of $ c_1 $ 
and $ d_2 $). Thus this rotation of $\mathbb{R}^4$ restricts to a rotation of the 
$3$-dimensional Euclidean space $\{(0,x_2,x_3,x_4)|x_j\in \mathbb{R} \}$.  

\begin{figure}[phbt]
\begin{center}
\includegraphics[width=0.4\linewidth]{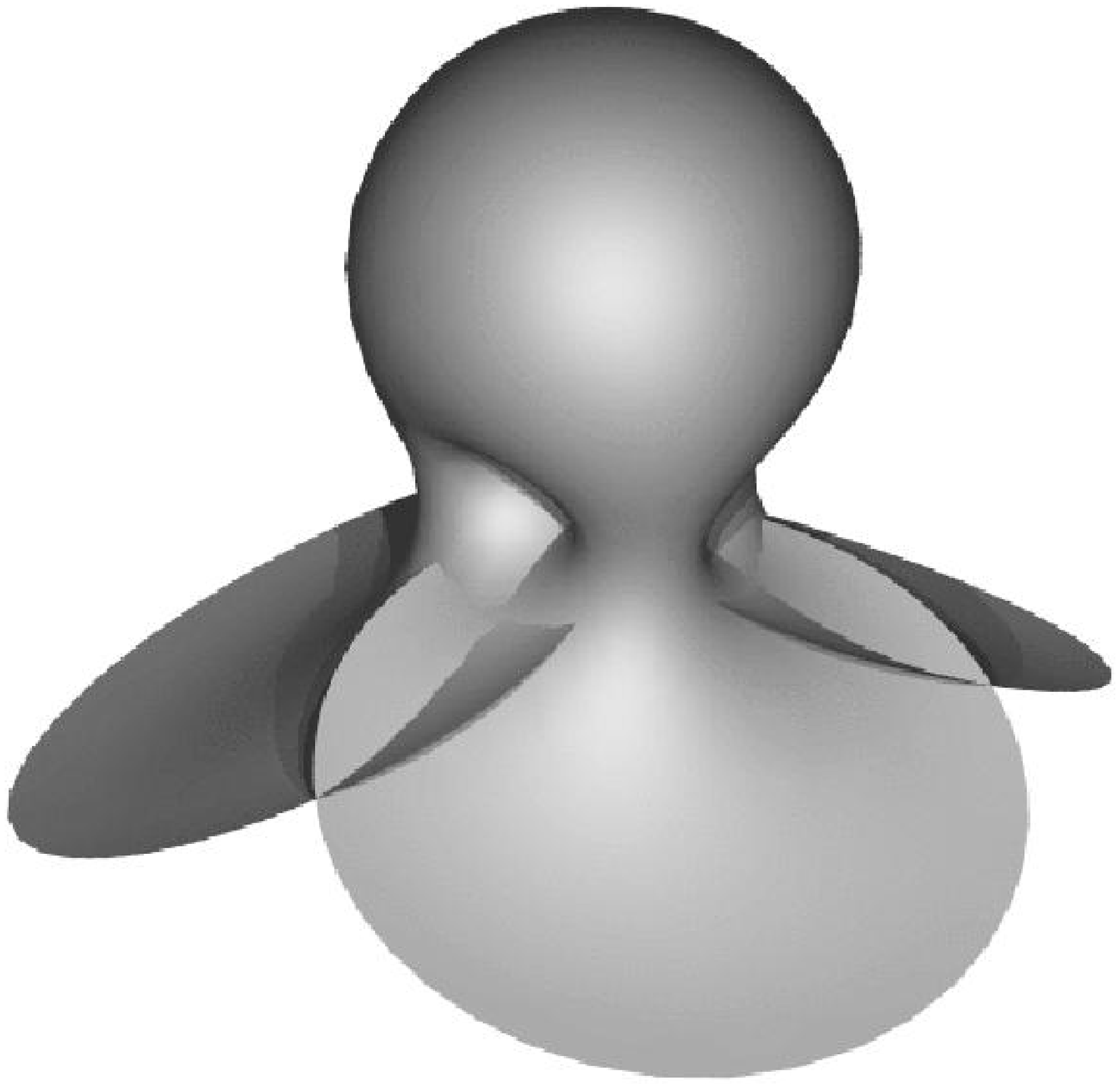}
\includegraphics[width=0.4\linewidth]{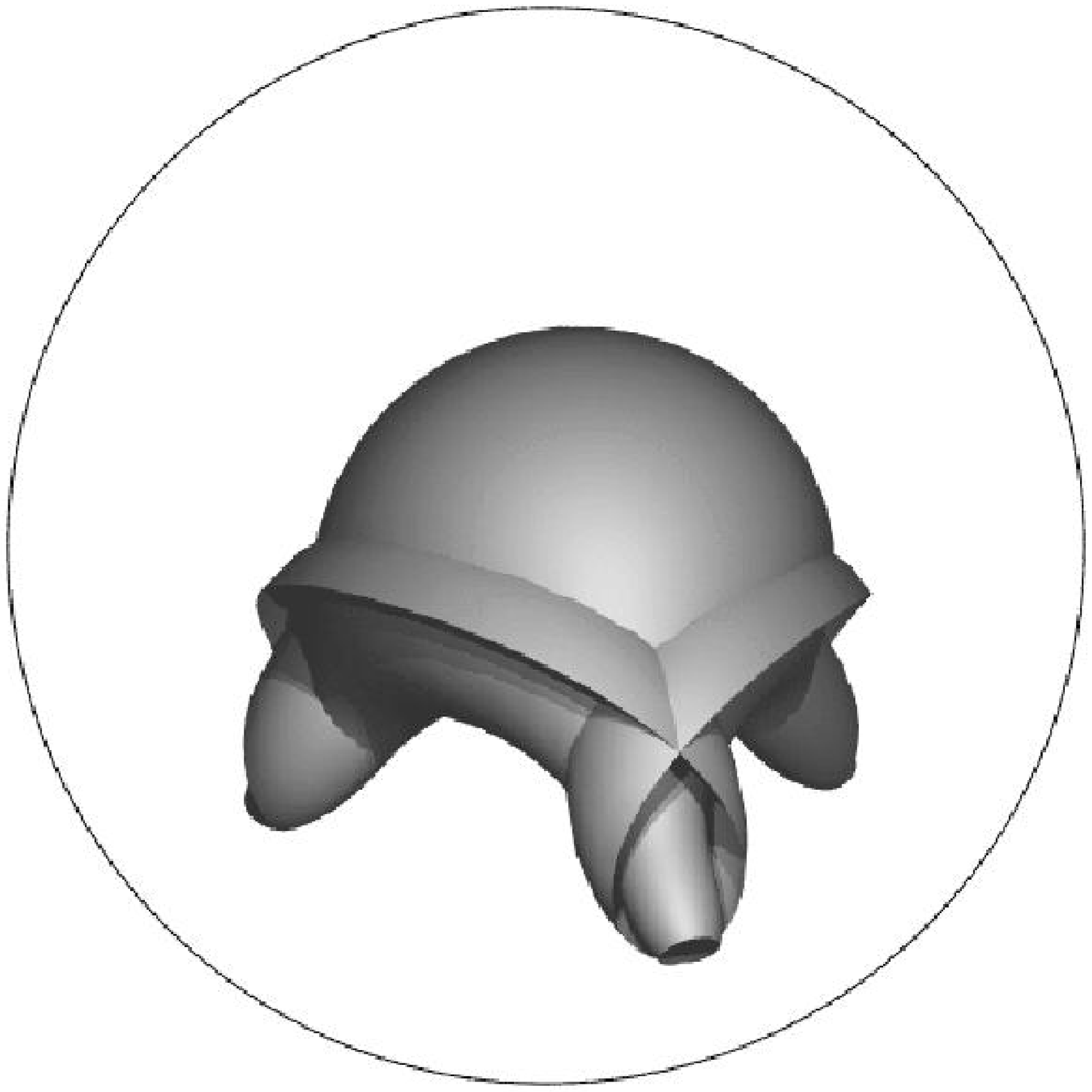}
\end{center}
\caption{Smyth surfaces in $\mathbb{S}^3$ and $\mathbb{H}^3$. 
To see a surface in $\mathbb{S}^3$, we use a stereographic projection
$$\hspace*{-60pt}
\mathbb{R}^4\supset\mathbb{S}^3\ni (x_1,x_2,x_3,x_4)\mapsto
\dfrac{1}{1-x_4}\left(x_1,x_2,x_3\right)\in\mathbb{R}^3\cup\{\infty\} \; .
$$
Also, to view the surface in $\mathbb{H}^3$, we use the Poincare model.  
(These graphics were made using N. Schmitt's software CMCLab \cite{Sch:cmclab}.)}
\label{fg:SmythS3H3}
\end{figure}

\section{The Lax pair for CMC $H>1$ surfaces in $\mathbb{H}^3$}
\label{section12}

As described in Section \ref{hyperbolic3space}, 
the Minkowski model for hyperbolic $3$-space $\mathbb{H}^3$ is the space 
\[ \{ (x_1,x_2,x_3,x_0) \in \mathbb{R}^{3,1} \, | \, x_1^2+x_2^2+x_3^2-x_0^2=
-1 \, , \, x_0>0 \} \; , \]  
with metric induced from the $4$-dimensional Lorentz space 
$\mathbb{R}^{3,1}= \{ (x_1,x_2,x_3,x_0) \in \mathbb{R}^4 \}$ 
with the Lorentz metric 
\[ \langle (x_1,x_2,x_3,x_0), (y_1,y_2,y_3,y_0) \rangle_{\mathbb{R}^{3,1}} = 
x_1y_1+x_2y_2+x_3y_3-x_0y_0 \; . \]  
($\mathbb{H}^3$ is a space form because it is the unique simply-connected 
complete $3$-dimensional Riemannian manifold with constant sectional curvature $-1$, 
see Lemma \ref{H3hascurvatureminus1}.) 
Let $f:\Sigma\to\mathbb{H}^3$ be a conformally immersed 
surface in $\mathbb{H}^3$, where $\Sigma$ is a simply-connected Riemann 
surface with complex coordinate $w$.  
(By Theorem \ref{conformalityispossible}, any immersion $f:\Sigma\to\mathbb{H}^3$ can be 
reparametrized so that it is 
conformal, so without loss of generality we assume here that $f$ is conformal.)
Then $\langle f , f \rangle_{\mathbb{R}^{3,1}} =-1$.  
Let $N \in \mathbb{R}^{3,1}$ be a unit normal of $f$ in $\mathbb{H}^3$, then 
$N$ actually lies in de Sitter space, because it has norm $+1$.  
De Sitter $3$-space $\mathbb{S}^{2,1}$ is the Lorentzian manifold that is 
the $1$-sheeted hyperboloid in $\mathbb{R}^{3,1}$, 
\[ \mathbb{S}^{2,1} = \left\{ (x_1,x_2,x_3,x_0) \in \mathbb{R}^{3,1} \, \left| \, 
\sum_{j=1}^{3} x_j^2 - x_0^2 = 1 \right. \right\} \] with 
the metric $g$ induced on its tangent 
spaces by the restriction of the metric from the Minkowski space $\mathbb{R}^{3,1}$.  
It is homeomorphic to $\mathbb{S}^2 \times \mathbb{R}$, 
so is simply-connected (since both $\mathbb{S}^2$ and $\mathbb{R}$ are individually 
simply-connected), and has constant sectional curvature $1$.  

Then, by the definition 
of the tangent space $T_f \mathbb{H}^3$, and by the definition of $N$, and by the 
conformality of $f$, and because $\langle f,f \rangle_{\mathbb{R}^{3,1}}=-1$, we have 
\begin{eqnarray*}
&  0=\langle N,f\rangle =\langle N,f_w\rangle =\langle N,f_{\bar w}\rangle 
    =\langle f_w,f_w\rangle =\langle f_{\bar w},f_{\bar w}\rangle 
    =\langle f,f_w\rangle =\langle f,f_{\bar w}\rangle \; , & \\
&  \langle N,N\rangle =1 \; , \;\;\; 
\langle f_w,f_{\bar w}\rangle =2e^{2 u} \; , & 
\end{eqnarray*} 
where $u:\Sigma\to\mathbb{R}$ is defined this way, and 
$\langle \cdot, \cdot \rangle$ is the bilinear extension of 
$\langle \cdot, \cdot \rangle_{\mathbb{R}^{3,1}}$.  
Let $\mathcal{F}=(f,f_w,f_{\bar w},N)$ be the frame for $f$, 
and define the Hopf differential $A$ and mean curvature $H$ of $f$ by 
\[ A = \langle f_{ww},N\rangle \; , \;\;\; 
2 H e^{2 u} =\langle f_{w\bar w},N\rangle \; . \] 
(We again reserve the symbol $Q$ for a slightly different purpose later.)  
Similar to the previous section and Section \ref{section6}, we have 
\[ \mathcal{F}_w = \mathcal{F}{\mathcal U} \; , \;\;\; 
\mathcal{F}_{\bar w} = \mathcal{F}{\mathcal V} \; , \] where 
\[ {\mathcal U}=\begin{pmatrix}
0 & 0 & 2 e^{2 u} & 0 \\
1 & 2 u_w & 0 & -H \\
0 & 0 & 0 & -A e^{-2 u}/2 \\
0 & A & 2 H e^{2 u} & 0 
\end{pmatrix} \; , \;\;\; 
{\mathcal V}=\begin{pmatrix}
0 & 2 e^{2 u} & 0 & 0 \\
0 & 0 & 0 & -\bar A e^{-2 u}/2 \\
1 & 0 & 2 u_{\bar w} & -H \\
0 & 2 H e^{2 u} & \bar A & 0 
\end{pmatrix} \; . \] 
The compatibility condition is 
\[ {\mathcal U}_{\bar w}-{\mathcal V}_w-[{\mathcal U},{\mathcal V}]=0 \; , \] 
which implies that 
\begin{equation}\label{originalH3MC} 
2 u_{w\bar w}+2 e^{2 u} (H^2-1) -\frac{1}{2} A \bar A e^{-2 u} = 0 
\; , \;\;\; A_{\bar w}= 2H_we^{2u} \; . \end{equation} 
Now suppose that $H>1$ is constant, and thus $A_{\bar w}=0$.  
Making the change of variables $z=2\sqrt{H^2-1} w$ and $A=2\sqrt{H^2-1} 
e^{2i\psi} Q$ for a real constant $\psi$, we again have Equation 
\eqref{maurer-cartan}, which we again note is identical to Equation 
\eqref{firstMC}, with the $H$ there replaced by $1/2$.  
We conclude that a CMC $H$ conformal immersion into $\mathbb{H}^3$ with 
$H>1$ has a 
conformal parameter $z$ such that $Q=2 \sqrt{H^2-1} e^{-2 i \psi} 
\langle f_{zz},N\rangle$ is holomorphic and 
$u$ satisfies the first equation in \eqref{maurer-cartan}.  (The factor 
$e^{-2i\psi}$ in $Q$ will again play the role of the spectral parameter.)

We now consider how to get a surface from given $u$, $Q$, and $H$.  
We can consider $\mathbb{R}^{3,1}$ to be the $2 \times 2$ self-adjoint 
matrices ($\bar X^t = X$), by the identification 
\begin{equation} 
X=(x_1,x_2,x_3,x_0) \in \mathbb{R}^{3,1} \leftrightarrow 
X=x_0 \cdot \id + \sum_{k=1}^3 x_k \sigma_k 
 =\begin{pmatrix} x_0+x_3 & x_1-i x_2 \\ 
                  x_1+i x_2 & x_0-x_3 \end{pmatrix}\; ,
\label{eq:L4identifyHerm2}
\end{equation}
where $\sigma_1$, $\sigma_2$, $\sigma_3$ are the Pauli matrices as in 
Section \ref{2by2laxpair}.  See also Section \ref{hyperbolic3space}.  

Since, for any $2 \times 2$ matrix $Y$, when $\det Y \ne 0$, 
\begin{equation}\label{littleYrelation}
\sigma_2 Y^t \sigma_2 = Y^{-1} \det Y \; , \end{equation} and 
since, for $X,Y \in \mathbb{R}^{3,1}$ (i.e. $\bar X^t=X$, $\bar Y^t=Y$), we have 
\begin{equation}\label{littleYrelationA}
\langle X,Y\rangle_{\mathbb{R}^{3,1}} = (-1/2) \text{trace} 
(X \sigma_2 Y^t \sigma_2) \; , \end{equation} and we have 
\begin{equation}\label{littleYrelationB}
\langle X,X\rangle_{\mathbb{R}^{3,1}} = (-1/2) \text{trace} 
\begin{pmatrix}
\det X & 0 \\ 0 & \det X
\end{pmatrix} = - \det X \; . \end{equation}

Then $\mathbb{H}^3$ is 
\[ \mathbb{H}^3 = \{ X \, | \, \bar X^t = X \, , 
\langle X,X\rangle_{\mathbb{R}^{3,1}} =-1 \, , \, \text{trace} (X) > 0 \} \; . \] 

\begin{lemma}\label{lm:SL2Cform}
$\mathbb{H}^3$ can be written as 
$\{ F \overline{F}^t \, | \, F \in \SL_2\!\mathbb{C} \}$ (see also 
Section \ref{hyperbolic3space}).  
\end{lemma}

\begin{proof}
We must show 
\[
  \{ X \, | \, \overline{X}^t = X \, , \det X =1 
            \, , \, \text{trace} (X) > 0 \}
= \{ F \overline{F}^t \, | \, F \in \SL_2\!\mathbb{C} \} \; .
\]
Any element of the right-hand side is obviously an element of the left-hand 
side.  We show the converse.  For any $X$ in the left-hand side, 
let $\lambda$ and $v=(v_1, v_2)^t$ be an eigenvalue and corresponding 
eigenvector of $X$, with $|v_1|^2+|v_2|^2=1$.  
Because $X=\overline{X}^t$, $\lambda$ must be real.  
Using $X=\overline{X}^t$ and $\det X=1$, one can see that $1/\lambda$ and 
$\tilde v=(-\overline{v_2}, \overline{v_1})^t$ are another 
eigenvalue-eigenvector pair of X.  
Let $E$ be the $2\times 2$ matrix $E=(v, \tilde v)$, then 
\[
XE=E\begin{pmatrix}\lambda&0\\0&1/\lambda\end{pmatrix} \; ,
\]
so $X=F\overline{F}^t$, where 
\[
F = E\begin{pmatrix}\sqrt{\lambda}&0\\0&1/\sqrt{\lambda}\end{pmatrix}
\in \SL_2\!\mathbb{C} \; ,
\]
because $\lambda >0$ (since $\text{trace}X>0$).  
\end{proof}

Just as was the case for Theorem \ref{S3formulation} when the ambient space 
was $\mathbb{S}^3$, Theorem \ref{H3formulation} gives us a method for 
constructing CMC $H$ surfaces in $\mathbb{H}^3$ from given data 
$u$ and $Q$.  A CMC $H>1$ surface in $\mathbb{H}^3$ is produced, conformally 
parametrized by $z$ and with metric and Hopf differential that are again 
nonzero constant multiples of $4 e^{2 u}dzd\bar z$ and $Q$, respectively.  

And once again we have the advantage that the Lax pair 
\eqref{LaxPairS3}-\eqref{LaxPairMatricesS3} is precisely the same as that 
in \eqref{LaxPair}-\eqref{eq:UandV} (with $H$ chosen to be $1/2$ there).  
Therefore we can again use the DPW method to produce general solutions $F$, 
and then Theorem \ref{H3formulation} gives all CMC 
$H$ surfaces in $\mathbb{H}^3$ for any $H>1$.  Together with the case of 
$\mathbb{S}^3$, we will later see how this 
approach can be applied to make Delaunay surfaces in $\mathbb{H}^3$ and 
also to formulate period problems for non-simply-connected CMC $H>1$ surfaces 
in $\mathbb{H}^3$.  

\begin{theorem}\label{H3formulation}
Let $\Sigma$ be a simply-connected domain in $\mathbb{C}$ with complex 
coordinate $z$.  
Choose $\psi \in \mathbb{R}$ and $q \in \mathbb{R} \setminus \{ 0 \}$.  
Let $u$ and $Q$ solve \eqref{maurer-cartan} and let $F=F(z,\bar z,
\lambda)$ be a solution of the system 
\eqref{LaxPairS3}-\eqref{LaxPairMatricesS3}.  
Suppose that $det (F) = 1$ for all $\lambda$ and $z$.  
Set $F_0=F|_{\lambda=e^{q/2} e^{i \psi}}$.  Then defining 
\begin{equation}\label{H3Symformula} 
f(z,\bar z)=\tilde{F} \overline{\tilde{F}}^t \; , \;\;\; \tilde{F} := 
F_0 \begin{pmatrix} e^{q/4} & 0 
\\ 0 & e^{-q/4} \end{pmatrix} \; , \;\;\; 
N=\tilde{F} \sigma_3 \overline{\tilde{F}}^t \; , \end{equation} 
$f$ is a CMC $H=\coth (-q)$ surface in $\mathbb{H}^3$ with normal $N$.  
\end{theorem}

\begin{proof}
By properties \eqref{littleYrelation} and \eqref{littleYrelationA} and 
\eqref{littleYrelationB}, we have 
\[ \langle N, f\rangle =\langle N, f_z\rangle =\langle N, f_{\bar z}\rangle 
=\langle f_z,f_z\rangle =\langle f_{\bar z},f_{\bar z}\rangle =0 \]  
and 
\[ \langle N,N\rangle =1 \; , \;\;\; 
\langle f_z,f_{\bar z}\rangle = (1/2) e^{2 u} \sinh^2 (-q) \; . \]
Also, 
\[ \langle f_{zz},N\rangle = \frac{1}{2} Q e^{-2i\psi} \sinh (-q) \; , \]
and 
\[ f_{z\bar z}= \frac{1}{2} e^{2 u} \sinh^2 (-q) f + 
\frac{1}{2} e^{2 u} \sinh (-q) \cosh (-q) N \; . \] 
With $z=2\sqrt{H^2-1} w$, $A=2 \sqrt{H^2-1} e^{-2i\psi} Q$, and 
$\sinh (-q) = 1/\sqrt{H^2-1}$, then $A=\langle f_{ww},N\rangle$ and 
$f_{w \bar w} = 2e^{2 u} f + 2 e^{2 u} \coth (-q) N$.  So with $H=\coth (-q)$, 
then with the above and similar computations, we get that 
$\mathcal{F}_w=\mathcal{F}{\mathcal U}$, 
$\mathcal{F}_{\bar w}=\mathcal{F}{\mathcal V}$ hold, proving the theorem.  
\end{proof}

\subsection{Rigid motions in $\mathbb{H}^3$}
\label{subsectionH3}

In preparation for Sections \ref{section13} and \ref{spaceformdelaunay}, 
we give a description of rigid motions of $\mathbb{H}^3$ in terms of 
$2 \times 2$ matrices here.
For points $X \in \mathbb{R}^{3,1}$, we claim that 
\begin{equation}\label{H3rigidmotions}
X \to \bmath{A} \cdot X \cdot \overline {\bmath{A}}^{t}
\end{equation} 
represents a general rotation of $\mathbb{R}^{3,1}$ fixing the origin, or 
equivalently a general rigid motion of $\mathbb{H}^3$, where
\[
\bmath{A}=\begin{pmatrix}a&b\\c&d\end{pmatrix}\in\SL_2\!\mathbb{C}
\]
for $a=a_1+i a_2$, $b=b_1+i b_2$, $c=c_1+i c_2$, $d=d_1+i d_2$ with 
$a_j,b_j,c_j,d_j \in \mathbb{R}$.  
To see this, writing the point $X=(x_1,x_2,x_3,x_0)^t\in \mathbb{R}^{3,1}$ 
in vector form, the map 
$X\to \bmath{A}\cdot X\cdot\overline{\bmath{A}}^{t}$ translates into 
$X\to R\cdot X$ in the vector formulation for $\mathbb{R}^{3,1}$, where
\[
R=\begin{pmatrix}
    \Re (\bar{a} d+\bar{b} c) & \Im (\bar{b} c-\bar{a} d) & 
    \Re(\bar{a} c-\bar{b} d) & \Re(\bar{a} c+\bar{b} d) \\
    \Im (\bar{b} c+\bar{a} d)& \Re (\bar{a} d-\bar{b} c) & 
    \Im(\bar{a} c-\bar{b} d) & \Im(\bar{a} c+\bar{b} d) \\
    \Re(\bar{a} b-\bar{c} d) & \Im(a \bar{b}+\bar{c} d) & 
    (a \bar{a}-b \bar{b}-c \bar{c}+d \bar{d})/2 & 
    (a \bar{a}+b \bar{b}-c \bar{c}-d \bar{d})/2 \\
    \Re (\bar{a} b+\bar{c} d) & \Im(a \bar{b}+c \bar{d}) & 
    (a \bar{a}-b \bar{b}+c \bar{c}-d \bar{d})/2 & 
    (a \bar{a}+b \bar{b}+c \bar{c}+d \bar{d})/2
  \end{pmatrix} \; .
\]

Recalling that 
\[
\text{SO}_{3,1}^+:=\left\{R \in\SL_4\mathbb {R} \;|\;
R^t\cdot I_{3,1}\cdot R=I_{3,1}\;\;,\;\;
\underline{e_0}^t\cdot R\cdot\underline{e_0}>0\right\}\; , 
\]
where 
\[
I_{3,1}=\begin{pmatrix}1 & 0 & 0 & 0 \\ 0 & 1 & 0 & 0 \\
                       0 & 0 & 1 & 0 \\ 0 & 0 & 0 & -1\end{pmatrix}
\;\;\;\text{and}\;\;\;
\underline{e_0}=\begin{pmatrix}0\\0\\0\\1\end{pmatrix} \; , 
\] 
one can check that $R \in\text{SO}_{3,1}^+$, implying that 
\eqref{H3rigidmotions} does in fact represent a 
rotation of $\mathbb R^{3,1}$.  

Another way to argue that \eqref{H3rigidmotions} is an isometry is by 
computing directly with $2 \times 2$ matrices: 
\[ \langle \bmath{A} X \overline{\bmath{A}}^{t} , 
\bmath{A} Y \overline{\bmath{A}}^{t} \rangle_{\mathbb{R}^{3,1}} = 
- \tfrac{1}{2} \text{trace} (\bmath{A} X \overline{\bmath{A}}^{t} 
\sigma_2 (\bmath{A} Y \overline{\bmath{A}}^{t})^t \sigma_2) \] 
can be directly computed to give 
\[ \langle \bmath{A} X \overline{\bmath{A}}^{t} , 
\bmath{A} Y \overline{\bmath{A}}^{t} \rangle_{\mathbb{R}^{3,1}} = 
\langle X,Y \rangle_{\mathbb{R}^{3,1}} \; , \] 
implying that \eqref{H3rigidmotions} is an isometry.  

As an example, in the special case that $a=e^{\theta/2}$, $b=0$, $c=0$ and 
$d=e^{-\theta/2}$ for $\theta\in\mathbb{R}$, we have
\[
  R=\begin{pmatrix}
    1 & 0 & 0 & 0 \\
    0 & 1 & 0 & 0 \\ 
    0 & 0 & \cosh{\theta} & \sinh{\theta} \\ 
    0 & 0 & \sinh{\theta} & \cosh{\theta}
    \end{pmatrix} \; .  
\]
The lower-right $2 \times 2$ submatrix 
\[
\begin{pmatrix}
\cosh \theta & \sinh \theta \\ \sinh \theta & \cosh \theta 
\end{pmatrix}
\]
is the standard $2 \times 2$ Lorentz transformation (or velocity boost) 
in the theory of special relativity (see \cite{Callahan}, for example). 

As another example, if $c=-\bar{b}$ and $d=\bar{a}$, then 
\[
R = \begin{pmatrix}
        a_1^2-a_2^2-b_1^2+b_2^2 & 2(a_1a_2+b_1b_2) & -2a_1b_1+2a_2b_2 & 0 \\ 
        -2a_1a_2+2b_1b_2 & a_1^2-a_2^2+b_1^2-b_2^2 & 2(a_2b_1+a_1b_2) & 0 \\ 
         2(a_1b_1+a_2b_2) & 2a_2b_1-2a_1b_2 & a_1^2+a_2^2-b_1^2-b_2^2 & 0 \\
         0 & 0 & 0 & 1
    \end{pmatrix} \; .  
\]
Then the transformation $X\to R\cdot X$ fixes the axis 
$\{(0,0,0,r)|r\in \mathbb{R} \}$ in $\mathbb{R}^{3,1}$.  
Furthermore, the upper-left $3 \times 3$ cofactor matrix is in $\text{SO}_3$ 
and is the same as the matrix in the proof of Lemma \ref{lm:rigid_motion} 
(up to a sign change of $c_1$ and $d_2$).  

\section{Period problems in $\mathbb{R}^3$, $\mathbb{S}^3$, $\mathbb{H}^3$}
\label{section13}

When one starts with a CMC $H$ conformal immersion $f$ into $\mathbb{R}^3$ or 
$\mathbb{S}^3$ or $\mathbb{H}^3$ defined on a simply-connected domain $\Sigma$, 
and then extends $f$ to a conformal CMC immersion $\hat f$ on a larger 
non-simply-connected domain $\hat \Sigma$ (that is, $\Sigma \subset 
\hat \Sigma$ and $\hat f |_{\Sigma} = f$), the extension $\hat f$ 
will be unique.  However, it is not necessarily true that $\hat f$ is 
well-defined on $\hat \Sigma$.  The extended immersion $\hat f$ being well-defined 
on $\hat \Sigma$ is equivalent to $\hat f$ being well-defined on every closed 
loop $\delta$ in $\hat \Sigma$, and this property can be studied using the DPW method.  
Here we consider some conditions that imply $\hat f$ is 
well-defined on $\hat \Sigma$.  

We first assume the following: 

(1) $\hat f$ is constructed via the DPW method.  

(2) $\delta:[0,1] \to \hat \Sigma$ is a closed loop, and $\tau$ is the 
deck transformation associated to 
$\delta$ on the universal cover $\widetilde{\Sigma}$ of $\hat \Sigma$.  

(3) The holomorphic potential $\xi$ used to construct $\hat f$ is chosen to 
be well-defined on $\hat \Sigma$.  

(4) The holomorphic potential $\xi$ and solution $\phi$ on 
$\widetilde{\Sigma}$ of $d\phi=\phi \xi$ used to construct $\hat f$ are chosen 
so that the monodromy matrix $M_\delta$ of $\phi$, 
\begin{equation}\label{Mdelta} 
M_\delta = \phi(\tau(z),\lambda) \cdot (\phi(z,\lambda))^{-1} \; , \end{equation} 
satisfies $M_\delta \in \Lambda \SU_2$ 
(here $z$ is a local coordinate on $\hat \Sigma$).  

Note that although $M_\delta=M_\delta(\lambda)$ depends on $\lambda$, 
it is independent of $z$, since 
\[ dM_\delta = d(\phi(\tau(z),\lambda) (\phi(z,\lambda))^{-1}) = 0 \; , \] 
which holds because $\xi$ is well-defined on $\hat \Sigma$ and both 
$\phi(\tau(z),\lambda)$ and $\phi(z,\lambda)$ are solutions of the equation 
$d\phi=\phi \xi$.  

The monodromy $M_\delta$ is independent of the choice of loop within the 
homotopy class of $\delta$.  Also, changing the initial condition used to 
determine $\phi$ results in a conjugation of $M_\delta$.  

The first two conditions above are simply definitions.  The third condition 
can be satisfied by any noncompact CMC immersion, as we saw in 
Lemma \ref{lm:well-definedness}.  
The fourth condition is a rather strong condition, but there are many 
examples for which this condition can hold.  We have in fact already seen 
one such example, when we constructed Delaunay surfaces in 
Section \ref{section5}.  

\begin{remark}
If $\hat \Sigma$ is compact, then the third assumption cannot be satisfied,
i.e. holomorphic potentials cannot be chosen to be well-defined 
on $\hat \Sigma$. Therefore, in general, we cannot define a monodromy 
matrix $M_{\delta}$ of $\phi$ simply as in Equation \eqref{Mdelta}.  
However, there is a more general way to define a monodromy matrix of 
$\phi$.  Detailed arguments and definitions can be found in \cite{Dorf-Kil}.
In general, the fourth assumption is also not satisfied, 
i.e. $M_{\delta}$ is not in $\Lambda \SU_{2}$ in general.  
However, there are several methods for forcing a 
monodromy matrix to be in $\Lambda \SU_{2}$. Detailed arguments 
can be found in \cite{Dorf-Kil}, \cite{KilKRS} and \cite{DorH:cyl}.  
\end{remark}

Iwasawa decomposing $\phi$ into $\phi=F B$, note that $F$ as well is generally 
only defined on $\widetilde{\Sigma}$, and the fourth condition above 
implies that the monodromy matrix of $F$ is also $M_\delta$: 
\[ M_\delta = 
F(\tau(z),\overline{\tau(z)},\lambda) \cdot (F(z,\bar{z},\lambda))^{-1} \; . 
\]  
Now we consider what happens to the immersion $\hat f$ along $\delta$ in 
each of the three space forms.  
We will find necessary and sufficient conditions for $\hat f=\hat f(z,\bar z)$ 
to be well-defined on the loop $\delta$, that is, 
\[
\hat f(\tau (z),\overline{\tau (z)})=\hat f(z,\bar z) \; .
\]
This is something referred to as the ``period problem'' or 
``closing condition'' of $\hat f$ about $\delta$.  

{\bf $\mathbb{R}^3$ case, $H \neq 0$:} 
In the case of $\mathbb{R}^3$, applying the DPW method, $\hat f$ is 
produced from $F$ by the Sym-Bobenko formula \eqref{Sym-Bobenko}.  In this 
formula, the necessary and sufficient condition for $\hat f=\hat f(z,\bar z)$ 
to be well-defined along the loop $\delta$, that is, for 
${\hat f}(\tau(z),\overline{\tau(z)}) = {\hat f}(z,\bar{z})$, is that 
\begin{equation}\label{eq:monR3}
\left[{\hat f}(z,\bar{z})\right]_{\lambda =1}
 = \left[ M_\delta {\hat f}(z,\bar{z}) M_\delta^{-1} 
  -\frac{i\lambda}{2H} (\partial_\lambda M_\delta) 
   M_\delta^{-1} \right]_{\lambda=1} \; . 
\end{equation}
Since $\hat f$ depends on $z$ while $M_\delta$ does not, we show in 
Lemma \ref{lm:well-definednessR3} that this holds if and only if 
\begin{equation}\label{periodR3} 
M_\delta |_{\lambda=1}=\pm \id \;\;\;\text{and}\;\;\; 
\partial_\lambda M_\delta |_{\lambda=1}= 0 \; , 
\end{equation} 
that is, we need $(M_\delta,\partial_\lambda M_\delta )|_{\lambda=1} \in 
\SU_2 \times \su_2$ to be the identity element (up to sign).  
Forcing this to happen (with respect to one given homotopy class of loops 
represented by $\delta$) is a six-dimensional period problem, as 
the real dimension of the space $\SU_2 \times \su_2$ is six.  

For the cases of $\mathbb{S}^3$ and $\mathbb{H}^3$, we will see below that again 
the period problem is six-dimensional with respect to each homotopy class of 
loops.  

\begin{lemma}\label{lm:well-definednessR3}
The immersion $\hat f$ into $\mathbb{R}^3$ is well-defined along $\delta$ if 
and only if \eqref{periodR3} holds.
\end{lemma}

\begin{proof}
If \eqref{periodR3} holds, then clearly \eqref{eq:monR3} holds.  
Conversely, suppose \eqref{eq:monR3} holds.  Then $[\hat f]_{\lambda =1}$ and 
\[
\left[ M_\delta\cdot\hat f\cdot M^{-1}_\delta 
-\frac{i\lambda}{2H}\partial_\lambda M_\delta\cdot M^{-1}_\delta
\right]_{\lambda =1}
\]
have the same frames, which are represented by $F|_{\lambda =1}$ and 
$M_\delta F|_{\lambda =1}$ (and these representatives are determined uniquely 
up to sign).  Hence 
$F|_{\lambda =1}=(\pm M_\delta F)|_{\lambda =1}$ and so 
$M_\delta |_{\lambda =1}=\pm\id$.  Using $M_\delta |_{\lambda =1}=\pm\id$, 
\eqref{eq:monR3} then implies 
\[
\left.\left( \frac{i\lambda}{2H}\partial_\lambda M_\delta\cdot M^{-1}_\delta
\right)\right|_{\lambda=1}=0 \; , 
\]
so $\partial_\lambda M_\delta |_{\lambda =1}=0$.  
\end{proof}

{\bf $\mathbb{S}^3$ case:} Travelling about the loop $\delta$, 
$\hat f={\hat f}(z,\bar z)$ changes to 
\begin{equation}\label{eq:monS3}
{\hat f}(\tau(z),\overline{\tau(z)}) = M_\delta|_{\lambda=e^{i \gamma_1}} 
\cdot{\hat f}(z,\bar z)\cdot M_\delta^{-1}|_{\lambda=e^{i \gamma_2}} 
\end{equation}
in the Sym-Bobenko-type formula \eqref{SymS3}.  
So ${\hat f}(\tau(z),\overline{\tau(z)})={\hat f}(z,\bar z)$ is 
equivalent to 
\begin{equation}\label{periodS3} 
M_\delta|_{\lambda=e^{i \gamma_1}} = M_\delta|_{\lambda=e^{i \gamma_2}} = 
\pm \id \; , \end{equation} 
by Lemma \ref{lm:well-definednessS3}.  
Thus, to close the surface $\hat f$ about the loop $\delta$, we need 
\[
(M_\delta|_{\lambda=e^{i \gamma_1}}, 
M_\delta|_{\lambda=e^{i \gamma_2}}) \in \SU_2 \times \SU_2
\]
to become the identity element (up to sign).  
As $\SU_2 \times \SU_2$ is 
six-dimensional, so is the period problem.  

\begin{lemma}\label{lm:well-definednessS3}
The immersion $\hat f$ into $\mathbb{S}^3$ is well-defined along $\delta$ if 
and only if \eqref{periodS3} holds.
\end{lemma}

\begin{proof}
If \eqref{periodS3} holds, then clearly \eqref{eq:monS3} holds.  
Conversely, suppose \eqref{eq:monS3} holds.  
Then by Subsection \ref{subsectionS3}, $X\to M_\delta|_{\lambda =e^{i\gamma_1}}
\cdot X\cdot M_\delta^{-1}|_{\lambda =e^{i\gamma_2}}$ is a rigid motion of 
$\mathbb{S}^3$, and it must be the identity rigid motion, since $\hat f$ and 
$M_\delta|_{\lambda =e^{i\gamma_1}}\cdot\hat f\cdot
 M_\delta^{-1}|_{\lambda =e^{i\gamma_2}}$ are equal.  
By Subsection \ref{subsectionS3}, 
this is so if and only if \eqref{periodS3} holds.  
\end{proof}

{\bf $\mathbb{H}^3$ case, $H>1$:} Travelling about the loop $\delta$, 
$\hat f={\hat f}(z,\bar z)$ changes to 
\begin{equation}\label{eq:monH3} 
{\hat f}(\tau(z),\overline{\tau(z)}) = M_\delta|_{\lambda=e^{q/2}e^{i \psi}}
\cdot\hat f(z,\bar z)\cdot\overline{M_\delta|_{\lambda=e^{q/2}e^{i \psi}}}^t 
\end{equation}
in the Sym-Bobenko-type formula \eqref{H3Symformula}.  
So ${\hat f}(\tau(z),\overline{\tau(z)})={\hat f}(z,\bar z)$ is 
equivalent to 
\begin{equation}\label{periodH3} 
M_\delta|_{\lambda=e^{q/2}e^{i \psi}} = \pm \id \; , \end{equation}
by Lemma \ref{lm:well-definednessH3}.  
Thus, to close the surface $\hat f$ about the loop $\delta$, 
we need 
\[
M_\delta|_{\lambda=e^{q/2}e^{i \psi}} \in \SL_2\!\mathbb{C}
\]
to become the identity element (up to sign).  (In the case of $\mathbb{H}^3$, 
it is possible that $M_\delta|_{\lambda=e^{q/2}e^{i \psi}} 
\not\in \SU_2$, since $e^{q/2}e^{i \psi} \not\in \mathbb{S}^1$.)  
As $\SL_2\!\mathbb{C}$ is six-dimensional, so is the period problem.  

Just like as for Lemma \ref{lm:well-definednessS3}, 
we can use Subsection \ref{subsectionH3} to prove:

\begin{lemma}\label{lm:well-definednessH3}
The immersion $\hat f$ into $\mathbb{H}^3$ is well-defined along $\delta$ if 
and only if \eqref{periodH3} holds.
\end{lemma}

\begin{remark}
As noted in Remark \ref{re:monodromy}, the dressing $\hat\phi =h_+\cdot\phi$ 
of $\phi$ generally results in a highly nontrivial relation between the 
unitary parts of the Iwasawa splittings of $\phi$ and $\hat\phi$.  
However, it is simpler to understand how the monodromy matrices 
$M_\delta =\phi (\tau (z),\lambda)\cdot (\phi (z,\lambda))^{-1}$ and 
$\hat M_\delta =\hat\phi (\tau (z),\lambda)\cdot (\hat\phi (z,\lambda))^{-1}$ 
of $\phi$ and $\hat\phi$ are related by $h_+$.  In fact, 
\begin{equation}\label{eq:mono-dress}
\hat M_\delta = h_+ M_\delta h_+^{-1} \; . 
\end{equation} 
Then, using the closing conditions \eqref{periodR3} or \eqref{periodS3} or 
\eqref{periodH3}, there are situations where one can apply the relation 
\eqref{eq:mono-dress} to solve period problems of the immersion resulting from 
the dressing $\phi\to h_+\cdot\phi$.  
\end{remark}

\section{Delaunay surfaces in $\mathbb{R}^3$, $\mathbb{S}^3$, $\mathbb{H}^3$}
\label{spaceformdelaunay}

We define Delaunay surfaces in all three space forms 
$\mathbb{R}^3$ and $\mathbb{S}^3$ and $\mathbb{H}^3$ to be CMC surfaces of revolution.  
(For $\mathbb{R}^3$ we assume $H\ne 0$, and for $\mathbb{H}^3$ we assume 
$|H|>1$.)  More exactly, Delaunay surfaces are CMC surfaces of revolution about a 
geodesic line in the ambient space form.  That geodesic line must also lie in the 
space form itself.  (In some definitions of Delaunay surfaces, the geodesic line is 
allowed to lie only in some larger outer space of at least four 
dimensions that contains the space form, such as the $4$-dimensional Minkowski 
space for $\mathbb{H}^3$, but we will not consider such cases here).  
Cylinders in each of the space forms are one special limiting 
case of Delaunay surfaces.  

We will show here that Delaunay surfaces in all three space forms can be 
constructed using the potential $\xi$ and solution $\phi$ of $d\phi =\phi\xi$ 
in Section \ref{section5}.  

We define $\xi$ and $\phi \in \Lambda \SL_2\!\mathbb{C}$ and 
$F \in \Lambda \SU_2$ and $B \in \Lambda_{+} \SL_2\!\mathbb{C}$ as in 
Section \ref{section5}, all defined on the punctured $z$-plane 
$\Sigma =\mathbb{C}\setminus\{0\}$.  
For a loop $\delta$ circling counterclockwise once about $z=0$ in $\Sigma$, 
the monodromy matrix $M_\delta$ of both $\phi$ and $F$ is 
\[ M_\delta = \exp(2 \pi i D) = \begin{pmatrix} \cos (2 \pi  \kappa) + i r 
\kappa^{-1} \sin (2 \pi  \kappa) & i \kappa^{-1} \sin 
(2 \pi  \kappa) (s \lambda^{-1} + \bar{t} \lambda) \\
i \kappa^{-1} \sin 
(2 \pi \kappa) (\bar{s} \lambda + t \lambda^{-1}) & 
\cos (2 \pi  \kappa) - i r \kappa^{-1} \sin (2 \pi  \kappa)
\end{pmatrix} \; , \]  where 
\[ \kappa =\sqrt{r^2+|s+\bar t|^2+st(\lambda -\lambda^{-1})^2} \; . \]

{\bf Closing the period:}  
The above $\xi$ and $\phi$ produce CMC surfaces $f$ via the Sym-Bobenko 
formulas \eqref{Sym-Bobenko}, \eqref{SymS3} and \eqref{H3Symformula}.  
Now we consider the period closing conditions for $f$ in each of the three 
space forms $\mathbb{R}^3$, $\mathbb{S}^3$ and $\mathbb{H}^3$:

For $\mathbb{R}^3$, we choose $r$, $s$, $t$ so that 
\begin{equation}\label{delperiodR3}
r^2+|s+\bar{t}|^2 = \frac{1}{4} \;\; \text{and} \;\; st \in \mathbb{R} \; , 
\end{equation}
and we choose $f$ as in \eqref{Sym-Bobenko}.  
Then \eqref{periodR3} is satisfied, so the surface $f$ will close to become 
homeomorphic to a cylinder.  

For $\mathbb{S}^3$, we choose $r$, $s$, $t$ so that 
\begin{equation}\label{delperiodS3}
r^2+|s+\bar{t}|^2 - 4st  \sin^2(\gamma)      
= \frac{1}{4} \;\; \text{and} \;\; st \in \mathbb{R} \; , 
\end{equation}
and we choose $f$ as in \eqref{SymS3}.  
Then \eqref{periodS3} is satisfied, where we choose 
$\gamma_1 = -\gamma$ and $\gamma_2 = \gamma$ in Theorem \ref{S3formulation} 
for some real number $\gamma$, so the surface $f$ closes to become 
homeomorphic to a cylinder.  

For $\mathbb{H}^3$, we choose $r$, $s$, $t$ so that 
\begin{equation}\label{delperiodH3}
r^2+|s+\bar{t}|^2 + 4st \sinh^2\left(\frac{q}{2}\right) 
= \frac{1}{4} \;\; \text{and} \;\; st \in \mathbb{R} \; , 
\end{equation}
and we choose $f$ as in \eqref{H3Symformula}.  
Then \eqref{periodH3} is satisfied, where we choose 
$q\in\mathbb{R}\setminus\{0\}$ and $\psi =0$ in Theorem \ref{H3formulation}, 
so again the surface $f$ will close to become homeomorphic to a cylinder.  

\begin{remark}\label{re:realrst}
Because $s\cdot t \in \mathbb{R}$, we have $s=\hat s\cdot e^{2i\alpha}$ and 
$t=\hat t\cdot e^{-2i\alpha}$ for some $\hat s,\hat t,\alpha\in\mathbb{R}$.  
If we apply the gauge \[\phi \to \phi\cdot g \] for 
\[ g = \begin{pmatrix} e^{i \alpha} & 0 \\ 0 & e^{-i \alpha} 
\end{pmatrix} \; , \] 
then $\xi$ changes to $\hat \xi =g^{-1}\cdot\xi\cdot g$, and the Iwasawa 
splitting $\phi =F\cdot B$ changes to $\hat\phi =\hat F\cdot\hat B$, where 
$\hat F=F\cdot g$ and $\hat B=g^{-1}\cdot B\cdot g$.  
$F$ and $\hat F$ produce the same surface in the Sym-Bobenko formulas 
\eqref{Sym-Bobenko}, \eqref{SymS3}, and \eqref{H3Symformula}, so without loss 
of generality we may use $\hat\xi$ and $\hat\phi$ instead of $\xi$ and $\phi$. 
The advantage of $\hat\xi$ is that 
\[
\hat\xi =\begin{pmatrix}r&\hat s\lambda^{-1}+\hat t\lambda \\
                        \hat s\lambda +\hat t\lambda^{-1}&-r
         \end{pmatrix}\frac{dz}{z} \; , 
\]
where now $r$, $\hat s$, $\hat t$ are all real numbers.  
So without loss of generality, we now assume in the remainder of this section 
that $r,s,t\in\mathbb{R}$.  
\end{remark}

{\bf Showing that surfaces of revolution are produced:} 
Here we show that the surfaces produced by the above 
choice of $\xi$ and $\phi$ are surfaces of revolution.  
By Remark \ref{re:realrst}, it is sufficient to do this for real $s$ and $t$, 
so we assume $s,t \in \mathbb{R}$.  

As we saw in Section \ref{section5}, under the rotation of the domain 
\[ z \to R_{\theta_0}(z)=e^{i \theta_0} z \; , \]  
the following transformations occur: 
\[ F \to M_{\theta_0} F \; , \;\;\;\;\; B \to B \; , \]
where 
\[ M_{\theta_0}=\exp (i \theta_0 D) \in \Lambda \SU_2 \; . \]

When the ambient space is $\mathbb{R}^3$, under the mapping $z \to R_{\theta_0}(z)$, 
$f$ as in \eqref{Sym-Bobenko} changes, like in \eqref{eq:monR3}, by 
\begin{equation}\label{delaunayR3} 
f \to (M_{\theta_0}|_{\lambda=1}) 
f (M_{\theta_0}^{-1}|_{\lambda=1}) -\frac{i}{2H} 
(\partial_\lambda M_{\theta_0}|_{\lambda=1}) 
(M_{\theta_0}^{-1}|_{\lambda=1}) \; . \end{equation}  
One can check, using Section \ref{section4}, that Equation \eqref{delaunayR3} 
represents a rotation of angle $\theta_0$ about the line 
\begin{equation}\label{DelAxisR3} 
\left\{\left.\!x \cdot (-s-t,0,r)+\left(\frac{s-t}{r},0,0\right) \, \right| \, 
x \in \mathbb{R}\right\} \; , 
\end{equation}
when $r\ne 0$.  (We also computed this in Section \ref{section5}.)  Therefore 
$f$ is a surface of revolution (since the line \eqref{DelAxisR3} does not 
depend on $\theta_0$), and hence a Delaunay surface in $\mathbb{R}^3$.  

When the ambient space is $\mathbb{S}^3$, 
under the mapping $z \to R_{\theta_0}(z)$, 
$f$ as in \eqref{SymS3} changes, like in \eqref{eq:monS3}, by 
\begin{equation}\label{delaunayS3} 
f \to (M_{\theta_0}|_{\lambda=e^{-i \gamma}}) 
f (M_{\theta_0}^{-1}|_{\lambda=e^{i \gamma}}) \; . 
\end{equation}  
One can check, using Subsection \ref{subsectionS3}, that Equation \eqref{delaunayS3} 
represents a rotation 
of angle $\theta_0$ about the geodesic line 
\begin{equation}\label{DelAxisS3} \{ (x_1,x_2,0,x_4) \in \mathbb{S}^3 \; | \; 
\sin(\gamma) (s - t) x_1 + r x_2 - \cos(\gamma) (s + t) x_4 = 0 \} 
\; . \end{equation}  So we have a surface of revolution in this case also 
(since the geodesic line \eqref{DelAxisS3} does not depend on $\theta_0$), 
and hence a Delaunay surface in $\mathbb{S}^3$.  

When the ambient space is $\mathbb{H}^3$, 
under the mapping $z \to R_{\theta_0}(z)$, 
$f$ as in \eqref{H3Symformula} changes, like in \eqref{eq:monH3}, by 
\begin{equation}\label{delaunayH3} 
f \to (M_{\theta_0}|_{\lambda=e^{q/2}}) 
f (\overline{M_{\theta_0}}^t|_{\lambda=e^{q/2}}) \; . 
\end{equation}  
One can check, using Subsection \ref{subsectionH3}, that Equation \eqref{delaunayH3} 
represents a rotation 
of angle $\theta_0$ about the geodesic line 
\begin{equation}\label{DelAxisH3} 
\{ (x_1,0,x_3,x_0) \in \mathbb{H}^3 \; | \; 
\sinh (q/2) (s-t) x_0 - r x_1 + \cosh (q/2) (s+t) x_3 = 0 \} \; . \end{equation}  
Therefore $f$ is a surface of revolution 
(since the geodesic line \eqref{DelAxisH3} does not depend on $\theta_0$), 
and hence a Delaunay surface in $\mathbb{H}^3$.  

\subsection{The weight of Delaunay surfaces} 
\label{subsec:weightDel}
We now describe the weight of Delaunay surfaces in all three space 
forms $\mathbb{R}^3$ and $\mathbb{S}^3$ and $\mathbb{H}^3$.  

Let $\delta$ be an oriented loop about an annular end of a CMC $H$ surface in 
$\mathbb{R}^3$ or $\mathbb{S}^3$ or $\mathbb{H}^3$, and let $\mathcal{Q}$ be an immersed 
disk with 
boundary $\delta$.  Let $\eta$ be the unit conormal of the surface along 
$\delta$ and let $\nu$ be the unit normal of $\mathcal{Q}$, the signs 
of both of them determined by the orientation of $\delta$.  
Then the weight of the end 
with respect to a Killing vector field $Y$ (in $\mathbb{R}^3$ or $\mathbb{S}^3$ or 
$\mathbb{H}^3$) is 
\[ w(Y) = \int_\delta \langle \eta , Y \rangle - 2 H \int_\mathcal{Q} \langle \nu , 
Y \rangle \; . \]  
In the case that the end is asymptotic to a Delaunay surface with axis 
$\ell$ and $Y$ is the Killing vector field associated to unit translation 
along the direction of $\ell$, we abbreviate $w(Y)$ to $w$ and say that 
$w$ is the {\em weight}, or the {\em mass} or {\em flux}, of the end.  This weight 
changes sign when the 
orientation of $\delta$ is switched, but otherwise it is independent of the 
choices of $\delta$ and $\mathcal{Q}$.  In other words, as shown in 
\cite{KorKS} and \cite{KorKMS}, it is a homology invariant.  (Since the mean 
curvature in \cite{KorKS} and \cite{KorKMS} is 
defined as the sum of the principal curvatures, rather than the average, 
we must replace $H$ by $2 H$ in the formulas for the weights there.)  

We have $F |_{z \in \mathbb{S}^1} \in \SU_2$ for all nonzero $\lambda$, 
hence \[ F |_{z \in \mathbb{S}^1} = \exp (i \theta D) \; , \] 
where $z=e^{i \theta}$.  

For the case of $\mathbb{R}^3$, \eqref{Sym-Bobenko} and \eqref{eq:normalR3} imply 
\[ f|_{z=1}=-N|_{z=1}=\frac{-i}{2} \sigma_3 \cong (0,0,-1) \; . \]  
From this and \eqref{DelAxisR3}, it follows that $\mathbb{S}^1$ in the 
$z$-plane is mapped to a circle of radius $2|s|/|H|$.  
Furthermore, if $\beta$ is the angle between the conormal $\eta$ of the 
surface along the circle and the normal $\nu$ to the plane containing the 
circle, then 
\[ \cos \beta = 2 (s+t) \]  
and the weight of the Delaunay surface in $\mathbb{R}^3$ is 
\[ w = \frac{8 \pi s t}{|H|} \; . \]  

An unduloid is produced when $w>0$, and a nodoid is produced when $w<0$, 
and for the limiting singular case of a chain of spheres, $w=0$.  
A cylinder is produced when $w$ attains the maximal value $w=\pi /(2|H|)$, 
that is, when $r=0$ and $s=t=1/4$.  

For the case of $\mathbb{S}^3$, \eqref{SymS3} and \eqref{eq:normalS3} imply 
\[ f|_{z=1}=\begin{pmatrix} e^{i \gamma} & 0 \\ 0 & e^{-i \gamma} 
\end{pmatrix} \; , \;\;\; N|_{z=1}=
\begin{pmatrix} i e^{i \gamma} & 0 \\ 0 & -i e^{-i \gamma} 
\end{pmatrix} \; . \]  
The circle $\mathbb{S}^1$ in the $z$-plane is mapped to a circle contained in 
the surface in $\mathbb{S}^3$.  Let $\ell$ denote the radius of this image 
circle in $\mathbb{S}^3$, and let $\beta$ denote the angle between the 
conormal $\eta$ of the surface along this circle and the normal $\nu$ to the 
geodesic plane in $\mathbb{S}^3$ containing this circle.  
Then, using \eqref{DelAxisS3}, 
\[ \sin \ell = 2 t \sin(2 \gamma) \; , \;\;\; \cos \beta = 
2 \frac{s+t \cos(2 \gamma)}{\cos(\ell)} \; . \]  
It follows that the weight of the Delaunay surface in $\mathbb{S}^3$ is 
\[ w = \frac{8 \pi s t}{\sqrt{H^2 +1}} \; . \]  

For the case of $\mathbb{H}^3$, \eqref{H3Symformula} implies 
\[ f|_{z=1}=\begin{pmatrix} e^{-q/2} & 0 \\ 0 & e^{q/2} 
\end{pmatrix} \; , \;\;\; N|_{z=1}=
\begin{pmatrix} e^{-q/2} & 0 \\ 0 & -e^{q/2} 
\end{pmatrix} \; . \]  
The circle $\mathbb{S}^1$ in the $z$-plane is mapped to a circle contained in 
the surface in $\mathbb{H}^3$.  Let $\ell$ denote the radius of this image 
circle in $\mathbb{H}^3$, and again let $\beta$ denote the angle between the 
conormal $\eta$ of the surface along this circle and the normal $\nu$ to the 
geodesic plane in $\mathbb{H}^3$ containing this circle.  
Then, using \eqref{DelAxisH3}, 
\[ \sinh \ell = 2 s \sinh (q) \; , \;\;\; \cos \beta = 
2 \frac{t+s \cosh (q)}{\cosh (\ell)} \; . \]  
Thus the weight of the Delaunay surface in $\mathbb{H}^3$ is 
\[ w = \frac{8 \pi s t}{\sqrt{H^2 -1}} \; . \]  

\begin{figure}[phbt]
\begin{center}
\includegraphics[width=0.4\linewidth]{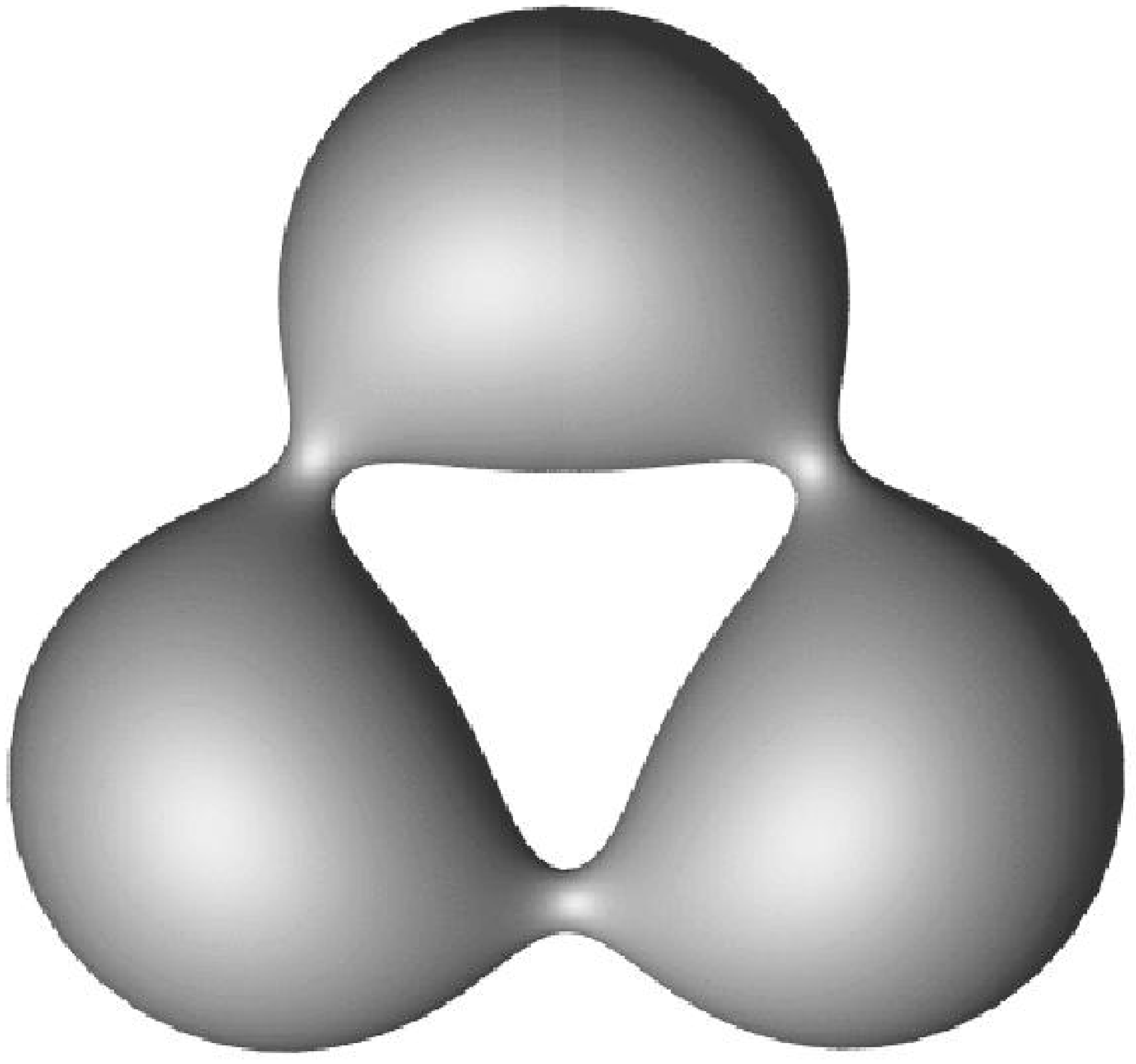}
\includegraphics[width=0.4\linewidth]{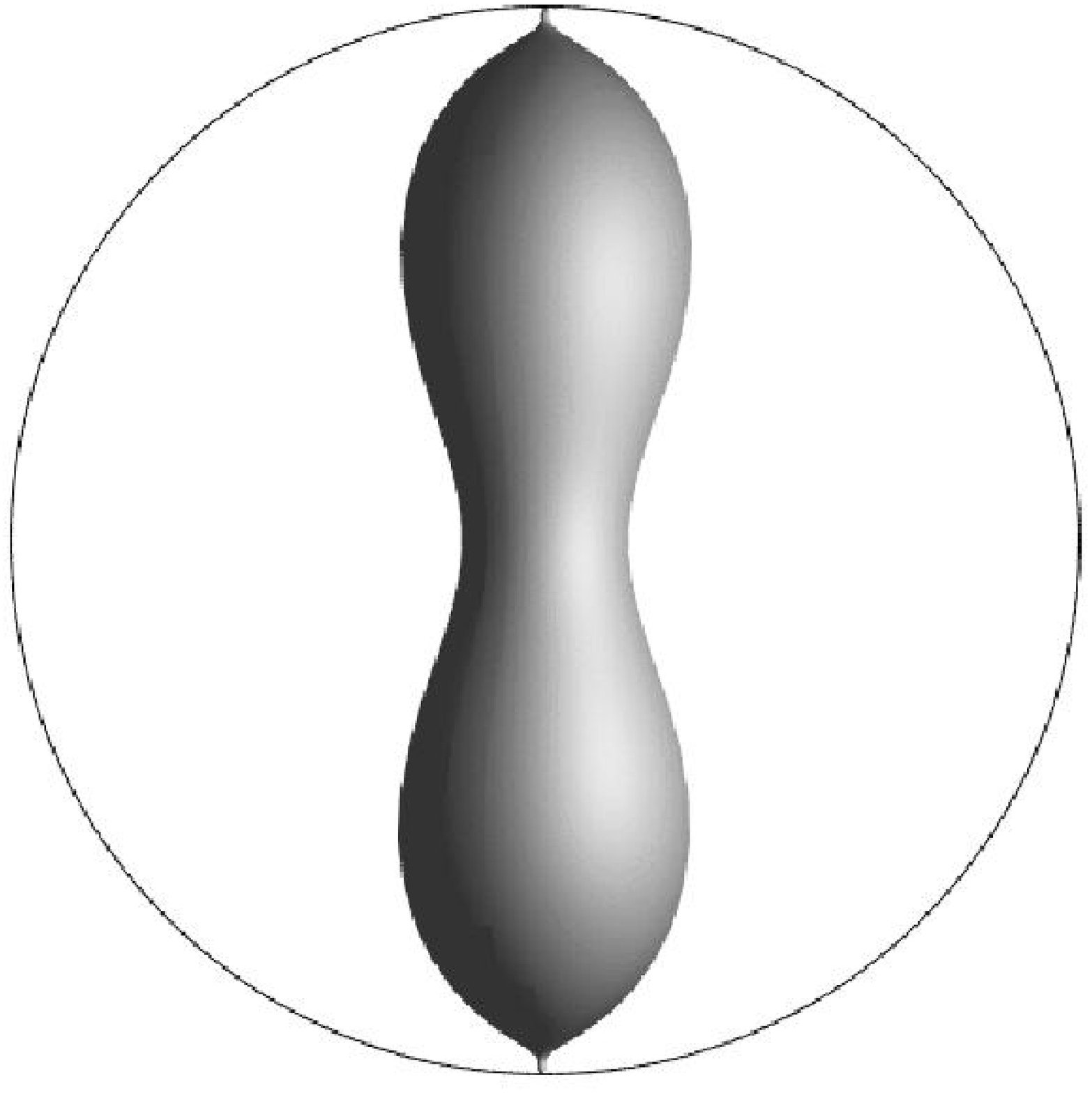}
\end{center}
\caption{Delaunay surfaces in $\mathbb{S}^3$ and $\mathbb{H}^3$.  
(These graphics were made using N. Schmitt's software CMCLab \cite{Sch:cmclab}.)}
\label{fg:DelaunayS3H3}
\end{figure}

\section{The representation of Bryant}
\label{bryantrep}

In this section we show a representation due to Bryant for CMC $1$ surfaces in 
$\mathbb{H}^3$ that is similar to the Weierstrass representation for minimal 
surfaces in $\mathbb{R}^3$.  It is natural to expect the existence of this 
analogous representation for CMC $1$ surfaces in $\mathbb{H}^3$, since CMC $1$ 
surfaces in $\mathbb{H}^3$ are related to minimal surfaces by a correspondence 
that is often called the Lawson correspondence.  
This correspondence was undoubtedly known before Lawson's time, but Lawson 
described it explicitly in the literature \cite{L}.  
In the first two subsections here, we describe the Lawson correspondence.  
The third subsection gives the Bryant representation, and could be read 
independently of the first two subsections.  

\subsection{The Lawson correspondence in a general setting}
\label{LawsonSection}

Regarding the Lawson correspondence between CMC surfaces in $\mathbb{R}^3$ 
and $\mathbb{S}^3$ and $\mathbb{H}^3$, the essential ingredient is the 
fundamental theorem of surface theory, telling us that on a simply-connected 
Riemann surface $\Sigma$ there exists an immersion with first and second 
fundamental forms $g$ and $b$ if and only if 
$g$ and $b$ satisfy the Gauss and Codazzi equations.  (This is true 
regardless of whether the ambient space is $\mathbb{R}^3$ or $\mathbb{S}^3$ or 
$\mathbb{H}^3$, but the Gauss and Codazzi equations are not the same in the 
three cases.)  So if $f : \Sigma \to \mathbb{R}^3$ (resp. $\mathbb{S}^3$) is 
a conformal CMC immersion, then its fundamental forms $g$ and $b$ satisfy 
the Gauss and Codazzi equations for surfaces in $\mathbb{R}^3$ 
(resp. $\mathbb{S}^3$), and one can then check 
that $\tilde{g}=g$ and $\tilde{b}=c \cdot g+b$ satisfy the Gauss and 
Codazzi equations for surfaces in $\mathbb{H}^3$ (resp. $\mathbb{R}^3$) 
for correctly chosen $c$ depending on the mean curvature $H$ of $f$, so 
there exists an immersion $\tilde{f} : \Sigma \rightarrow 
\mathbb{H}^3$ (resp. $\mathbb{R}^3$) with fundamental forms $\tilde{g}$ and 
$\tilde{b}$.  Since $\mbox{trace}(g^{-1}b)=2H$, we have 
$\mbox{trace}(\tilde{g}^{-1} \tilde{b})=2(H+c)$, so 
$\tilde{f}$ is CMC ($H+c$).  And since $g = \tilde{g}$, $\tilde{f}$ is 
also conformal, and $f(\Sigma)$ and $\tilde{f}(\Sigma)$ are 
isometric.  This is Lawson's correspondence.  

We now explain the Lawson correspondence in more detail: Let $M^3(\bar{K})$ be 
the unique simply-connected complete $3$-dimensional space form with constant 
sectional curvature $\bar{K}$ 
(e.g. $M^3(0)=\mathbb{R}^3$, $M^3(-1)=\mathbb{H}^3$, $M^3(1)=\mathbb{S}^3$).  
For an immersion $f:\Sigma \to M^3(\bar{K})$ with induced metric 
$\langle \cdot , \cdot \rangle$ and Levi-Civita connection $\nabla$ and 
Gaussian curvature $K$ and shape operator $\mathcal{A}$, the Gauss and Codazzi 
equations are satisfied:  
\begin{equation}\label{eq:GCequations}
K - \bar{K} = \det(\mathcal{A}) \; , \; \; \; 
\langle \mathcal{A}([X,Y]),Z \rangle = 
   \langle \nabla_X \mathcal{A}(Y),Z \rangle
 - \langle \nabla_Y \mathcal{A}(X),Z \rangle 
\end{equation}  
for all smooth vector fields $X$, $Y$, and $Z$ in the tangent space of 
$\Sigma$.  
Assume $f$ is CMC $H$, so $H = (1/2)\mbox{trace}(\mathcal{A})$ is constant.  
Now choose any $c \in \mathbb{R}$ and define 
\[ \tilde{\mathcal{A}} = \mathcal{A}+c \cdot \id \; , \; \; \tilde{K} = 
\bar{K} - c\cdot \mbox{trace}(\mathcal{A}) - c^2 \; . \] 
Then Equations \eqref{eq:GCequations} imply that the Codazzi equation still 
holds when $\mathcal{A}$ is replaced by 
$\tilde{\mathcal{A}}$:  
\begin{eqnarray*}
\langle \tilde{\mathcal{A}}([X,Y]),Z \rangle 
&=& \langle \mathcal{A}([X,Y]),Z \rangle + c \langle [X,Y],Z \rangle \\
&=& \langle \nabla_X \mathcal{A}(Y),Z \rangle
   -\langle \nabla_Y \mathcal{A}(X),Z \rangle + c \langle \nabla_X Y,Z \rangle 
   -c \langle \nabla_Y X,Z \rangle \\
&=& \langle \nabla_X \tilde{\mathcal{A}}(Y),Z \rangle 
   -\langle \nabla_Y \tilde{\mathcal{A}}(X),Z \rangle \; . 
\end{eqnarray*} 
Equations \eqref{eq:GCequations} also imply that the Gauss equation in $M^3(\tilde{K})$ 
also holds with $\tilde{\mathcal{A}}$ 
(note that $K$ is intrinsic and does not change): 
\begin{eqnarray*}
K-\tilde{K}
 &=& K - (\bar{K} - c\cdot\mbox{trace}(\mathcal{A}) - c^2) \\
 &=& \det(\mathcal{A}) + c\cdot\mbox{trace}(\mathcal{A}) + c^2 
  =  \det(\mathcal{A}+c\cdot\id) = \det(\tilde{\mathcal{A}}) \; .
\end{eqnarray*}
Therefore there exists 
an immersion $\tilde{f}: \Sigma \to M^3(\tilde{K})$ with metric 
$\langle \cdot , \cdot \rangle$ and shape operator $\tilde{\mathcal{A}}$, and 
$\tilde{f}(\Sigma)$ is isometric to $f(\Sigma)$.  
As the mean curvature of $\tilde{f}(\Sigma)$ is 
\[ \tilde{H} = \frac{1}{2}\mbox{trace}(\tilde{\mathcal{A}}) = 
\frac{1}{2}\mbox{trace}(\mathcal{A}) + c = H + c \; , \] this 
demonstrates the Lawson correspondence between a CMC $H$ surface in 
$M^3(\bar{K})$ and a CMC $(H+c)$ surface in $M^3(\bar{K} - 2 c H - c^2)$.  

In particular, when $H=\bar{K}=0$ and $c=1$, we have the 
correspondence between minimal surfaces in $\mathbb{R}^3$ and CMC $1$ surfaces in 
$\mathbb{H}^3$.  And when $H=0$ and $\bar{K}=1$ and $c=1$, we have the 
correspondence between minimal surfaces in $\mathbb{S}^3$ and CMC $1$ surfaces 
in $\mathbb{R}^3$.  

\subsection{The Lawson correspondence in our setting.}
\label{LawsonSection2}

Now we consider what the Lawson correspondence means in our situation.  
It means that the metric function $u$ in the metric $4e^{2u}dzd\bar{z}$ and 
the Hopf differential $Q=\langle f_{zz},N\rangle$ are unchanged, 
and that the constant mean curvature $H$ is changed to $H+c$.  
(In Sections \ref{section11} and \ref{section12}, $w$ and $A$ were used to 
denote the conformal parameter and the Hopf differential, respectively.  
However, here we revert to the notations $z$ and $Q$, i.e. we rename 
$w$ and $A$ to $z$ and $Q$, respectively.)  
When changing from $\mathbb{S}^3$ (resp. $\mathbb{R}^3$) to $\mathbb{R}^3$ 
(resp. $\mathbb{H}^3$), the real number $c$ is $c=-H \pm \sqrt{H^2+1}$.  So if 
$H_r$ and $H_s$ and $H_h$ represent the 
mean curvatures in $\mathbb{R}^3$ and $\mathbb{S}^3$ 
and $\mathbb{H}^3$, respectively, then we have 
\begin{equation}\label{Lawsonrelation} 
H_s^2+1=H_r^2 \; , \;\;\; H_r^2+1=H_h^2 \; . \end{equation}  
Looking back at Equations \eqref{firstMC}, \eqref{originalS3MC} and 
\eqref{originalH3MC}, we have that the Gauss and Codazzi equations in the 
three cases are:  
\[ \begin{array}{c@{\;:\;\;\;}c@{\;,\;\;\;}c}
\mathbb{R}^3 \text{ case} & 
4 u_{z\bar z}+4 e^{2 u} (H_r^2) - Q \bar Q e^{-2 u} = 0 & 
Q_{\bar z}= 0 \; . \\
\mathbb{S}^3 \text{ case} & 
4 u_{z\bar z}+4 e^{2 u} (H_s^2+1) - Q \bar Q e^{-2 u} = 0 & 
Q_{\bar z}= 0 \; . \\
\mathbb{H}^3 \text{ case} & 
4 u_{z\bar z}+4 e^{2 u} (H_h^2-1) - Q \bar Q e^{-2 u} = 0 & 
Q_{\bar z}= 0 \; . \end{array} \]
The last two equations have had their 
notations changed so that $w$ is now $z$ and $A$ is now $Q$.  

Thus these three sets of Gauss-Codazzi equations are the same in 
$\mathbb{R}^3$ and $\mathbb{S}^3$ and $\mathbb{H}^3$, by 
Equation \eqref{Lawsonrelation}. This demonstrates the Lawson correspondence. 

\subsection{The Bryant representation}
\label{sub:B-rep}

Recall Lorentz $4$-space $\mathbb{R}^{3,1}$ with Lorentz metric 
$\langle\cdot ,\cdot\rangle_{\mathbb{R}^{3,1}}$ as in Section \ref{section12}. 
Recall also that hyperbolic $3$-space $\mathbb{H}^3$ can be written as in 
Lemma \ref{lm:SL2Cform}.  
We have the following fact from Subsection \ref{subsectionH3}: 
\begin{lemma}
For any $F \in \SL_2\!\mathbb{C}$,
\[
\langle X,Y\rangle = \langle F X \bar{F}^t, F Y \bar{F}^t 
\rangle \;\;\; \text{for all}\;\; X,Y \in \mathbb{R}^{3,1} \; . 
\]  
Hence the map $X \to F X \bar{F}^t$ is a rigid motion of $\mathbb{R}^{3,1}$, 
and therefore also of $\mathbb{H}^{3}$.  
\end{lemma}  

Let \[ f:\Sigma \to \mathbb{H}^3 \] be a conformal immersion, where $\Sigma$ is a 
simply-connected Riemann surface with complex coordinate $w$.  
(By Theorem \ref{conformalityispossible}, without loss of generality we may assume 
$f$ is conformal.)  We also have: 

\begin{lemma}\label{rigidmotionlemma}
There exists an $F \in \SL_2\!\mathbb{C}$ (unique up to sign $\pm F$) so that 
\[ f = F \bar F^t \; , \;\; 
 e_1 := \frac{f_x}{|f_x|} = F \sigma_1 \bar{F}^t \; , \;\; 
 e_2 := \frac{f_y}{|f_y|} = F \sigma_2 \bar{F}^t \; , \;\; 
 N := F \sigma_3 \bar{F}^t \; , \] where $N$ is the unique unit vector in 
the tangent space of $\mathbb{H}^3$ at the point $F \bar F^t$ such that $N$ is 
perpendicular to both $e_1$ and $e_2$ and $\{ e_1,e_2,N\}$ is positively 
oriented.  
\end{lemma}

\begin{proof}
There exists an $F\in\SL_2\!\mathbb{C}$ so that $f=F\bar F^t$ and this $F$ is 
unique up to the choice $F \cdot B$ for any $B \in \SU_2$.  Choose $B$ so 
that the above equations for $e_1$, $e_2$, and $N$ hold.  This $B$ can be 
found because $\sigma_1$, $\sigma_2$ and $\sigma_3$ are orthonormal in the 
usual Euclidean $\mathbb{R}^3$ sense (see Sections \ref{section0} and 
\ref{2by2laxpair}).  Furthermore, this $B$ is unique up to sign.  
\end{proof}

Define the metric factor $u$ and Hopf differential $A$ by 
\[ \langle f_w, f_{\bar w}\rangle = 2 e^{2 u} \; , \;\;\; 
   A = \langle f_{w w},N\rangle \; . \]  
As we saw in Section \ref{section12}, it follows that 
\begin{eqnarray*} 
& f_{ww} = 2 u_w f_w + A N \; , \;\;\; A_{\bar w}=2 H_w e^{2u} \; , & \\
& f_{w \bar w} = 2 e^{2 u} f + 2 H e^{2 u} N \; , \;\;\; 
N_w = - H f_w - \frac{1}{2} A e^{-2 u} f_{\bar w} \; , & 
\end{eqnarray*}
and the Gauss equation is 
\begin{equation}\label{eq:GaussforCMC1inH3}
 2 u_{w \bar w} + 2e^{2u}(H^2-1) - \frac{1}{2} A \bar A e^{-2 u} = 0 \; . 
\end{equation}

Choosing $F$ as in Lemma \ref{rigidmotionlemma}, where $w=x+iy$ for 
$x,y\in\mathbb{R}$, we have 
$e_1 = F \sigma_1 \bar F^t$ and $e_2 = F \sigma_2 \bar F^t$, and so 
\[ f_w = 2 e^u F \begin{pmatrix} 0 & 0 \\ 1 & 0 \end{pmatrix} \bar F^t 
\; , \;\;\; f_{\bar w} = 2 e^u F \begin{pmatrix} 0 & 1 \\ 0 & 0 
\end{pmatrix} \bar F^t \; . \]  
We define $U=(U_{ij})_{i,j=1,2}$ and $V=(V_{ij})_{i,j=1,2}$ by 
\[ F_w = F \cdot U \; , \;\;\; F_{\bar w} = F \cdot V \; . \]  
Because $f_{w \bar w}=f_{\bar w w}= 2 e^{2 u} f + 2 H e^{2 u} N$, we 
have $U_{21}=(1-H)e^u$ and $V_{12}=(1+H) e^u$ and $u_{\bar w} + V_{22} + \bar U_{11} 
= 0$.  Because $\det F = 1$, we have that $U$ and $V$ are trace free, and so 
$U_{11}=-U_{22}$ and $V_{11}=-V_{22}$.  Because 
$f_{ww} = 2 u_w f_w + A N$, we have $U_{12}=(1/2) e^{-u} A$ and 
$V_{21} = -(1/2) e^{-u} \bar A$.  Because $N_w = - Hf_w - (1/2) A e^{-2 u} 
f_{\bar w}$, we have $U_{11}=-\bar V_{11}$ and $U_{22}=-\bar V_{22}$.  
Therefore 
\[ U = \frac{1}{2} \begin{pmatrix}
-u_w & e^{-u} A \\ 2(1-H)e^u & u_w 
\end{pmatrix} \; , \;\;\; V = \frac{1}{2} \begin{pmatrix}
u_{\bar w} & 2(1+H)e^u \\ -e^{-u} \bar A & -u_{\bar w} 
\end{pmatrix} \; . \]  
Because $U \neq 0$, $F$ is not antiholomorphic, but we now show that if $f$ 
has constant mean curvature $H=1$, then we can 
change $F$ to $F B$ for some $B=B(w,\bar w) \in \SU_2$ so that $(F B)_w=0$.  
(Note that multiplying by $B$ will not change the immersion $f$, as 
$f=F\overline{F}^t=(FB)\overline{(FB)}^t$, even though the last three 
properties in Lemma \ref{rigidmotionlemma} will no longer hold.)  
To accomplish this, we first assume 
\[ H=1 \]
and then we need $B$ to satisfy 
\begin{equation}\label{BLax1} B_w = - U B \; . \end{equation}  
Consider also the equation 
\begin{equation}\label{BLax2} B_{\bar w} = W B \; , \end{equation}  
which we take to be the definition of $W$.  
Then the Lax pair \eqref{BLax1}-\eqref{BLax2} is equivalent to 
\begin{equation}\label{LaxPairforB} 
B_x = (W-U) B \; , \;\;\; B_y =-i (W+U) B \; . \end{equation} 
A simple computation shows that 
\[ 
W- U \, , \; -i (W+U) \in \su_2 \;\;\;\text{if and only if}\;\;\;W=\bar U^t , 
\] 
so we set $W=\bar U^t$.  Then we have that 
\[ U_{\bar w} + \bar U_w^t + [U,\bar U^t] = 2(H-1)e^{2u}\sigma_3 \; . \]
So if $H\equiv 1$, the compatibility condition 
$U_{\bar w} + \bar U_w^t + [U,\bar U^t] = 0$ for the Lax pair 
\eqref{LaxPairforB} (or equivalently, for the Lax pair 
\eqref{BLax1}-\eqref{BLax2}) holds.  Thus an analog of 
Proposition \ref{prop:MCeqn}, with $\slg_n\!\mathbb{C}$ and $\SL_n\!\mathbb{C}$ 
replaced by $\su_2$ and $\SU_2$, implies that there exists a solution 
$B \in \SU_2$ of the Lax pair 
\eqref{BLax1}-\eqref{BLax2}.  In particular, \eqref{BLax1} has a 
solution $B \in \SU_2$.  

Now, since $f$ has constant mean curvature $H=1$, writing that solution $B$ as 
\[ B = \begin{pmatrix} p & q \\ -\bar q & \bar p \end{pmatrix} \; , \;\;\;
p\bar p+q\bar q=1 \; ,
\]  
we have \[ (F B)^{-1} (F B)_{\bar w} = B^{-1} (V+\bar U^t) B = 
\begin{pmatrix} -2 e^u \bar p \bar q & 2 e^u \bar p^2 \\ 
-2 e^u \bar q^2 & 2 e^u \bar p \bar q \end{pmatrix} \; , \] which must be 
antiholomorphic, because $F B$ is.  Define 
\[ \omega = -2 e^u q^2 \; , \;\;\; g = \frac{p}{q} \; , \;\;\; 
\hat F = \bar F \cdot \bar B \; . \]  
Then $\hat F$ is holomorphic in $w$ and 
\begin{equation}\label{eq:Bryant-rep}
   \hat F^{-1} \hat F_{w}
 = \begin{pmatrix} g & -g^2 \\ 1 & -g \end{pmatrix} \omega \; , 
\end{equation}
as in the 
Bryant representation in \cite{uy1} (modified there from \cite{Br}, and 
proven in \cite{Br}).  
\begin{quote}
{\bf The Bryant representation}: Any conformal CMC $1$ immersion from a 
simply-connected domain $\Sigma$ into $\mathbb{H}^3$ can be written as 
$\hat F\overline{\hat F}^t$ for some solution $\hat F\in\SL_2\!\mathbb{C}$ of 
an equation of the form \eqref{eq:Bryant-rep}.  
\end{quote}
Also, 
\[ (1+|g|^2)^2 |\omega|^2 = 4 e^{2 u} \; , \] so we also have the metric 
determined from the $g$ and $\omega$ used in the Bryant representation.  
Note that changing $F$ to $\hat F$ changes $f$ to $\hat f=\bar f$, 
which is only a reflection of $\mathbb{H}^3$, so $F$ and $\hat F$ produce 
the same surface up to a reflection.  

Conversely, starting with any $\hat F$ solving \eqref{eq:Bryant-rep}, 
one can check that a CMC $1$ surface $\hat F\overline{\hat F}^t$ is obtained.  

\begin{remark}
We note that choosing $H=1$ was essential to proving the Bryant representation 
above.  
\end{remark}

\begin{remark}\label{re:KforCMCinH3}
Equations \eqref{eq:egregium} and \eqref{eq:GaussforCMC1inH3} imply that the Gaussian 
curvature for a CMC 1 surface in $\mathbb{H}^3$ is 
always non-positive, just as was the case for minimal surfaces in 
$\mathbb{R}^3$ (Remark \ref{re:KformininR3}).  
\end{remark}

\subsection{The Gauss map}
\label{sub:gaussmap}

In the Bryant representation, the map $g:\Sigma\to\mathbb{C}\cup\{\infty\}$ 
does not represent a Gauss map (as it did in Section \ref{Wrepsection}).  
The actual Gauss map can be represented as follows: The normal vector 
$\hat N=\bar N=\bar F\sigma_3 F^t$ is both perpendicular to the surface 
$\hat f=\bar FF^t=\hat F\overline{\hat F}^t$ in $\mathbb{H}^3$ and tangent to 
the space $\mathbb{H}^3$ in $\mathbb{R}^{3,1}$ at each point 
$\hat f\in\mathbb{H}^3$.  
Let $P$ be the unique $2$-dimensional plane in $\mathbb{R}^{3,1}$ containing 
the three points $(0,0,0,0)$ and $\hat f$ and $\hat f+\hat N$.  
Then $P$ contains the geodesic $\delta$ in $\mathbb{H}^3$ that starts at 
$\hat f$ and extends in the direction of $\hat N$ (this follows from the fact 
that all isometries of $\mathbb{H}^3$ are of the form in Subsection 
\ref{subsectionH3}, and there is a rotation of $\mathbb{R}^{3,1}$ that moves 
$P\cap\mathbb{H}^3$ to the geodesic 
$\{(\sinh t,0,0,\cosh t)\in\mathbb{R}^{3,1}\;|\;t\in\mathbb{R}\}$ 
that lies in a plane).  
This geodesic $\delta$ has a limiting direction in the upper half-cone 
$\mathcal{L}^+=\{(x_1,x_2,x_3,x_0)\in\mathbb{R}^{3,1}\;|\;
                 x_0^2=x_1^2+x_2^2+x_3^2,\;x_0\ge 0\}$ 
of the light cone of $\mathbb{R}^{3,1}$.  
This limiting direction, in the direction of $\hat f+\hat N$, 
is the Gauss map of $\hat f$.  We identify the set of limiting directions with 
$\mathbb{S}^2$ by associating the direction 
$[(x_1,x_2,x_3,x_0)]:=\{\alpha (x_1,x_2,x_3,x_0)\;|\;\alpha\in\mathbb{R}^+\}$ 
for $(x_1,x_2,x_3,x_0)\in\mathcal{L}^+$ with the point 
$(x_1/x_0,-x_2/x_0,x_3/x_0)\in\mathbb{S}^2$; that is, we identify
\[
\begin{pmatrix}a&b\\\bar b&c\end{pmatrix}\in\mathcal{L}^+\;\;\;\text{with}
\;\;\;\left(\frac{b+\bar b}{a+c},i\frac{\bar b-b}{a+c},\frac{a-c}{a+c}\right)
\in\mathbb{S}^2 \; .
\]
Moreover, composing with stereographic projection
\[
\mathbb{S}^2\ni(\xi_1,\xi_2,\xi_3)\mapsto\left\{
\begin{array}{@{\,}ll}
(\xi_1+i\xi_2)/(1-\xi_3)\in\mathbb{C} & \text{if}\;\;\;\xi_3\ne 1 \; , \\
\infty                        & \text{if}\;\;\;\xi_3  = 1 \; ,
\end{array}\right.
\]
we can then identify 
\[
\begin{pmatrix}a&b\\\bar b&c\end{pmatrix}\in\mathcal{L}^+\;\;\;\text{with}
\;\;\;\frac{b}{c}\in\mathbb{C}\cup\{\infty\} \; .
\]

\begin{remark}
When identifying the limiting direction $[(x_1,x_2,x_3,x_0)]$ with the point 
\[ (x_1/x_0,-x_2/x_0,x_3/x_0)\in\mathbb{S}^2 \; , \] the minus sign in the second 
coordinate of $(x_1/x_0,-x_2/x_0,x_3/x_0)$ at first appears rather arbitrary.  
But it is an acceptable identification, and we need it this way to make the 
Gauss map holomorphic (rather than antiholomorphic), as we will see below.  
Our identification here is in fact the same as the identification used in 
\cite{Br} and \cite{uy1}.  The reason that the minus sign appears here 
(but not in \cite{Br} and \cite{uy1}) is that our choice of matrix 
representation \eqref{eq:L4identifyHerm2} for $\mathbb{R}^{3,1}$ is slightly 
different than that of \cite{Br} and \cite{uy1}.  
\end{remark}

Since 
\begin{eqnarray*}
 \left[\hat f+\hat N\right]
&=& \left[\bar F\begin{pmatrix}1&0\\0&0\end{pmatrix}F^t\right]
 =  \left[\hat FB^t\begin{pmatrix}1&0\\0&0\end{pmatrix}
          \bar B\overline{\hat F}^t\right] \\
&=& \left[\begin{pmatrix}(F_{11}p+F_{12}q)\overline{(F_{11}p+F_{12}q)}&
                         (F_{11}p+F_{12}q)\overline{(F_{21}p+F_{22}q)}\\
                         (F_{21}p+F_{22}q)\overline{(F_{11}p+F_{12}q)}&
                         (F_{21}p+F_{22}q)\overline{(F_{21}p+F_{22}q)}
          \end{pmatrix}\right] \; ,
\end{eqnarray*}
where 
\[
\hat F=\begin{pmatrix}F_{11}&F_{12}\\F_{21}&F_{22}\end{pmatrix} \; , 
       \;\;\;\text{and}\;\;\;
     B=\begin{pmatrix}p&q\\-\bar q&\bar p\end{pmatrix} \; ,
\]
the Gauss map is represented by $G:\Sigma\to\mathbb{S}^2$ with 
\[
G=\frac{F_{11}p+F_{12}q}{F_{21}p+F_{22}q} 
 =\frac{F_{11}g+F_{12}}{F_{21}g+F_{22}}
 =\frac{dF_{11}}{dF_{21}}
 =\frac{dF_{12}}{dF_{22}} \; , 
\]
by \eqref{eq:Bryant-rep}.  So we find that the Gauss map of $\hat f$ is 
represented by $G:\Sigma\to\mathbb{S}^2$ with 
\[
G=\frac{dF_{11}}{dF_{21}}=\frac{dF_{12}}{dF_{22}} \; .  
\]

\begin{remark}
If we change $\hat{F}$ to $\hat{F}^{-1}$ in the Bryant representation we get 
the CMC 1 ``dual'' surface $f^\sharp=(\hat{F}^{-1})(\overline{\hat{F}^{-1}}^t)$ 
in $\mathbb{H}^3$.  Note that changing $f$ to $f^\sharp$ switches $g$ and $G$. 
\end{remark}

\begin{figure}[phbt]
\begin{center}
\includegraphics[width=0.25\linewidth]{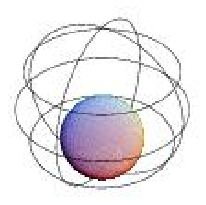}
\includegraphics[width=0.35\linewidth]{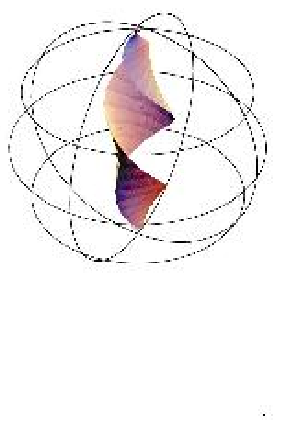}
\includegraphics[width=0.35\linewidth]{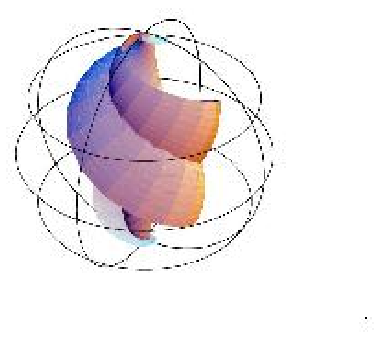}
\end{center}
\vspace{2in}
\caption{The horosphere (on the left) and the helicoid cousin (in the center).  
The horosphere is given by the Bryant representation with 
$(g,\eta)=(c_1,c_2dw)$, $c_1\in\mathbb{C}$, $c_2\in\mathbb{C}\setminus\{0\}$ 
on $\Sigma=\mathbb{C}$, like the data for a plane in $\mathbb{R}^3$.  
The helicoid cousin is given by the Bryant representation with 
$(g,\eta)=(e^w,cie^{-w}dw)$, $c\in\mathbb{R}\setminus\{0\}$ 
on $\Sigma=\mathbb{C}$, like the data for a helicoid in $\mathbb{R}^3$.  
The figure on the right shows a larger portion of the helicoid cousin, so we 
can see that this surface is not embedded.}
\label{fg:horoheli}
\end{figure}
\begin{figure}[phbt]
\begin{center}
\includegraphics[width=0.4\linewidth]{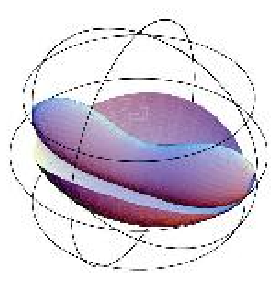}
\includegraphics[width=0.4\linewidth]{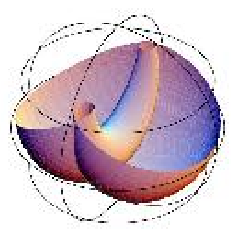}
\end{center}
\vspace{2in}
\caption{An Enneper cousin (on the left) and its dual surface (on the 
right).  These surfaces can be given by the Bryant representation with 
$(g,\eta)=(w,cdw)$, $c\in\mathbb{R}\setminus\{0\}$ on $\Sigma=\mathbb{C}$, 
like the data for Enneper's surface in $\mathbb{R}^3$.  The value of $c$ is 
important, as these surfaces for different $c$ do not differ by only a 
dilation, as would happen in the Weierstrass representation for minimal 
surfaces.  We note that the Enneper cousin dual has infinite total curvature 
even though the Enneper cousin has finite total curvature.}
\label{fg:enn_cousins}
\end{figure}
\begin{figure}[phbt]
\begin{center}
\includegraphics[width=0.4\linewidth]{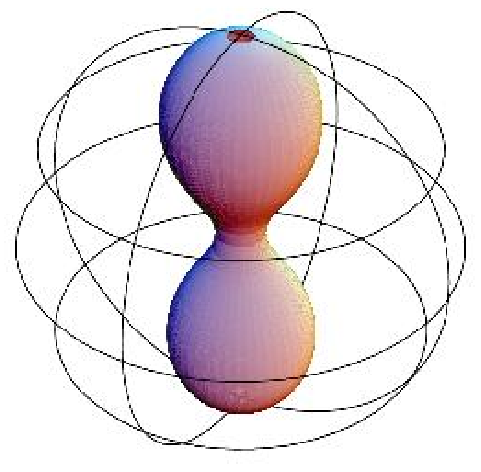}
\includegraphics[width=0.4\linewidth]{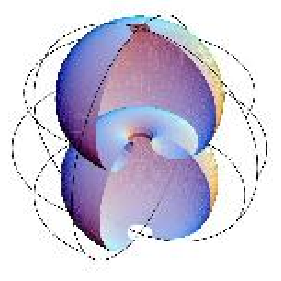}
\end{center}
\vspace{2in}
\caption{Catenoid cousins can be given by the Bryant representation with 
$(g,\eta)=(w^{\mu},(1-\mu^2)dw/(4\mu w^{\mu +1}))$, 
$\mu\in\mathbb{R}\setminus\{1\}$ on $\Sigma=\mathbb{C}\setminus\{0\}$.  
The left-hand side is drawn with $\mu =0.8$ and the right-hand side is drawn 
with $\mu =1.2$.  
$f$ is well-defined on $\Sigma$, even though $g$ (and $\hat F$) is not.  
The catenoid cousin is locally isometric to a minimal catenoid in 
$\mathbb{R}^3$, if one makes the coordinate transformation $z=w^\mu$ and an 
appropriate homothety of the minimal catenoid.}
\label{fg:cmc1cat}
\end{figure}
\begin{figure}[phbt]
\begin{center}
\includegraphics[width=0.4\linewidth]{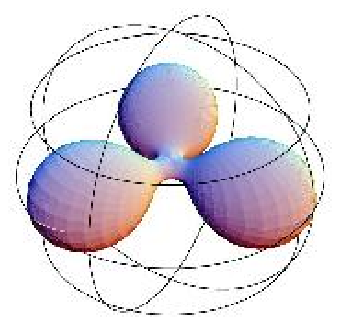}
\includegraphics[width=0.4\linewidth]{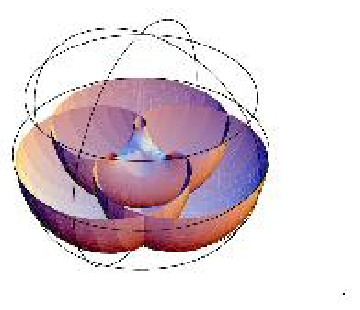}
\end{center}
\vspace{2in}
\caption{Trinoid cousin duals can be given by the Bryant representation with 
$(g,\eta)=(w^2,cdw/(w^3-1)^2)$, $c\in\mathbb{R}\setminus\{0\}$ on 
$\Sigma=(\mathbb{C}\cup\{\infty\}\setminus\{w\in\mathbb{C}|w^3=1\}$, like the 
data for the Jorge-Meeks 3-noid.  The left-hand side is drawn with $c=-0.3$ and the 
right-hand side is drawn with $c=0.8$.  Although numerical experiments show 
that the surface on the left-hand side is embedded, none have yet been proven to be 
embedded.}
\label{fg:cmc1tri}
\end{figure}

\section{Flat surfaces in $\mathbb{H}^3$}
\label{sectionGMM}

A surface is {\em flat} if its intrinsic Gaussian curvature $K$ is identically 
zero.  In this section we describe a representation for flat surfaces in 
$\mathbb{H}^3$ that is due to G\'alvez, Mart\'{\i}nez and Mil\'an \cite{GMM}, 
and is similar in spirit to both the Weierstrass representation for minimal 
surfaces in $\mathbb{R}^3$ and the representation of Bryant for CMC $1$ 
surfaces in $\mathbb{H}^3$.  

Let $\Sigma$ be a simply-connected Riemann surface with 
global coordinate $w=x+i y$.  Let $f: 
\Sigma \to \mathbb{H}^3$ be a smooth flat conformal 
immersion.  (By Theorem \ref{conformalityispossible}, without loss of generality we 
may assume $f$ is conformal.)  
Because $f$ is flat, the first fundamental form is of the form 
\[ ds^2=\rho^2 (dx^2+dy^2) \; , \;\;\; \rho=2 e^u \; , \;\;\; u : 
\Sigma \to \mathbb{R} \; , \] 
with \[ K = - \rho^{-2} \cdot \triangle \log \rho = 0 \; , \]  
where $\triangle=\partial_x \partial_x + 
\partial_y \partial_y$ denotes the Laplacian.  This implies that 
\[ \triangle u = 0 \; . \]  So there exists a function $v:\Sigma \to 
R$ so that $u+i v$ is holomorphic in $w$.  Now define 
\[ z = 2 \int e^{u+i v} dw \; , \] so $z$ is holomorphic in $w$, and 
\[ \frac{dz}{dw} = e^{u+i v} \neq 0 \; , \] and 
\[ ds^2 = \frac{4e^{2u}}{\left| \frac{dz}{dw} \right|^2} dz d\bar z = 4 
dz d\bar z \; . \] So without loss of generality, we may assume that $u$ is 
a constant and is 
\[ u = \frac{1}{2} \log \left( \frac{1}{4} \right) \; . \]  Thus 
$4 e^{2u}=1$ and 
\[ ds^2=dx^2+dy^2 \] with respect to the parameter $w=x+i y$.  

With the usual definitions for the Hopf differential $A=
\langle f_{ww}, N \rangle$ and mean curvature 
$H=(1/2) e^{-2u} \langle f_{w\bar w}, N \rangle 
  =2\langle f_{w\bar w}, N \rangle$, 
the Gauss-Weingarten equations are 
\[ f_{ww} = A \cdot N \; , \;\;\; f_{w\bar w} = \frac{1}{2} f + 
\frac{1}{2} H N \; , \;\;\; N_w = -H f_w - 2 A f_{\bar w} \; , \] 
and the Gauss-Codazzi equations are 
\[ H^2-1=4 A \bar A \; , \;\;\; H_{\bar w} = 2 \bar A_w \; . \]  
So either 
\[ H = +\sqrt{1+4A\bar A} \geq 1 \;\;\; \text{ or } \;\;\; 
H = -\sqrt{1+4A\bar A} \leq -1 \] holds.  Since the surface $f$ is 
smooth, $H$ is continuous, and 
precisely one of these two cases will hold on all of $\Sigma$.  
Let us assume that the first case holds, that is, $H \geq 1$.  
(The second case $H \leq -1$ can be dealt with in similar fashion to the 
following computations.)  

Choose a solution $F \in \SL_2\!\mathbb{C}$ of the Lax pair 
\[ F_w = F \cdot U \; , \;\;\; F_{\bar w} = F \cdot V \; , \] 
with 
\begin{equation}
U = \begin{pmatrix} 0 & (1-H)/2 \\
A & 0 \end{pmatrix} \; , \;\;\; 
V = \begin{pmatrix} 0 & - \bar A \\
(1+H)/2 & 0 \end{pmatrix} \; , \end{equation}
and defining an immersion $\tilde{f}:\Sigma \to \mathbb{H}^3$ via 
the Sym-Bobenko type formula as in Lemma \ref{lm:SL2Cform}, 
\begin{equation}
\tilde f=F \overline{F}^t \; , \end{equation} we find that $\tilde f$ 
is also a flat surface with the same mean curvature $H$ and Hopf differential 
$A$ as the immersion $f$.  Thus $f$ and $\tilde f$ are the same surface, 
up to a rigid motion of $\mathbb{H}^3$.  Hence we may simply 
assume $f = F \overline{F}^t$.  

We can rewrite the above Lax pair as 
\begin{equation}\label{eq:flatLax}
F^{-1} dF = \begin{pmatrix} 0 & \frac{1}{2}(1-H)dw - \bar A d\bar w \\
Adw+\frac{1}{2}(1+H)d\bar w & 0 \end{pmatrix} \; . 
\end{equation}
We see that $F$ is not holomorphic in $w$, because $d\bar w$ terms exist 
in the above equation (or equivalently, because $V$ is not zero).  
Our goal in this section is to show that we can make a change of parameter 
so that $F$ is holomorphic with respect to the new parameter.  

The Codazzi equation implies that there exist functions 
$z=z(w):\Sigma \to \mathbb{C}$ and $\zeta=\zeta(w):\Sigma \to \mathbb{C}$ so that 
\[ F^{-1} dF = \begin{pmatrix} 0 & d\zeta \\
dz & 0 \end{pmatrix} \; . \]  

\begin{lemma}\label{lm:461}
The function $z$ is a global coordinate of $f$.  
\end{lemma}

\begin{proof}
The function $z$ will be global coordinate if the Jacobian 
$z_w\bar z_{\bar w}-z_{\bar w}\bar z_w
= \partial_wz\cdot\partial_{\bar w}\bar z
 -\partial_{\bar w}z\cdot\partial_w\bar z$ is never zero, and since
\[ z_w \bar z_{\bar w}-z_{\bar w} \bar z_w = 
\frac{-1}{2} (1+\sqrt{1+4A\bar A}) \leq -1 \; , \]  
this Jacobian is never zero.  
\end{proof}

\begin{remark}
Note that the coordinate $z$ is {\em not} a conformal coordinate for $f$, 
since $z_{\bar w}=(1+H)/2$ is not zero.  
\end{remark}

\begin{lemma}\label{lm:463}
$\zeta_{\bar z}=\partial_{\bar z}\zeta =0$.  
\end{lemma}  

\begin{proof}
Computing that 
\[ \zeta_{\bar z} = \zeta_w w_{\bar z}+\zeta_{\bar w} \bar w_{\bar z} = 
(1-\bar z_w) w_{\bar z} - \bar z_{\bar w} \bar w_{\bar z}=w_{\bar z} - 1 \; , 
\] 
we see that we need only show $w_{\bar z}=1$.  We have 
\[ w_{\bar z} \bar z_w+w_z z_w= 1 \; , \;\;\; 
w_{\bar z} \bar z_{\bar w}+w_z z_{\bar w}= 0 \; , \] and this implies that 
\[ \begin{pmatrix} w_{\bar z} \\
w_z \end{pmatrix} = \frac{1}{z_{\bar w} \bar z_w-z_w \bar z_{\bar w}}
\begin{pmatrix} z_{\bar w} & -z_{w} \\
-\bar z_{\bar w} & \bar z_{w} \end{pmatrix} 
\begin{pmatrix} 1 \\ 0 \end{pmatrix} \; . \]  
Since $z_{w}$, $z_{\bar w}$, $\bar z_{w}$, $\bar z_{\bar w}$ are all known 
in terms of $A$ and $H$, we can compute that 
\[ 
w_{\bar z}=\frac{2(1+H)}{(1+H)^2-4A\bar A} = 1 \; . 
\]  
\end{proof}

So with this new coordinate $z$, we have that 
\[
f= F \overline{F}^t \; , \;\;\text{where} \;\;\;
dF=F \begin{pmatrix} 0 & h(z) \\
                     1 & 0 \end{pmatrix} dz \] 
for some holomorphic function $h(z)$.  

If we had assumed the second case $H \leq -1$, then $\zeta$ would 
have been the global coordinate and $z$ would have been a holomorphic function 
of $\zeta$.  
Thus the analogous result would hold, and the final relation would have been 
\[ f= F \overline{F}^t \; , \;\;\text{where} \;\;\;
dF=F \begin{pmatrix} 0 & 1 \\
h(\zeta) & 0 \end{pmatrix} d\zeta \] 
for some holomorphic function $h(\zeta)$.  

\begin{figure}[phbt]
\begin{center}
\includegraphics[width=0.4\linewidth]{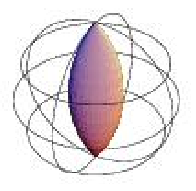}
\includegraphics[width=0.4\linewidth]{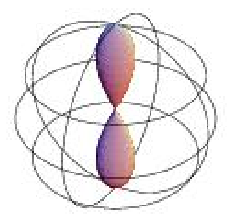}
\end{center}
\vspace{2in}
\caption{Flat surfaces of revolution for $\alpha =-1$ and 
$\alpha <0$ ($\alpha\ne -1$).  (These graphics were made 
by Hiroya Shimizu using Mathematica.)}
\label{fg:flatH3}
\end{figure}

\begin{figure}[phbt]
\begin{center}
\includegraphics[width=0.4\linewidth]{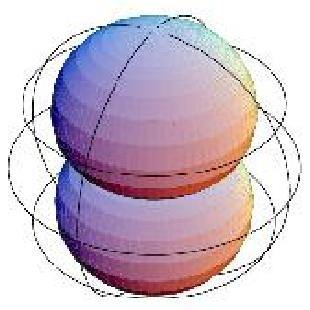}
\includegraphics[width=0.4\linewidth]{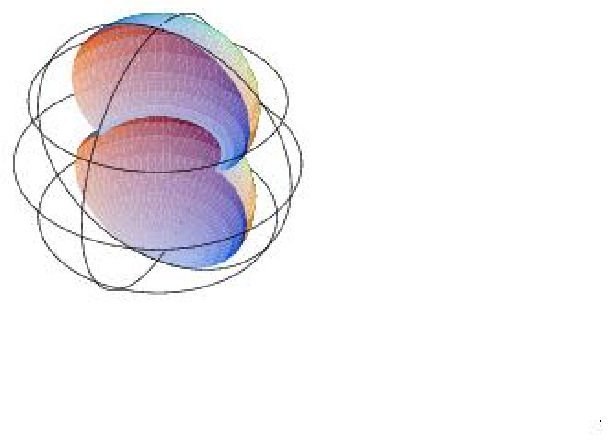}
\end{center}
\vspace{2in}
\caption{Flat surfaces of revolution for 
$\alpha >0$.  (These graphics were made 
by Hiroya Shimizu using Mathematica.)}
\label{fg:flatH3b}
\end{figure}

To include both these cases, we state the conclusion as follows:

\begin{quote}
{\bf The representation of G\'alvez, Mart\'{\i}nez and Mil\'an \cite{GMM}}:
A flat immersion $f: \Sigma \to \mathbb{H}^3$ can be described as 
\[ f= F \overline{F}^t \; , \;\; \text{where}\;\;\; 
dF=F \begin{pmatrix} 0 & h(z) \\
g(z) & 0 \end{pmatrix} dz \] 
for some holomorphic functions $h(z)$ and $g(z)$ defined on $\Sigma$ with 
coordinate $z$.  
\end{quote}

\begin{remark}
This method above will not work for $K=0$ surfaces in $\mathbb{R}^3$ or 
$\mathbb{S}^3$, for the following reason: In the case of $\mathbb{R}^3$, 
Equation \eqref{eq:flatLax} becomes 
\[ F^{-1} dF = \begin{pmatrix} 0 & \frac{-1}{2}Hdw - \bar A d\bar w \\
Adw+\frac{1}{2}Hd\bar w & 0 \end{pmatrix} = \begin{pmatrix} 0 & d\zeta \\
dz & 0 \end{pmatrix} \; , \]  with $H^2=4A\bar A$ and $H_{\bar w} = 
2 \bar A_w$.  In this case we find that 
$\zeta_w \bar \zeta_{\bar w}-\zeta_{\bar w} \bar \zeta_w = 
z_w \bar z_{\bar w}-z_{\bar w} \bar z_w = 0$, and thus neither $z$ nor 
$\zeta$ can become the new coordinate.  

In the case of $\mathbb{S}^3$, Equation \eqref{eq:flatLax} becomes 
\[ F^{-1} dF = \begin{pmatrix} 0 & \frac{1}{2}(-i-H)dw - \bar A d\bar w \\
Adw+\frac{1}{2}(-i+H)d\bar w & 0 \end{pmatrix} = \begin{pmatrix} 0 & d\zeta \\
dz & 0 \end{pmatrix} \; , \]   with $H^2+1=4A\bar A$ and $H_{\bar w} = 
2 \bar A_w$.  Then, just as in the case of $\mathbb{R}^3$, neither $z$ nor 
$\zeta$ can become the new coordinate.  
\end{remark}

\begin{example}
Let $\alpha\in\mathbb{R}\setminus\{0,1\}$ be a constant and define holomorphic 
functions $g$, $h$ on $\Sigma=\mathbb{C}\setminus\{0\}$ as follows:
\[
g=-\frac{1}{c^2}z^{-2/(1-\alpha)} \; \, \quad
h=\frac{c^2\alpha}{(1-\alpha)^2}z^{2\alpha /(1-\alpha)} 
\]
for some constant $c$.  
Then we have a hyperbolic cylinder if $\alpha =-1$ and an hourglass if 
$\alpha <0$ ($\alpha\ne -1$) and a snowman if $\alpha >0$.  
If we set $\alpha =0$ and replace $\Sigma$ by $\mathbb{C}$, then we have a 
horosphere as shown in Figure \ref{fg:horoheli}.  See \cite{KoUY}, \cite{krsuy}, 
\cite{kruy}.  
\end{example}


\end{document}